\documentclass[11pt,reqno]{article}
\usepackage{amsfonts}
\usepackage{graphicx}
\usepackage{imakeidx}
\usepackage{amsmath}
\usepackage{amssymb}
\usepackage{xcolor}

\begin{document}

\noindent
{\Large\bf Gibbs measures for contact Anosov flows are exponentially mixing}

\bigskip

\noindent
{\bf Luchezar Stoyanov\footnote{luchezar.stoyanov@uwa.edu.au;
Department of Mathematics and Statistics,
University of Western Australia, 35 Stirling Hwy, Crawley WA 6009, Australia}}

\bigskip

\def\R{{\mathbb R}}
\def\T{{\mathbb T}}
\def\S{{\mathbb S}}
\def\C{{\mathbb C}}
\def\Z{{\mathbb Z}}
\def\N{{\mathbb N}}
\def\H{{\mathbb H}}
\def\B{{\mathbb B}}
\def\diam{\mbox{\rm diam}}
\def\sn{\S^{n-1}}
\def\rr{{\cal R}}
\def\mt{{\Lambda}}
\def\e{\emptyset}
\def\dQ{\partial Q}
\def\dk{\partial K}
\def\endofproof{{\rule{6pt}{6pt}}}
\def\di{\displaystyle}
\def\dist{\mbox{\rm dist}}
\def\sa+{\Sigma_A^+}
\def\du{\frac{\partial}{\partial u}}
\def\dv{\frac{\partial}{\partial v}}
\def\dt{\frac{d}{d t}}
\def\dx{\frac{\partial}{\partial x}}
\def\con{\mbox{\rm const }}
\def\mat{\mm_{at}}
\def\ma{\mm_{a}}
\def\lab{L_{ab}}
\def\labt{L_{abt}}
\def\mabn{\mm_{a}^N}
\def\man{\mm_a^N}
\def\labn{L_{ab}^N}
\def\fa{f^{(a)}}
\def\i{{\bf i}}
\def\gge{{\mathcal G}_\epsilon}
\def\gej{\chi^{(j)}_\mu}
\def\ge{\chi_\epsilon}
\def\chio{\chi^{(1)}}
\def\chit{\chi^{(2)}}
\def\chii{\chi^{(i)}}
\def\chil{\chi^{(\ell)}}
\def\gett{\chi^{(2)}_{\mu}}
\def\geol{\chi^{(1)}_{\ell}}
\def\getl{\chi^{(2)}_{\ell}}
\def\geil{\chi^{(i)}_{\ell}}
\def\gee{\chi_{\ell}}
\def\wloc{W_{\epsilon}}
\def\Int{\mbox{\rm Int}}
\def\dist{\mbox{\rm dist}}
\def\pr{\mbox{\rm pr}}
\def\tpp{\widetilde{\pp}}
\def\supp{\mbox{\rm supp}}
\def\Arg{\mbox{\rm Arg}}
\def\In{\mbox{\rm Int}}
\def\con{\mbox{\rm const}\;}
\def\Re{\mbox{\rm Re}}
\def\li{\mbox{\rm li}} 
\def\Seo{S^*_\epsilon(\Omega)}
\def\sdk{S^*_{\dk}(\Omega)}
\def\lae{\Lambda_{\epsilon}}
\def\ep{\epsilon}
\def\be{\begin{equation}}
\def\ee{\end{equation}}
\def\beqn{\begin{eqnarray*}}
\def\eeqn{\end{eqnarray*}}
\def\Pr{\mbox{\rm Pr}}

\def\gi{\gamma^{(i)}}
\def\ii{{\imath }}
\def\jj{{\jmath }}
\def\II{{\mathcal I}}
\def\ccij{ \cc_{i'_0,j'_0}[\eta]}
\def\la{\langle}
\def\ra{\rangle}
\def\bs{\bigskip}
\def\xio{\xi^{(0)}}
\def\xo{x^{(0)}}
\def\zo{z^{(0)}}
\def\Con{\mbox{\rm Const}\;}
\def\do{\partial \Omega}
\def\dk{\partial K}
\def\dl{\partial L}
\def\pr{{\rm pr}}
\def\dist{{\rm dist}}
\def\dds{\frac{d}{ds}}
\def\con{{\rm const}\;}
\def\Con{{\rm Const}\;}
\def\di{\displaystyle}
\def\oo{\mbox{\rm O}}
\def\hess{\mbox{\rm Hess}}
\def\gi{\gamma^{(i)}}
\def\endofproof{{\rule{6pt}{6pt}}}
\def\xm{x^{(m)}}
\def\vm{\varphi^{(m)}}
\def\km{k^{(m)}}
\def\dm{d^{(m)}}
\def\kam{\kappa^{(m)}}
\def\dem{\delta^{(m)}}
\def\xim{\xi^{(m)}}
\def\ep{\epsilon}
\def\ms{\medskip}
\def\ex{\mbox{\rm extd}}

\def\clip{C^{\mbox{\footnotesize \rm Lip}}}
\def\wlocs{W^s_{\mbox{\footnote\rm loc}}}
\def\Lip{\mbox{\rm Lip}}

\def\Xr{X^{(r)}}
\def\lip{\mbox{{\footnotesize\rm Lip}}}
\def\Vol{\mbox{\rm Vol}}

\def\naf{\nabla f(z)}
\def\so{\sigma_0}
\def\Xo{X^{(0)}}
\def\z1{z^{(1)}}
\def\Vo{V^{(0)}}
\def\Yo{Y{(0)}}

\def\uo{u^{(0)}}
\def\vo{v^{(0)}}
\def\no{\nu^{(0)}}
\def\psa{\partial^{(s)}_a}
\def\hcd{\hc^{(\delta)}}
\def\Md{M^{(\delta)}}
\def\Uo{U^{(1)}}
\def\Ut{U^{(2)}}
\def\Uj{U^{(j)}}
\def\no{n^{(1)}}
\def\nt{n^{(2)}}
\def\nj{n^{(j)}}
\def\ccm{\cc^{(m)}}

\def\f0{f^{(0)}}

\def\gl{\gamma_\ell}
\def\id{\mbox{\rm id}}
\def\piU{\pi^{(U)}}

\def\cca{C^{(a)}}
\def\bba{B^{(a)}}
\def\saa{\Sigma^+_A}
\def\sa{\Sigma_A}

\def\Int{\mbox{\rm Int}}
\def\epo{\ep^{(0)}}
\def\pH{\partial \H^{n+1}}
\def\sh{S^*(\H^{n+1})}
\def\zoo{z^{(1)}}
\def\yoo{y^{(1)}}
\def\xoo{x^{(1)}}


\def\supp{\mbox{\rm supp}}
\def\Arg{\mbox{\rm Arg}}
\def\In{\mbox{\rm Int}}
\def\diam{\mbox{\rm diam}}
\def\e{\emptyset}
\def\endofproof{{\rule{6pt}{6pt}}}
\def\di{\displaystyle}
\def\dist{\mbox{\rm dist}}
\def\con{\mbox{\rm const }}
\def\Box{\spadesuit}
\def\Int{\mbox{\rm Int}}
\def\dist{\mbox{\rm dist}}
\def\pr{\mbox{\rm pr}}
\def\be{\begin{equation}}
\def\ee{\end{equation}}
\def\beqn{\begin{eqnarray*}}
\def\eeqn{\end{eqnarray*}}
\def\la{\langle}
\def\ra{\rangle}
\def\bs{\bigskip}
\def\Con{\mbox{\rm Const}\;}
\def\clip{C^{\alpha}}
\def\wlocs{W^s_{\mbox{\footnote\rm loc}}}
\def\Lip{\mbox{\rm Lip}}
\def\lip{\mbox{\footnotesize\rm Lip}}
\def\Re{\mbox{\rm Re}}
\def\li{\mbox{\rm li}} 
\def\ep{\epsilon}
\def\ms{\medskip}
\def\dds{\frac{d}{ds}}
\def\oo{\mbox{\rm O}}
\def\hess{\mbox{\rm Hess}}
\def\id{\mbox{\rm id}}
\def\ii{{\imath }}
\def\jj{{\jmath }}
\def\graph{\mbox{\rm graph}}
\def\span{\mbox{\rm span}}
\def\Intu{\Int^u}
\def\Ints{\Int^s}

\def\i{{\bf i}}

\def\G{{\mathcal G}}
\def\nn{{\mathcal N}}
\def\mm{{\mathcal M}}
\def\kk{{\mathcal K}}
\def\ll{{\mathcal L}}
\def\vv{{\mathcal V}}
\def\bb{{\mathcal B}}
\def\ff{{\mathcal F}}
\def\tt{{\mathcal T}}
\def\uu{{\mathcal U}}
\def\pp{{\mathcal P}}
\def\qq{{\mathcal Q}}
\def\aa{{\mathcal A}}
\def\cc{{\mathcal C}}
\def\oo{{\mathcal O}}
\def\dd{{\mathcal D}}
\def\ss{{\mathcal S}}
\def\Ee{{\mathcal E}}
\def\tEe{{\widetilde{\cal E}}}
\def\rr{{\mathcal R}}
\def\hh{{\mathcal H}}
\def\II{{\mathcal I}}
\def\jj{{\mathcal J}}

\def\hs{\hat{s}}
\def\hz{\hat{z}}
\def\hL{\widehat{L}}
\def\hl{\hat{l}}
\def\hl{\hat{l}}
\def\hc{\hat{\cc}}
\def\hbb{\widehat{\cal B}}
\def\hu{\hat{u}}
\def\hX{\hat{X}}
\def\hx{\hat{x}}
\def\hu{\hat{u}}
\def\hv{\hat{v}}
\def\hQ{\hat{Q}}
\def\hC{\widehat{C}}
\def\hF{\hat{F}}
\def\hf{\hat{\varphi}}
\def\hvarphi{\hat{\varphi}}
\def\hii{\hat{\ii}}
\def\hr{\hat{r}}
\def\hq{\hat{q}}
\def\hy{\hat{y}}
\def\hZ{\widehat{Z}}
\def\hz{\hat{z}}
\def\hE{\widehat{E}}
\def\hR{\widehat{R}}
\def\hell{\hat{\ell}}
\def\hs{\hat{s}}
\def\hW{\widehat{W}}
\def\hS{\widehat{S}}
\def\hV{\widehat{V}}
\def\hB{\widehat{B}}
\def\hhh{\widehat{\cal H}}
\def\hK{\widehat{K}}
\def\hU{\widehat{U}}
\def\hhh{\widehat{\hh}}
\def\hdd{\widehat{\dd}}
\def\hZ{\widehat{Z}}
\def\hGa{\widehat{\Gamma}}

\def\hal{\hat{\alpha}}
\def\hbe{\hat{\beta}}
\def\hg{\hat{\gamma}}
\def\hrho{\hat{\rho}}
\def\hd{\hat{\delta}}
\def\hphi{\hat{\phi}}
\def\hmu{\hat{\mu}}
\def\hnu{\hat{\nu}}
\def\hsi{\hat{\sigma}}
\def\htau{\hat{\tau}}
\def\hpi{\hat{\pi}}
\def\hep{\hat{\epsilon}}
\def\hxi{\hat{\xi}}
\def\hLa{\widehat{\Lambda}^u}
\def\hPhi{\widehat{\Phi}}
\def\hPsi{\widehat{\Psi}}
\def\hPhii{\widehat{\Phi}^{(i)}}
\def\hath{\hat{h}}

\def\tc{\tilde{c}}
\def\tC{\widetilde{C}}
\def\tg{\tilde{\gamma}}  
\def\tV{\widetilde{V}}
\def\tb{\tilde{b}}
\def\tt{\tilde{t}}
\def\tx{\tilde{x}}
\def\tp{\tilde{p}}
\def\tz{\tilde{Z}}
\def\tZ{\tilde{Z}}
\def\tF{\tilde{F}}
\def\tf{\tilde{\varphi}}
\def\tvarphi{\tilde{\varphi}}
\def\tp{\tilde{p}}
\def\te{\tilde{e}}
\def\tv{\tilde{v}}
\def\tu{\tilde{u}}
\def\tw{\tilde{w}}
\def\ts{\tilde{\sigma}}
\def\talpha{\tilde{\alpha}}
\def\tr{\tilde{r}}
\def\tU{\widetilde{U}}
\def\tS{\tilde{S}}
\def\tP{\widetilde{P}}
\def\ttau{\tilde{\tau}}
\def\tLip{\widetilde{\Lip}}
\def\tz{\tilde{z}}
\def\tS{\tilde{S}}
\def\tts{\tilde{\sigma}}
\def\tVl{\widetilde{V}^{(\ell)}}
\def\tVj{\widetilde{V}^{(j)}}
\def\tVo{\widetilde{V}^{(1)}}
\def\tVj{\widetilde{V}^{(j)}}
\def\tPsi{\tilde{\Psi}}
 \def\tp{\tilde{p}}
 \def\tVjo{\widetilde{V}^{(j_0)}}
\def\tvj{\tilde{v}^{(j)}}
\def\tVjj{\widetilde{V}^{(j+1)}}
\def\tvl{\tilde{v}^{(\ell)}}
\def\tVt{\widetilde{V}^{(2)}}
\def\tR{\widetilde{R}}
\def\tQ{\widetilde{Q}}
\def\oL{\tilde{\Lambda}}
\def\tq{\tilde{q}}
\def\tk{\tilde{k}}
\def\tx{\tilde{x}}
\def\ty{\tilde{y}}
\def\tz{\tilde{z}}
\def\txo{\tilde{x}^{(0)}}
\def\tso{\tilde{\sigma}_0}
\def\tmt{\tilde{\Lambda}}
\def\tg{\tilde{g}}
\def\tsi{\tilde{\sigma}}
\def\tell{\tilde{\ell}}
\def\trho{\tilde{\rho}}
\def\ts{\tilde{s}}
\def\tB{\widetilde{B}}
\def\thh{\widetilde{\cal H}}
\def\tV{\widetilde{V}}
\def\trr{\tilde{r}}
\def\te{\tilde{e}}
\def\tv{\tilde{v}}
\def\tu{\tilde{u}}
\def\tw{\tilde{w}}
\def\trho{\tilde{\rho}}
\def\tell{\tilde{\ell}}
\def\tz{\tilde{Z}}
\def\tF{\tilde{F}}
\def\tf{\tilde{f}}
\def\tp{\tilde{p}}
\def\ttau{\tilde{\tau}}
\def\tz{\tilde{z}}
\def\tg{\tilde{\gamma}}  
\def\tV{\widetilde{V}}
\def\tCC{\widetilde{\cc}}
\def\tLa{\widetilde{\Lambda}^u}
\def\tR{\widetilde{R}}
\def\tr{\tilde{r}}
\def\tD{\widetilde{D}}
\def\tt{\tilde{t}}
\def\tp{\tilde{p}}
\def\tS{\tilde{S}}
\def\tts{\tilde{\sigma}}
\def\tZ{\widetilde{Z}}
\def\tdelta{\tilde{\delta}}
\def\th{\tilde{h}}
\def\tB{\widetilde{B}}
\def\thh{\widetilde{\hh}}
\def\tep{\tilde{\ep}}
\def\tE{\widetilde{E}}
\def\tu{\tilde{u}}
\def\txi{\tilde{\xi}}
\def\teta{\tilde{\eta}}
\def\tRR{\widetilde{\rr}}

\def\sr{{\sc r}}
\def\mt{{\Lambda}}
\def\do{\partial \Omega}
\def\dk{\partial K}
\def\dl{\partial L}
\def\wloc{W_{\epsilon}}
\def\piU{\pi^{(U)}}
\def\Rio{\R_{i_0}}
\def\Ri{\R_{i}}
\def\Rii{\R^{(i)}}
\def\Riii{\R^{(i-1)}}
\def\hRii{\widehat{\R}_i}
\def\hRiio{\widehat{\R}_{(i_0)}}
\def\Eii{E^{(i)}}
\def\Eio{E^{(i_0)}}
\def\Rj{\R_{j}}
\def\Vio{{\cal V}^{i_0}}
\def\Vi{{\cal V}^{i}}
\def\Wio{W^{i_0}}
\def\Wioo{W^{i_0-1}}
\def\hi{h^{(i)}}
\def\Psii{\Psi^{(i)}}
\def\pii{\pi^{(i)}}
\def\piii{\pi^{(i-1)}}
\def\gxyii{g_{x,y}^{i-1}}
\def\span{\mbox{\rm span}}
\def\Jac{\mbox{\rm Jac}}
\def\Vol{\mbox{\rm Vol}}
\def\limp{\lim_{p\to\infty}}

\def\saa{\Sigma^+_{\aa}}
\def\sa{\Sigma_{\aa}}
\def\mta{\Lambda({\aa}, \tau)}
\def\mtaa{\Lambda^+({\aa}, \tau)}

\def\xijl{X_{i,j}^{(\ell)}}
\def\xij{X_{i,j}}
\def\hyijl{\widehat{Y}_{i,j}^{(\ell)}}
\def\hxijl{\widehat{X}_{i,j}^{(\ell)}}
\def\hxij{\widehat{X}_{i,j}}
\def\eijl{\omega_{i,j}^{(\ell)}}
\def\eij{\omega_{i,j}}
\def\Gl{\Gamma_\ell}

\def\cB{\check{B}}
\def\tpi{\tilde{\pi}}
\def\J{{\sf J}}
\def\bJ{{\mathbb J}}

\def\hcc{\widehat{\cc}}
\def\hpp{\widehat{\pp}}
\def\ttP{\widetilde{\pp}}
\def\tP{\widetilde{P}}
\def\hP{\widehat{P}}
\def\hY{\widehat{Y}}

\def\diamtef{{\footnotesize \diam_\theta}}


\def\tg{\tilde{\gamma}}  
\def\tV{\widetilde{V}}
\def\tW{\widetilde{W}}
\def\tCC{\widetilde{\cc}}
\def\tKo{\widetilde{K_0}}
\def\tUKo{\widetilde{U\setminus K_0}}

\def\wo{w^{(1)}}
\def\vo{v^{(1)}}
\def\uo{u^{(1)}}
\def\wt{w^{(2)}}
\def\xio{\xi^{(1)}}
\def\hxio{\hat{\xi}^{(1)}}
\def\xit{\xi^{(2)}}
\def\etao{\eta^{(1)}}
\def\etat{\eta^{(2)}}
\def\zetao{\zeta^{(1)}}
\def\zetat{\zeta^{(2)}}
\def\vt{v^{(2)}}
\def\ut{u^{(2)}}
\def\Wo{W^{(1)}}
\def\Vo{V^{(1)}}
\def\Uo{U^{(1)}}
\def\Wt{W^{(2)}}
\def\Vt{V^{(2)}}
\def\Ut{U^{(2)}}
\def\tmu{\tilde{\mu}}
\def\tla{\tilde{\lambda}}
\def\diamf{{\rm\footnotesize diam}}
\def\Intu{\mbox{\rm Int}^u}
\def\Ints{\mbox{\rm Int}^s}

\def\Bmt{\overline{B_{\ep_0}(\mt)}}
\def\Lye{L_{y,\eta}}
\def\Lyep{L^{(p)}_{y,\eta}}
\def\Fyp{F^{(p)}_y}
\def\Fxp{F^{(p)}_x}
\def\Lxx{L_{x,\xi}}
\def\Lxxp{L^{(p)}_{x,\xi}}

\def\tc{\tilde{C}}
\def\tg{\tilde{\gamma}}  
\def\tV{\widetilde{V}}
\def\tW{\widetilde{W}}
\def\tC{\widetilde{\cc}}
\def\tKo{\widetilde{K_0}}
\def\tUKo{\widetilde{U\setminus K_0}}

\def\wo{w^{(1)}}
\def\vo{v^{(1)}}
\def\uo{u^{(1)}}
\def\wt{w^{(2)}}
\def\xio{\xi^{(1)}}
\def\xit{\xi^{(2)}}
\def\etao{\eta^{(1)}}
\def\etat{\eta^{(2)}}
\def\zetao{\zeta^{(1)}}
\def\zetat{\zeta^{(2)}}
\def\vt{v^{(2)}}
\def\ut{u^{(2)}}
\def\ui{u^{(i)}}
\def\Wo{W^{(1)}}
\def\Vo{V^{(1)}}
\def\Uo{U^{(1)}}
\def\Wt{W^{(2)}}
\def\Vt{V^{(2)}}
\def\Ut{U^{(2)}}
\def\tmu{\tilde{\mu}}
\def\tla{\tilde{\lambda}}
\def\diamf{{\rm\footnotesize diam}}
\def\Intu{\mbox{\rm Int}^u}
\def\Ints{\mbox{\rm Int}^s}

\def\Wuo{W^{u,1}}
\def\Wui{W^{u,i}}
\def\Wuj{W^{u,j}}
\def\Wut{\tW^{u,2}}
\def\Wuk{W^{u,k}}
\def\Wuh{\hW^{u}}
\def\tWuo{\tW^{u,1}}
\def\tWui{\tW^{u,i}}
\def\tWuj{\tW^{u,j}}
\def\tWuk{\tW^{u,k}}
\def\hWuo{\hW^{u,1}}
\def\hWui{\hW^{u,i}}
\def\hWuj{\hW^{u,j}}
\def\hWuk{\hW^{u,k}}
\def\dj{\delta^{(j)}}
\def\do{\delta^{(1)}}
\def\epj{\ep^{(j)}}
\def\epo{\ep^{(1)}}
\def\hSj{\widehat{S}^{(j)}}
\def\hSo{\widehat{S}^{(1)}}

\def\hGammam{\widehat{\Gamma}^{(m)}}

\def\Gammam{\Gamma^{(m)}}

\def\tmu{\tilde{\mu}}
\def\tla{\tilde{\lambda}}
\def\hE{\widehat{E}}


\def\tPsi{\widetilde{\Psi}}
\def\chU{\check{U}}


\def\Ulo{U^{(\ell_0)}}
\def\dte{D_\theta}
\def\diamte{\mbox{\rm diam}_{\theta}}
\def\Ial{I^{(\alpha)}}
\def\uml{u_m^{(\ell)}}
\def\yl{y^{(\ell)}}
\def\tyl{\tilde{y}^{(\ell)}}
\def\ool{\oo^{(\ell)}}
\def\fl{f^{(\ell)}}
\def\hep{\hat{\ep}}
\def\dl{d^{(\ell)}}
\def\dli{d_{\ell,i}}
\def\dlo{d_{\ell,1}}
\def\dlt{d_{\ell,2}}
\def\Lipt{{\Lip_\theta}}
\def\lipt{{\footnotesize \Lip_\theta}}
\def\tm{\tilde{m}}
\def\tj{\tilde{j}}
\def\lengthf{\mbox{\rm\footnotesize length}}
\def\length{\mbox{\rm length}}


\def\Xijl{X^{(\ell)}_{i,j}}
\def\hXijl{\widehat{X}^{(\ell)}_{i,j}}
\def\Wl{W^{(\ell)}}
\def\omijl{\omega^{(\ell)}_{i,j}}

\def\hXitl{\widehat{X}^{(\ell)}_{i,t}}
\def\Vl{V^{(\ell)}}
\def\omitl{\omega^{(\ell)}_{i,t}}
\def\Xisl{X^{(\ell)}_{i,s}}
\def\hXisl{\widehat{X}^{(\ell)}_{i,s}}
\def\omisl{\omega^{(\ell)}_{i,s}}
\def\hGa{\widehat{\Gamma}}
\def\hOm{\widehat{\Omega}}
\def\tGa{\widetilde{\Gamma}}
\def\hA{\widehat{A}}
\def\tnu{\tilde{\nu}}
\def\tX{\widetilde{X}}

\def\ww{{\mathcal W}}
\def\Zl{Z^{(\ell)}}
\def\hpp{\widehat{\pp}}
\def\tnn{\widetilde{\nn}}

\def\ftt{f^{(t)}}
\def\f0{f^{(0)}}
\def\fat{f^{(at)}}
\def\Fat{F^{(at)}}
\def\Fa{F^{(a)}}
\def\F0{F^{(0)}}
\def\tu{\tilde{u}}
\def\tD{\widetilde{D}}
\def\tchi{\tilde{\chi}}
\def\hC{\widehat{C}}
\def\hQ{\widehat{Q}}
\def\hF{\widehat{F}}
\def\hD{\widehat{D}}
\def\hr{\hat{r}}
\def\psid{\psi^\dag}
\def\taud{\tau^\dag}
\def\Omn{\Omega^{(n)}}
\def\Omm{\Omega^{(m)}}
\def\Omk{\Omega^{(k)}}
\def\Conf{{\mbox{\footnotesize\rm Const}}}
\def\hp{\hat{p}}

\def\tj{t^{(j)}}
\def\tyj{\tilde{y}^{(j)}}
\def\tyjo{\tilde{y}_{j,1}}
\def\tyjt{\tilde{y}_{j,2}}
\def\tyji{\tilde{y}_{j,i}}
\def\yjo{y_{j,1}}
\def\yjt{y_{j,2}}
\def\yji{y_{j,i}}
\def\tylo{\tilde{y}_{\ell,1}}
\def\tylt{\tilde{y}_{\ell,2}}
\def\tyli{\tilde{y}_{\ell,i}}
\def\ylo{y_{\ell,1}}
\def\ylt{y_{\ell,2}}
\def\yli{y_{\ell,i}}

\def\tdlo{\tilde{d}_{\ell,1}}
\def\tdlt{\tilde{d}_{\ell,2}}
\def\tdli{\tilde{d}_{\ell,i}}
\def\dlo{d_{\ell,1}}
\def\dlt{d_{\ell,2}}
\def\dli{d_{\ell,i}}
\def\wjo{w_{j,1}}
\def\wjt{w_{j,2}}
\def\wji{w_{j,i}}
\def\sj{s^{(j)}}
\def\Yj{Y^{(j)}}
\def\Vj{V^{(j)}}
\def\Zj{Z^{(j)}}
\def\vj{v^{(j)}}
\def\wj{w^{(j)}}
\def\twj{\tilde{w}^{(j)}}
\def\gj{g^{(j)}}
\def\tgj{\tilde{g}^{(j)}}
\def\tg{\tilde{g}}
\def\hn{\hat{n}}

\def\hbeta{\hat{\beta}}
\def\hmu{\hat{\mu}}
\def\piS{\pi^{(S)}}
\def\hb{\hat{b}}
\def\shP{P^\sharp}
\def\tshP{\widetilde{P}^\sharp}
\def\T{\mathcal T}
\def\tut{\tu^{(2)}}
\def\twt{\tw^{(2)}}
\def\piS{\pi^{(S)}}
\def\tsigma{\tilde{\sigma}}
\def\td{\tilde{d}}
\def\m{{\sf m}}
\def\tXi{\widetilde{\Xi}}

\def\Omb{\Omega^{(\hb)}}
\def\Xib{\Xi^{(L\hb)}}

\def\tomega{\tilde{\omega}}
\def\vlim{v^{(\ell)}_{i,m}}
\def\tvlim{\tilde{v}^{(\ell)}_{i,m}}
\def\thetlim{\theta^{(\ell)}_{i,m}}
\def\lamblim{\lambda^{(\ell)}_{i,m}}
\def\vl{v^{(\ell)}}
\def\vlik{v^{(\ell)}_{i,k}}
\def\vlom{v^{(\ell)}_{1,m}}
\def\vltm{v^{(\ell)}_{2,m}}
\def\tvlik{\tilde{v}^{(\ell)}_{i,k}}
\def\tvlom{\tilde{v}^{(\ell)}_{1,m}}
\def\tvltm{\tilde{v}^{(\ell)}_{2,m}}
\def\Cum{\cc^{(u_m)}}
\def\Cu{\cc^{(u)}}
\def\gao{\gamma^{(1)}}
\def\gat{\gamma^{(2)}}
\def\diamtef{{\footnotesize\mbox{\rm  diam}_\theta}}
\def\halpha{\hat{\alpha}}
\def\hgamma{\hat{\gamma}}
\def\tXijl{\widetilde{X}^{(\ell)}_{i.j}}
\def\tdd{\widetilde{\dd}}
\def\hcc{\widehat{\cc}}
\def\lambdam{\lambda^{(m)}}
\def\thetam{\theta^{(m)}}

\def\tvlo{\tilde{v}^{(\ell)}_{1}}
\def\tvlt{\tilde{v}^{(\ell)}_{2}}
\def\tvli{\tilde{v}^{(\ell)}_{i}}

\def\omegak{\omega^{(k)}}
\def\nnk{\nn^{(k)}}
\def\omegao{\omega^{(o)}}
\def\nno{\nn^{(o)}}
\def\htheta{\hat{\theta}}

\def\vlo{v^{(\ell)}_{1}}
\def\vlt{v^{(\ell)}_2}
\def\vli{v^{(\ell)}_{i}}
\def\vloij{v^{(\ell,0)}_{i,j}}
\def\vlooj{v^{(\ell,0)}_{1,j}}
\def\vlotj{v^{(\ell,0)}_{2,j}}
\def\vlkij{v^{(\ell,k)}_{i,j}}
\def\vlkoj{v^{(\ell,k)}_{1,j}}
\def\vlktj{v^{(\ell,k)}_{2,j}}

\def\ho{h^{(1)}}
\def\Ho{H^{(1)}}
\def\tcc{\widetilde{\cc}}
\def\utk{u^{(\tk)}}
\def\hXi{\widehat{\Xi}}
\def\chR{\check{R}}
\def\llm{\ll^{(m)}}
\def\llk{\ll^{(k)}}

\def\hB{\widehat{B}}
\def\Omb{\Omega^{(\hb)}}
\def\tnn{\widetilde{\nn}}

\def\Xim{\Xi^{(m)}}
\def\Omm{\Omega^{(m)}}
\def\Xin{\Xi^{(n)}}
\def\Omn{\Omega^{(n)}}
\def\hP{\widehat{P}}

\def\halpha{\hat{\alpha}}
\def\hgamma{\hat{\gamma}}

\def\E{\mathcal E}
\def\tp{\tilde{p}}
\def\tq{\tilde{q}}

\def\tyj{\tilde{y}^{(j)}}
\def\yl{y^{(\ell)}}
\def\wo{w^{(1)}}
\def\wi{w^{(i)}}
\def\wj{w^{(j)}}
\def\vo{v^{(1)}}
\def\vi{v^{(i)}}
\def\vj{v^{(j)}}
\def\etao{\eta^{(1)}}
\def\etai{\eta^{(i)}}
\def\etaj{\eta^{(j)}}
\def\zetao{\zeta^{(1)}}
\def\zetai{\zeta^{(i)}}
\def\zetaj{\zeta^{(j)}}
\def\etak{v^{(\tk)}}
\def\uo{u^{(1)}}
\def\ui{u^{(i)}}
\def\uj{u^{(j)}}
\def\huo{\hat{u}^{(1)}}
\def\wt{w^{(2)}}
\def\xio{\xi^{(1)}}
\def\xit{\xi^{(2)}}
\def\xii{\xi^{(i)}}
\def\xij{\xi^{(j)}}

\def\ao{a^{(1)}}
\def\bo{b^{(1)}}
\def\zi{z^{(i)}}
\def\hGa{\widehat{\Gamma}}

\def\tvo{\tv^{(1)}}
\def\tetao{\teta^{(1)}}

\def\Bmt{\overline{B_{\ep_0}(\mt)}}
\def\chBo{\check{B}^{u,1}}
\def\tBo{\tB^{u,1}}
\def\hBo{\hB^{u,1}}
\def\hpi{\hat{\pi}}
\def\tU{\widetilde{U}}
\def\tr{\tilde{r}}
\def\chBo{\check{B}^{u,1}}
\def\tBo{\tB^{u,1}}
\def\hBo{\hB^{u,1}}
\def\hpi{\hat{\pi}}
\def\tL{\widetilde{L}}

\def\tde{\tilde{\delta}}
\def\tSt{\widetilde{S}^{(2)}}
\def\So{S^{(1)}}
\def\hr{\hat{r}}
\def\li{\mbox{\rm li}}

\def\vlo{v^{(\ell)}_{1}}
\def\vlt{v^{(\ell)}_2}
\def\vli{v^{(\ell)}_{i}}
\def\vloij{v^{(\ell,0)}_{i,j}}
\def\vlooj{v^{(\ell,0)}_{1,j}}
\def\vlotj{v^{(\ell,0)}_{2,j}}
\def\vlkij{v^{(\ell,k)}_{i,j}}
\def\vlkoj{v^{(\ell,k)}_{1,j}}
\def\vlktj{v^{(\ell,k)}_{2,j}}

\def\pij{\pi^{(j)}}
\def\wo{w^{(1)}}
\def\vo{v^{(1)}}
\def\uo{u^{(1)}}
\def\wt{w^{(2)}}
\def\vt{v^{(2)}}
\def\twt{\tw^{(2)}}
\def\tvt{\tv^{(2)}}
\def\xio{\xi^{(1)}}
\def\xit{\xi^{(2)}}
\def\txio{\txi^{(1)}}
\def\txit{\txi^{(2)}}
\def\ut{u^{(2)}}
\def\tut{\tu^{(2)}}
\def\tuo{\tu^{(1)}}
\def\tut{\tu^{(2)}}
\def\Wo{W^{(1)}}
\def\Vo{V^{(1)}}
\def\Uo{U^{(1)}}
\def\tUo{\tU^{(1)}}
\def\Wt{W^{(2)}}
\def\Vt{V^{(2)}}
\def\Ut{U^{(2)}}
\def\tWt{\tW^{(2)}}
\def\Vt{\tV^{(2)}}
\def\tUt{\tU^{(2)}}
\def\MM{{\mathcal M}}
\def\tv{\tilde{v}}
\def\Xij{X_{i,j}}

\def\chu{\check{u}}
\def\tkk{\widetilde{\kk}}
\def\hGammam{\widehat{\Gamma}^{(m)}}

\def\varphis{\varphi^{(s)}}
\def\varphit{\varphi^{(t)}}
\def\varphisp{\varphi^{(s')}}
\def\Qus{Q^{(s)}}
\def\Qalpha{Q^{(\alpha)}}
\def\Qbeta{Q^{(\beta)}}
\def\Qusp{Q^{(s')}}
\def\tvarphia{\tilde{\varphi}^{(\alpha)}}
\def\Lalpha{L^{(\alpha)}}
\def\Lbeta{L^{(\beta)}}
\def\Phialpha{\Phi^{(\alpha)}}
\def\Phibeta{\Phi^{(\beta)}}


\noindent
{\large\bf Abstract}

\medskip

\noindent
In this work we study strong spectral properties of Ruelle  transfer operators related to 
Gibbs measures for contact Anosov flows. As a consequence we establish exponential decay of
correlations for H\"older observables with respect to any Gibbs measure.
The approach invented in 1997 by Dolgopyat, and further developed in our papers in 2011 and 2023, is substantially enhanced here, 
allowing to deal with the general case of arbitrary  contact Anosov flows and arbitrary Gibbs measures. 
The results obtained here naturally apply to geodesic flows on compact Riemannian manifolds.

As is now well-known, the strong spectral estimates for Ruelle operators and a well-established technique by Dolgopyat  
lead to exponential decay of correlations for H\"older continuous potentials.
Other immediate consequences are: (a) existence of a non-zero analytic continuation of the 
Ruelle zeta function  with a pole at the entropy in a vertical strip containing  the entropy in its interior; 
(b) a Prime Orbit Theorem with an exponentially small error.

\bigskip



\section{Introduction}
\renewcommand{\theequation}{\arabic{section}.\arabic{equation}}

More than 20 years ago Liverani \cite{L} proved exponential decay of correlations for $C^4$ contact Anosov flows
for the Sinai-Bowen-Ruelle measure determined by the Riemann volume.
In this work, as a consequence of the main result, we derive exponential decay of correlations 
for $C^5$ contact Anosov flows on compact Riemannian manifolds $M$ with respect to 
any  Gibbs measure on $M$.  In \cite{St4} this was done for Gibbs measures admitting a Pesin set with exponentially
small tails\footnote{I.e. a Pesin set whose pre-images along the flow have measures decaying exponentially fast.}.
Here we succeed to establish exponential decay without any additional assumptions. As a consequence, for geodesic flows
on compact Riemannian manifolds every Gibbs measure is exponentially mixing.

As is now well-known, the study of statistical properties of continuous dynamical systems proved to be significantly more
difficult than the one for discrete systems. In particular the study of rates of correlation decay for H\"older continuous potentials
turned out to be highly non-trivial. After the extensive work of Sinai, Bowen and Ruelle in the 70's on statistical properties of Anosov
diffeomorphisms, and the discovery in the 80's by Ruelle \cite{R1} and Pollicott \cite{Po} that for Axiom A flows on basic sets the decay
of correlations for H\"older potentials could be arbitrarily slow, there was a period of more than 15 years when it appeared no major results
in this area had been established.  Significant  breakthrough was achieved at the end of the 90's.
First, it was Chernov \cite{Ch1} who proved sub-exponential decay of correlations for Anosov flows on 3D 
Riemannian manifolds  (with respect 
to the Sinai-Bowen-Ruelle measure). Then Dolgopyat \cite{D1} proved exponential decay of correlations for H\"older continuous 
potentials in two major cases: (i) geodesic flows on compact surfaces of negative curvature (with respect to any 
Gibbs measure);  (ii) transitive Anosov flows  on compact Riemannian manifolds with $C^1$ jointly non-integrable local 
stable and unstable foliations (with respect to the Sinai-Bowen-Ruelle measure). The results in \cite{D1} on decay of correlations were 
derived as a consequence of some very strong spectral estimates for Ruelle transfer operators defined by means of a Markov family for the
flow. Dolgopyat's paper \cite{D1} nowadays is regarded as fundamental, not just for the results that it established 
but also for the general framework 
created there.  The latter, or parts of it, has been used by various people to establish some very significant results: see e.g. \cite{AGY},
\cite{BaV}, \cite{PoS}, \cite{St1},  \cite{N}, \cite{OWi}, \cite{PeS1},  \cite{PeS2}, \cite{DMS},  just to name a few of these.
Although Liverani studied a different Ruelle operator in \cite{L}, at some stage he used, as he said, "Dolgopyat's cancellation mechanism".
In our work \cite{St2} we developed a modification of Dolgopyat's  approach to establish strong spectral properties for Ruelle transfer
operators for Axiom A flows on basic sets satisfying certain regularity conditions, and as a consequence established exponential decay of
correlations for arbitrary Gibbs measures for such flows. A much more sophisticated modification of 
Dolgopyat's framework was done in \cite{St4},
where we proved  strong spectral properties for Ruelle transfer operators for $C^3$ contact Anosov flows with respect to any Gibbs measure
admitting a Pesin set with exponentially small tails, and used this to establish exponential decay of correlations for such measures.
Recently Tsujii and Zhang \cite{TZ} proved exponential decay of correlations for 
arbitrary mixing measures for transitive Anosov flows on 3D compact manifolds, also using a modification of Dolgopyat's appoach.

Another significant phase in the study of statistical properties of continuous dynamical systems originated from the works of Young
\cite{Y1}, \cite{Y2} where she developed her so called "tower method". This was also a major event that prompted and facilitated 
significant research activities and it turned out to be very useful in the study of both uniformly and non-uniformly hyperbolic systems --
see e.g. \cite{Mel} and the references there. Various other approaches in studying decay of correlations
have been developed as well -- see e.g. \cite{D2}, \cite{D3}, \cite{MelV} and the references there. 
Recently very sophisticated tools from PDE's involving microlocal analysis have been 
used in studying various properties of hyperboilic flows -- decay of correlations, dynamical zeta functions, 
distribution of Ruelle-Pollicott resonances --
see e.g. \cite{NZ}, \cite{DyG}, \cite{DyZ}, \cite{FaSj}, \cite{FaT} and the references there, 
just to mention a few of the large number of publications in
this area. And speaking about decay of correlations for hyperbolic systems we have to mention here the major result in \cite{BaDL} 
about exponential decay of correlations for 2D Sinai billiards.

Let $\phi_t : M \longrightarrow M$ be a $C^2$ contact Anosov flow on a $C^2$ compact Riemannian manifold $M$.
Let $\varphi = \phi_1$ be the time-one map of the flow, and let $\m$ be an $\varphi$-invariant probability measure on $M$.
Given $\alpha > 0$ denote by $C^{\alpha} (M)$  {\it the space of all $\alpha$-H\"older 
complex-valued functions} on $M$, i.e. functions $h : M \longrightarrow \C$ for which there exists
$L \geq 0$ with  $|h(x) - h(y)| \leq L\, (d(x,y))^\alpha$ for all $x,y \in M$. For such $h$,  let $|h|_\alpha$ be the smallest possible
choice for $L$. Set $\|h\|_0 = \sup_{x\in M} |h(x)|$, and $\|h\|_\alpha = \|h\|_0 + |h|_\alpha$.

The main result in this paper is the following.

\bs

\noindent
{\bf Theorem 1.1.} {\it Let $\phi_t : M \longrightarrow M$ be a $C^5$ contact Anosov flow,
let $F_0$ be a H\"older continuous function on $M$ and let $\m$ be the Gibbs measure
determined by $F_0$ on $M$. 
Then for every $\alpha > 0$ there exist constants $C = C(\alpha) > 0$ and $c = c(\alpha) > 0$ such that 
$$\left| \int_{M} A(x) B(\phi_t(x)) d\m (x) -  \left( \int_{M} A(x)  d\m (x)\right)\left(\int_{M} B(x)  d\m (x)\right)\right|
\leq C e^{-ct} \|A\|_\alpha  \|B\|_\alpha \;$$
for any two functions $A, B\in C^\alpha(M)$.}

\bs

We obtain this as a consequence of Theorem 1.2 below  and the procedure described in \cite{D1}.
The assumption that the flow is $C^5$ is made so that one can apply the procedure in \cite{D1}.
In particular this is essential when estimating the Laplace transform of the correlation function
$$\rho(t) = \int_{M} A(x) B(\phi_t(x))\; d\m (x) - 
\left( \int_{M} A(x)\; d\m (x)\right)\left(\int_{M} B(x) \; d\m (x)\right) $$
(see part VI in Sect. 4 in \cite{D1}).

To our knowledge there are several results known so far on
exponential decay of correlations for  general Gibbs potentials: that 
of Dolgopyat \cite{D1} for  geodesic flows on compact surfaces, the one in \cite{St2} for Axiom A flows on basic sets
 (under some additional  assumptions); the recent result of Tsujii and Zhang for Anosov flows on 3D manifolds, and
 the one in \cite{St4} for contact Anosov flows on arbitrary compact manifolds but only for Gibbs measures admitting a Pesin
 set with exponentially small tails.
 
Let $\rr = \{R_i\}_{i=1}^{k_0}$ be a (pseudo-) Markov partition for $\phi_t$  consisting of 
rectangles $R_i = [U_i ,S_i ]$, where $U_i$ (resp. $S_i$) are (admissible) subsets of  $W^u_{\ep}(z_i)$
(resp. $W^s_{\ep}(z_i)$) for some $\ep > 0$ and $z_i\in M$ (cf. Sect. 2 for details). 
The first return time function 
$$\tau : R = \cup_{i=1}^{k_0} R_i  \longrightarrow [0,\infty)$$
is essentially $\alpha_1$-H\"older continuous on $R$  for some $\alpha_1 > 0$,
i.e. there exists a constant $L > 0$ such that if $x,y \in R_i \cap \pp^{-1}(R_j)$  for some $i,j$, 
where $\pp: R \longrightarrow R$ is the 
standard Poincar\'e map,  then $|\tau(x) - \tau(y)| \leq L\, (d(x,y))^{\alpha_1}$.
The {\it shift map} 
$$\sigma : U = \cup_{i=1}^{k_0} U_i \longrightarrow U$$
is defined by $\sigma = \piU\circ \pp$, where 
$\piU : R \longrightarrow U$ is the projection along the leaves of local stable manifolds.
Let $\hU$ be the set of all $x \in U$ whose orbits do not have common points with the boundary of $R$.
Given $\theta \in (0,1)$, as in \cite{St4}, define the {\it metric $\dte$} on $\hU$ 
by $\dte(x,y) = 0$ if $x = y$, $\dte(x,y) = 1$ if $x,y$ belong to different $U_i$'s and $\dte(x,y) = \theta^N$ if 
$\pp^j(x)$ and $\pp^j(y)$ belong to the same rectangle $R_{i_j}$ for all $j = 0,1, \ldots,N-1$, and 
$N$ is the largest integer  with this property. 

Denote by $\ff_\theta(\hU)$ {\it the space of all bounded functions $h : \hU \longrightarrow \C$ 
with Lipschitz constants }
$$|h|_\theta = \sup \left\{ \frac{|h(x) - h(y)|}{\dte(x,y)} : x\neq y; \ ; x,y \in \hU\right\} < \infty .$$
Apart from the standard norm $\|h\|_\theta = \|h\|_0 + |h|_\theta$ on $\ff_\theta (\hU)$, where $\|h\|_0 = \sup_{x\in \hU} |h(x)|$,
as in \cite{D1}, for $b\in \R$, $|b| > 0$, we will also use the special norm $\|.\|_{\theta,b}$  defined by  
$$\di \| h\|_{\theta,b} = \|h\|_0 + \frac{|h|_{\theta}}{|b|} .$$

Consider a real-valued function $f \in \ff_{\theta_0} (\hU)$ for some small $\theta_0 > 0$. Then $f \in \ff_\theta(\hU)$
for all $\theta \in [\theta_0 , 1)$. Set $g = g_f = f - P_f\tau$, where   $P_f\in \R$ 
is the unique  number such that the topological pressure $\Pr_\sigma(g)$ of $g$ with respect to $\sigma$  is zero (cf. \cite{PP}). 


We say  that  {\it Ruelle transfer operators related to $f$ are eventually contracting }
on $\ff_\theta(\hU)$ if there exist constants $0 < \rho < 1$, $a_0 > 0$, 
$b_0 \geq 1$, $T \geq 1$ and  $C > 0$ such that if $a,b\in \R$  satisfy $|a| \leq a_0$ and $|b| \geq b_0$,  then 
$$\|L_{f -(P_f+a+ \i b)\tau}^m h \|_{\theta,b} \leq C \;\rho^m \; \| h\|_{\theta,b}$$
for any integer $m \geq T \log |b|$ and any  $h\in \ff_\theta (\hU)$.

This condition implies that the spectral radius  of $L_{f-(P_f+ a+\i b)\tau}$ on 
$\ff_\theta (\hU)$  does not exceed  $\rho$. It is also easy to see that it implies the following\footnote{Which is the way we defined 
eventual contraction of Ruelle transfer operators in \cite{St2}, and it agrees with the way 
the main result in \cite{D1} is stated.}: for every $\epsilon > 0$ there exist  constants $0 < \rho < 1$, $a_0 > 0$, 
$b_0 \geq 1$ and  $C > 0$ such that if $a,b\in \R$  satisfy $|a| \leq a_0$ and $|b| \geq b_0$,  then
$$\|L_{f -(P_f+a+ \i b)\tau}^m h \|_{\theta,b} \leq C \;\rho^m \; |b|^{\ep}\; \| h\|_{\theta,b}$$
for any integer $m \geq 0$ and any  $h\in \ff_\theta (\hU)$.

The central result in this paper is the following.

\bs

\noindent
{\bf Theorem 1.2.} {\it  Let $\phi_t : M \longrightarrow M$ be a $C^2$ contact Anosov flow 
on a $C^2$ compact  Riemannian manifold $M$, let  $\rr = \{R_i\}_{i=1}^{k_0}$ be a (pseudo-) Markov 
partition for $\phi_t$ as above and let $\sigma : U \longrightarrow U$ be the corresponding shift map. 
There exists constants  $0 <  \theta_1 < \theta < 1$ such that for any  real-valued function $f\in \ff_{\theta_1}(\hU)$
the Ruelle transfer operators related to $f$ are eventually contracting on $\ff_{\theta}(\hU)$}.

\bs

\noindent
{\bf Remark.} The constant $\theta$ is defined by using some regularity properties of the flow
-- see Sect. 6.1 below for details.
A slightly more general result can be established. Namely, one can prove the same result as above for every
$f \in \ff_{\theta}$. This can be done repeating the lengthy approximation procedure in Sect. 8 in \cite{St4}. 
However the main consequences -- Corollary 1.3 below, as well as Theorems 1.1 and 1.3 -- will be 
the same\footnote{Except for a better constant $\alpha_0 > 0$ in Corollary 1.3.}.


\bs


A similar result for H\"older continuous functions (with respect to the Riemannian metric) looks a 
bit more complicated, since in general Ruelle transfer operators do not preserve any of the spaces 
$C^\alpha(\hU)$. However, they preserve a certain `filtration' $\cup_{0 < \alpha \leq \alpha_0} \clip (\hU)$.
For $b\in \R$, $b \neq 0$, define the norm $\|.\|_{\alpha,b}$ on $\clip (\hU)$ by  $\| h\|_{\alpha,b} = \|h\|_0 + \frac{|h|_{\alpha}}{|b|}$.

\bs

\noindent
{\bf Corollary 1.3.} {\it  Under the assumptions of Theorem 1.2,  there exists a constant $\alpha_0 \in (0,1]$ 
such that for any real-valued function  $f \in C^{\alpha_0}(\hU)$ the  Ruelle transfer operators related to  $f$  are 
eventually contracting on $\cup_{0 < \alpha \leq \alpha_0} \clip (\hU)$. More precisely, there 
exists a constant $\hbeta\in (0,1]$
and for each $\epsilon > 0$ there exist constants $0 < \rho < 1$, $a_0 > 0$, $b_0 \geq 1$,
$C > 0$ and $M > 0$ such that if $a,b\in \R$  
satisfy $|a| \leq a_0$ and $|b| \geq b_0$, then for every integer $m  \geq M\, \log |b|$ and every 
$\alpha \in (0,\alpha_0]$ the operator
$\di L_{f -(P_f+a+ \i b)\tau}^m : C^\alpha (\hU) \longrightarrow C^{\alpha \hbeta}  (\hU)$
is well-defined and
$$\|L_{f -(P_f+a+ \i b)\tau}^m h \|_{\alpha\hbeta,b} \leq C \;\rho^m \; \| h\|_{\alpha,b}$$ 
for every  $h\in C^\alpha  (\hU)$.}

\ms

The maximal constant $\alpha_0 \in (0,1]$ that one can choose above 
is related to the regularity of the local stable/unstable foliations. Estimates for this constant 
can be derived from certain bunching 
conditions concerning the rates of expansion/contraction of the  flow along local unstable/stable 
manifolds (see \cite{Ha1},  \cite{Ha2}, \cite{PSW}).

The above  was first proved by Dolgopyat (\cite{D1}) in the case of geodesic flows on compact 
surfaces of negative curvature with $\alpha_0 = 1$ (then one can choose $\hbeta = 1$ as well).
The second main result in \cite{D1} concerns transitive Anosov flows  on compact Riemannian
manifolds with $C^1$ jointly  non-integrable local stable and unstable foliations. For such 
flows Dolgopyat proved that the conclusion of Corollary 1.3 with 
$\alpha_0 = 1$ holds for the Sinai-Bowen-Ruelle  potential  $F_0 = \log \det (d\phi_\tau)_{|E^u}$. 
More general results were proved in \cite{St2} for mixing Axiom A flows on basic sets
(again for $\alpha_0 = 1$) under some additional regularity assumptions\footnote{As mentioned earlier,
the results apply e.g. to  $C^2$ mixing Axiom A flows on basic sets satisfying a certain pinching condition (similar to the
$1/4$-pinching condition for geodesic flows on manifolds of negative curvature).}, and more recently
in \cite{St4} for Gibbs measures for contact Anosov flows admitting a Pesin set with exponentially small tails.

Next, consider the {\it Ruelle zeta function} 
$$\zeta(s) = \prod_{\gamma} (1- e^{-s\ell(\gamma)})^{-1} \quad, \quad s\in \C ,$$
where $\gamma$ runs over the set of primitive  closed orbits of $\phi_t: M \longrightarrow M$
and $\ell(\gamma)$ is the least period of $\gamma$.  Denote by $h_T$ the  {\it topological entropy} 
of $\phi_t$ on $M$.

Using Theorem 1.2 above  and an argument of Pollicott and Sharp \cite{PoS}, one derives the 
following\footnote{As remarked in \cite{St4}, instead of using the norm $\|\cdot \|_{1,b}$ as in \cite{PoS}, in the present 
case one has to work with $\|\cdot \|_{\theta,b}$ for some  $\theta \in (0,1)$, 
and then one has to use the so called Ruelle's Lemma in the form proved in  \cite{W}. 
This is enough to prove the estimate (2.3) 
for $\zeta(s)$ in \cite{PoS}, and from there the arguments are the same.}.

\bs

\noindent
{\bf Theorem 1.4.} {\it Let $\phi_t : M \longrightarrow M$ be a $C^2$ contact Anosov flow on a $C^2$ compact  
Riemannian manifold $M$.  Then:} 

\ms

(a) {\it The Ruele zeta function $\zeta(s)$ of the flow  $\phi_t: M \longrightarrow M$  has an analytic  and non-vanishing 
continuation in a half-plane $\Re(s) > c$ for some $c < h_T$ except for a  simple pole at $s = h_T$.  }

\ms

(b) {\it There exists $c \in (0, h_T)$ such that
$$\di \pi(\lambda) = \# \{ \gamma : \ell(\gamma) \leq \lambda\} = \li(e^{h_T \lambda}) + O(e^{c\lambda})$$
as $\lambda\to \infty$, where 
$\di \li(x) = \int_2^x \frac{du}{\log u} \sim \frac{x}{\log x}$
as  $x \to \infty$.}

\bigskip

Parts (a) and (b) were first established by Pollicott and Sharp \cite{PoS} for geodesic flows on 
compact surfaces of negative curvature (using \cite{D1}), and then similar results were proved in 
\cite{St2} for mixing Axiom A flows on basic sets satisfying certain additional assumptions 
(as mentioned above). More recently, using different methods, it was proved in \cite{GLP} that: 
(i) for volume preserving three dimensional Anosov flows (a) holds, and moreover, in the case of 
$C^\infty$ flows,  the Ruelle zeta function $\zeta(s)$ is meromorphic in $\C$ and  $\zeta(s) \neq 0$ for $\Re(s) > 0$; 
(ii) (b) holds for geodesic flows on $\frac{1}{9}$-pinched compact Riemannian manifolds of negative 
curvature. These were obtained as consequences of more general results in \cite{GLP}.

In Sect. 2 we collect some preliminary information used later on. Sect. 3 contains basic facts from Pesin's theory
of Lyapunov exponents, some notation and three lemmas from \cite{St3} and \cite{St4} that we need in subsequent sections.
Sect. 4 deals with some technical results concerning cylinders defined by means of Markov families for general Anosov flows.
The most important fact in Sect. 4 is 
Lemma 4.4 where for every unstable cylinder $\cc$ intersecting the specially chosen in Sect. 4.2 Pesin set $P_0$ we construct two families of 
sub-cylinders, each with a significant total measure\footnote{All considerations involve a given Gibbs measure on $R$
and the related Gibbs measure $\nu$ on $U$.}, and a certain
pairing between elements of the two families that is used later in Sects. 5 and 6 to develop a cancelations procedure for the so called
contraction operators $\nn_J$. The construction in Lemma 4.4 is enhanced in Lemma 5.5, where we use significantly the contact form,
and establish a certain strong non-integrability property of the flow.  

In Sect. 6 we  define for a given parameter $b \in \R$, a family of cylinders $\cc_m$ with
lengths of size $\Con \log |b|$ covering the Pesin set $P_0$, and then use the constructions in Sects. 4 and 5 to define
nice families of sub-cylinders of these with appropriate pairings between their elements. The so called contraction operators are
then defined in a rather different way to what was done in previous papers. 
As in \cite{St4}, the contraction operators $\nn_J$ only contract in the vicinity of the Pesin set $P_0$ -- this is established in Sect. 7
using the families of sub-cylinders of the cylinders $\cc_m$ from Lemmas 4.4 and 5.5 
and the particular pairing between these provided by the lemmas.

In Sect. 8 we succeed to obtain global contraction properties of the contraction operators. 
Then using an appropriate modification of arguments from \cite{D1} and \cite{St4} we prove the main results in the paper.

\def\tp{\tilde{p}}

\newpage

\section{Preliminaries}
\setcounter{equation}{0}

Let $M$ be a $C^2$ compact Riemannian manifold,  and  let $\phi_t : M \longrightarrow M$ 
 ($t\in \R$) be a  $C^2$  transitive {\it Anosov flow}  on $M$. 
This means that there exist constants $C > 0$ and $0 < \lambda < 1$ such that there exists a $d\phi_t$-invariant decomposition  
$T_xM = E^0(x) \oplus E^u(x) \oplus E^s(x)$ of $T_xM$ ($x \in M$) 
into a direct  sum of non-zero linear subspaces, where $E^0(x)$ is the one-dimensional subspace determined by the direction of the flow at $x$, 
$\| d\phi_t(u)\| \leq C\, \lambda^t\, \|u\|$ for all  $u\in E^s(x)$ and $t\geq 0$,  and
$\| d\phi_t(u)\| \leq C\, \lambda^{-t}\, \|u\|$ for all $u\in E^u(x)$ and $ t\leq 0$.
Throughout we denote by $\|\cdot \|$ the {\it norm determined by the Riemannian metric on $M$.} Transitivity means that the flow has a dense orbit in $M$.

Given $x\in M$ and a sufficiently small $\epsilon > 0$ the (strong) {\it stable} and {\it unstable manifolds} $\wloc^s(x)$ and $\wloc^u(x)$ of size $\ep$
are defined in the usual way. The corresponding tangent bundles are then
$E^u(x) = T_x \wloc^u(x)$ and $E^s(x) = T_x \wloc^s(x)$.  Set 
$E^u(x;\delta) = \{ u\in E^u(x) : \|u\| \leq \delta\}$ for any $\delta > 0$ and define
$E^s(x;\delta)$ similarly. Let
$\exp^u_x : E^u(x;\ep) \longrightarrow \wloc^u(x)$ and $\exp^s_x : E^s(x;\ep) \longrightarrow \wloc^s(x)$
be the corresponding {\it exponential maps}.

The so called {\it temporal distance function} $\Delta(x,y)$ is defined as follows. Given a sufficiently small $\ep_0 > 0$,
it follows from the hyperbolicity of the flow on $M$ that there exists $\ep_1 \in (0,\ep_0)$ such that if $x,y\in M$ and $d (x,y) < \ep_1$, 
then there exists a unique point
$[x,y] = W^s_{\ep}(x) \cap \phi_{[-\ep_0,\ep_0]} (W^u_{\ep_0}(y)) $
(cf. \cite{KH}). Hence $\phi_t([x,y]) \in W^u_{\ep_0}(y)$ for some $t\in [-\ep_0, \ep_0]$ 
and we set\footnote{A different definition for $\Delta$ is given in \cite{D1} and \cite{L}, however
in the only case considered in this paper when $x\in W^u_\ep(z)$ and 
$y \in W^s_\ep(z)$ for some $z \in M$,  these definitions coincide with  the present one.}
 $\Delta(x,y) = t$ (\cite{KB},\cite{D1}, \cite{L}). For $x, y\in M$ with $d (x,y) < \ep_1$, define $\pi_y(x) = [x,y]$.
In this way for every $y \in M$, on a small open neighbourhood $W$ of $y$ in $M$ we get 
a {\it projection} $\pi_y : W \longrightarrow \phi_{[-\ep_0,\ep_0]} (W^u_{\ep_0}(y))$  along local stable manifolds.
The map
$\pi_y: \phi_{[-\ep_1,\ep_1]} (W^u_{\ep_1}(x)) \longrightarrow \phi_{[-\ep_0,\ep_0]} (W^u_{\ep_0}(y))$
is called a {\it local stable holonomy map}. In a similar way one defines
holonomy maps between any two sufficiently close local transversals to stable laminations (see e.g. \cite{PSW}).
Combining such a map with a shift along the flow we get another  local stable holonomy  map 
$\hh_x^y : W^u_{\ep_1}(x)  \longrightarrow W^u_{\ep_0}(y) .$
In a similar way one defines local holonomy maps along unstable laminations.

For convenience of the reader we will now provide the definition of a Markov family for the flow (see e.g. Sect. 4 in \cite{Ch2}
for details). 
Given a submanifold $D$ of $M$ of codimension one with $\diam(D) \leq \ep_0$ which is transversal to the flow, the projection 
$\pr_D : \phi_{[-\ep,\ep]}(D) \longrightarrow D$ along the flow is well-defined and smooth. For $x,y\in D$, set $\la x, y\ra_D = \pr_D([x,y])$. 
A subset $\tR$ of $D$ is called a {\it rectangle} if $\la x, y\ra_D \in \tR$ for all $x,y\in \tR$.  A {\it proper rectangle} is a rectangle $\tR$ 
that coincides with the closure of its interior in $D$. The stable and unstable leaves through $x\in \tR$ are defined by 
$W^s_{\tR}(x) = \phi_{[-\ep,\ep]}(W^s_\ep(x)) \cap \tR$ and $W^u_{\tR}(x) = \phi_{[-\ep,\ep]}(W^s_\ep(x)) \cap \tR$.
By $\Int_D(A)$ we denote the {\it interior} of a subset $A$ of $D$ in $D$. 

Let $\tRR = \{ \tR_i\}_{i=1}^{k_0}$  be a family of proper rectangles, where each $\tR_i$ is contained
in a submanifold $D_i$ of $M$ of codimension one and has the form 
$\tR_i = \la U_i  , S_i \ra_{D_i} = \{ \la x,y\ra_{D_i} : x\in U_i, y\in S_i\}$,
where $U_i \subset \wloc^u(z_i)$ and $S_i \subset \wloc^s(z_i)$, respectively, 
for some $z_i\in M$.  
Set 
$\di \tR =  \cup_{i=1}^{k_0} \tR_i .$
The family $\tRR$ is called {\it complete} if  there exists a constant $\chi > 0$ such that for every $x \in M$, $\phi_{t}(x) \in \tR$ for some  
$t \in (0,\chi]$.  The {\it Poincar\'e map} $\tpp: \tR \longrightarrow \tR$
related to a complete family $\tRR$ is defined by $\tpp(x) = \phi_{\ttau(x)}(x) \in \tR$, where
$\ttau(x) > 0$ is the smallest positive time with $\phi_{\ttau(x)}(x) \in \tR$. The function $\ttau$ 
is called the {\it first return time}  associated with $\tRR$. 
A complete family $\tRR = \{ \tR_i\}_{i=1}^{k_0}$ of rectangles in $M$ is called a  {\it Markov family} 
of size $\chi > 0$ for the  flow $\phi_t$ if: 
                                         
(a) $\diam(\tR_i) < \chi$ for all $i$; 

(b)  for any $i\neq j$ and any $x\in \Int_{D}(\tR_i) \cap \tpp^{-1}(\Int_{D}(\tR_j))$ we have   
$$W_{\tR_i}^s(x) \subset \overline{ \tpp^{-1}(W_{\tR_j}^s(\tpp(x)))} \quad , \quad                                                         
\overline{\tpp(W_{\tR_i}^u(x))} \supset W_{\tR_j}^u(\tpp(x)) ;$$

(c) for any $i\neq j$ at least one of the sets $\tR_i \cap \phi_{[0,\chi]}(\tR_j)$ and $\tR_j \cap \phi_{[0,\chi]}(\tR_i)$ is empty. 

The existence of a Markov family $\tRR$ of an arbitrarily small size $\chi > 0$ for $\phi_t$ follows from the construction of Bowen \cite{B1}.

As in  \cite{R1} and \cite{D1}, we can slightly change the Markov family $\tRR$ to a  {\it pseudo-Markov family}  
$\rr = \{ R_i\}_{i=1}^{k_0}$  of {\it pseudo-rectangles} 
$R_i = [U_i  , S_i ] =  \{ [x,y] : x\in U_i, y\in S_i\}$,
where $U_i$ and $S_i$ are as above. Set 
$\di R =  \cup_{i=1}^{k_0} R_i .$
Notice that  $\pr_{D_i} (R_i) = \tR_i$ for all $i$. For any $\xi = [x,y] \in R_i$ set
$$W^u_{R_i}(\xi) = [U_i ,y] = \{ [x',y] : x'\in U_i\} \quad , \quad 
W^s_{R_i}(\xi) = [x,S_i] = \{[x,y'] : y'\in S_i\} \subset W^s_{\ep_0}(x) .$$
The corresponding {\it Poincar\'e map} $\pp: R \longrightarrow R$ is defined by  
$\pp(x) = \phi_{\tau(x)}(x) \in R$, where $\tau(x) > 0$ 
is the smallest positive time with $\phi_{\tau(x)}(x) \in R$.  
The {\it interior} $\Int(R_i)$ of a rectangle $R_i$ is defined by  $\pr_D(\Int(R_i)) = \Int_D(\tR_i)$. In a similar way one defines
$\Intu(A)$ for a subset $A$ of some  $W^u_{R_i}(x)$ and $\Ints(A)$ for a subset of  $W^s_{R_i}(x)$.
The family $\rr = \{ R_i\}_{i=1}^{k_0}$ has the same properties as $\tRR$: 

$(a')$ $\diam(R_i) < \chi$ for all $i$; 

$(b')$ for any $i\neq j$ and any  $x\in \Int(R_i) \cap \pp^{-1}(\Int(R_j))$ we have  
$$\pp(\Int (W_{R_i}^s(x)) ) \subset \Ints (W_{R_j}^s(\pp(x))) \quad , \quad   \pp(\Int(W_{R_i}^u(x))) \supset \Int(W_{R_j}^u(\pp(x))) ;$$

$(c')$ for any $i\neq j$ at least one of the sets $R_i \cap \phi_{[0,\chi]}(R_j)$ and  $R_j \cap \phi_{[0,\chi]}(R_i)$ is empty. 

Define the matrix $\aa = (\aa_{ij})_{i,j=1}^{k_0}$  by $\aa_{ij} = 1$ if $\pp(\Int (R_i)) \cap \Int(R_j) \neq  \e$ 
and $\aa_{ij} = 0$ otherwise.  According to Sect. 2 in  \cite{BR}, we may assume that $\rr$ 
is chosen in such a way that $\aa^{\tp_0} > 0$ (all entries of the $\tp_0$-fold product of $\aa$  by itself are 
positive) for some integer $\tp_0 > 0$.  In what follows we assume that the matrix $\aa$ has this property.

One should remark here that while in general $\pp$ and $\tau$ are only (essentially) H\"older continuous,
the map $\tpp$ is (essentially) Lipschitz; see (2.1) below.

From now on we will assume that $\tRR = \{ \tR_i\}_{i=1}^{k_0}$ is a {\bf fixed Markov family} for  the flow
$\phi_t$ of size $\chi < \ep_1/2 < \ep_0/2 < 1$ and that  $\rr = \{ R_i\}_{i=1}^{k_0}$ is the related pseudo-Markov family. Set  
${\bf U = \cup_{i=1}^{k_0} U_i }$
and $\di \Intu (U) = \cup_{j=1}^{k_0} \Intu(U_j)$. It follows from the hyperbolicity of the flow that there exist
constants $d_0 \in (0,1]$ and $\gamma_1 > \gamma > 1$ such that
\be
d_0 \gamma^m\; d (x,y) \leq  d (\tpp^m(x)), \tpp^m(y)) \leq \frac{\gamma_1^m}{d_0} d (x,y)
\ee
for all $x,y\in \tR$ such that $\tpp^j(x), \tpp^j(y)$ belong to the same $\tR_{i_j}$ for all $j = 0,1, \ldots, m$. 
Fix constants $d_0$ and $\gamma_1 > \gamma > 1$ with these properties; these will be used throughout the whole paper.

Using the {\it projection} $\piU : R \longrightarrow U$ along stable leaves we define the
{\it shift map} $\sigma : U   \longrightarrow U$ by $\sigma  = \piU \circ \pp$.
We will use the same notation for the related projection $\piU : \tR \longrightarrow U$ along stable leaves which defines
another {\it shift map} $\tsigma : U   \longrightarrow U$ such that $\tsigma  = \piU \circ \tpp$. The two shift maps
$\sigma$ and $\tsigma$ are naturally related. Notice that our definitions of $R$ and $\tR$ are so that both contain $U$
and both projections $\piU : R \longrightarrow U$ and $\piU : \tR \longrightarrow U$ along stable leaves take values in $U$.

Notice that  $\tau$ is constant on each stable leaf $W_{R_i}^s(x) = W^s_{\ep_0}(x) \cap R_i$ (however not
on the corresponding stable leaves $W_{\tR_i}^s(x) = W^s_{\ep_0}(x) \cap \tR_i$).
The shift map $\sigma$ is naturally conjugate to the Bernoulli shift map $\sigma_{\aa} : \Sigma_{\aa} \longrightarrow \Sigma_{\aa}$
on the  symbol space 
$$\di \Sigma_{\aa} = \{  (i_j)_{j=-\infty}^\infty : 1\leq i_j \leq k_0 , \aa_{i_j\; i_{j+1}} = 1 \:\: \mbox{ \rm for all } \: j\; \},$$ 
given by $\sigma_{\aa} ( (i_j)) = ( (i'_j))$, where $i'_j = i_{j+1}$ for all $j$.
There exists a natural surjection $\pi : \sa \longrightarrow \tR$ 
such that $\pi \circ \sigma_A = \pp\circ \pi$ on a residual subset of $\tR$ (see e.g. \cite{B1} or Sect. 4 in \cite{Ch2}).
Denoting by $\tR^*$ be the set of those $x\in R$ such that $\phi_t(x) \notin \Int(R_j)$ for any $t \in \R$ and any $j$,
we have $\pi \circ \sigma_A = \pp\circ \pi$ on $\pi^{-1}(\tR^*)$.
Moreover $\pi$ is Lipschitz on $\pi^{-1}(\tR^*)$ if the latter is considered with  the {\it metric} $d_\theta$ for some  
$\theta\in (0,1)$, defined by  $d_\theta(\xi,\eta) = 0$ if $\xi = \eta$ and $d_\theta(\xi,\eta) = \theta^m$ if
$\xi_i = \eta_i$ for $|i| \leq m$ and $m$ is maximal with this property. Notice that
$\htau = \tau \circ \pi$ defines a Lipschitz function on $\pi^{-1}(\tR^*)$,
so it has a Lipschitz extension $\htau  : \sa \longrightarrow \R_+$ (\cite{B1}, \cite{Ch2}).
The  {\it space of Lipschitz functions} on $\Sigma_{\aa}$ with respect to the
metric $d_\theta$  will be denoted by $C_\theta(\Sigma_{\aa})$ and will be considered with the norm
$\|h\|_\theta = \|h\|_0 + |h|_\theta$.

The {\it shift map} $\sigma_{\aa} : \saa \longrightarrow \saa$ on the one-sided subshift of finite type
$$\saa = \{  (i_j)_{j=0}^\infty : 1\leq i_j \leq k, \aa_{i_j\; i_{j+1}} = 1 \:\: \mbox{ \rm for all } \: j \geq 0\; \},$$
is defined similarly.
Notice that $\htau(\xi) = \tau(\pi(\xi))$ depends only on the
forward coordinates of $\xi\in \Sigma_{\aa}$. 
In particular we can consider $\htau$ as a function on
$\Sigma_A^+$ such that $\htau =  \tau\circ \pi$ on a residual subset of $\saa$. 
The metric $d_\theta$ on $\saa$ and the  space of Lipschitz functions $C_\theta(\Sigma_{\aa}^+)$
are defined as for $\sa$.
If $\hpi : \Sigma_{\aa} \longrightarrow \Sigma_{\aa}^+$ is the {\it natural projection}, one shows easily that 
there exists a continuous surjection $ \pi^+ : \saa \longrightarrow U$ such that
then $\pi^+\circ \hpi = \piU\circ \pi$. Moreover, $\sigma \circ \pi^+ = \pi^+ \circ \sigma^+_{\aa}$.

As stated in Sect. 1, we will denote by $\widehat{U}$ the set  of those $x\in U$ 
such that  $\pp^m(x) \in \Int(R) = \cup_{i=1}^k \Int(R_i)$  for all $m \in \Z$. This is a residual subset 
of $U$ and has full measure with respect to any Gibbs measure on $U$ (see e.g.  \cite{B2}).
Set $\hU_i = U_i \cap \hU$.

Let $B(U)$ be the {\it space of  bounded functions} $g : U \longrightarrow \C$ with its standard norm  $\|g\|_0 = \sup_{x\in U} |g(x)|$.
Given a function $g \in B(\hU)$, the  {\it Ruelle transfer operator }  $L_g : B(U) \longrightarrow B(U)$  is defined by 
$$\di (L_gh)(u) = \sum_{\sigma(v) = u} e^{g(v)} h(v) .$$
Via the natural projection  $ \pi^+ : \saa \longrightarrow U$, the above corresponds to the well-known definition of a Ruelle transfer
operator on $\saa$ (see \cite{Ba}, \cite{B2} or \cite{PP}).

{\bf Fix constants  $0 < \ttau_0 < \tau_0 < \chi$ so that }
$\ttau_0 \leq \tau(x) \leq \tau_0$ for all $x\in R$.
We will assume that $\ttau$ satisfies the same estimates, namely $\ttau_0 \leq \ttau(x) \leq \tau_0$ for all $x \in \tR$.
As we remarked earlier, in general $\ttau$ is not constant on the stable leaves $W_{\tR_i}^s(x) = W^s_{\ep_0}(x) \cap \tR_i$.
However, since $\tpp$ is Lipschitz, we can assume the constant $d_0 \in (0,1]$ is chosen so that
$|\ttau(\tx) - \ttau(\ty)| \leq \frac{1}{d_0} d(\tx,\ty)$
for all $\tx, \ty\in \tR_i$ and all $i = 1,\ldots, k$. 

Throughout this paper $\alpha_1 \in (0,1]$ denotes a {\it fixed constant } such that  $\tau\in C^{\alpha_1}(\hU)$ and the 
local stable/unstable holonomy maps are uniformly $\alpha_1$-H\"older. Choosing the constant $d_0 \in (0,1]$ from (2.1)
sufficiently small, we will assume that for any $x, y\in R$ the  local stable holonomy  map
$\hh_x^y : W^u_{\ep_1}(x)  \longrightarrow W^u_{\ep_0}(y) $
is so that 
$$\di d(\hh_x^y(x'), \hh_x^y(x'')) \leq \frac{1}{d_0} (d (x', x''))^{\alpha_1}$$
for all $x', x'' \in W^u_{\ep_1}(x)$. 
We will also assume that a similar condition is satisfied by the corresponding local unstable holonomy maps.
Furthermore, we will assume that $\alpha_1$ is chosen so that
the shift $\tPsi : R \longrightarrow \tR$ along the flow is $\alpha_1$-H\"older and satisfies a condition similar to the one about the
local holonomy maps.


\def\chvt{\check{v}^{(2)}}
\def\chVt{\Check{V}^{(2)}}
\def\chUt{\Check{U}^{(2)}}
 
\def\hq{\hat{q}}
\def\txio{\tilde{\xi}^{(1)}}
\def\tP{\widetilde{P}}
\def\tQ{\widetilde{Q}}
\def\hcc{\widehat{\cc}}
\def\hdd{\widehat{\dd}}

\section{Lyapunov exponents and Lyapunov regularity functions}
\setcounter{equation}{0}

Throughout this paper $M$ denotes a $C^2$ compact Riemannian manifold, $\phi_t$ is a $C^2$ Anosov flow
on $M$ and  $\varphi = \phi_1$.  Let $F_0$ be a H\"older continuous real-valued function on $M$ and let 
{\it $\m$  be the Gibbs measure  generated by $F_0$ on $M$.} 
The  {\it Oseledets Multiplicative Ergodic Theorem} (\cite{Os}) implies
that in the situation considered here  there exists a $\phi_t$-invariant subset $\ll$ of $M$ with $\m(\ll) = 1$ 
such that for every  $x \in \ll$  there exist numbers
$0 < \chi_1  < \chi_2 < \ldots < \chi_{\tk}$
and a $d\phi_t$-invariant decomposition  
$\di E^u(x) = E^u_1(x) \oplus E^u_2(x) \oplus \ldots \oplus E^u_{\tk}(x)$
of $E^u(x)$ into subspaces of constant dimensions such that
\begin{equation}
\lim_{t\to \infty} \frac{1}{t}\, \log \|d \phi_t(x)\cdot v\| = \chi_i \quad , \quad v \in E^u_i(x)\setminus \{0\} ,
\end{equation}
for all $i = 1, \ldots,\tk$ (see \cite{BP}, \cite{Ar}, \cite{V} or \cite{KH}).  
Moreover $\ll$ is so that $\mu(\ll) = 1$ for every $\phi_t$-invariant Borel probability measure $\mu$ on $M$.

Here and in what follows we denote by $\|\cdot \|$ the {\it norm on the tangent spaces
$T_xM$ ($x\in M$) induced by the Riemannian metric} on $M$.

The numbers $\chi_i > 0$ are called the the {\it positive Lyapunov exponents} of $\phi_t$.
In our case, the dimension $n_i$ of $E^u_i(x)$ is constant on $\ll$, and  clearly 
$n_1 + n_2 + \ldots + n_{\tk} = n^u = \dim(E^u(x))$ for all $x\in \ll$.
For $E^s(x)$, $x\in \ll$, we have a similar decomposition involving the
corresponding negative Lyapunov exponents. For contact flows we have 
$n^s = \dim(E^s(x))  = n^u$ for all $x \in \ll$, and the negative Lyapunov exponents are $-\chi_i$, $i = 1, \ldots, \tk$.

It follows from (3.1) that for every $\ep > 0$, every $x\in \ll$ and every $i = 1, \ldots, \tk$ we have
$$\di \lim_{n\to\infty} \frac{\|d\varphi^n(x)_{|E^u_i(x)}\|}{e^{(\chi_i+\ep) n}} = 0 ,$$
therefore
$$R_0(x) = \max_{1\leq i \leq \tk} \: \sup_{n \geq 0} \frac{\|d\varphi^n(x)_{|E^u_i(x)}\|}{e^{(\chi_i+\ep) n}} < \infty .$$
The function $R_0(x)$ just defined is an example of a Lyapunov $\ep$-regularity function.

More generally, a Borel function $R_\ep: \ll \longrightarrow (1,\infty)$ such that
\be
\frac{1}{R_\ep (x)} \leq \frac{\|d\varphi^n(x)\cdot v\|}{e^{(\chi_i + \ep)n}\, \|v\|} 
\leq R_\ep (x) \quad , \quad x\in \ll \;, \; v\in E^u_i(x)\setminus \{0\} \;, \; n \geq 0 ,
\ee
for all $i = 1, \ldots, \tk$, and
\be
e^{-\ep} \leq \frac{R_\ep (\varphi(x))}{R_\ep (x)} \leq  e^{\ep} \quad , \quad x\in \ll ,
\ee
is called a {\it Lyapunov $\epsilon$-regularity function}.  
As in \cite{PS}, by an {\it $\ep$-slowly varying radius function}
we mean a function of the form $r_\ep(x) = 1/R_\ep(x)$, $x \in \ll$, where $R_\ep$ is a Lyapunov $\ep$-regularity function on $\ll$.
For such  $r_\ep$ and $x \in \ll$, 
the linear map $d\varphi^n (x)$ behaves on the ball $B^u(x, r_\ep(x))$  as in the case of an uniformly hyperbolic flow --
see the relations (3.11) - (3.14) below. 

In this paper by a {\it Pesin set} in $M$ we mean a compact subset $P$ of $\ll$ with $\m(P) > 0$ such that
there exist constants $\ep > 0$ and $C > 0$ and a Lyapunov
$\ep$-regularity function $R_\ep(x)$ on $\ll$ with $R_\ep(x) \leq C$ for all $x \in P$. In a similar way we define Pesin sets in $R$
(i.e. in $R\cap \ll$) with respect to the induced measure $\mu_0$ on $R$ (see Sect. 4.2).
Clearly Pesin sets always exist.

In \cite{St4} we considered Gibbs measures admitting
{\it Pesin sets with exponentially small tails}, i.e. Pesin sets whose preimages along the flow have measures decaying exponentially fast.
Existence of such Pesin sets for a variety of Gibbs measures for hyperbolic diffeomorphisms and flows was established in \cite{GSt}.
However, to our knowledge, there are no general results of this kind.

\ms

Set  $\lambda_i = e^{\chi_i}$  for all $ i = 1, \ldots,\tk$.  Fix an arbitrary constant $\beta \in (0,1]$  such that 
\be
\lambda_j^{1+\beta} < \lambda_{j+1}
\ee
for all $1 \leq j < \tk .$
Take $\hep > 0$ so small that 
$\di e^{8\hep} < \lambda_1$ and $ e^{8\hep} < \lambda_{j}/\lambda_{j-1}$ for all $j = 2, \ldots,\tk$.
Some further assumptions about $\hep$ will be made later.
Set
\be
1 < \nu_0 = \lambda_1 e^{-8\hep} <  \mu_{j} = \lambda_{j} e^{-\hep} < \lambda_{j} < \nu_{j} = \lambda_{j} e^{\hep}
\ee
for all $j = 1, \ldots,\tk$. 
{\bf Fix $\hep > 0$ with the above properties.}

We will now provide some basic definitions and set-up from Sect. 3 in \cite{St4}.
For $x \in \ll$ and $1 \leq j \leq \tk$ set 
$$\hE^u_j(x) = E^u_1(x) \oplus \ldots \oplus E^u_{j-1}(x) \quad , \quad \tE^u_{j}  = E^u_j(x) \oplus \ldots \oplus E^u_{\tk}(x) .$$
Also set $\hE^u_1(x) = \{0\}$ and $\hE^u_{\tk+1}(x) = E^u(x)$.
For any $x\in \ll$ and any $u\in E^u(x)$ we will write $u = (\uo,\ut , \ldots, \utk)$,  where $\ui \in E^u_i(x)$ for all $i$.

According to results in the theory of non-uniformly hyperbolic systems (see \cite{P}, \cite{BP}) 
for any $j = 1, \ldots,\tk$ the invariant bundle  $\{\tE^u_{j}(x)\}_{x\in \ll}$ is uniquely integrable over $\ll$, i.e. 
there exists a measurable $\varphi$-invariant family $\{ \tWuj_{\tr(x)}(x)\}_{x\in \ll}$  
of $C^2$ submanifolds $\tWuj(x) = \tWuj_{\tr(x)}(x)$ of $M$ tangent to the bundle $\tE^u_{j}$ for some 
 {\it $\hep$-slowly varying radius function} $\tr = \tr_{\hep} : \ll \longrightarrow (0,1)$.
Moreover, with $\beta \in (0,1]$ as above,  it follows from Theorem 6.6 in 
\cite{PS} and (3.4) that there exists a $\varphi$-invariant family $\{ \hWuj_{\tr(x)}(x)\}_{x\in \ll}$ of
$C^{1+\beta}$ submanifolds $\hWuj(x) = \hWuj_{\tr(x)}(x)$ 
of $M$ tangent to the bundle $\hE^u_j$ for every $j > 1$.
(However this family is not unique in general.)  Fix
an $\varphi$-invariant family $\{ \hWuj_{\tr(x)}(x)\}_{x\in \ll}$ with the latter properties for all $x\in \ll$ and  $j = 2, \ldots, \tk$.
Then there exist an {\it $\hep$-slowly varying radius function}  $r = r_{\hep}: \ll \longrightarrow (0,1)$ 
and for any $x\in \ll$ a $C^{1+\beta}$ diffeomorphism
$$\Phi_x^u: E^u (x; r(x)) \longrightarrow \Phi_x (E^u (x; r (x)) \subset W^u_{\tr (x)}(x)$$
such that 
\be
\Phi^u_x(\hE^u_j(x; r(x))) \subset \hWuj_{\tr(x)}(x)\quad , \quad
\Phi^u_x(\tE^u_{j}(x; r(x))) \subset \tWuj_{\tr (x)}(x) 
\ee
for all $x \in \ll$ and $j = 2, \ldots , \tk$. Moreover, since  the submanifolds
$ \hWuj_{r(x)}(x)$ and $\exp^u_x(\hE^u_j(x;r(x)))$ of $W^u_{\tr(x)}(x)$ are tangent at $x$ of order 
$1+\beta$ for each $j > 1$, we can choose  $\Phi^u_x$ so that the diffeomorphism
$$\Psi^u_x =  (\exp^u_x)^{-1} \circ \Phi^u_x   : E^u(x: r(x)) \longrightarrow \Psi^u_x (E^u(x: r(x))) \subset E^u(x; \tr(x))$$ 
is $C^{1+\beta}$-close to identity. So, we can choose a Lyapunov $\hep$-regularity function $R(x) = R_{\hep}(x)$ such that
\be
\| \Psi^u_x(u) - u\| \leq R(x) \|u\|^{1+\beta} \quad , \quad \| (\Psi^u_x)^{-1} (u) - u\| \leq R(x) \|u\|^{1+\beta}
\ee
for all $ x\in \ll$ and  $u \in E^u(x; r(x))$, and also  that 
\be
\|d \Phi^u_x(u)\| \leq R(x) \quad ,\quad  \|(d\Phi^u_x(u))^{-1}\| \leq R(x) \quad,  
\quad x\in \ll \:, \: u \in E^u (x ; r(x)) .
\ee
More precisely, as in \cite{LY2} (see (v) in (8.1) there), we will assume that there exists a global constant $R_0 > 0$
such that for all $u,v \in E^u(x;r(x))$ we have
\be
\frac{1}{R_0} d(\Phi^u_x(u), \Phi^u_x(v)) \leq \|u-v\| \leq R(x)\, d(\Phi^u_x(u), \Phi^u_x(v)) .
\ee

For any $x\in \ll$ consider the $C^{1+\beta}$ map 
$$\di \hf_x = (\Phi^u_{\varphi(x)})^{-1} \circ \varphi \circ \Phi^u_x : E^u (x)  \longrightarrow E^u (\varphi(x))$$
(defined locally near $0$).
We then have the relations
$$\hf^{-1}_x(\hE^u_j (\varphi(x) ; r(\varphi(x))) \subset \hE^u_j (x; r(x)) \quad , \quad 
\hf^{-1}_x(\tE^u_j (\varphi (x) ; r(\varphi (x))) \subset \tE^u_j (x; r(x)) \;$$
for all $x \in \ll$ and $j > 1$.

Given $y \in \ll$ and any integer $j \geq 1$ we will use the notation
$$\hf_y^j = \hf_{\varphi^{j-1}(y)} \circ \ldots \circ \hf_{\varphi (y)} \circ \hf_y\quad,
\quad \hf_y^{-j} = (\hf_{\varphi^{-j}(y)})^{-1} \circ \ldots \circ (\hf_{\varphi^{-2}(y)})^{-1}  \circ (\hf_{\varphi^{-1}(y)})^{-1} \;,$$
at any point where these sequences of maps are well-defined.

It follows from well known results (see  e.g. the Appendix in \cite{LY1}, Sect. 8 in \cite{LY2} or Sect. 3 in \cite{PS})  
that there exist a 
Lyapunov $\hep$-regularity function $\Gamma = \Gamma_{\hep} : \ll \longrightarrow [1,\infty)$ and an $\hep$-slowly varying
radius function $r = r_{\hep} : \ll \longrightarrow (0,1)$ (we will assume that it is the same as the one chosen above)
and  for each $x\in \ll$ a norm $\| \cdot \|'_x$  on $T_xM$ such that
\be
\|v\| \leq \|v\|'_x \leq \Gamma (x) \|v\| \quad ,\quad x\in \ll\:,\: v \in T_xM ,
\ee
and for any $x \in \ll$ and any integer $m \geq 0$, assuming 
$\hf_x^j(u), \hf_x^j(v) \in E^u(\varphi^j(x), r (\varphi^j(x)))$ are well-defined for all  $j =1, \ldots,m$, the following hold:
\be
\mu^m_j\, \|u-v\|'_x  \leq \|\hf_x^m(u) - \hf_x^m(v)\|'_{\varphi^m(x)} \leq \nu^m_{\tk}\, \|u-v\|'_x 
\quad , \quad u,v \in \tE^u_j(x; r(x)) ,
\ee
\be
\mu^m_1\, \|u-v\|'_x  \leq \|\hf_x^m(u) - \hf_x^m(v)\|'_{\varphi^m(x)}  \leq \nu_{\tk}^m \, \|u-v\|'_x\quad , \quad u,v \in E^u(x; r(x))  ,
\ee
\be 
\mu_1^m \, \|v\|'_x \leq \|d\hf_x^m(u)\cdot v\|'_{\varphi^m(x)} \leq \nu_{\tk}^m \, \|v\|'_x
\quad , \quad u \in E^u(x;r(x))\;,\; v \in E^u(x) ,
\ee
\be 
\mu_j^m \, \|v\|'_x \leq \|d\hf_x^m(0)\cdot v\|'_{\varphi^m(x)} \leq \nu_j^m \, \|v\|'_x
\quad , \quad v \in E^u_j(x) .
\ee

Another useful norm  is given by
$|u| = \max\{ \|\ui\| : 1 \leq i \leq \tk\}$, which is easily related to $\|\cdot\|$. Clearly,
$ \|u\| 
\leq \sum_{i=1}^{\tk} \|\ui\| \leq \tk \, |u| .$
Taking the regularity function $\Gamma(x)$ appropriately,   we have $|u|\leq \Gamma(x) \|u\|$, so
$\di \frac{1}{\tk}\, \|u\| \leq |u| \leq \Gamma (x) \|u\|$ for all $ x\in \ll$ and $u\in E^u(x)$.

Taylor's formula (see also Sect. 3 in \cite{PS}) implies that there exists 
a Lyapunov $\hep$-regularity function  $G = G_{\hep} : \ll \longrightarrow  [1,\infty)$ such that  for any $i = \pm 1$ 
and any $x \in \ll$ we have
$$\quad\|\hf_x^i(v) - \hf_x^i(u) - d\hf^i_x(u) \cdot (v-u)\| \leq G (x) \, \|v-u\|^{1+\beta} \:\:, \: u,v \in E^u (x ; r(x)) ,$$
and
$\|d\hf^i_x(u) - d\hf^i_x(0)\| \leq G (x)\, \|u\|^\beta$ for all  $u \in E^u (x; r(x))$.
{\bf Fix a global constant $\beta > 0$} with the above properties. We will assume that $\beta$
is chosen so that it satisfies (3.4) as well.

Next, for convenience of the reader we state three lemmas from \cite{St3} and \cite{St4} which will be used  in Sect. 5 below.

\bs

\noindent
{\bf Lemma 3.1.} (Lemma 3.3 in \cite{St3}) {\it There exist a Lyapunov $\hep$-regularity function
$L : \ll \longrightarrow [1,\infty)$ and an $\hep$-slowly varying radius function 
$r: \ll \longrightarrow (0,1)$
such that for any $x\in \ll$, any integer $p \geq 1$ and any $v\in E^u(z, r (z))$ with  $\|\hf_z^p(v)\|\leq r(x)$,  where $z = \varphi^{-p}(x)$, we have
$$\|\wo_p - \vo_p\| \leq L(x) |v_p|^{1+\beta} ,$$
where $v_p = \hf^p_z(v) \in E^u (x)$ and $w_p = d\hf_z^p(0)\cdot v \in E^u (x)$. Moreover, if $|v_p| = \|\vo_p\| \neq 0$, then}
$1/2 \leq \|\wo_p\|/\|\vo_p\| \leq 2$.

\bs

\noindent
{\bf Remark.} Notice that if $v\in E^u_1(z, r (z))$ in the above lemma, then 
$v_p, w_p \in E^u_1(x)$, so $\|w_p - v_p \| \leq L(x)\, \|v_p\|^{1+\beta}$.

\bs

We will now state some consequences from Sect. 10 in \cite{St4} that apply to every Anosov flow.
For any $v = \vo + \vt + \ldots+ v^{(\tk)} \in E^u(x)$ with $v^{(j)} \in E^u_j(x)$, set
$\chvt =  \vt + \ldots+ v^{(\tk)} \in \tE_2^u(x) .$

Repeating the proof of Lemma 3.5 in \cite{St3} (with the above slightly different choices of $\hmu_1$, $\hnu_1$, $\hmu_2$
and $\nu_2$; see also Lemma 10.1 in \cite{St4}) we get the following.

\bs

\noindent
{\bf Lemma 3.2.}  {\it  Choosing the constant $\hep' \in (0,\hep)$ sufficiently small, there exists an $\hep'$-slowly 
varying radius function  $\hr(x) \leq r(x)$ on $\ll$ such that for any $x\in \ll$ and any 
$V = \Vo + \chVt \in E^u(x; \hr(x))$,
setting $y = \varphi^{-1}(x)$ and $U = \hf_x^{-1}(V)$, we have}
$\di \|\chUt\|'_{y} \leq \frac{\|\chVt\|'_x}{\mu_2}$ and
$\di \|\Uo\|'_{y} \geq \frac{\|\Vo\|'_x}{\nu_1}$.

\bs

\noindent
{\bf Lemma 3.3.}   (Lemma 10.7(b) in \cite{St4})  
 {\it There exist a $\frac{5\hep}{\beta}$-slowly varying radius  functions $\hr(x) \leq r(x)$ 
and a Lyapunov $4 \hep$-regularity function $L(x)$, $x\in \ll$, such that
for any $x \in \ll$ and any integer $p \geq 1$, setting $z = \varphi^{-p}(x)$,  the map
$$F^p_x = d \hf^p_z(0) \circ (\hf^p_x)^{-1}:  E^u (x; \hr (x)) \longrightarrow E^u (x; \hr(x))$$
satisfies
$$\left\| \left[ (F^p_x (a))^{(1)} - (F^p_x (b))^{(1)} \right] - [\ao - \bo ]\right\|  
\leq L (x) \, \left( \|a - b \|^{1+\beta} + \| b \|^\beta \cdot \|a - b\| \right)  $$
for all $a,b \in E^u (x ; \hr (x))$. Moreover,
$$\di \frac{1}{2} \| a - b\| \leq \left\| d\hf^p_z(0)\cdot \left[ (\hf^p_x)^{-1}(a) - (\hf^p_x)^{-1}(b)\right] \right\|\leq 2  \| a - b\| $$
for all $a,b \in E^u_1 (x ;  \hr (x))$. }

\def\tct{\tilde{c}_2}
\def\hX{\widehat{X}}
\def\hY{\widehat{Y}}
\def\hZ{\widehat{Z}}
\def\diamteh{\diam_{\htheta}}
\def\dteh{D_{\htheta}}
\def\talpha{\tilde{\alpha}}
\def\ttheta{\tilde{\theta}}
\def\diamtet{\diam_{\ttheta}}
\def\dtet{D_{\ttheta}}
\def\diamteo{\diam_{\theta_1}}
\def\dteo{D_{\theta_1}}

\section{Estimates and constructions involving cylinders}
\setcounter{equation}{0}

\subsection{Cylinders defined by the Markov family}

Let again $M$ be a $C^2$ compact  Riemannian manifold and let $\phi_t$  be a  $C^2$ transitive Anosov flow on $M$.
{\bf  Here we do not assume that the flow is contact.}

Choose a small $\hep > 0$, as in Sect. 3,  but we may need to make it smaller later.
Throughout we will again assume that 
$R(x)$, $\Gamma(x)$, $D(x)$ and $L(x)$ are Lyapunov $\hep$-regularity functions, while $r(x)$ is
an $\hep$-slowly varying radius function so that it satisfies (3.6) -- (3.14)  and the conclusions of  Lemma 3.1. 
Replacing $r(x)$ with the smaller regularity function $\hr(x)$, 
without loss of generality we will assume that the conclusions of Lemmas 3.2 and 3.3 hold with $\hr(x)$ replaced by $r(x)$.

In what follows we assume that $\tRR = \{ \tR_i\}_{i=1}^{k_0}$  is a fixed Markov falimiy for 
$\phi_t$ on $M$ of size $\chi < \ep_1/2 < \ep_0/2$ and $\rr = \{ R_i\}_{i=1}^{k_0}$  is the related pseudo-Markov 
family as in Sect. 2. We will use the notation associated with these from Sect. 2, and
we will assume that for any $i = 1, \ldots,k_0$ the point $z_i$ is chosen so that $z_i \in \Intu(W^u_{R_i}(z_i))$.
For any $x\in R$, any $y \in \tR$ and $\delta > 0$ set 
$$B^u(x,\delta) = \{ y\in W^u_{R_i}(x) : d(x,y) < \delta\}  \quad , \quad \tB^u(y,\delta) = \{ z\in W^u_{\tR_i}(z) : d(z,y) < \delta\} .$$
In a similar way define $B^s(x,\delta)$. The {\it open ball} with centre $x$ and radius $r > 0$ on $M$ with respect to
the Riemannian metric will be denoted simply by $B(x,r)$.

Given an unstable leaf $W = W^u_{R_i}(z)$ in some rectangle $R_i$ 
and an admissible sequence\\ $\ii = (i_0, \ldots,i_m)$ of integers $i_j \in \{ 1,\ldots, k_0\}$, the set
$$\cc_W[\ii] = \{ x\in W : \pp^j(x) \in R_{i_j} \:, \: j = 0,1,\ldots,m\}$$
will be called a {\it cylinder of length} $m$ in $W$ (or an {\it unstable cylinder} in $R$ in general). 
When $W = U_i$ we will simply write $\cc[\ii]$. 
In a similar way one defines cylinders $\tcc_V[\ii]$, 
where $V = W^u_{\tR_i}(z)$ is an unstable leaf in some rectangle $\tR_i$.

Let
$\pr_D : \cup_{i=1}^{k_0} \phi_{[-\ep,\ep]}(D_i) \longrightarrow \cup_{i=1}^{k_0} D_i$
be the {\it projection along the flow}, i.e. for all $i = 1,\ldots,k_0$ and all 
$x\in \phi_{[-\ep,\ep]}(D_i) $ we have $\pr_D(x) = \pr_{D_i}(x)$ (see Sect. 2). 
The shift along the flow determines bi-H\"older continuous bijections 
$$\T_z :  W^u_{\tR}(z) \longrightarrow \T_z(W^u_{\tR}(z)) \subset W^u_{\ep_0}(z) \quad , 
\quad \tPsi : R \longrightarrow \tR ,$$
where both $\T_z$ and $\tPsi : W^u_{R}(z) \longrightarrow W^u_{\tR}(z)$ are bijection for all $z\in S_i$ and all $i$.  

Given a cylinder $\cc = \cc_W[\ii]$ in some $R_i$ we will frequently use the notation
$\tcc = \tcc_W[\ii] = \tPsi(\cc)$. This is then a cylinder of the same length in $\tR_i$.
Sometimes it will be more convenient to work with the projection of the cylinders on actual
unstable manifolds, and we will use the notation $\hcc = \T_z(\tcc)$
for any cylinder $\tcc$ in some
$\tR_i$ and any $z \in W^u_{\tR_i}$. Then $\hcc \subset W^u_\ep(z)$. The map 
$\T_z : \tcc \longrightarrow \hcc$ is uniformly Lipschitz.
Clearly there exist global constants $0 < \tilde{c}_1 < \tilde{c}_2$, independent of $\cc$ and $z$  such that 
\be
 \tilde{c}_1\, \diam(\tcc) \leq \diam(\hcc) \leq \tilde{c}_2\, \diam (\tcc) .
\ee

\def\ttc{\tilde{c}}

Although the rectangles
$\tR_i$ could have complicated structure\footnote{They are not connected in general except in 3D, as Chernov points out in Sect. 9
in \cite{Ch1}.}, and could be rather "fragmented", each of them contains a non-empty open subset\footnote{In fact $\tR_i$ is the closure 
of such an open subset of $D_i$.} of the corresponding submanifold $D_i$.
It is rather easy to show that there exists a constant $r_0 > 0$ such that for every $i = 1, \ldots, k_0$ and every 
$x \in \tR_i$ there exists $y \in W^u_{\tR_i}(x)$ such that $\dist(y, \partial \tR_i) > r_0$ and 
$B^u(y,r_0) \cap \tR_i \subset W^u_{\tR_i}(x)$.
From now on we will assume that the constant $r_0 > 0$ is chosen so that it has the property just described.
A few more restrictions on $r_0$ will be imposed later.


\def\thetaoo{\theta_0}
\def\Kn{K^{(n)}}
\def\Knm{K^{(n_m)}_m}
\def\Knp{K^{(n_p)}_p}
\def\tK{\widetilde{K}}
\def\tn{\tilde{n}}
\def\tKn{\tK^{(\tn)}}
\def\KnK{K^{(n(K))}}

\subsection{The Gibbs measure $\nu$, the Ruelle operators $\lab$ and the metric $D_\theta$}

Throughout this paper we will be using some basic concepts and facts from the so called\\ Thermodynamics Formalism (see e.g.
\cite{PP}, \cite{R2}, \cite{Wa} or \cite{K}).

Let $F_0 : M \longrightarrow \R$ be a H\"older continuous function  and let 
$\m$ be the  {\it Gibbs measure} determined by $F_0$ on $M$ defined on the set $\ll$ of Lyapunov
regular points (\cite{Si}, \cite{B2}, \cite{Ch2}, \cite{PP}).
It induces a {\it Gibbs measure} $\mu_0$  on $R$ (with respect to the Poincar\'e map $\pp$) for the function
$$\di F(x) = \int_0^{\tau(x)} F_0 (\phi_s(x))\, ds \quad, \quad x\in \R .$$
The latter is H\"older and, using Sinai's Lemma, it is cohomologous to a {\it H\"older function $f: R \longrightarrow \R$  
which is constant on stable leaves in rectangles} $R_i$ in $R$. 
Setting $g = f - P_f \tau$,  where $P_f \in \R$ is chosen so that the topological pressure of $g$ with respect to the Poincar\'e map 
$\pp : R \longrightarrow R$  is $0$, we get a function on $R$ that
depends on forward coordinates only, 
so it can be considered as a function on $U$, i.e. on $\saa$. 

From now on in this paper we will assume that {\bf $f \in \ff_{\theta_0} (\hU)$ is a fixed real-valued function 
for some small constant $\thetaoo > 0$
and $g = f - P_f\, \tau$}, where $P_f\in \R$ is such that $\Pr_{\sigma} (g) = 0$. 
Notice that $\ff_{\theta_0}(\hU) \subset \ff_\theta(\hU)$ for $0 < \theta_0 \leq \theta < 1$.
In what follows we will consider parameters $\theta \in [\theta_0, 1)$.
Set $\Fa = f - (P_f + a) \tau$.
By Ruelle-Perron-Frobenius' Theorem (see e.g. Theorem 2.2 in \cite{PP}) for any real number $a$  
with $|a|$ sufficiently small, as an operator on the space $C(U)$ of continuous functions $h : U \longrightarrow \R$ with the sup-norm
(which we identify with $C(\saa)$ with the sup-norm), $L_{\Fa}$ has a  {\it largest eigenvalue}  $\lambda_{a}$ and 
there exists a (unique) regular probability measure $\hnu_a$ on $U$ with  $L_{\Fa}^*\hnu_a = \lambda_a\, \hnu_a$, i.e.  
$$\int_U L_{\Fa} H \, d\hnu_a = \lambda_a\, \int_U H\ d\hnu_a$$
for every $H \in C(U)$. The corresponding eigenfunctions belong to $\ff_{\thetaoo}(\hU)$.
Fix a corresponding (positive) eigenfunction $h_{a} \in \ff_{\thetaoo}(\hU)$ such that $\int h_{a} \, d\hnu_a = 1$. 
Then $d\nu = h_0\, d\hnu_0$ defines a  $\sigma$-invariant  probability measure  $\nu$ on $U$, called the {\it  Gibbs measure}
determined by the function $F^{(0)}$. 
{\bf This is the measure on $U$ that we will use throughout this paper.}
Since $\Pr_\sigma (f- P_f\tau) = 0$, it follows from the main properties of pressure (cf. e.g. chapter 3 in \cite{PP} or \cite{R2}) 
that $|\Pr_\sigma(\Fa)| \leq \|\tau \|_0 \, |a|$.  Moreover, for small $|a|$ the maximal eigenvalue
$\lambda_{a}$ and the eigenfunction $h_{a}$ are Lipschitz in  $a$, so
there exist constants $a_0 > 0$ and $C > 0$ such that $|h_{a} - h_0| \leq C\, |a|$ on $\hU$ and
$|\lambda_a - 1| \leq C |a|$ for  $|a| \leq a_0$.

We will identify $\mu_0$ with the measure on $\tR$ defined by $\mu_0(\tPsi(A)) = \mu_0(A)$ for every Borel subset $A$ of $R$.
Apart from that we will use the measure $\nu$ on subsets of $\tPsi(U)$ simply by setting $\nu(\tPsi(A)) = \nu(A)$
for every measurable subset $A$ of $U$. The same will apply to subsets of $W^u_{\tR}(x)$ for $x \in \tR$, identifying these with
subsets of $U$ using projections along stable leaves in $\tR$.

We should also remark that the properties of the set $\ll$ in $R$ imply that $\mu(U \cap \ll) = 1$ for every
$\sigma$-invariant Borel probability measure $\mu$ on $U$. In particular, for every such measure $\mu$ on $U$ 
and every Borel subset $A$ of $U$ we have $\mu(A\setminus \ll) = 0$, so $\mu(A) = \mu(A \cap \ll)$.


For $a,b\in \R$, $|a|\leq a_0$ and $|b|\geq 1$,  as in \cite{D1}, consider the function
$$\fa(u) = f (u) - (P_ f+ a) \tau(u) + \ln h_{a}(u) -  \ln h_{a}(\sigma(u)) - \ln \lambda_{a}$$
and the operators 
$$\lab = L_{\fa - \i\,b\tau} : C(U) \longrightarrow C (U)\:\:\: , \:\:\:  \ma = L_{\fa} : C (U) \longrightarrow C(U) .$$
Then $\ma \; 1 = 1$.
It is easy to see that
$$\di |(\lab^m h)(u)| \leq (\ma^m |h|)(u)$$
 for all $u \in U$,
$h\in C (U)$ and  $m \geq 0$. Moreover, $L_{f^{(0)}}^*\nu = \nu$,  i.e.  
$$\di \int_U L_{f^{(0)}} H \, d\nu = \int_U H\, d\nu \quad, \quad H \in C (U) .$$

For any integer $m \geq 1$ and any function $h : U \longrightarrow \C$ define $h_m : U \longrightarrow \C$ by
$$h_m(u) = h(u) + h(\sigma(u)) + \ldots + h(\sigma^{m-1}(u)) .$$


Since $g$ has zero topological pressure with respect to the shift map $\sigma : U \longrightarrow U$, 
there exist constants $0 < c_1 \leq c_2$  
such that for any cylinder $\cc = \cc^u[i_0, \ldots, i_{m-1}]$  of length $m$ in $U$ we have
\be\label{eq:Gibbs}
 c_1 \leq  \frac{\nu(\cc)}{e^{g_m(y) }} \leq c_2 \quad , \quad y \in \cc ,
\ee
(see e.g. \cite{PP} or \cite{R2}). Moreover, we may assume that $g_0 = \max g < 0$ and there exist constants 
$0 < \rho_1 \leq \rho_2 < 1$ such that for some constants $c_1, c_2 > 0$ as above we have 
\be\label{eq:exp}
c_1 \, \rho_1^m \leq  \nu(\cc)  \leq c_2 \, \rho_2^m  ,
\ee
for every cylinder $\cc = \cc^u[i_0, \ldots, i_{m-1}]$  of length $m$ in $U$
(see e.g. Proposition 2 in \cite{Po} or pp. 54-55 in \cite{Ch2}).

As in \cite{St4}, here we will make a substantial use of the {\it metric $\dte$ on $U$} defined in Sect. 1 above.
For a non-empty subset  $A$ of $U$ (or some $W^u_R(x)$) let $\diamte(A)$ be the {\it diameter}  of $A$ with respect to $\dte$.

Let the constants $d_0 > 0$, $1 < \gamma < \gamma_1$ and $\alpha_1 > 0$  be as in Sect. 2.
{\bf Choose the constants $\alpha_2 \in (0,1)$ and $\theta \in [\thetaoo ,1)$ so that }
$$\frac{1}{\gamma^{\alpha_1 \beta}} \leq \theta \leq \frac{1}{(\gamma_1)^{\alpha_2}}  < 1 ,$$
where $\beta \in (0,1)$ is as in (3.4). 
In particular $\alpha_2 \in (0,1)$ satisfies 
$\di 0 < \alpha_2 \leq \frac{ \alpha_1 \beta\, \log \gamma}{\log \gamma_1} .$

The following lemma is similar to Lemma 5.1 in \cite{St4} and its proof is the same.

\bs

\noindent
{\bf Lemma 4.1.}  {\it Let $\theta \in [\theta_0 ,1)$ be as above}.

\ms

(a) {\it For any cylinder $\cc$ in $U$ the characteristic function $\chi_{\cc}$ of
$\cc$ on $U$ is Lipschitz with respect to $\dte$ and $\Lip_{\theta}(\chi_\cc) \leq 1/\diamte(\cc)$.}

\ms

(b) {\it There exists a constant $C_1 > 0$ such that if $x,y\in \hU_i$ for some $i$,  then 
$|\tau(x) - \tau(y) |  \leq  C_1\, \dte(x,y) .$
That is, $\tau\in \ff_{\theta}(\hU)$. }

\ms

(c) {\it There exists a  constant $C_1 > 0$ such that for any $z\in R$, any cylinder 
$\cc$ in $W^u_R(z)$ and any $x,y\in \cc$ we have 
$d(\tPsi(x),\tPsi(y))  \leq  C_1\, \dte(x,y)$ and
$\dte (x,y) \leq C_1 (d (\tPsi(x), \tPsi(y)))^{\alpha_2}$.
Therefore, for $\tcc = \tPsi(\cc)$ we have }
$$\di  \diam (\tcc) \leq C_1 \, \diamte(\cc)  \quad, \quad
\diamte (\cc) \leq C_1 (\diam(\tcc))^{\alpha_2} .$$

\ms


Using the analytical dependence of $h_a$ and $\lambda_a$ on $a$ and assuming that the 
constant $a_0 > 0$ is sufficiently small, there exists $T_0 = T_0(a_0) > 0$ such that
\begin{equation}\label{eq:T0-cond}
T_0 \geq \max \{ \, \|\fa \|_0 \, , \, |\fa_{|\hU}|_{\theta} \, , |\tau_{|\hU}|_{\theta} \, \}\;
\end{equation}
for all $|a| \leq a_0$. We will assume from now on that $a_0 > 0$ and $T_0 > 0$ have these properties.
Taking the constant $T_0 > 0$ sufficiently large, we have  $\|\fa - f^{(0)}\|_0 \leq T_0\, |a|$ 
on $\hU$ for $|a| \leq a_0$.

As in \cite{D1}  we have the following Lasota-Yorke type inequality (see Lemma 5.2 and its proof in the Appendix in \cite{St4}). 

\bs

\noindent
{\bf Lemma 4.2.}  {\it There exists a constant $A_0 > 0$, depending on $\|f\|_\theta$, such that for all $a\in \R$ with $|a|\leq a_0$ 
the following hold: If the functions $h$ and  $H$ on $\hU$ and the constant $B > 0$ are such that $H > 0$ on $\hU$ and 
$|h(v) - h(v')| \leq B\, H(v')\, \dte (v,v')$
for any $i$ and any $v,v'\in \hU_i$, 
then for any $b\in \R$ with  $|b|\geq 1$ and any integer $m \geq 1$ we have 
$$ |\lab^m h(u) - \lab^m h(u')| \leq  A_0\,\left[ B\,\theta^m \, (\ma^m H)(u')  + |b|\, (\ma^m |h| )(u')\right]\, \dte (u,u')$$
whenever $u,u'\in \hU_i$ for some $i = 1, \ldots,k_0$.}

\bs

By a { \bf Pesin set  in $R$} we mean a compact subset $P_0$ of $R \cap \ll$ with $\mu_0(P_0) > 0$ such that
the Lyapunov $\hep$-regularity function $R(x)$ is bounded on $P_0$. 

For later convenience we will choose a Pesin set $P_0$ which has an additional property. The construction can be
done so that we work with compact subsets of a given Pesin set, so we do not have to worry about the
Lyapunov regularity function $R(x)$ being bounded on the new set. 
It is more convenient to work with projections $K_0 = \piU(P_0)$ of Pesin sets in $U$.

In what follows we always consider topological pressure $\Pr = \Pr_\sigma$ with respect to the
shift map $\sigma : U \longrightarrow U$.

Given a compact subset $K$ of $U$ and an integer $n \geq 1$, denote by $\Kn$ the union of all cylinders
$C$ of length $n$ in $U$ such that $C \cap K \neq \e$. Clearly $\Kn$ is a finite union of cylinders, so it
is both an open and closed subset of $U$ and its characteristic function $\chi_{\Kn}$ is continuous. Thus,
$g - \chi_{\Kn}$ has a well-defined topological pressure.

\bs

\noindent
{\bf Remark.} Let $V$ be a non-empty subset of $U$ which is a finite union of cylinders. Let $\mu$ the equilibrium state
of $g - \chi_V$ in $U$ (see e.g. \cite{Wa} or \cite{K}). Then
$$\Pr(g - \chi_V) = h_{\mu}(\sigma) + \int_U (g- \chi_V)\, d\mu
\geq h_{\hmu}(\sigma) + \int_U (g- \chi_V)\, d\hmu$$
for any $\sigma$-invariant Borel probability measure $\hmu$ on $U$. Here $h_\mu(\sigma)$ is the entropy of $\sigma$
with respect to $\mu$. What we claim here is that, although the function $g- \chi_{V\cap \ll}$ is not continuous in general,
it has an equilibrium state and it is the same measure $\mu$. Indeed, let $\hmu$ be an arbitrary
$\sigma$-invariant Borel probability measure on $U$. Then, as above,
$$h_{\mu}(\sigma) + \int_U (g- \chi_V)\, d\mu \geq h_{\hmu}(\sigma) + \int_U (g- \chi_V)\, d\hmu .$$
However, the properties of $\ll$ imply $\mu(V) = \mu(V\cap \ll) + \mu(V\setminus \ll) = \mu(V\cap \ll)$, since
$\mu(V\setminus \ll) = 0$. Thus,
$$\int_U \chi_V\, d\mu = \mu(V) = \mu(V\cap \ll) = \int_U \chi_{V \cap \ll}\, d\mu .$$
Similarly, $\di \int_U \chi_V\, d\hmu = \int_U \chi_{V \cap \ll}\, d\hmu$. It follows now that
\begin{eqnarray*}
h_{\mu}(\sigma) + \int_U (g- \chi_{V\cap \ll})\, d\mu
& = & h_{\mu}(\sigma) + \int_U (g- \chi_V)\, d\mu \geq h_{\hmu}(\sigma) + \int_U (g- \chi_V)\, d\hmu\\
& = & h_{\hmu}(\sigma) + \int_U (g- \chi_{V\cap \ll})\, d\hmu .
\end{eqnarray*}
This proves that $\mu$ is an equilibrium state for $g - \chi_{V\cap \ll}$ and therefore the latter has a well-defined
topological pressure $\Pr (g - \chi_{V\cap\ll}) = \Pr (g - \chi_{V})$. 

\ms

The following lemma is derived from simple properties of topological pressure and the above Remark.


\bs

\noindent
{\bf Lemma 4.3.} {\it We can choose the Pesin set $P_0$ in $R\cap \ll$ so that $\mu_0(P_0) > 0$ and its projection $K_0 = \piU(P_0)$ satisfies
$\nu(K_0) > 0$ and}
\be\label{eq:K-cond}
\lim_{n\to\infty} \Pr (g - \chi_{\Kn_0}) < 0 .
\ee

\noindent
{\it Proof.} Notice that $\Pr (g - \chi_{\Kn_0})$ is an increasing function of $n$, so the above limit exists.

Assume that there is no Pesin set $P_0$ in $R\cap \ll$ so that its projection $K_0 = \piU(P_0)$ has $\nu(K_0) > 0$ and satisfies (4.5).
Then for every compact subset $K$ of $U\cap \ll$ with $\nu(K) > 0$ which is the projection of a Pesin set $P$ with $\mu_0(P) > 0$
and every $\ep > 0$ there exists an integer $n \geq 1$ so that $\Pr (g - \chi_{\Kn}) > - \ep$.

Fix an arbitrary compact subset $\tK$ of $U \cap \ll$ with $\nu(\tK) > 0$ which is the projection of 
Pesin set $\tP$ with $\mu_0(\tP) > 0$ and an integer $\tn \geq 1$ so that $- 1 < \Pr (g - \chi_{\tKn}) .$
Then
\be
- 1 < \Pr (g - \chi_{\tKn}) = - c
\ee
for some constant $c > 0$. 
For every point $x\in \tKn \cap \ll$ consider a compact subset $K = K(x)$ of $\ll$ containing $x$
with $\nu(K) > 0$ which is the projection of a Pesin set $P = P(x)$ with $\mu_0(P) > 0$. Then
take a large integer $n = n(x) > 1/c$ such that $\di - \frac{1}{n(x)} < \Pr (g - \chi_{K^{(n(x)}}) .$
The open sets $K^{(n(x))}(x)$ obtained in this way form
an open cover of $\tKn \cap \ll$. Since $\saa$ and so $U$ are separable metric spaces, we can choose
a countable sub-cover. Ordering the sets $K = K(x)$ involved in this sub-cover in an increasing order, we
obtain a sequence
$$K_1 \subset K_2 \subset \ldots \subset K_m \subset \ldots \subset U \cap \ll$$
of compact subsets of $U \cap \ll$ and a corresponding sequence of integers $n_m > 1/c$ such that
$$\tKn \cap \ll \subset \bigcup_{m=1}^\infty \Knm  .$$ 
Let us emphasise here that the compact sets $K_m$ and the integers $n_m > 1/c$ are chosen so that
\be
- \frac{1}{n_m} < \Pr (g - \chi_{\Knm}) .
\ee
Moreover, $K_m$ is the projection of a Pesin set $P_m$ with $\mu(P_m) > 0$
and $\nu(K_m) > 0$ for all $m \geq 1$.

It follows now that
$$\chi_{\tKn \cap \ll} \leq \sup_{m \geq 1} \chi_{\Knm \cap \ll} .$$
Hence
$$g - \chi_{\tKn \cap \ll} \geq \inf_{m\geq 1} (g- \chi_{\Knm \cap \ll} ), $$
and therefore, using the Remark before Lemma 4.3 and (4.6), we have
$$- c = \Pr\left(g - \chi_{\tKn}\right) = \Pr\left(g - \chi_{\tKn \cap \ll}\right)  \geq \inf_{m\geq 1} \Pr \left(g- \chi_{\Knm \cap \ll} \right)
= \inf_{m\geq 1} \Pr \left(g- \chi_{\Knm} \right) . $$
However by the choice of the integers $n_m$ we have $\frac{1}{n_m} < c$ for all $m \geq 1$, so (4.7) implies
$$\Pr \left(g- \chi_{\Knm} \right) \geq - \frac{1}{n_m} > - c $$
for all $m \geq 1$. Now the above yields
$$- c = \Pr\left(g - \chi_{\tKn}\right)  \geq \inf_{m \geq 1} \Pr \left(g- \chi_{\Knm} \right) > - c ,$$
which is a contradiction.

This proves  that there exists a Pesin set $P_0$ in $R\cap \ll$ with $\mu_0(P_0) > 0$ so that its projection $K_0 = \piU(P_0)$ 
has $\nu(K_0) > 0$ and satisfies (4.5).
\endofproof

\bs

{ \bf Fix an arbitrary Pesin set $P_0$ in $R\cap \ll$} such that $\mu_0(P_0) > 0$ and its projection $K_0 = \piU(P_0)$ has
$\nu(K_0) > 0$ and satisfies (4.5).
This will be the Pesin set we will use throughout the rest of the paper.

Since the Lyapunov $\hep$-regularity function $R(x)$ is bounded on $P_0$, the functions 
$\Gamma(x)$, $D(x)$ and $L(x)$ are also bounded above by some constants on $P_0$.
Similarly, the $\hep$-slowly varying radius function $r(x)$ 
is bounded  below by some constant on $P_0$. Thus, we may assume that
$$R (x) \leq R_0 \:\:, \:\: r(x) \geq r_0 \:\: , \:\: \Gamma(x) \leq \Gamma_0 \:\: , \:\:
L(x) \leq L_0 \: \: , \:\: D(x) \leq D_0 $$
for all $x\in P_0$ {\bf for some positive constants $R_0, \Gamma_0, L_0, D_0 \geq 1$ and $r_0 > 0$.  We fix $r_0 > 0$}
so that $r_0 \leq \frac{1}{R_0}$.
We will also use the Pesin set {\bf $\tP_0 = \tPsi(P_0)$} in $\tR$, and we will assume that the 
functions $R(x)$, $r(x)$, etc. satisfy the same bounds as above on $\tP_0$.

\def\hgamma{\hat{\gamma}}
\def\tH{\widetilde{H}}
\def\Uo{U^{(1)}}
\def\Wo{W^{(1)}}
\def\hH{\widehat{H}}
\def\hG{\widehat{G}}

\subsection{Technical lemmas on sizes of cylinders}

We continue with the notation and assumptions in Sects. 4.1 and 4.2. In addition we will assume that
the Markov family $\rr$ is chosen so that  $6\tau_0 < \beta ,$
where $\tau_0$ is the constant from Sect. 2 and $\beta > 0$ is the constant with (3.4), fixed in Sect. 3.

The following technical lemma will be used significantly later on.

\bs

\noindent
{\bf Lemma 4.4.}  {\it There exist a global constant $C_3 > 0$ and constants $0 < \hep_2 < \hep_1$ 
with $\hep_i \leq \con \hep$ that can be made arbitrarily small with $\hep$, such that if $\tcc$ is a cylinder of length $m$ 
in $\tR$  with $\diam(\tcc) < r_0$ and $z_0 \in \tcc\cap \tP_0$, then:}

\ms

(a) {\it There exists an  integer $k$ with
$m \hep_2 \leq k \leq m \hep_1$
such that $\di \tpp^{m-k}(\tcc) \subset B^u(z', r(z'))$, where $z' = \tpp^{m-k}(z_0)$.}

\ms

(b) {\it For $p = [\ttau_{m-k}(z_0)]$ with $k$ as above,  we have
$\di \hvarphi_{z_0}^p(\hcc) \subset B^u(z_p, r(z_p))$, where $z_p = \varphi^p(z_0)$ and $\hcc = \T_{z_0}(\tcc)$.
Moreover,
$$\di \diam(\varphi^p(\tcc)) \leq \frac{3C_3}{\gamma^{\alpha_1 k}} \leq r_0 e^{-p \hep} \leq r(z_p) ,$$
where $\alpha_1 > 0$ is the H\"older constant from Sect. 2.}

\ms

  (c) {\it For every unstable cylinder $\cc$ of length $m$  in $\tR$ and $z_0 \in \tcc  \cap \tP_0$, 
we have
\be\label{eq:diamC}
\frac{e^{-m \hep_7/\ttau_0}}{C_3 \lambda_1^{m\tau_0} } \leq \frac{e^{-q \hep_7}}{C_3 \lambda_1^q } \leq \diam(\tcc) \leq 
\frac{C_3\, e^{q \hep_7}}{\lambda_1^{q}} \leq \frac{C_3 e^{m \tau_0\, \hep_7}}{\lambda_1^{m\ttau_0}},
\ee
where $q = [\ttau_m(z_0)]$ and $\hep_7$  is a constant with $0 < \hep_7 < \con\, \hep$.
Moreover, there exist $\hx_0 \in \T_{z_0}(\tcc)$ and a constant $0 < \hep_5 < \con\, \hep$ such that 
$u_0 = (\Phi^u_{z_0})^{-1}(\hx_0) \in E^u(z_0 ,\ep_1)$ and we have
\be
\|\uo_0\| \geq \frac{e^{-q \hep_7}}{C \lambda_1^{q} }
\ee
where $C =  2 R^2_0\Gamma_0/c > 0$ for some global constant $c > 0$. }

\bs

\noindent
{\it Proof of Lemma.} 4.4. We use some bits from the proof of Lemma 4.2(a) in \cite{St4}, however we need a lot
more precision and details.

Let $\tcc$ be a cylinder of length $m$ in  $\tRR$. 
Fix an arbitrary $z_0 \in \tcc \cap \tP_0$. 
Since $m$ is the length of  $\tcc$, $\tpp^m(\tcc)$ contains a whole unstable leaf of a proper rectangle $\tR_{j_0}$.  
Let $Z_m  = \tpp^m(z_0) \in \tR_{j_0}$. Set $Z_j = \tpp^j(z_0)$ and $z_j = \varphi^j(z_0)$ for $j \geq 1$.
By the choice of the constant $r_0 > 0$ (see Sect. 4.1) there exists $y_0 \in W^u_{\tR_{j_0}}(Z_m)$ such that
$\dist(y_0, \partial \tR_{j_0}) > r_0$ and  $B^u(y_0,r_0) \cap \tR_{j_0} \subset W^u_{\tR_{j_0}}(Z_m)$.
In particular, for every point $b'\in B^u(y_0,r_0)$
there exists $b \in \tcc$ with $\tpp^{m}(b) = b'$. Also, $\diam (\tR_{j_0}) \leq \chi$, so $\tR_{j_0} \subset B^u(y_0, 2\chi)$
and $\tR_{j_0} \subset B^u(Z_m, 2\chi)$. Thus, $\tpp^m(\tcc) \subset B^u(Z_m, 2\chi)$.

We will now choose $k$ with $0 < k \leq m$ so that 
\be
\tpp^{m-k}(\tcc) \subset B^u(Z_{m-k}, r(Z_{m-k})) ,
\ee
where  $Z_{m-k} = \tpp^{m-k}(z_0) = \tpp^{-k}(Z_m) .$
Since $r(\cdot)$ is a Lyapunov $\hep$-regularity function (see Sect. 3), we have 
$r(Z_{m-k}) \geq r(z_0) e^{-(m-k)\hep} \geq r_0 e^{-(m-k)\hep}$.
For every integer $0 \leq k \leq  m$, by (2.1)  we have 
$\diam (\tpp^{-k}(B^u(Z_{m-k}, 2\chi))) \leq \frac{2\chi}{d_0 \gamma^k}$. Thus, (4.10) will be satisfied if
\be
\frac{2\chi}{d_0 \gamma^k} \leq r_0\, e^{-(m-k) \hep} .
\ee

Let $k\geq 0$ be a number with (4.11). 
Then 
$\di  \gamma^k e^{-(m-k) \hep} \geq \frac{2 \chi}{r_0 d_0 }$,
so
$e^{-m \hep} (\gamma e^{\hep})^k \geq \con$, i.e. $(\gamma e^{\hep})^k \geq \con e^{m \hep}$.
Setting $\hgamma = \gamma e^{\hep}$, we get
$k \log \hgamma \geq m \hep + \con $
for some positive global constant $\con$.

It follows from all the above that
\begin{equation}
\tpp^{-k}(B^u(y_0, r_0)) \subset \tpp^{m-k}(\tcc) \subset \tpp^{-k}(B^u(Z_{m} , 2 \chi)) 
 \subset B^u(Z_{m-k}, 2\chi/(d_0 \gamma^k)) \subset B^u(Z_{m-k}, r(Z_{m-k})).
\end{equation}

Set $T = \ttau_{m-k}(z_0)$ and $p = [T]$, so that $p \leq T < p+1$.

It follows from (2.1) and the $\alpha_1$-H\"older continuity of weak unstable manifolds of the flow $\phi_t$ that
for every $z \in \tcc$ we have
\be
|\ttau_{m-k}(z) - \ttau_{m-k}(z_0)| \leq \Con \, d(\tpp^{m-k}(z), \tpp^{m-k}(z_0)) \leq \Con \left(\frac{\chi}{d_0\gamma^k}\right)^{\alpha_1}
\leq \frac{C}{\gamma^{\alpha_1 k}}
\ee
for some global constant $C > 0$. We will show now that
\be
d(\varphi^p(z), \varphi^p(z_0)) \leq \frac{3 C}{\gamma^{\alpha_1 k}}
\ee
for all $z \in \tcc$. 

\ms

\noindent
{\it Proof of} (4.14): Given $z \in \tcc$, there are two cases to consider for $t = \ttau_{m-k}(z)$ and $T = \ttau_{m-k}(z_0)$.

\ms

\noindent
{\bf Case 1.} $t  < p$. Then $t < p \leq T$ and by (4.13), $p-t \leq T-t \leq \frac{C}{\gamma^{\alpha_1 k}}$
and $T - p \leq \frac{C}{\gamma^{\alpha_1 k}}$. Thus,
\begin{eqnarray*}
d(\varphi^p(z), \varphi^p(z_0)) 
& =     & d(\phi_p(z), \phi_p(z_0) )\\
& \leq & d(\phi_p(z), \tpp^{m-k}(z)) + d(\tpp^{m-k}(z), \tpp^{m-k}(z_0)) + d(\tpp^{m-k}(z_0), \phi_p(z_0))\\
& \leq & |p- t | + \frac{C}{\gamma^{\alpha_1 k}} + |T - p| \leq \frac{3 C}{\gamma^{\alpha_1 k}} .
\end{eqnarray*}

\ms

\noindent
{\bf Case 2.} $p \leq t$. First, assume that $t \leq T$. Then, using (4.13),
\begin{eqnarray*}
d(\varphi^p(z), \varphi^p(z_0)) 
& =     & d(\phi_p(z), \phi_p(z_0)) \leq d (\phi_t(z), \phi_t(z_0))
 \leq  d(\phi_t(z), \phi_T(z_0)) + d(\phi_T(z_0), \phi_t(z_0))\\
& =     & d(\tpp^{m-k}(z), \tpp^{m-k}(z_0)) + |T - t| \leq \frac{2C}{\gamma^{\alpha_1 k}} .
\end{eqnarray*}
The other case to consider is $ t > T$. Then $t > T \geq p$, and as above we get
\begin{eqnarray*}
d(\varphi^p(z), \varphi^p(z_0)) 
& =     & d(\phi_p(z), \phi_p(z_0)) \leq d (\phi_T(z), \phi_T(z_0))
 \leq  d(\phi_t(z), \phi_T(z_0)) + d(\phi_t(z), \phi_T(z))\\
& =     & d(\tpp^{m-k}(z), \tpp^{m-k}(z_0)) + |T - t| \leq \frac{2C}{\gamma^{\alpha_1 k}} .
\end{eqnarray*}
This proves (4.14). It implies that
\be
\diam(\varphi^p(\tcc)) \leq \frac{3C}{\gamma^{\alpha_1 k}} .
\ee

For $z_p = \varphi^{p}(z_0)$ we have 
$r_p = r(z_p) \geq r_0\, e^{- p \hep}$. We need to have 
\be
r(z_p)  > \diam(\varphi^p(\tcc))\quad \mbox{\rm and }\quad r(z_p) \geq \diam(\hvarphi^p(\hcc)) .
\ee
By (4.15) and (4.1), for this it would be enough to have
$\di 3C \ttc_2 \, e^{-\alpha_1 k\, \log \gamma} \leq  r_0 e^{-p \hep }$,
that is $-( \alpha_1 \log \gamma) k \leq \log(r_0/(3C \ttc_2)) - p \hep$ which is equivalent to
$\di p \hep \leq \log \frac{r_0}{3C \ttc_2} + (\alpha_1 \log \gamma) k$.
Since $p = [\ttau_{m-k}(z_0)] \leq (m-k) \tau_0$, we have
$p \hep \leq (m-k) \tau_0\hep$, and so the above holds if
$$(m-k) \tau_0 \hep \leq (\alpha_1 \log \gamma) k - D_0 ,$$
with $D_0 = |\log \frac{r_0}{3C \ttc_2}|$, that is if
\be
m \tau_0 \hep \leq  (\tau_0 \hep + \alpha_1  \log \gamma) k - D_0 .
\ee

We will now assume that
\be\label{eq:k-cond}
m \hep_2 \leq k \leq m \hep_1,
\ee
where $\hep_2 = \frac{\hep \tau_0}{\alpha_1 \log \gamma}$ is a small number 
(can be made arbitrarily small choosing the initial $\hep$ small), and (there is a lot freedom in this choice) e.g.
$\hep_1 = 2 \hep_2$. Then for $m \geq m_0$ sufficiently large, $m \hep_2 \leq k$ implies
$$m \tau_0 \hep  \leq (\alpha_1 \log \gamma) k  < (\tau_0 \hep + \alpha_1 \log \gamma) k - D_0 ,$$
assuming $k > 1$ is sufficiently large,
so (4.17) holds and therefore (4.16) holds as well. The latter yields
\be
\varphi^p(\tcc) \subset B^u(z_p, r(z_p)) \quad \mbox{\rm and} \quad \hvarphi^p_{z_0}(\hcc) \subset B^u(z_p, r(z_p)) .
\ee
This completes the proofs of parts (a) and (b).

\ms
\noindent
{\it Proof of part} (c). We will continue to use the notation introduced above.

Since $p = [\ttau_{m-k}(z_0)]$, we have $\ttau_{m-k}(z_0) = p +t$ for some $t \in [0,1)$.
For $y' = \tpp^{-k}(y_0)$, by (2.1) and (4.12), it follows that
$$B^u(y', d_0r_0/\gamma_1^k) \subset \tpp^{-k}(B^u(y_0,r_0) ) \subset \tpp^{m-k}(\tcc) = \phi_{\ttau_{m-k}(z_0)} (\tcc) .$$
Setting $y'' = \phi_{-t}(y')$ and $\hy = \T_{z_0}(y'')$ and using (4.1), we have
$$B^u(y'', d_0r_0/\gamma_1^{k+1}) \subset \phi_{\ttau_{m-k}(z_0)-t} (\tcc) = \varphi^p(\tcc) 
\quad \mbox{\rm and} \quad B^u(\hy, d_0r_0 \ttc_1/\gamma_1^{k+1}) \subset \hvarphi_{z_0}^p(\hcc) .$$
By (4.17),  
$$d_0r_0 \ttc_1/\gamma_1^{k+1} = \frac{d_0 r_0 \ttc_1}{\gamma_1} e^{-k \log \gamma_1} 
\geq \frac{d_0 r_0 \ttc_1}{\gamma_1} e^{-m \hep_1 \, \log \gamma_1} = c'_3\, e^{-m\hep_3} ,$$
for some global constant $c'_3 > 0$ and $\hep_3 = (\log\gamma_1)\, \hep_1 > 0$. 
Hence, taking into account (4.19) as well, we obtain
\be
B^u(\hy, c'_3 \, e^{-m \hep_3}) \subset \hvarphi_{z_0}^p(\hcc) \subset B^u(z_p, r(z_p)) .
\ee
and
\be
B^u(y'', c''_3 \, e^{-m \hep_3}) \subset \varphi^p(\tcc) \subset B^u(z_p, r(z_p)) 
\ee
with $c''_3 = c'_3/ \ttc_1$.
Then for every $b'\in B^u (\hy, c'_3 \, e^{-m\hep_3})$ there exists $b \in \hcc$ with $\hvarphi_{z_0}^{p}(b) = b'$. 
Notice also that $r = c'_3 e^{-m\hep_3} \leq r(z_p)$.

To apply the map $(\Phi^u_{z_p})^{-1}$ to (4.20) we will use (3.9) and
$R(z_p) \leq R(z_0) e^{\hep\, p} \leq R_0 e^{\hep\, (m-k) \tau_0} \leq R_0 e^{\hep \, m \,\tau_0}$.
Setting $\xi = (\Phi^u_{z_p})^{-1}(\hy)$, $c_3 = c'_3/R_0$ and $\hep_4 = \hep_3 + \hep \tau_0$, 
it follows from (4.20) that
$$(\Phi^u_{z_p})^{-1} (B^u(\hy, c_3 \, e^{-m\hep_4}) ) \subset B(\xi, R(z_p) c_3 \, e^{-m \hep_4}) \subset B(\xi, r)
\subset E^u(z_p). $$
An elementary argument shows that in the normed
space $E^u(z_p)$ we can always find an element $\hxi$ in the ball $B(\xi,r)$ in $E^u(z_p)$ such that $\|\hxio\| \geq r/2$. Indeed, assume
e.g. $\xio \geq 0$ in the natural coordinates in $E^u(z_p)$. If $\xio \geq r/2$ just take $\hxi = \xi$. If $0 \leq \xio < r/2$, take
$\hxi = (\xio + r/2, \xi^{(2)}, \ldots, \xi^{(n_u)})$, where $n_u = \dim(E^u(z_p))$. 
We have
$$v_0 = \Phi^u_{z_p}(\hxi) \in \Phi^u_{z_p}(B(\xi,r)) \subset B^u(z_p, r(z_p)) \cap \hvarphi_{z_0}^p(\hcc) .$$
Set $u_0 = \hvarphi^{-p}_{z_{p}}(v_0) \in E^u(z_0)$. For $v_0$ the above gives $\|\vo_0\|_{z_p} \geq \|\hxio\| \geq r/2$, while
$\hx_0 = \Phi^u_{z_0}(u_0) \in \hcc$. Hence
$\di \diam(\hcc) \geq d(z_0, \hx_0) \geq \frac{\|u_0\|}{R_0} \geq \frac{ \|\uo_0\|'_{z_0}}{\Gamma_0 R_0}$, and
by (4.1), 
$$\di \diam(\tcc) \geq \frac{\ttc_1 }{\Gamma_0 R_0} \, \|\uo_0\|'_{z_0} .$$
Notice also that for $\hx_0 = \Phi^u_{z_0}(u_0)$ 
we have $\hvarphi^p_{z_0}(u_0) = v_0 \in \hvarphi_{z_0}(\hcc)$.
Thus, for $x_0 = (\T_{z_0})^{-1}(\hx_0) \in \tcc$ we have $\varphi^p(x_0) \in \varphi^p(\tcc)$. 

Recall that $p = [\ttau_{m-k}(z_0)]$.
It follows from $q = [\ttau_m(z_0)]$ that $m\, \ttau_0 - 1\leq q \leq m\tau_0$, while (4.18) yields
$$q \geq p + [\ttau_{k}(\tpp^p(z_0))]  \geq p + k \ttau_0 -1  \geq p + m \ttau_0 \hep_2 - 1 \geq p  + q \frac{\ttau_0 \hep_2}{\tau_0}  - 1.$$
Similarly, using (4.18) again, we get $q \leq p  + q \frac{\tau_0 \hep_1}{\ttau_0} + \frac{\tau_0 \hep_1}{\ttau_0} +1$,
therefore 
$$q(1 - \frac{\tau_0 \hep_1}{\ttau_0}) - 2 \leq p \leq q(1 - \frac{\ttau_0 \hep_2}{\tau_0}) +1 .$$
Thus, for sufficiently large $m$ (then $q$ is also large) we have
\be\label{q-cond}
q (1- \hep_6)  \leq p \leq  q (1- \hep_5) \leq q ,
\ee
for some small constants $0 < \hep_5  \leq \con \hep$ and  $0 < \hep_6 \leq \con \hep$.

Now $\nu_1 = \lambda_1 e^{\hep}$ and $p \leq q$ give 
\be
\|\uo_0\|'_{z_0} \geq \frac{\|\vo_0\|'_{z_p}}{\nu_1^{p}}
 \geq  \frac{r}{2 \nu_1^{p}}  
  \geq  \frac{c'_3\, e^{-m \hep_3}}{2 \nu_1^{q}}  
    \geq  \frac{c'_3 \, e^{-(q+1) \hep_3/\ttau_0 }}{2  \nu_1^{q}} 
     \geq  \frac{c'_3 \, e^{- (q+1)\hep_3/\ttau_0 - q \hep}}{2 \lambda_1^{q}} 
  \geq \frac{c'_3 \, e^{-q \hep_7 - \hep_3/\ttau_0}}{2 \lambda_1^{q} } ,
\ee
for some $0 < \hep_7 \leq \con \, \hep$.
Hence $\diam(\tcc) \geq  \frac{e^{-q \hep_7}}{C_3 \lambda_1^{q} }$, taking $C_3 \geq  2 R_0\Gamma_0/(\ttc_1 c'_3 e^{1/\ttau_0}) > 0$. 

This proves the left-hand-side inequality in (\ref{eq:diamC}).

We will now prove in a similar way the other inequality in (\ref{eq:diamC}).

It follows from (4.21) and (4.19) and the choice of $k$ and $p$ that
$\diam(\varphi^p(\tcc)) \leq r_0 e^{-p \hep} < r(z_p) .$
By (4.19) there exists $\xi \in E^u(z_p)$ with $\|\xi\| \leq r_0 e^{-p \hep}  < r(z_p)$ so that
$z = \Phi^u_{z_p} (\xi) \in \varphi^p(\tcc)$ and for $y = \varphi^{-p}(z) \in \tcc$ we have $d(y,z_0) \geq \frac{1}{2} \diam(\tcc)$.
Thus for $\eta = (\Phi^u_{z_0})^{-1}(y) \in E^u(z_0)$ we have 
$$\|\eta\|'_{z_0} \geq \|\eta\| \geq \frac{1}{\Gamma_0} d(y,z_0)  \geq \frac{1}{2 \Gamma_0} \diam(\tcc) .$$
Using (3.12), $R(z_p) \leq R_0 e^{p \hep}$, $\mu_1 = \lambda_1 e^{-\hep}$, $q(1-\hep_6) \leq p \leq q$, $r_0 < 1$   and the above we obtain
\begin{eqnarray*}
\diam(\tcc) 
& \leq &  2 \Gamma_0 \|\eta\|'_{z_0} \leq 2 \Gamma_0 \frac{ \|\xi\|'_{z_p}}{\mu_1^p} 
\leq 2 \Gamma_0 R_0 e^{p \hep} \frac{ \|\xi\|}{\mu_1^p} \leq 2 \Gamma_0 R_0 \frac{ r_0 e^{2p \hep}}{\lambda_1^p} \nonumber \\
& \leq &
\frac{2 \Gamma_0 R_0 \lambda_1^{q \hep_6} e^{2q \hep}}{\lambda_1^{q}} 
\leq \frac{2 \Gamma_0 R_0  e^{q (\hep_6 \log \lambda_1 + 2\hep)}}{\lambda_1^{q}} 
\leq  \frac{C_3 e^{q \hep_7}}{\lambda_1^{q}} ,
\end{eqnarray*}
assuming $\hep_7 \geq 2\hep + \hep_6 \log \lambda_1$ and $C_3 \geq 2 \Gamma_0 R_0$.
This proves the left-hand-side inequality in (\ref{eq:diamC}). 

The existence of $u_0$ with (4.9) follows from (4.23).
\endofproof

\def\tn{\tilde{n}}
\def\tGamma{\widetilde{\Gamma}}
\def\tLambda{\widetilde{\Lambda}}
\def\tX{\widetilde{X}}
\def\tY{\widetilde{Y}}
\def\tW{\widetilde{W}}
\def\tOmega{\widetilde{\Omega}}
\def\tw{\tilde{w}}
\def\Fp{F^{(p)}}
\def\F0{F^{(0)}}
\def\tTheta{\widetilde{\Theta}}

\subsection{Constructing two families of  sub-cylinders}

We will now get important consequences of some of the arguments in the proof of Lemma 4.3. 


In the next lemma we will use the following notation. Given a cylinder $\tcc$ in $\tR$ containing a 
point $z_0 \in \tcc \cap \tP_0$ and an integer $p \geq 0$, consider the  locally defined near $z_0$ map
$$(\F0)^{-1}\circ \hvarphi_{z_0}^{-p} \circ \Fp : \tcc \longrightarrow \tcc ,$$ 
where
$$\Fp = (\Phi^u_{z_p})^{-1} \circ \T_{z_p} \circ \varphi^p : \tcc \longrightarrow E^u (z_p) \quad, \quad
\F0 = (\Phi^u_{z_0})^{-1} \circ \T_{z_0} : \tcc \longrightarrow E^u (z_0) .$$
See Figure 1 below.

\bs


\noindent
{\bf Lemma 4.5.}  {\it Let $C > 1$ and $d > 1$ be given constants and let $\frac{1}{\gamma^{\alpha_1 \beta}} \leq \theta < 1$.}

{\it There exist global constants $q_0, q_1, q_2 \in \N$, $\tm_0 \in \N$ and $\hep_{12} > 0$,
$0 < \hep_{12} \leq \con \hep$,  which can be made arbitrarily small with $\hep$, such that $q_0 \leq q_1 \leq q_2$ and
 for every cylinder $\tcc$ in $\tR$ of length $m \geq \tm_0$ containing a point $z_0 \in \tcc \cap \tP_0$, setting 
$$ p = [\ttau_{m-k}(z_0)]  \quad \mbox{with}\quad k = \frac{4 d (m-k)  \hep}{|\log \theta|},$$
we have the following:}

\ms

(a) {\it  There exists a
sub-cylinder $\tGamma = \tGamma (\tcc)$ of $\tcc$ of co-length $q_0$ in $\tcc$ 
so that for every $x \in \tGamma$, for $u = \F0 (x) \in E^u(z_0)$ we have
\be
\|\uo\| \geq \kappa\, \diam(\tcc) 
\ee
with $\kappa = d_2 e^{- m \hep_{12} }$.
Moreover,  there exists a sub-cylinder $\tLambda$ of $\tcc$ of co-length $q_1$ in $\tcc$ 
such that
for every $x \in \tLambda$
we have
$\di C \hvarphi_{z_0}^{p} (\F0(x)) \in \Fp(\tGamma)$. That is,}
\be\label{eq:in-cond}
C\, \hvarphi_{z_0}^{p} (\F0 (\tLambda)) \subset \Fp (\tGamma) .
\ee

\ms

(b) {\it For the sub-cylinders $\tGamma$ and  $\tLambda$ of $\tcc$ as in part (a), there exists a sub-cylinder $\tTheta$
of $\tGamma$ with co-length $q_2$ in $\tcc$ such that }
\be\label{eq:out-cond}
 \Fp (\tTheta) \subset C\, \hvarphi_{z_0}^{p} (\F0 (\tLambda)) .
\ee

\bs

Notice that, since the sub-cylinders $\tGamma$ of $\tcc$ that appear in Lemma 4.5
have co-lengths bounded by global constants, there exists a global constant $d_1 > 0$ such that
\be
\nu(\tGamma') \geq d_1 \, \nu(\tcc') ,
\ee
where $\tGamma' = \piU(\tGamma)$, $\tcc' = \piU(\tcc)$. A similar estimate applies
to $\nu(\tTheta)$ and $\nu(\tLambda)$.

\bs

\noindent
{\it Proof.} 
We will use the notation and the assumptions from the beginning of the proof of Lemma 4.4,
with some small changes. Set $Z_m = \tpp^m(z_0)$, and more generally, $Z_j = \tpp^j(z_0)$, $1 \leq j \leq m$.
As before, let $r_0 > 0$ be a small constant and let $y_0$ be a point
with $B^u(y_0, 4r_0) \subset W^u_{\tR_{j_0}}(Z_m)$ and now {\bf we assume $d(y_0, Z_m) > 4r_0$}. Let $y \in \tcc$
be the point with $\tpp^m(y) = y_0$. 

Again we will assume that the integer  $k$ satisfies (4.11), however now we will impose a stronger condition on $k$. 
Namely  we want something similar to (4.17):
\be\label{eq:m-k-cond}
 |\log \theta|\, k =  4 d   (m-k) \hep .
\ee
That is  $(4 d \hep + |\log \theta|) k = 4  d m \hep$, which implies
\be
m \hep_8 = k \leq m \hep_9 ,
\ee
where $\di \hep_9 = \frac{4  d}{|\log\theta|}\, \hep \leq \con \hep$ and
$\di \hep_8 = \frac{4  d}{ 4  d \hep + |\log\theta|}\, \hep \leq \con \hep$.

Notice that, since $1/\gamma^{\alpha_1 \beta} \leq \theta < 1$,
we have $\alpha_1 \beta \log \gamma \geq |\log \theta|$, so (\ref{eq:m-k-cond}) implies 
$$(\alpha_1\beta \log \gamma) k \geq  4 d   (m-k) \hep .$$

As in the proof of Lemma 4.4, set $p = [\ttau_{m-k}(z_0)]$. Notice that with our new definition of $k$
and the related $p$, (4.13), (4.14) and (4.15) in the proof of Lemma 4.4 still hold, and therefore
(4.16) hold as well. Moreover, 
$$\diam(\varphi^p(\tcc)) \leq \frac{3C}{\gamma^{\alpha_1 k}}  \leq \frac{3C}{\gamma^{\alpha_1 \beta k}} .$$
Also, as in the proof of Lemma 4.4 we show that 
\be\label{eq:p,q-cond}
q (1 - \hep'_6) \leq p \leq q (1 - \hep'_5) \leq q
\ee
for some constants $0 < \hep'_5 \leq \con \hep$ and $0 < \hep'_6 \leq \con \hep$, which is the analogue
of (4.22) for our new definition of $p$.

Next, set  $t = \ttau_{m-k}(z_0) - p\in [0,1)$, $Z_{m-k} = \tpp^{-k}(Z_m) = \tpp^{m-k}(z_0)$,
$z_p = \varphi^p(z_0) = \phi_{-t}(Z_{m-k})$, $y' = \tpp^{-k}(y_0)$, $y'' = \phi_{-t}(y')$. 
Consider the map
$$F = \tpp^{k}\circ \phi_{t} : \varphi^p(\tcc) \longrightarrow W^u_{\tR_{j_0}}(Z_m) .$$
As before, by (2.1), we have
$$\di B^u(y', \frac{4r_0d_0}{\gamma_1^k}) \subset \tpp^{-k}(B^u(y_0, 4r_0))$$ 
and
$$\di B^u(y'', \frac{4r_0d_0}{\gamma_1^{k+1}}) \subset F^{-1}(B^u(y_0, 4r_0)) .$$
Set
$\hy = \T_{z_p}(y'') \in W^u_{\ep}(z_p)$ and $\eta = (\Phi^u_{z_p})^{-1}(\hy) \in E^u(z_p)$.
Since $d(y_0,Z_m) > 4r_0$, we have $d(y'', z_p) \geq \frac{4r_0d_0}{\gamma_1^{k+1}}$, and it follows from (4.1) that
$d(\hy, z_p) \geq \frac{4r_0d_0}{\tct\gamma_1^{k+1}}$ and therefore
$$\di \|\eta\| \geq \frac{d(\hy, z_p)}{R(z_p)} \geq \frac{4r_0d_0}{\ttc_2 \gamma_1^{k+1} R(z_p)} .$$


\begin{center}
\begin{picture}(200,180)       
\put(10,0){$\tcc$}
\put(30,4){\vector(1,0){40}}  
\put(50,10){$\varphi^p$}  
\put(80,0){$\varphi^p(\tcc)$} 
\put(120,5){\vector(1,0){35}} 
\put(160,0){$\tpp^{m-k}(\tcc)$}
\put(130,10){$\phi_t$}
\put(130,90){$\phi_t$}
\put(0,160){$E^u(z_0)$}
\put(80,160){$E^u(z_p)$}
\put(160,160){$E^u(Z_{m-k})$}

\put(0,80){$W^u_{\ep}(z_0)$}
\put(40,83){\vector(1,0){35}}  
\put(40,163){\vector(1,0){35}}  
\put(120,163){\vector(1,0){35}}  
\put(50,90){$\varphi^p$}  
\put(50,170){$\hvarphi^p_{z_0}$}  

\put(80,80){$W^u_{\ep}(z_p)$} 
\put(120,83){\vector(1,0){35}} 
\put(160,80){$W^u_{\ep}(Z_{m-k})$}
\put(130,90){$\phi_t$}
\put(130,170){$\hat{\phi}_t$}
\put(15,15){\vector(0,1){60}}
\put(15,95){\vector(0,1){60}}
\put(95,95){\vector(0,1){60}}
\put(175,95){\vector(0,1){60}}
\put(20,45){$\T_{z_0}$}
\put(20,125){$(\Phi^u_{z_0})^{-1}$}
\put(100,125){$(\Phi^u_{z_p})^{-1}$}
\put(180,125){$(\Phi^u_{Z_{m-k}})^{-1}$}
\put(95,15){\vector(0,1){60}}
\put(102,45){$\T_{z_p}$}
\put(175,15){\vector(0,1){60}}
\put(182,45){$\T_{Z_{m-k}}$} 
\end{picture}
\end{center}


\begin{center}
Figure 1
\end{center}

\bs

In what follows
it is more convenient to deal with $\tpp^{m-k} (\tcc)$ instead of $\varphi^p(\tcc)$, and the space $E^u(Z_{m-k})$, rather then $E^u(z_p)$.
As before, for $z_p = \varphi^p(z_0)$ we have $r(z_p) \geq r_0 e^{-p \hep}$.
It follows from (4.15), which holds again with the present choice of $k$ as remarked earlier, that we have
$$\diam(\varphi^p(\tcc)) \leq \frac{3C }{\gamma^{\alpha_1  k}}
\leq \frac{3C }{\gamma^{\alpha_1 \beta k}} = 3C e^{-k\, \alpha_1 \beta \log \gamma} ,$$
therefore
$$\diam(\tpp^{m-k}(\tcc)) \leq \frac{3C \gamma_1}{d_0} \, e^{-k\, \alpha_1 \beta \log \gamma} ,$$
which, combined with $(\alpha_1 \beta \log \gamma) k \geq  4 d   (m-k) \hep$, implies
\begin{eqnarray*}
 \diam((\Phi^u_{Z_{m-k}})^{-1}(\T_{Z_{m-k}}(\tpp^{m-k}(\tcc))) 
& \leq & \frac{3C \gamma_1 \ttc_2 R(Z_{m-k})}{d_0} \, e^{-k\, \alpha_1\beta \log \gamma} \nonumber\\
& \leq & \frac{3C \gamma_1 \ttc_2 R_0}{d_0} \,e^{(m-k) \hep}\, e^{-  4 d   (m-k) \hep} 
 \leq   e^{-  3 d (m-k) \hep}   ,
\end{eqnarray*}
since  $d> 1$, so
$\di \frac{3C \gamma_1 \ttc_2 R_0}{d_0} \,e^{(m-k) \hep} \leq e^{ d  (m-k) \hep}$, assuming $k$ is sufficiently large.
Also,  we have $e^{- d  (m-k) \hep} \leq r_0 e^{ - (m-k)\hep} \leq r(Z_{m-k})$, again for large $k$.
The above gives
\be
 (\Phi^u_{Z_{m-k}})^{-1}(\T_{Z_{m-k}}(\tpp^{m-k}(\tcc)) \subset
E = \left\{ u \in E^u(Z_{m-k}) : \|u\| \leq e^{- 3 d   (m-k) \hep} \right\} .
\ee

Let $e_1,e_2, \ldots, e_{n_u}$ be an orthonormal basis in $E^u(Z_{m-k})$, so that $\{e_1, \ldots, e_{n_1}\}$ is a
basis in $E^u_1(Z_{m-k})$. 
Consider the subsets
$$H = \{ u = (u_1, \ldots, u_{n_u}) \in E^u(Z_{m-k}) : \|u\| < e^{- 3 d  (m-k) \hep} \:, \: u_1 = 0  \} $$
and
$$H_1 = \left\{ u = (u_1, \ldots, u_{n_u}) \in E^u(Z_{m-k}) : \|u\| < e^{- 3 d  (m-k) \hep} \:, \: 
\|u_1\| < \frac{1}{4} e^{- 4 d  (m-k) \hep}  \right\} $$
of $E$. Clearly $H \subset H_1$. We will also need $\hH_1 = (\hphi_t)^{-1}(H_1) \subset E^u (z_p) .$


Set
\be\label{eq:ep}
\ep = \frac{e^{-4d (m-k) \hep}}{4} .
\ee
Consider a rectangular box $\Delta$ contained entirely  in $E\setminus H_1$ with sides of length $\ep$ parallel to the coordinate axes.
To construct $\Delta$, {\bf fix a point $w \in H$ with }
$$\|w_1\| = \frac{1}{2} e^{-4d (m-k) \hep} ,$$
and then {\bf consider the cube $\Delta$ in $E$ with centre $w$ and sides of length $\ep$. }
Then, applying the map $(\hphi_t)^{-1}$ to the cube $\Delta$, we get that its 
image is contained in $E(z_p) \setminus H_2$, where
$$H_2 = \left\{ v = (v_1, \ldots, v_{n_u}) \in E^u(z_p) : \|v\| <  e^{- 3  d (m-k) \hep} \:, \: 
\|v_1\| < \frac{d_0 \ttc_1}{ 4 \gamma_1} e^{- 4 d  (m-k) \hep}  \right\} .$$

Consider the largest cylinder $\tX$ in $\tcc_1 = \tpp^{m-k}(\tcc)$
such that $\Delta \supset (\Phi^u_{Z_{m-k}})^{-1}\circ \T_{Z_{m-k}} (\tX)$.
Denote by $n$ the length of the cylinder $\tX$. 
Then $\tX \subset \T^{-1}_{Z_{m-k}}(\Phi^u_{Z_{m-k}}(\Delta))$, however
the cylinder $\tY$ of length $n-1$ containing $\tX$ is not contained in
$\T^{-1}_{Z_{m-k}}(\Phi^u_{Z_{m-k}}(\Delta)))$, and using (3.9) we get
$$\diam(\tY)  \geq \diam(\T^{-1}_{Z_{m-k}}(\Phi^u_{Z_{m-k}}(\Delta))) \geq \frac{\ttc_1\, \ep}{R_0} .$$
By Lemma 4.1(c),
$\diam(\tY) \leq C_1 \diamte(\tY) = C_1 \theta^{n-1} .$
Hence $C_1 \theta^{n-1} \geq \ttc_1\, \ep/R_0$. So, we have
$\theta^{n} \geq \con \ep $ for some global constant $\con > 0$. This gives
$n \log \theta \geq \con + \log \ep$, hence  $n |\log \theta| \leq c' - \log \ep $, and therefore
$$n \leq \frac{c' -  \log \ep}{|\log\theta|} = c'' - \frac{\log \ep}{|\log \theta|}$$
for some global constants $c', c'' \in \R$.


Consider an arbitrary point $x' \in \tX' = \piU(\tX)$. Since the length of $\tcc_1 = \tpp^{m-k}(\tcc)$ is $k$,
it follows from (4.2) and (4.3) that
$$\frac{\nu(\tX')}{\nu(\tcc_1')} \geq \frac{c_1 e^{g_{n}(x')}}{c_2 e^{g_{k}(x')} } = \frac{c_1}{c_2} e^{g_{n-k}(\sigma^k(x'))}
\geq \frac{c_1^2}{c_2} \rho_1^{n-k}
= \frac{c_1^2}{c_2} e^{(n-k) \log \rho_1} = \frac{c_1^2}{c_2} \left( e^{- n} e^{k}\right)^{ |\log \rho_1|} . $$
From the estimate above,
$\di e^{-n} \geq e^{-c''} \, e^{\frac{\log \ep}{|\log\theta|}} = c'''\,  \ep^{1/|\log\theta| } $
for some global constant $c''' > 0$.
On the other hand, it follows from (\ref{eq:m-k-cond}) that
$\di e^k = e^{4d (m-k) \hep/|\log \theta|}$.
Therefore, using (\ref{eq:ep}),
$$\frac{\nu(\tX')}{\nu(\tcc_1')} \geq \frac{c_1^2}{c_2} \left(c''' \, \ep\, e^{4d(m-k) \hep}\right)^{\frac{|\log \rho_1|}{|\log\theta|}} 
= \frac{c_1^2}{c_2} \left(c'''/4\right)^{\frac{|\log \rho_1|}{|\log\theta|}} = d'_1 .$$
Thus, there exists a global constant $d'_1 > 0$ so that $\nu(\tX') \geq d'_1\, \nu(\tcc_1')$.

In this way we have constructed a cylinder $\tX$ in $\tcc_1$ so that
$(\Phi^u_{Z_{m-k}})^{-1}\circ \T_{Z_{m-k}} (\tX) \subset E \setminus H_1$,
which implies
\be
(\hphi_t)^{-1}\circ (\Phi^u_{Z_{m-k}})^{-1}\circ \T_{Z_{m-k}} (\tX) \subset E^u(z_p) \setminus H_2 ,
\ee
and $\nu(\tX') \geq d'_1\, \nu(\tcc_1')$.

Let $\tGamma$ be the sub-cylinder of $\tcc$ such that $\tpp^{m-k}(\tGamma) = \tX$.
Set $\tGamma' = \piU(\tGamma) \subset U$.

To estimate $\nu(\tGamma')/\nu(\tcc')$ we will use (4.2). Consider an arbitrary $z'\in \tGamma'$.
Since $n$ is the length of the cylinder $X$, the length of $\Gamma$ is $n + m-k$. It follows from (4.2) that
$$\frac{\nu(\tGamma')}{\nu(\tcc')} \geq \frac{c_1 e^{g_{n+ m-k}(z')}}{c_2 e^{g_{m}(z')} } = \frac{c_1}{c_2} e^{g_{n-k}(\sigma^m(z'))} .$$
Then exactly as in the estimate for $\di \frac{\nu(\tX)}{\nu(\tcc_1)}$ we obtain 
\be
\nu(\tGamma') \geq d'_1\, \nu(\tcc') .
\ee
It follows from (4.33) and $\tpp^{m-k}(\tGamma) = \tX$, that
$$(\Phi^u_{z_p})^{-1} ( \T_{z_p} (\varphi^p(\tGamma))) = (\hphi_t)^{-1}\circ 
(\Phi^u_{Z_{m-k}})^{-1}\circ \T_{Z_{m-k}} (\tpp^{m-k}(\tGamma)) ,$$
hence
\be\label{eq:Gamma}
(\Phi^u_{z_p})^{-1} ( \T_{z_p} (\varphi^p(\tGamma))) \subset \hH_1 \setminus H_2 .
\ee

\ms

{\it Proof of } (4.24):
Consider the map
$$\hat{\phi}_t = (\Phi^u_{Z_{m-k}})^{-1} \circ \phi_t \circ \Phi^u_{z_p} : E^u(z_p, r(z_p)) \longrightarrow E^u(Z_{m-k}, r(Z_{m-k})) .$$
Given $v \in \hvarphi_{z_0}^{m-k}(\tGamma)$, it follows from (\ref{eq:Gamma}) and the definition of $H_2$ that
$$ \|\vo\| \geq \|v_1\| \geq \frac{d_0 \ttc_1}{4 \gamma_1} e^{- 4d   (m-k) \hep} .$$
Set $\xi = \hat{\phi}_{t}^{-1} \cdot v\in E^u(z_p)$. Then 
$ \|\xio\| \geq  \frac{d_0 \ttc_1}{4 \gamma_1} e^{- 4 d  (m-k) \hep} $.
Consequently for $u = \hvarphi_{z_p}^{-p}(\xi)$, using
$\nu_1^p = \lambda_1^p e^{p \hep}$ and $(m-k) \ttau_0 \leq p  \leq q $ 
as in the proof of Lemma 4.4, it follows from (3.12) that
$$\|\uo\|'_{z_0} \geq \frac{\|\xio\|'_{z_p}}{\nu_1^{p}} \geq \frac{\|\xio\|}{\lambda_1^{p} e^{2 p\hep} R_0} 
  \geq  \frac{d_0 \ttc_1 e^{- \frac{4 d q \hep}{\ttau_0}}}{4 \gamma_1 R_0 e^{2q\hep}\lambda_1^{q }}  
    \geq  \frac{d_0 \ttc_1 e^{- q \hep_{11} }}{4 \gamma_1  R_0 \lambda_1^{q}} ,$$
where $\hep_{11} = \frac{4 d \hep}{ \ttau_0} + 2 \hep  \leq \con \, \hep$.

On the other hand, (\ref{eq:diamC}) gives
$$\diam(\tcc) \leq  \frac{C_3 \, e^{q \hep_7}}{\lambda_1^{q}} 
 \leq \frac{e^{-q \hep_{11}}}{e^{-q \hep_{11}}}\, \frac{C_3 e^{q\hep_7}}{\lambda_1^{q}} 
\leq \|\uo_1\|'_{z_0} \, \frac{4  \gamma_1 R_0 C_3}{d_0\ttc_1} \,  e^{q(\hep_7 + \hep_{11})} .$$
Thus, for every $x \in \tGamma'$, for $u = (\Phi^u_{z_0})^{-1}(\T_{z_0}(x)) \in E^u(z_0)$
we have $\|\uo\| \geq \frac{\|\uo\|'_{z_0}}{R_0} \geq \kappa\, \diam(\tcc)$,
where $\kappa = d_2 e^{-q \hep_{12}}$ for some $0 < \hep_{12} \leq \con \, \hep$,
where $d_2 = \frac{d_0 \ttc_1}{4 \gamma_1 R_0^2 C_3}$.
This proves (4.24).

\ms


We will now repeat the above for a different subset of $\tcc$.
Choose a sufficiently large integer $\tq \geq 1$ (a global constant determined by $C$) 
so that the cylinder $\tdd_2$ in $W^u_{\tR_{j_0}}(Z_m)$ of length
$\tq$ containing the point $Z_m$ is so that $\tdd_2 \subset B^u(Z_m, r_0/C)$
for the given constant $C > 1$. Let $\tdd$ be the sub-cylinder
of $\tcc$  of co-length $\tq$ so that $\tpp^m(\tdd) = \tdd_2$. Then for the map $F$ defined by
$F = \tpp^k\circ \phi_t : W^u_\ep (z_p) \longrightarrow W^u_{\tR_{j_0}}(Z_m)$ we have
$$\varphi^p(\tdd) = F^{-1}(\tdd_2) \subset F^{-1}(B^u(Z_m, r_0/C)) .$$

Set
$$\tH = \frac{1}{C} H = \left\{ u = (u_1, \ldots, u_{n_u}) \in E^u(Z_{m-k}) : \|u\| < \frac{1}{C} e^{- 3 d (m-k) \hep } \: , \: u_1 = 0 \right\} ,$$
$$\tH_1 = \frac{1}{C} H_1 = \left\{ u = (u_1, \ldots, u_{n_u}) \in E^u(Z_{m-k}) : \|u\| < \frac{1}{C} e^{- 3 d (m-k) \hep} \: , 
\: \|u_1\| < \frac{1}{4 C} e^{- 4 d (m-k) \hep} \right\} ,$$
$$\tH_2 =  \frac{1}{C} H_2 = 
\left\{ v = (v_1, \ldots, v_{n_u}) \in E^u(z_p) : \|v\| < \frac{1}{C } e^{- 3 d (m-k)\hep } \:, \: 
\|v_1\| < \frac{d_0 \ttc_1}{4\gamma_1 C } e^{- 4 d (m-k)\hep}  \right\} ,$$
$$H_1' = (\hphi_t)^{-1}(\tH_1) =  (\hphi_t)^{-1}(\frac{1}{C} H_1) = \frac{1}{C} \hH_1 \subset E^u(z_p) ,$$
$$\ep' = \frac{e^{-4d (m-k) \hep}}{4C} = \frac{\ep}{C} .$$

Recall the fixed point $w \in H$. 
{\bf Consider now the point $\di \tw = \frac{w}{C} \in \tH$.} Then
$$\|\tw_1\| = \frac{1}{2C} e^{-4d (m-k) \hep} .$$
Let $\Delta'$ be the cube in $E$ with centre $\tw$ and sides of length $\ep'$ so that
$(\hphi_t)^{-1}(\Delta') \subset E(z_p) \setminus \tH_2$.
Consider the largest cylinder $\tY$ in $\tdd_1 = \tpp^{m-k}(\tdd)$
such that $\Delta' \supset (\Phi^u_{Z_{m-k}})^{-1}\circ \T_{Z_{m-k}} (\tY)$.
If $\tn$ is the {\it length of the cylinder} $\tY$, as in the previous case we get
$$\tn \leq \frac{c' - \log \ep'}{|\log\theta|} = c'' - \frac{\log \ep'}{|\log \theta|}$$
for some global constants $c', c'' \in \R$. 

Since $\tdd_1 = \tpp^{m-k}(\tdd)$ is a sub-cylinder of $\tcc_1$ of length $k+ \tq$ and $\tY$ is 
a sub-cylinder of $\tcc_1$ of length $\tn$, for $\tY' = \piU(\tY)$ and $\tdd_1' = \piU(\tdd_1)$,
given an arbitrary point $z' \in \tY'$, it follows from (4.2) that
$$\frac{\nu(\tY')}{\nu(\tdd'_1)} \geq \frac{c_1 e^{g_{\tn}(z')}}{c_2 e^{g_{k+ \tq}(z')} } = \frac{c_1}{c_2} e^{g_{\tn-k- \tq}(\tpp^{k+ \tq}(z'))}
\geq \frac{c_1^2}{c_2} \rho_1^{\tn- (k+ \tq)}
= \frac{c_1^2}{c_2} \left( e^{- \tn} e^{k+ \tq}\right)^{ |\log \rho_1|} . $$
As before $\di e^k = e^{4d (m-k) \hep/|\log \theta|}$, and
$\di e^{-\tn} \geq e^{-c''} \, e^{\frac{\log \ep'}{|\log\theta|}} = c'''\,  (\ep')^{1/|\log\theta| } $
for some global constant $c''' > 0$.
Since $\tq$ is a global constant, it follows that
$$\frac{\nu(\tY')}{\nu(\tdd'_1)} \geq \frac{c_1^2}{c_2} \left(\ttc'\, \ep'\, e^{4d(m-k) \hep}\right)^{\frac{|\log \rho_1|}{|\log\theta|}} 
= \frac{c_1^2}{c_2} \left(\frac{\ttc'}{4C}\right)^{\frac{|\log \rho_1|}{|\log\theta|}} = d''_1 $$
for some  constants $\ttc' > 0$ and  $d''_1 = d''_1(C) > 0$.

For the cylinder $\tY$ in $\tdd_1$ just constructed we have
$(\Phi^u_{Z_{m-k}})^{-1}\circ \T_{Z_{m-k}} (\tY) \subset E\setminus \tH_1$,
so
\be
(\hphi_t)^{-1}\circ (\Phi^u_{Z_{m-k}})^{-1}\circ \T_{Z_{m-k}} (\tY) \subset E(z_p) \setminus \tH_2 ,
\ee
and $\nu(\tY') \geq d''_1\, \nu(\tdd'_1)$.

Let $\tLambda$ be the sub-cylinder of $\tdd$ such that $\tpp^{m-k}(\tLambda) = \tY$; then the length of $\tLambda$ is $\tn+ m-k$.
Using (4.2) again, we get
$$\frac{\nu(\tLambda')}{\nu(\tdd')} \geq \frac{c_1 e^{g_{\tn+ m-k}(z)}}{c_2 e^{g_{m+ \tq}(z)} } 
= \frac{c_1}{c_2} e^{g_{\tn- (k+q_1)}(\tpp^{m+q_1}(z))}
\geq \frac{c_1}{c_2} \rho_1^{\tn - (k+q_1)} \geq d''_1 ,$$
repeating the argument from the estimate of $\di \frac{\nu(\tY')}{\nu(\tdd'_1)}$. Thus,
\be
\nu(\tLambda') \geq d''_1\, \nu(\tdd') .
\ee
It follows from (4.36) and $\tpp^{m-k}(\tLambda) = \tY$ that 
$$(\Phi^u_{z_p})^{-1} ( \T_{z_p} (\varphi^p(\tLambda))) = (\hphi_t)^{-1}\circ 
(\Phi^u_{Z_{m-k}})^{-1}\circ \T_{Z_{m-k}} (\tpp^{m-k}(\tdd)) \subset \hG\setminus \tH_2 ,$$
where 
$$\hG = (\Phi^u_{z_p})^{-1}(\T_{z_p}(\varphi^p(\tdd))) \subset H'_1 = (\hphi_t)^{-1}(\tH_1) = \frac{1}{C}\, \hH_1  .$$
Therefore
\be\label{tW}
(\Phi^u_{z_p})^{-1} ( \T_{z_p} (\varphi^p(\tLambda))) \subset
(\Phi^u_{z_p})^{-1} ( \T_{z_p} (\varphi^p(\tdd))) \subset H'_1 \setminus \tH_2 ,
\ee
Now (4.37) implies
$\frac{\nu(\tLambda')}{\nu(\dd')} \geq  d''_1$
with the same constant $d''_1 = d''_1(C) > 0$. Since $\tdd$ has co-length $\tq$ in $\tcc$,
it follows that $\frac{\nu(\tLambda')}{\nu(\cc')} \geq  d'''_1$ for some global constant $d'''_1 > 0$.
Thus, there exists a global constant $q_1 \geq 1$ so that the co-length of $\tLambda$ in $\tcc$
is $\leq q_1$.

It follows from (4.38) that if $z \in \tLambda$, then for $w = (\Phi^u_{z_0})^{-1}(\T_{z_0}(z)) \in E^u(z_0)$ we have
$\xi = \hvarphi_{z_0}^{p}(w) \in H'_1 \setminus  \tH_2$, so $\|\xi_1\| \geq \frac{d_0\ttc_1}{4\gamma_1 C} e^{- 4 d (m-k)\hep}$.
This shows that $C \xi  \in H'_1\setminus H_2$.

So, for the maps
$$\Fp = (\Phi^u_{z_p})^{-1} \circ \T_{z_p} \circ \varphi^p : \tcc \longrightarrow E^u (z_p) \quad, \quad
\F0 = (\Phi^u_{z_0})^{-1} \circ \T_{z_0} : \tcc \longrightarrow E^u (z_0)$$
we have 
$\Fp (x) = C\, \hvarphi_{z_0}^p(\F0(x)) $
for all $x \in \tLambda$. While it appears that in general the locally defined near $z_0$ map
$$(\F0)^{-1}\circ \hvarphi_{z_0}^{-p} \circ \Fp : \tcc \longrightarrow \tcc ,$$
does not have to send a sub-cylinder to a sub-cylinder, it follows from our construction that, assuming that
global constant $q_1 \geq 1$ is sufficiently large, for every sub-cylinder $\tGamma$ of $\tcc$
of co-length $q_0$ there exists a sub-cylinder $\tLambda_1$ of $\tcc$ of co-length $q_1$ such that
$$\Fp(\tLambda'_1) \subset C\, \hvarphi_{z_0}^p(\F0(\tGamma')) .$$
Also,  for every sub-cylinder 
$\tLambda$ of $\tcc$ of co-length $q_0$ there exists a sub-cylinder $\tGamma_1$ of $\tcc$ of co-length $q_2$ such that
$$\Fp(\tLambda') \supset C\, \hvarphi_{z_0}^p(\F0(\tGamma_1')) .$$

This proves the lemma.
\endofproof


\section{Contact Anosov flows}
\setcounter{equation}{0}

\def\tyj{\tilde{y}^{(j)}}
\def\cvj{\check{v}_j}
\def\heta{\hat{\eta}}

\def\etao{\eta^{(1)}}
\def\etai{\eta^{(i)}}
\def\zetao{\zeta^{(1)}}
\def\hbeta{\hat{\beta}}
\def\tdelta{\tilde{\delta}}
\def\cxi{\check{\xi}}
\def\cxio{\cxi^{(1)}}
\def\cxit{\cxi^{(2)}}
\def\cet{\check{\eta}}
\def\ceto{\cet^{(1)}}
\def\cett{\cet^{(2)}}
\def\cv{\check{v}}
\def\cvo{\cv^{(1)}}
\def\cvt{\cv^{(2)}}
\def\cu{\check{u}}
\def\cuo{\cu^{(1)}}
\def\cut{\cu^{(2)}}
\def\cj{c^{(j)}}
\def\fj{f^{(j)}}
\def\gji{g^{(j,i)}}
\def\tvarphi{\tilde{\varphi}}
\def\hvarphi{\hat{\varphi}}

\def\hvo{\hat{v}^{(1)}}
\def\huo{\hat{u}^{(1)}}
\def\Uo{U^{(1)}}
\def\hw{\hat{w}}
\def\ha{\hat{a}}
\def\hb{\hat{b}}
\def\hc{\hat{c}}
\def\hk{\hat{k}}

\def\tn{\tilde{n}}
\def\omij{\omega_{i,j}}
\def\chep{\check{\epsilon}}

\def\hw{\hat{w}}
\def\vj{v^{(j)}}
\def\tbeta{\tilde{\beta}}

\subsection{Temporal distance function vs contact form}

We continue here with the assumptions and notation from Sect. 4. However now we assume that 
$\phi_t$ is a {\bf $C^2$ contact Anosov flow} on the compact Riemannian manifold $M$
with a $C^2$ invariant  contact form $\omega$.

For any $x\in \ll$ consider the $C^{2}$ map 
$$\tvarphi_x = (\exp^u_{\varphi(x)})^{-1} \circ \varphi \circ \exp^u_x : E^u (x; r(x))  \longrightarrow E^u (\varphi(x), \tr (\varphi(x)))\;.$$
It is well-defined assuming that the $\hep$-slowly varying radius function $r(x)$ and the $\hep/2$-slowly varying radius function $\tr(x)$
are chosen appropriately as in Sect. 3.
As with the maps $\hvarphi_x$, for $y \in \ll$ and an integer $j \geq 1$ we will use the notation
$$\tvarphi_y^j = \tvarphi_{\varphi^{j-1}(y)} \circ \ldots \circ \tvarphi_{\varphi(y)} \circ \tvarphi_y\quad,
\quad \tvarphi_y^{-j} = (\tvarphi_{\varphi^{-j}(y)})^{-1} \circ \ldots \circ (\tvarphi_{\varphi^{-2}(y)})^{-1}  \circ (\tvarphi_{\varphi^{-1}(y)})^{-1} \;,$$
at any point where these sequences of maps are well-defined.
In a similar way one defines the maps $\tvarphi_x$ and their iterations on $E^s(x;r(x))$.
It follows from the definitions that $\hvarphi_x = (\Psi^u_{\varphi(x)})^{-1} \circ \tvarphi_x \circ \Psi^u_x$.


The main ingredient in this section is the following lemma of Liverani (Lemma B.7 in \cite{L})
which significantly strengthens a lemma of Katok and Burns (\cite{KB}).

\bs

\noindent
{\bf Lemma 5.1.} (\cite{L}) {\it Let $\phi_t$ be a $C^2$ contact Anosov flow on $M$ with a 
$C^2$ contact form $\omega$.
Then there exist constants $C_0  > 0$, $\beta > 0$ and $\ep_0 > 0$ such that
for any $z\in M$, any $x\in  W^u_{\ep_0}(z)$ and any  $y\in W^s_{\ep_0}(z)$ we have
\be
|\Delta(x,y)  - d\omega_z(u,v)| \leq C_0\,\left[ \|u\|^2\, \|v\|^{\beta} + \|u\|^{\beta} \|v\|^2 \right]\;,
\ee
where $u \in E^u(z)$ and $v \in E^s(z)$ are such that $\exp^u_z(u) = x$ and $\exp^s_z(v) = y$.}

\bs

Replacing the constant $\beta > 0$ from Sect. 3 with a smaller one if necessary, we will assume
that {\bf the constant $\beta > 0$ above is the same as the one in Sect. 3}. {\bf Fix a constant $\ep_0$}
with the above property. We will also assume that $0 < \ep_1 < \ep_0$ satisfy the assumptions
in Sect. 2 and $\ep_0 < r_0$. 

The two-form $d\omega$ is $C^1$, so there exists a constant $C_0 > 0$ such that
\be
|d\omega_x(u,v)| \leq C_0 \|u\|\, \|v\| \quad , \quad u, v \in T_xM\:, \: x\in M .
\ee
Moreover, there exists a constant $c_0 > 0$ such that for any $x\in M$ and any $u\in E^u(x)$ with $\|u\| = 1$
there exists $v\in E^s(x)$ with $\|v\| = 1$ such that $|d\omega_x(u,v)| \geq 2 c_0$. {\bf Fix constants $C_0$ and $c_0$ }
with these properties and (5.1) as well.

\bs

\noindent
{\bf Corollary 5.2.} {\it Under the assumptions in Lemma {\rm 5.1}, we can choose the constant
$C_0  > 0$
so that for any $\hz\in M$, any $x , z\in  W^u_{\ep_0}(\hz)$ and any  $y\in W^s_{\ep_0}(z)$ we have
\be
|\Delta(x,y) | \leq C_0\,\|\hu - \hw\|^{\beta}\, \|v\|^{\beta} ,
\ee
where $\hu, \hw \in E^u(\hz)$ and $v \in E^s(z)$ are such that $\exp^u_{\hz}(\hu) = x$, $\exp^u_{\hz}(\hw) = z$ and $\exp^s_{z}(v) = y$.
Thus, we can choose the constant $C_0 > 0$ and $\beta > 0$ so that 
$\di |\Delta(x,y)| \leq C_0 \, (d(x,z))^{\beta} \, (d(z, y))^{\beta} $
under the above assumptions about $x$ and $y$.}

\bs

\noindent
{\it Proof of Corollary} 5.2. 
From the assumptions about $M$, $W^u_{\ep_0}(\hz)$ is a $C^2$ local submanifold of $M$. 

Given $\hz\in M$, $x\in  W^u_{\ep_0}(\hz)$ and  $y\in W^s_{\ep_0}(z)$, let 
$\hu, \hw \in E^u(\hz)$ and $v \in E^s(z)$ be such that $\exp^u_{\hz}(\hu) = x$, $\exp^u_{\hz}(\hw) = z$ and $\exp^s_{z}(v) = y$.
It follows from Lemma 5.1 used for the points $x \in W^u_{\ep_0}(z)$ and $y\in W^s_{\ep_0}(z)$ that if 
$u' \in E^u(z)$ is such that $\exp^u_z(u') = x$, then 
\be
|\Delta(x,y)  - d\omega_z(u' , v)| \leq C_0\,\left[ \|u'\|^2\, \|v\|^{\beta} + \|u'\|^{\beta} \|v\|^2 \right] .
\ee

Consider the map
$\psi  = ( \exp^u_{\hz})^{-1} \circ \exp^u_z : E^u(z, \ep) \longrightarrow E^u(\hz, \ep_0)$, defined for
appropriately chosen small $0 < \ep \leq \ep_0$. It follows from general properties of normal neighbourhoods
on Riemannian manifolds that $\frac{1}{C} \|\xi -\eta\|\leq \|\psi(\xi ) - \psi( \eta)\| \leq C\, \| \xi - \eta\|$ for all $\xi , \eta \in E^u(z,\ep)$ 
for some global constant $C > 0$. Since $x = \exp^u_z(u') = \exp^u_{\hz}(\hu)$, so $\hu = (\exp^u_{\hz})^{-1} (\exp^u_z(u'))$,
and similarly $\hw = (\exp^u_{\hz})^{-1} (z) = (\exp^u_{\hz})^{-1} (\exp^u_z(0))$, we get $\psi(u') - \psi(0) = \hu - \hw$, so
$ \|u'\| \leq C\, \|\hu - \hw\| .$
Now using (5.4) we get
\begin{eqnarray*}
|\Delta(x,y) |
& \leq & | d\omega_z(u' , v)| + C_0\,\left[ \|u'\|^2\, \|v\|^{\beta} + \|u'\|^{\beta} \|v\|^2 \right] \\
& \leq  & C_0\, \|u'\| \, \|v\| + 2C_0 \, \|u'\|^{\beta}\, \|v\|^{\beta} \leq 3C C_0 \|\hu- \hw\|^{\beta}\, \|v\|^{\beta} .
\end{eqnarray*}
Thus (5.3) holds, taking an appropriate larger constant $C_0 > 0$.
\endofproof

\bs

As is well-known the contact form $\omega$ vanishes on every stable/unstable manifold of a point on $M$,
while $d\omega$ vanishes on every weak stable/unstable manifold (see e.g.
\cite{KH} or Appendix B in \cite{L}). Moreover,
for Lyapunov regular points $x\in \ll$ and every
$u = (\uo, \ldots, u^{(\tk)}) \in E^u(x;r(x))$ and $v = (\vo, \ldots, v^{(\tk)})\in E^s(x;r(x))$
we have\footnote{See e.g. Lemma 9.1 in \cite{St4}, where this is formally proved.}
$\di d\omega_x(u,v) = \sum_{i=1}^{\tk} d\omega_x(\ui,\vi) .$

In the proof of the following lemma we will use some arguments from Sect. 9 in \cite{St4}.

\bs

\noindent
{\bf Lemma 5.3.} 
 {\it  There exist constants $C_4 > 0$ and $\tbeta$ with $0 < \tbeta \leq \beta$ such that:}
 
 \ms
 
 (a) {\it For any unstable cylinder $\tcc$ 
 in $\tR$ with $\tcc \cap \tP_0 \neq \e$,
any $x_0,z_0 \in \tcc$, and any $y_0, b_0 \in W^s_{\ep_0}(z_0)$ we have
\be
|\Delta(x_0 ,y_0)  - \Delta(x_0, b_0)| \leq C_4 \, \diam(\tcc) \, (d(y_0,b_0))^{\beta} .
\ee
In particular,
$$|\Delta(x_0 , y_0)| \leq C_4 \, \diam(\tcc)\, (d(y_0,z_0))^{\beta}  \leq C_4 \, \diam(\tcc) .$$
More precisely we have
\begin{eqnarray}
 |\Delta(x_0, y_0) - \Delta (x_0, b_0)|
& \leq &  |d\omega_{z_0} (\uo_0, \vo_0 - \etao_0)|
       + C_0 L_0 \|\uo_0\|\, \|v_0 - \eta_0\| \, (\|v_0-\eta_0\|^\beta + \|\eta_0\|^\beta) \nonumber\\
&        &                      + C_4 \, (\diam(\tcc))^{1+ \tbeta} ,
\end{eqnarray}
and
\begin{eqnarray}
|\Delta(x_0,y_0)|
 \leq   |d\omega_{z_0} (\uo_0, \vo_0 )| + C_0 L_0 \|\uo_0\| \, \|v_0\|^{1+\beta} + C_4 \, (\diam(\tcc))^{1+ \tbeta}  ,
\end{eqnarray}
where $u_0 \in E^u(z_0)$ and $v_0, \eta_0 \in E^s(z_0)$ are such that $\exp^u_{z_0}(u_0) = \hx_0 = \T_{z_0}(x_0)$, 
$\exp^s_{z_0}(v_0) = y_0$ and $\exp^s_{z_0}(\eta_0) = b_0$.}

\ms

(b) {\it For any unstable cylinder $\tcc$ in $\tR$, any $z_0 \in \tcc \cap  \tP_0$, 
any $x_0 \in W^u_R(z_0)$ with $\hx_0 = \T_{z_0}(x_0) \in \T_{z_0}(\cc)$, and any $y_0, b_0 \in W^s_{\ep_0}(z_0)$ we have 
\begin{eqnarray}
 |\Delta(x_0, y_0) - \Delta (x_0, b_0)|
& \geq & |d\omega_{z_0} (\uo_0, \vo_0 - \etao_0)|  
- C_0 L_0 \|\uo_0\|\, \|v_0 - \eta_0\| \, (\|v_0-\eta_0\|^{\beta} + \|\eta_0\|^{\beta})  \nonumber\\
&         & - C_4 \, (\diam(\tcc))^{1+ \tbeta} ,
\end{eqnarray}
and
\begin{eqnarray}
|\Delta(x_0, y_0) | \geq |d\omega_{z_0} (\uo_0, \vo_0 )|  
- C_0 L_0 \|\uo_0\|\, \|v_0\|^{1+\beta} - C_4 \, (\diam(\tcc))^{1+ \tbeta} ,
\end{eqnarray}
where $u_0$, $v_0$ and $\eta_0$ are as in part} (a).

\bs

\noindent
{\it Proof.} 
Let $\tcc$ be a cylinder  of length $m$ in $\tR$ and let $z_0 \in \tcc \cap \tP_0$.
Let  $x_0, z_0 \in \tcc$, $y_0, b_0 \in W^s_{\ep_0}(z_0)$.  
Then $R(z_0) \leq R_0$, $r(z_0) \geq r_0$.
We have
$x_0 = \Phi^u_{z_0}(u_0) = \exp^u_{z_0}(\tu_0)$
for some $u_0, \tu_0 \in E^u(z_0)$ with $\tu_0 = \Psi^u_{z_0} (u_0)$.  
Then 
$\|u_0\|, \|\tu_0\| \leq R_0 \; \diam(\tcc) .$
Similarly, write 
$$y_0 = \exp^s_{z_0}(\tv_0) = \Phi^s_{z_0}(v_0) \quad \mbox{\rm and} \quad b_0 = \exp^s_{z_0}(\teta_0) = \Phi^s_{z_0}(\eta_0)$$
for some 
$v_0, \tv_0, \eta_0, \teta_0 \in E^s(z_0)$ with $\tv_0 =  \Psi^s_{z_0} (v_0)$ and
$\teta_0 =  \Psi^s_{z_0} (\eta_0)$. It follows from  (3.7) that
\be
\|\tv_0 - v_0\| \leq R_0 \|v_0\|^{1+\beta} \:\: , \:\: \|\tu_0 - u_0\|\leq R_0\|u_0\|^{1+\beta} \:\:, \:\:
 \|\teta_0 - \eta_0\| \leq R_0\| \eta_0\|^{1+\beta} .
\ee
In particular, assuming $R_0 r_0^\beta \leq 2$, we have 
$\|\tv_0\| \leq 2\|v_0\|$, $\|\teta_0\| \leq 2 \|\eta_0\|$ and $\|\tu_0\| \leq 2\|u_0\| \leq 2 R_0 \diam(\tcc)$.
For $j \geq 0$ define
$$z_j = \varphi^j(z_0) \quad, \quad x_j = \varphi^j(x_0) \quad , \quad y_j = \varphi^j(y_0) \quad , \quad
u_j =  \hvarphi^j_{z_0}(u_0) \quad, \quad \tu_j = \tvarphi^j_{z_0}(\tu_0) ,$$
$$\hu_j =  d\hvarphi^j_{z_0}(0)\cdot u_0 \quad , \quad
\hv_j =  d\hvarphi^j_{z_0}(0)\cdot v_0 \quad  , \quad v_j = \hf^j_{z_0}(v_0) \quad , \quad \tv_j = \tvarphi^j_{z_0}(\tv_0) ,$$
$$b_j = \varphi^j(b_0) \quad  , \quad \heta_j = \hf^j_{z_0}(0)\cdot \eta_0 \quad , \quad \eta_j = \hf^j_{z_0}(\eta_0) 
 \quad , \quad \teta_j = \tvarphi^j_{z_0}(\teta_0).$$

Notice that $\tu_j = \Psi^u_{z_j} (u_j)$, $\tv_j = \Psi^s_{z_j} (v_j)$,  $\teta_j = \Psi^s_{z_j} (\eta_j),$
so it follows from (3.7) that
\be
\|u_j - \tu_j \| \leq R(z_j) \|u_j\|^{1+\beta} \: , \: \|v_j - \tv_j\| \leq R(z_j) \|v_j\|^{1+\beta} 
\: , \: \|\eta_j - \teta_j\| \leq R(z_j) \|\eta_j\|^{1+\beta} .
\ee
Moreover, 
$\exp^u_{z_j}(\tu_j) = \varphi^j(\exp^u_{z_0}(\tu_0)) = \varphi^j(x_0) = x_j $, $\exp^s_{z_j}(\tv_j) = y_j $ and
$\exp^s_{z_j}(\teta_j) = b_j$,
so Lemma 5.1 implies
\be
|\Delta(x_j,y_j)  - d\omega_{z_j}(\tu_j,\tv_j)| 
\leq  C_0 \,\left[ \|\tu_j\|^2\, \|\tv_j\|^\beta + \|\tu_j\|^\beta \|\tv_j\|^2 \right] 
\ee
and
\be
|\Delta(x_j,b_j)  - d\omega_{z_j}(\tu_j,\teta_j)| 
\leq  C_0 \,\left[ \|\tu_j\|^2\, \|\teta_j\|^\beta + \|\tu_j\|^\beta \|\teta_j\|^2 \right] .
\ee
We will use these a bit later.

Let the constants $0 < \hep_2 \leq \hep_1$ be as in Lemma 4.4. By the latter there exists an integer $k < m$ so that (4.11) holds,
and moreover for $p = [\ttau_{m-k}(z_0)]$, (4.15) and (4.16) hold, in particular we have
$\di \varphi^p(\tcc) \subset B^u(z_p, r(z_p))$,
so $x_p = \varphi^p(x_0) \in B^u(z_p, r(z_p))$ and $\|u_p\| = \|(\Phi^u_{z_p})^{-1}(x_p)\| \leq R(z_p) \chi \leq R_0 e^{p\hep}$.
As before, set $q = [\ttau_m(z_0)]$.
Then, as in the proof of Lemma 4.4 we have (4.22).

We will now use an argument from \cite{St4} with $\ell = p/2$, 
assuming for simplicity that $p$ is an even number (the other case is similar).
We have $r(z_j) \geq r_0 e^{- j \hep}$, so for all $0 \leq j \leq p$ 
it follows from (3.12) that
\be
\|u_j\| \leq \|u_j\|'_{z_j} \leq \frac{\|u_p\|'_{z_p}}{\mu_1^{p-j}} \leq
\Gamma(z_p) \frac{\|u_p\|}{\mu_1^{p-j}}  \leq \Gamma_0 e^{p \hep} \frac{R_0  e^{p \hep}}{\mu_1^{p-j}} 
= \Gamma_0 R_0  \frac{e^{2p\hep}}{\mu_1^{p-j}}\leq  r(z_j) 
\ee
for $j \leq \ell = p/2$. Indeed, for such $j$ we have
$\di \Gamma_0 R_0  \frac{e^{2p \hep} }{\mu_1^{p-j}} e^{j \hep} \leq  \Gamma_0 R_0  \frac{e^{3p\hep}}{\mu_1^{p/2}}
= \Gamma_0 R_0  \left(\frac{e^{3\hep}}{\mu_1}\right)^{p/2} \leq r_0 $,
assuming that $p$ is  sufficiently large. Thus (5.14) holds for all $0 \leq j \leq \ell$. 

The above also shows that $\|u_j\| \leq r_0$ for all  $0 \leq j \leq \ell$.
Essentially repeating the above estimate, we get
\be
\|u_{\ell}\| \leq \|u_\ell\|'_{z_\ell} \leq \frac{\|u_p\|'_{z_p}}{\mu_1^{p-\ell}}
\leq  \frac{\Gamma(z_p) e^{\ell \hep} \|u_p\|}{\lambda_1^{\ell}} \leq \Gamma_0 R_0 \chi \frac{e^{\ell \hep} e^{2p\hep}}{\lambda_1^\ell}
 \leq  \frac{R_0 \Gamma_0 e^{5\ell \hep}}{\lambda_1^{\ell}} \leq r_0 < 1 ,
\ee
assuming $\ell$ is sufficiently large.
Using (3.12) again (on stable manifolds) and assuming $\|v_0\| < 1$, we get
\be
\|v_{\ell}\| = \|v_\ell\|'_{z_\ell} \leq  \frac{\|v_0\|'_{z_0}}{\mu_1^\ell}  \leq \frac{\Gamma_0 e^{\ell\hep} \|v_0\|}{\lambda_1^\ell} 
\leq \frac{\Gamma_0 e^{\ell\hep}}{\lambda_1^\ell} .
\ee
Similarly, $\|\eta_\ell\| \leq  \frac{\Gamma_0 e^{\ell\hep}}{\lambda_1^\ell}$, assuming $\|\eta_0\| < 1$.

It follows from (5.10) and (5.11), repeating yet again some of the above estimates, that for $0 \leq j \leq \ell$ we have 
$\|\tu_j\| \leq \|u_j\| (1+ R(z_j) \|u_j\|^\beta) \leq \|u_j\| (1 + R_0 r_0)  \leq 2 \|u_j\|$.
Also,
$\|\tv_j\| \leq 2 \| v_j\|$ for all $j$ since the flow is contracting on stable manifolds.
Using these, it follows from (5.12) that
\begin{eqnarray}
|\Delta(x_j,y_j)  - d\omega_{z_j}(u_j,v_j)| 
& \leq &   4 C_0 R(z_j)\, \|u_j\|\, \|v_j\| (\|u_j\|^\beta + \|v_j\|^\beta) \nonumber\\
&        & + 8 C_0\,\left[ \|u_j\|^2\, \|v_j\|^\beta + \|u_j\|^\beta \|v_j\|^2 \right] 
\end{eqnarray}
for $0 \leq j \leq \ell$. Similarly (5.13) implies
\begin{eqnarray}
|\Delta(x_j,b_j)  - d\omega_{z_j}(u_j, \eta_j)| 
& \leq &  4 C_0 R(z_j)\, \|u_j\|\, \|\eta_j\| (\|u_j\|^\beta + \|\eta_j\|^\beta)\nonumber\\
&        & + 8 C_0\,\left[ \|u_j\|^2\, \|\eta_j\|^\beta + \|u_j\|^\beta \|\eta_j\|^2 \right] 
\end{eqnarray}
for $0 \leq j \leq \ell$. 

Next, we will be estimating 
$|\Delta(x_0,y_0) - d\omega_{z_0}(u_0,v_0)|  .$
Since $\Delta$ is $\varphi$-invariant and $d\omega$ is $d\varphi$-invariant  we have
$\Delta(x_0,y_0) = \Delta (x_j,y_j)$ and $d\omega_{z_0}(u_0,v_0) = d\omega_{z_j}(\hu_j , \hv_j) $
for all $j$. (Notice that $d\hf_{x}(0) = d\varphi(x)$ for all $x\in M$.)
With $j = \ell$, it follows from Lemma 3.1 and $z_0 \in \tP_0$  that
\be
\|\huo_\ell - \uo_\ell\| \leq L(z_\ell) \, \|u_\ell\|^{1+\beta}  \leq L_0 e^{\ell \hep} \, \|u_\ell\|^{1+\beta}.
\ee
Using Lemma 10.7(b) in \cite{St4} (see Lemma 3.3 above), backwards on stable manifolds, with \\
$a =  d\hf^{-\ell}_{z_\ell}(0) \cdot (\vo_\ell) \in E^s(z_0)$ and $b = d\hf^{-\ell}_{z_\ell}(0) \cdot (\etao_\ell)  \in E^s(z_0)$, 
since
$v_0 = \hf^{-\ell}_{z_\ell}(v_\ell)$ and $\eta_0 = \hf^{-\ell}_{z_\ell}(\eta_\ell)$, it follows that
$$\|(\ao - \bo) - (\vo_0 - \etao_0)\| \leq L_0 \left[ \|v_0 - \eta_0\|^{1+\beta}  + \|\eta_0\|^\beta \|v_0-\eta_0\|\right]
\leq 2L_0 \| v_0 - \eta_0\| .$$
Thus,
\be
\| d\hf^{-\ell}_{z_\ell}(0) \cdot (\vo_\ell - \etao_\ell) - (\vo_0 - \etao_0)\| \leq L_0 \|v_0 - \eta_0\| \, (\|v_0-\eta_0\|^\beta + \|\eta_0\|^\beta) .
\ee


Proceeding as in Sect. 9 in \cite{St4} and using (5.19), we obtain
\begin{eqnarray*}
& &        |d\omega_{z_\ell} (u_\ell, v_\ell - \eta_\ell)| 
 \leq  |d\omega_{z_\ell} (\uo_\ell, \vo_\ell - \etao_\ell)| 
        + C_0 \sum_{i=2}^{\tk} \|u_\ell^{(i)}\| \, (\|v_\ell^{(i)}\| + \| \eta_\ell^{(i)}\|)\nonumber\\
& \leq &  |d\omega_{z_\ell} (\huo_\ell, \vo_\ell - \etao_\ell)| + C_0 L_0 e^{\ell \hep} \|u_\ell\|^{1+\beta}  \|\vo_\ell - \etao_\ell\| \nonumber
         +  C_0 \sum_{i=2}^{\tk} \|u_\ell^{(i)}\| \, (\|v_\ell^{(i)}\|  + \| \eta_\ell^{(i)}\|)\nonumber\\
& = & |d\omega_{z_\ell} (d \hvarphi^{\ell}_{z_0}(0) \cdot \uo_0, \vo_\ell - \etao_\ell)| 
+  C_0 L_0 e^{\ell \hep} \|u_\ell\|^{1+\beta} \|\vo_\ell - \etao_\ell\|\nonumber
         +  C_0 \sum_{i=2}^{\tk} \|u_\ell^{(i)}\| \,  (\|v_\ell^{(i)}\| + \| \eta_\ell^{(i)}\|)\nonumber\\
& =     & |d\omega_{z_0} (\uo_0, d\hvarphi^{-\ell}_{z_\ell}(0)\cdot (\vo_\ell - \etao_\ell))|
+ C_0 L_0  e^{\ell \hep} \|u_\ell\|^{1+\beta} \|\vo_\ell - \etao_\ell\| \nonumber
         +  C_0 \sum_{i=2}^{\tk} \|u_\ell^{(i)}\| \,  (\|v_\ell^{(i)}\| + \| \eta_\ell^{(i)}\|) .
\end{eqnarray*}
Now (5.20) implies
\begin{eqnarray*}
&        & |d\omega_{z_0} (\uo_0, d\hf^{-\ell}_{z_\ell}(0) \cdot (\vo_\ell - \etao_\ell)| \\
& \leq & |d\omega_{z_0} (\uo_0, \vo_0 - \etao_0)| 
+ |d\omega_{z_\ell} (\uo_0, d\hf^{-\ell}_{z_\ell}(0) \cdot (\vo_\ell - \etao_\ell) - (\vo_0 - \etao_0)|\\
& \leq & |d\omega_{z_0} (\uo_0, \vo_0 - \etao_0)| + C_0 L_0 \|\uo_0\| \, \|v_0 - \eta_0\| \, (\|v_0-\eta_0\|^\beta + \|\eta_0\|^\beta) .
\end{eqnarray*}
The latter and the previous estimate yield
\begin{eqnarray}
         |d\omega_{z_\ell} (u_\ell, v_\ell - \eta_\ell)| 
& \leq &  |d\omega_{z_0} (\uo_0, \vo_0 - \etao_0)| + C_0 L_0 \|\uo_0\|\, \|v_0 - \eta_0\| \, (\|v_0-\eta_0\|^\beta + \|\eta_0\|^\beta) \nonumber\\
&      & + C_0 L_0 e^{\ell \hep} \|u_\ell\|^{1+\beta} \|\vo_\ell - \etao_\ell\| +  C_0 \sum_{i=2}^{\tk} \|u_\ell^{(i)}\| \, (\|v_\ell^{(i)}\| + \| \eta_\ell^{(i)}\|) .
\end{eqnarray}
Consequently,
\begin{eqnarray}
 |d\omega_{z_\ell} (u_\ell, v_\ell - \eta_\ell)| 
& \leq & \Con \, \diam(\tcc) \, \|v_0 - \eta_0\|  + \Con e^{\ell \hep} \|u_\ell\|^{1+\beta} \|\vo_\ell - \etao_\ell\| \nonumber\\
&        & +  C_0 \sum_{i=2}^{\tk} \|u_\ell^{(i)}\| \, (\|v_\ell^{(i)}\| + \| \eta_\ell^{(i)}\|) .
\end{eqnarray}

Using similar estimates from below, we obtain
\begin{eqnarray}
           |d\omega_{z_\ell} (u_\ell, v_\ell - \eta_\ell)| 
& \geq &  |d\omega_{z_0} (\uo_0, \vo_0 - \etao_0)| - C_0 L_0 \|\uo_0\|\, \|v_0 - \eta_0\| \, (\|v_0-\eta_0\|^\beta + \|\eta_0\|^\beta) \nonumber\\
&         &  -  C_0 L_0  e^{\ell \hep} \|u_\ell\|^{1+\beta} \|\vo_\ell - \etao_\ell\| 
- C_0 \sum_{i=2}^{\tk} \|u_\ell^{(i)}\| \, (\|v_\ell^{(i)}\| + \| \eta_\ell^{(i)}\|) .
\end{eqnarray}

Next, it follows from (\ref{eq:diamC}) that 
$\di \diam(\tcc) \geq 
\frac{e^{-q \hep_7}}{C_3 \lambda_1^q} $
for some global constant $C_3 > 0$, so
$\di \lambda_1^q \, e^{q \hep_7} \geq \frac{1}{C_3\, \diam(\tcc)}$,
that is
$\di q \, \log \lambda_1 + q \hep_7 \geq \log \frac{1}{C_3\, \diam(\tcc)} > 0$.
Assuming $\hep$ and so $\hep_7$ is sufficiently small and using (4.22), this gives
$\frac{p}{1-\hep_6} \, (\log \lambda_1 +  \hep_7)\geq \log \frac{1}{C_3\, \diam(\tcc)}$,
therefore for $\ell = p/2$ we get
\be
\ell 
> \frac{1- \hep_6}{2(\log \lambda_1 + \hep_7)}\, \log \frac{1}{C_3\, \diam(\tcc)} .
\ee
Now (5.15) implies
\begin{eqnarray}
\|u_{\ell}\|
& \leq & R_0 \Gamma_0 (\lambda_1 e^{-5\hep})^{-\ell} = R_0 \Gamma_0 e^{- \ell \log (\lambda_1 e^{-5\hep})} 
 \leq  R_0 \Gamma_0 \, e^{-\frac{(1-\hep_6) \log (\lambda_1 e^{-5\hep})}{ 2 (\log \lambda_1 + \hep_7)}
\log \left(\frac{1}{C_3\, \diamf(\tcc)}\right) } \nonumber \\
& =    &  R_0 \Gamma_0 \,\left(\frac{ 1}{C_3\, \diam(\tcc)}\right)^{- \frac{(1-\hep_6) (\log\lambda_1 - 5 \hep)}{2 (\log \lambda_1 + \hep_7)}}
\leq R_0 \Gamma_0 C_3\, \left(\diam(\tcc)\right)^{\frac{(1-\hep_6) (\log\lambda_1 - 5 \hep)}{2  (\log \lambda_1 + \hep_7)}} .
\end{eqnarray}
(Using here $C_3 > 1$ and $\frac{(1-\hep_6) (\log\lambda_1 - 5 \hep)}{2( \log \lambda_1 + \hep_7)} < 1$.)
Similarly, (5.16) yields
\begin{eqnarray}
\|v_{\ell}\|
& \leq & \Gamma_0 (\lambda_1 e^{-\hep})^{-\ell}
 \leq  \Gamma_0 \, e^{-\frac{(1-\hep_6) \log (\lambda_1 e^{-\hep})}{2 ( \log \lambda_1 +\hep_7) }  \log \left(\frac{1}{C_3\, \diamf( \tcc)}\right) } \nonumber\\
& \leq & \Gamma_0 C_3 \, \left(\diam(\tcc)\right)^{\frac{(1-\hep_6) (\log\lambda_1 -  \hep)}{2( \log \lambda_1 + \hep_7)}} .
\end{eqnarray}
The same estimate holds for $\|\eta_\ell\|$.

Using these we get the following estimates for terms in (5.18) with $j = \ell$:
\begin{eqnarray*}
\|u_\ell\|\, \|v_\ell\| (\|u_\ell\|^\beta + \|v_\ell\|^\beta) 
\leq  C'_3 \left(\diam(\tcc)\right)^{(2+\beta)\frac{(1-\hep_6)(\log\lambda_1 -  5\hep)}{2 (\log \lambda_1+ \hep_7)}}  \leq  C'_3 (\diam(\tcc))^{1+\hbeta} ,
\end{eqnarray*}
where $C'_3 = 2 (\Gamma_0 R_0 C_3)^3$ and we choose 
\be\label{eq:hbeta}
0 < \hbeta  = \frac{1}{2} \min\left\{  \frac{\beta}{4}  \; ,\:
\frac{\log\lambda_2 - \log \lambda_1}{2  \log \lambda_1}  \right\} .
\ee
Then 
$$(2+\beta)\frac{(1-\hep_6)(\log \lambda_1 -  5\hep)}{2 ( \log \lambda_1 + \hep_7)} 
= (1+\beta/2)(1-\hep_6) \left(1- \frac{5 \hep + \hep_7}{\log\lambda_1 + \hep_7} \right) \geq 1+ \hbeta ,$$
assuming $\hep_6 > 0$, $\hep_7 > 0$  and $\hep > 0$ are sufficiently small.
Similarly, we obtain
\be
e^{\ell \hep} \|u_\ell\|^{1+\beta} \|v_\ell\| \leq C''_3 (\diam(\tcc))^{1+\hbeta} 
\ee
for some global constant $C''_3 > 0$.
To prove the latter we need an estimate from above similar to (5.24).
First, using (\ref{eq:diamC}) we get
$\di \diam(\tcc) \leq  \frac{C_3 e^{q \hep_7}}{\lambda_1^q}, $
so
$\di \lambda_1^q \, e^{-q \hep_7} \leq \frac{C_3}{\diam(\tcc)} ,$
that is
$$\di q \, \log \lambda_1 - q \hep_7 \leq \log \frac{C_3}{\diam(\tcc)} .$$
Now (4.22) yields
$\di p \, (\log \lambda_1 - \hep_7)\leq \log \frac{C_3}{ \diam(\tcc)} ,$
therefore for $\ell = p/2$ we get
$$\ell 
\leq \frac{1}{2(\log \lambda_1 - \hep_7)}\, \log \frac{C_3}{\diam(\tcc)} . $$
This, (5.15) and (5.16), using (5.24) again as well, imply
\begin{eqnarray*}
e^{\ell \hep} \|u_\ell\|^{1+\beta} \|v_\ell\|
& \leq & \left( \frac{C_3}{\diam(\tcc)}\right)^{\frac{\hep}{2(\log \lambda_1 - \hep_7)}}\, 
C'_3 \left(\diam(\tcc)\right)^{(2+\beta)\frac{(1-\hep_6)(\log\lambda_1 -  5\hep)}{2 (\log \lambda_1 + \hep_7)}}\\
& \leq & C_3 C'_3\, \left(\diam(\tcc)\right)^{(2+\beta)\frac{(1-\hep_6)(\log\lambda_1 -  5\hep)}{2 (\log \lambda_1 + \hep_7)} 
- \frac{\hep}{2(\log \lambda_1 - \hep_7)}}
 \leq  C''_3 \, (\diam(\tcc))^{1 + \hbeta} ,
\end{eqnarray*}
where $C''_3 = C_3 \, C'_3$. The above holds since $\frac{\hep}{2(\log \lambda_1 - \hep_7)} < \frac{\beta}{4}$,
assuming $\hep > 0$ is sufficiently small, and
\begin{eqnarray*}
(2+\beta)\frac{(1-\hep_6)(\log\lambda_1 -  5\hep)}{2 (\log \lambda_1 + \hep_7)}  - \frac{\beta}{4}
 \geq  (1 +\beta/2) (1-\hep_6)\left( 1  -  \frac{5\hep + \hep_7}{\log \lambda_1 + \hep_7} \right)  - \frac{\beta}{4}
 \geq  1 + \hbeta ,
\end{eqnarray*}
assuming that $\hep_6 > 0$ and $\hep > 0$ are sufficiently small.
This proves (5.28). We get a similar estimate replacing $v_\ell$ by $\eta_\ell$.

Next, for any $u = \uo + \ut + \ldots + u^{(\tk)} \in E^u(z)$ or $E^s(z)$, $z\in M$, set
$\cut = \ut + \ldots + u^{(\tk)}$,
so that $u = \uo + \cut$.
Using Lemma 3.5 in \cite{St3} (see Lemma 3.2 above), $p = 2\ell $, $u_p \in E^u(z_p, r(z_p))$ and $\|u_p\| \leq R_0 e^{p \hep}$, 
we get
$$\| \cut_{\ell}\|'_{z_\ell}  \leq \frac{\Gamma_0  e^{p \hep} \|\cut_{p}\|}{\mu_2^{\ell}} 
\leq  \frac{\Gamma_0  e^{p \hep} \|u_{p}\|}{\mu_2^{\ell}} \leq  \frac{\Gamma_0 R_0 e^{4\ell \hep}}{\mu_2^{\ell}} . $$
Similarly, using Lemma 3.2 above (backwards for the map $\varphi^{-1}$ on stable manifolds),
$z_0 \in \tP_0$, $v_0\in E^s(z_0,r_0)$ and the fact that  $\|v_0\| \leq \ep < 1$, we get
$$\| \cvt_{\ell}\|'_{z_\ell}  \leq \frac{\Gamma_0 \|v_0\| e^{2 \ell \hep}}{\mu_2^{\ell}} 
\leq  \frac{\Gamma_0 e^{2 \ell \hep} }{\mu_2^{\ell}} . $$
Hence for $i \geq 2$ we have
$$\|\ui_\ell\| \leq |\cut_\ell| \leq \Gamma(z_\ell) \|\cut_\ell\| \leq \frac{\Gamma^2_0 R_0 e^{5\ell \hep}}{\lambda_2^{\ell}} , $$
where we used $\mu_2 = \lambda_2 e^{-\hep}$. Similarly 
$\di \|\vi_\ell\| \leq \frac{\Gamma^2_0 e^{3\ell \hep}}{\lambda_2^\ell} .$
From these estimates, (5.24), and the assumptions about $\hep$, we get
\begin{eqnarray*}
\|\ui_\ell\|\, \| \vi_\ell\|
& \leq & \Gamma_0^4 R_0\,  (\lambda_2 e^{-5\hep})^{-  2\ell}  = \Gamma_0^4 R_0\, e^{-2\ell \log(\lambda_2e^{-5\hep})}\\
& \leq &  \Gamma_0^4 R_0\, e^{\frac{- (1-\hep_6)\log(\lambda_2e^{-5 \hep})}{\log \lambda_1 + \hep_7}\, \log\left(\frac{1}{C_3\, \diamf(\tcc)} \right)}\\
& \leq & \Gamma_0^4 R_0\,\left( C_3\, \diam(\tcc)\right)^{\frac{(1-\hep_6) (\log \lambda_2 - 5 \hep)}{\log\lambda_1 + \hep_7} } 
\leq C''_3 (\diam(\tcc))^{1+\hbeta} ,
\end{eqnarray*}
where $C''_3 = \Gamma_0^4 R_0 (C_3)^{\log\lambda_2/\log\lambda_1}$. Here we are assuming (5.27) and that
$\hep_6 > 0$ and $\hep_7 > 0$ are sufficiently small. Since
$$\frac{(1 - \hep_6)\log \lambda_2}{\log \lambda_1 + \hep_7} > \frac{\log \lambda_2 + \log \lambda_1}{2 \log \lambda_1} ,$$ 
for small $\hep_6 > 0$, $\hep_7 > 0$, we derive
$$\frac{(1-\hep_6) (\log \lambda_2 - 5\hep)}{\log\lambda_1 + \hep_7} =
 \frac{(1-\hep_6)\log \lambda_2}{\log\lambda_1 + \hep_7} - \frac{5 (1-\hep_6) \hep}{\log\lambda_1 + \hep_7}
> \frac{\log \lambda_2 + \log \lambda_1}{2 \log \lambda_1} - \frac{5  \hep}{\log\lambda_1} \geq 1 + \hbeta ,$$
assuming that $\hep > 0$ is sufficiently small. This proves the above estimate for $\|\ui_\ell\|\, \| \vi_\ell\|$.

Using  (5.22), (5.28) and the above estimates for $\|\ui_\ell\|$ and $\| \vi_\ell\|$,  we obtain
\begin{eqnarray*}
|d\omega_{z_\ell} (u_\ell, v_\ell - \eta_\ell)|
& \leq & \Con \, \diam(\tcc) \, \|v_0 - \eta_0\| + \Con \, (\diam(\tcc))^{1+ \hbeta} .
\end{eqnarray*}

It now follows from  (5.17) and (5.18) with $j = \ell$  and the previous estimates that
\begin{eqnarray}
&        & |\Delta(x_0, y_0) - \Delta (x_0, b_0)|
 =      |\Delta(x_\ell, y_\ell) - \Delta (x_\ell, b_\ell)| \nonumber\\
& \leq & |d\omega_{z_\ell}(u_\ell,v_\ell) - d\omega_{z_\ell}(u_\ell, \eta_\ell)| 
+   4 C_0 R(z_\ell)\, \|u_\ell\|\, \|v_\ell\| (\|u_\ell\|^\beta + \|v_\ell\|^\beta) \nonumber\\
&      &+ 8 C_0\,\left[ \|u_\ell\|^2\, \|v_\ell\|^\beta + \|u_\ell\|^\beta \|v_\ell\|^2 \right] \nonumber\\
&        & + 4 C_0 R(z_\ell)\, \|u_\ell\|\, \|\eta_\ell\| (\|u_\ell\|^\beta + \|\eta_\ell\|^\beta)
 + 8 C_0\,\left[ \|u_\ell\|^2\, \|\eta_\ell\|^\beta + \|u_\ell\|^\beta \|\eta_\ell\|^2 \right] \nonumber\\
& \leq & |d\omega_{z_\ell} (u_\ell, v_\ell - \eta_\ell)| + \Con \, (\diam(\tcc))^{1+ \hbeta}\nonumber\\
& \leq & \Con \, \diam(\tcc) \, \|v_0 - \eta_0\| + \Con \, (\diam(\tcc))^{1+ \hbeta} .
\end{eqnarray}
The  more precise estimate (5.6), for some global constant $C_4 > 0$, follows immediately from (5.21).
In particular, with $\eta_0 = 0$, the latter gives (5.7).

In a similar way, this time using (5.23)  we obtain (5.8), and with $\eta_0 = 0$ it gives (5.9).

To prove (5.5) we just repeat an argument from Sect. 9.3.3 in \cite{St4}.

\ms

\noindent
{\bf Case 1.} $\diam(\tcc) \leq \| v_0 -\eta_0\|^{\beta/2}$. Then by (5.29),
$$|\Delta(x_0, y_0) - \Delta (x_0, b_0)| \leq \Con \, \diam(\tcc) \, \|v_0 - \eta_0\|^{\hbeta \beta/2} .$$
So, (5.5) holds with $\beta$ replaced by $\hbeta \beta/2$.

\ms

\noindent
{\bf Case 2.} $\diam(\tcc) \geq \| v_0 -\eta_0\|^{\beta/2}$. 
Consider the point $X = \phi_{\Delta(x_0,y_0)}([x_0,y_0]) \in W^u_\ep(y_0)$.
It is easy to see that
$$\Delta(x_0, y_0) - \Delta (x_0, b_0) = \Delta(X, y_0) - \Delta(X, b_0) = - \Delta(X, b_0) .$$
We have $X = \exp^u_{y_{0}}(\tt)$ and $b_{0} = \exp^s_{y_{0}}(\ts)$ for some $\tt \in E^u(y_{0})$ and 
$\ts \in E^s(y_{0})$. Clearly $\|\tt\| \leq \Con$.
Using Lemma 5.1 we get
$$|\Delta(X,b_{0})| \leq C_0 [ |d\omega_{y_{0}}(\tt, \ts)|  + \|\tt\|^2 \|\ts\|^{\beta} + \|\tt\|^{\beta} \|\ts\|^2]
\leq \Con \, \|\ts\|^{\beta} .$$
However, $\|\ts\| \leq \Con d(y_{0},b_{0}) \leq \Con \|v_{0}-\eta_{0}\|$, so
$$|\Delta(X,b_{0})| \leq \Con \|v_0 - \eta_0\|^{\beta}  \leq \Con \diam(\tcc) \|v_0 - \eta_0\|^{\beta/2} .$$
Therefore
$\di |\Delta(x_0, y_0) - \Delta (x_0, b_0)|  \leq \Con \diam(\tcc) \|v_0 - \eta_0\|^{\beta/2}$,
so (5.5) holds with $\beta$ replaced by $\beta/2$.
This proves Lemma 5.3.
\endofproof

\def\tEE{{\mathcal E}}
\def\hwo{\hat{w}^{(1)}}
\def\yj{y^{(j)}}
\def\vj{v^{(j)}}
\def\Lj{L^{(j)}}
\def\Mj{M^{(j)}}
\def\dj{d^{(j)}}
\def\bj{b^{(j)}}
\def\hGamma{\widehat{\Gamma}}
\def\hLambda{\widehat{\Lambda}}

\subsection{Non-integrability of contact Anosov flows}

{\bf Fix constants} $C_4 > 0$, $\beta > 0$ and $\tbeta > 0$ with the properties in Lemma 5.3. 
Let the small constant $\hep > 0$ be as in Sect. 3 and $0 < \ep_1 < \ep_0$ as in Sect. 2.

Recall the projections $\T_z : W^u_{\tR}(z) \longrightarrow W^u_{\ep_0}(z)$ for $z\in R$ from Sect. 4.1
and the constants $C_0 > 0$ and $c_0 > 0$ introduced before Lemma 5.1. Set
\be\label{eq:delta0}
\delta_0 =  \left( \frac{c_0}{32 C_0 L_0 R_0} \right)^{1/\tbeta} .
\ee

The following lemma is derived from the non-integrability of the flow which stems from the fact that the
flow is contact.

\bs


\noindent
{\bf Lemma 5.4.} {\it  
There exist global constants $\tm_0 \geq 1$,  $\delta' \in (0,1)$,  $d_2 \in (0,1)$
and $\hep_{13} \leq \con \hep$ which can be made
arbitrarily small with $\hep$ such that for every integer  $m \geq \tm_0$ we have  the following:}

\ms

(a) {\it  For any $z_0 \in \tR \cap \tP_0$ and any cylinder $\tcc$ of length $m$ in $W^u_{\tR}(z_0)$ with $z_0 \in \tcc$,
if $x_0  \in \tcc $ is such that 
\begin{equation}\label{eq:u0}
u_0 = (\Phi^u_{z_0})^{-1}(\T_{z_0}(x_0))  \in E^u(z_0)
\end{equation}
satisfies
\be\label{eq:u0-cond}
\| \uo_0\|  \geq \kappa \, \diam(\tcc) ,
\ee
where  $\kappa = d_2 e^{- m \hep_{13}}$, 
then there exist  a point $y_0 = y_0(z_0, x_0)\in  B^s (z_0,\ep_1)$ such that we have
\be\label{eq:Delta-cond}
12 c_0 \delta_0 \, \|\uo_0\| \leq |\Delta(x_0, b_1)  - \Delta(x_0, b_2)| 
\ee
for any $b_1,b_2 \in W^s_R(z_0)$ with $d(z_0,b_1) < 2\delta'$ and $d(y_0, b_2) < 2\delta'$. Thus,
\be \label{eq:Delta1-cond}
12 c_0 \delta_0 \,\kappa\,   \diam (\tcc) \leq |\Delta(x_0, b_1)  - \Delta(x_0, b_2)| 
\ee
under the above conditions.}

\ms

(b)  {\it There exists an integer $N_0 \geq 1$ 
such that for any integer $N \geq N_0$,  any $z_0 \in \tR \cap \tP_0$, any cylinder $\tcc$ in $\tR$
of length $m \geq \tm_0$ in $W^u_{\tR}(z_0)$ with $z_0 \in \tcc$, and any $x_0  \in \tcc$ such that
$u_0 \in E^u(z_0)$ with {\rm (\ref{eq:u0})} satisfies {\rm (\ref{eq:u0-cond})}, there exist 
$y_0 = y_0(z_0, x_0)\in  B^s (z_0,\ep_1)$  as in part (a) and families of points 
$$y_1 = y_1(z_0,x_0) \in \tpp^N(B^u(z_0,\ep')) \cap B^s(z_0, \delta')  \:\: , \:\:
y_2 = y_2(z_0,x_0) \in \tpp^N(B^u(z_0,\ep')) \cap B^s(y_0, \delta') .$$
such that
{\rm (\ref{eq:Delta-cond})}  holds for any $b_1 \in B^s(y_1, \delta')$, any  $b_2 \in B^s(y_2, \delta')$.}

\bs

\noindent
{\it Proof of Lemma} 5.4. 
(a)
Fix for a moment $z_0\in \tR \cap \tP_0$. Let $\tcc$ be a cylinder of length $m \geq \tm_0$ in $\tR$ with $z_0 \in \tcc$,
and let $x_0$, $u_0$ and $\kappa$  satisfy (\ref{eq:u0}) -- (\ref{eq:u0-cond}).
By the choice of $c_0 > 0$ in Sect. 5.1, there exists a vector $\tv\in E^s_1(z_0)$ with $\|\tv\| = 1$ 
such that
\be\label{eq:tv}
d\omega_{z_0}(\uo_0 , \tv) \geq 2 c_0 \, \|\uo_0\|  .
\ee
{\bf Fix a vector $\tv$} with the above property and set
\be \label{eq:vo}
v_0 = 10 \delta_0\, \tv
\in E^s_1(z_0)\quad , 
\quad  y_0  = \Phi^s_{z_0}(v_0) \in W^s_{\ep_1}(z_0)  ,
\ee
assuming that $\delta_0 > 0$ is sufficiently small so that $10\delta_0 < \ep_1$.
Then 
$\frac{\|v_0\|}{R_0}  \leq d(z_0 , y_0 ) \leq  \|v_0\| $, 
and
$$|d\omega_{z_0}(u_0 , v_0 ) | = |d\omega_{z_0}(\uo_0 , v_0 ) | \geq 2 c_0   \|\uo_0\| \|v_0\| = 20 c_0 \delta_0 \|\uo_0\| ,$$
while (\ref{eq:u0}) and (\ref{eq:u0-cond}) imply
$ \kappa\, \diam(\tcc) \leq \|\uo_0\| \leq \|u_0\| \leq  R_0\, \diam(\tcc)$.
Notice that 
$$|d\omega_{z_0}(\uo_0 , v_0 ) | \leq C_0 \|\uo_0\|\, \|v_0\| = 10 C_0 \delta_0 \|\uo_0\| 
\leq 10 C_0\, \delta_0 \, R_0\, \diam(\tcc) .$$

We will now prove that
\be\label{eq:x0y0}
|\Delta(x_0,y_0)| \geq 14 c_0 \delta_0\, \|\uo_0\|  .
\ee

In what follows we use the notation from the proof of  Lemma 5.3 and also part (b) in Lemma 5.3 (see (5.9) there).
Using it with $b_0 = z_0$ (so $\eta_0 = 0$) and using $v_0 \in E^s_1(z_0)$ and (5.36) we get
\begin{eqnarray}
|\Delta(x_0, y_0)|
& \geq &  |d\omega_{z_0} (\uo_0, v_0 )|  - C_0L_0 \, \|\uo_0\| \,\|v_0\|^{1+\beta} - C_4 \, (\diam(\tcc))^{1+ \tbeta} \nonumber\\
& \geq &  2 c_0  \, \|\uo_0\|\, \|v_0\|  - C_0L_0 \, \|\uo_0\| \,\|v_0\|^{1+\beta} - C_4 \, (\diam(\tcc))^{1+ \tbeta}\nonumber \\
& =      &  \|\uo_0\| \|v_0\| (2 c_0  - C_0L_0 \|v_0\|^{\beta})  - C_4 \, (\diam(\tcc))^{1+ \tbeta}  \nonumber\\
& =      &  10 \delta_0 \|\uo_0\|  (2 c_0  - C_0L_0 (10 \delta_0)^{\beta})  - C_4 \, (\diam(\tcc))^{1+ \tbeta}  \nonumber\\
& \geq & 15 c_0 \delta_0 \|\uo_0\| - C_4 \, (\diam(\tcc))^{1+ \tbeta} ,
\end{eqnarray}
assuming $\delta_0 > 0$ is sufficiently small so that $\di C_0L_0 (10 \delta_0)^{\beta} < c_0/2$.
Since $\delta_0$,  $c_0$ and $C_4 > 0$ are global constants, using (\ref{eq:diamC}) and assuming that $m_0 \geq 1$ is sufficiently large 
and $m \geq m_0$ we have 
$$C_4\, (\diam(\tcc))^{\tbeta} \leq C_4 \left(\frac{C_3 e^{m\tau_0 \hep_7}}{\lambda_1^{m\ttau_0}}\right)^{\tbeta}
\leq \frac{c_0 \delta_0\, d_2}{8} e^{-m \hep_{13}} = \frac{c_0 \delta_0\, \kappa}{8} ,$$
choosing $\hep_{13} \leq \con \hep$. Thus,
\be\label{eq:C4-cond}
C_4 (\diam(\tcc))^{1+\tbeta} \leq  \frac{c_0 \delta_0\,\kappa}{8} \, \diam(\tcc) \leq   \frac{c_0 \delta_0}{8} \, \|\uo_0\| .
\ee
This and the above estimates imply (\ref{eq:x0y0}).

In a similar way, this time using (5.7), we get
\begin{eqnarray*}
|\Delta(x_0, y_0)|
& \leq &  |d\omega_{z_0} (\uo_0, \vo_0 )|  + C_0L_0 \, \|\uo_0\| \,\|v_0\|^{1+\beta} + C_4 \, (\diam(\tcc))^{1+ \tbeta} \nonumber\\
& \leq &  C_0  \, \|\uo_0\|\, \|v_0\|  + C_0L_0 \, \|\uo_0\| \,\|v_0\|^{1+\beta} + C_4 \, (\diam(\tcc))^{1+ \tbeta}\nonumber \\
& =      &  \|\uo_0\| \|v_0\| (C_0  + C_0L_0 \|v_0\|^{\tbeta})  + C_4 \, (\diam(\tcc))^{1+ \tbeta}  \nonumber\\
& \leq & 20 \delta_0 C_0 \|\uo_0\| + C_4 \, (\diam(\tcc))^{1+ \tbeta} .
\end{eqnarray*}
We will use this with $y_0$ replaced by $b_1 \in W^s_R(z_0)$ with $d(z_0, b_1) < 2\delta'$  for some small $\delta' > 0$
(to be determined later). Let $b_1 = \Phi^s_{z_0}(v_1)$ for some $v_1\in E^s(z_0)$. Then $\|v_1\| \leq R_0 d(z_0,b_1) < 2 R_0 \delta'$
and (\ref{eq:C4-cond}) give
\be
|\Delta(x_0, b_1)| \leq 2C_0 R_0 \delta' \|\uo_0\| + \frac{c_0 \delta_0}{8} \|\uo_0\|
\leq \frac{c_0 \delta_0}{4} \,  \|\uo_0\| ,
\ee
assuming that
$4C_0 R_0 \delta' \leq c_0 \delta_0/8$.
For this and later use define
\be\label{eq:delta'}
\delta' = \frac{c_0 \delta_0}{64 C_0 L_0 R_0}  = \frac{c_0}{64 C_0 L_0 R_0} \left(\frac{c_0}{32 C_0 L_0 \Gamma_0}\right)^{1/\tbeta} .
\ee

Assume also that $b_2 \in W^s_R(z_0)$ with $d(y_0, b_2) < 2 \delta'$. Using (5.6) with $b_0 = b_2$, so that
$\|v_0 - \eta_0\| \leq R_0 \, d(y_0,b_2) < 2 R_0 \delta'$, we get 
\begin{eqnarray*}
|\Delta(x_0 , y_0) - \Delta(x_0, b_2) | 
& \leq &\|\uo_0\|  2 R_0 \delta' (C_0+ 2C_0 L_0) + C_4\, (\diam(\tcc))^{1+\tbeta}\\
& \leq & \|\uo_0\| 2 \delta' 3 R_0 C_0 L_0 + \frac{c_0 \delta_0}{8} \, \|\uo_0\| \leq \frac{c_0 \delta_0}{4} \,  \|\uo_0\| ,
\end{eqnarray*}
using the choice of $\delta'$ in (\ref{eq:delta'}).

Now (\ref{eq:x0y0}) and the above two estimates yield
\begin{eqnarray}
|\Delta(x_0 , b_1) - \Delta(x_0, b_2) |
& \geq & |\Delta(x_0 , y_0)  | - |\Delta(x_0 , y_0) - \Delta(x_0, b_2) | - |\Delta(x_0, b_1) | \nonumber\\
& \geq & 14 c_0 \delta_0 \,  \|\uo_0\|  - c_0 \delta_0 \, \|\uo_0\| = 13 c_0 \delta_0 \, \|\uo_0\| ,
\end{eqnarray}
that is (\ref{eq:Delta-cond}) holds.
This proves part (a).

\bs

(b) As in  Sect. 6 in \cite{D1}, we will use the fact that $\tpp^N(U)$ fills in $\tR$ densely as $N \to \infty$.

Let $\delta_0 > 0$ and $\delta' > 0$ be as in part (a). 
Take the integer $\tm_0 \geq 1$ so large that for every unstable
cylinder $\tcc$ in $\tR$ of length $m \geq \tm_0$ we have $\diam(\tcc) < \ep' = \frac{\ep_1}{2\tilde{c}_2}$; see (4.1).

Using the symbolic coding provided by the Markov family $\{\tR_i\}$ it is easy to see that
there exists an integer $N_0 \geq 1$ such that for any integer $N \geq N_0$  we have
$\tpp^N(B^u_{\ep'}(z)) \cap B^s(z', \delta') \neq \e$ for any $z,z'\in \tR$. We choose $N_0$ such that
$1/\gamma^{N_0} < \delta' .$

Let $\tcc$ be a cylinder of length $m \geq \tm_0$ in $\tR$ and let $z_0 \in \tcc \cap \tP_0$.
Let $x_0 \in \tcc$ and $u_0 \in E^u(z_0)$ satisfy (\ref{eq:u0}), where $\kappa \in (0,1)$ satisfies (\ref{eq:u0-cond}).
Define $N_0$ as above, and let $N \geq N_0$. 
Choose $y_0 = y_0(z_0,x_0)\in B^s(z_0, \ep_1)$ 
as in part (a).

It follows from the choice of $N$ that for each $i = 1,2$  there exist
$$y_1(z_0,x_0) \in \tpp^N(B^u(z_0,\ep')) \cap B^s(z_0, \delta')  \:\: , \:\:
y_2(z_0,x_0) \in \tpp^N(B^u(z_0,\ep')) \cap B^s(y_0, \delta') .$$
Fix points $y_1 = y_1(z_0,x_0)$, $y_2 = y_2(z_0,x_0)$ with these properties; these are then points in  $W^s_{\ep_1}(z_0)$.  
Let $b_1, b_2 \in W^s_\ep(z_0)$ be so that $b_1 \in B^s(y_1, \delta')$ and  $b_2 \in B^s(y_2, \delta')$.
Then $b_1 \in B^s(z_0, 2\delta')$ and  $b_2 \in B^s(y_0, 2\delta')$, so (\ref{eq:Delta-cond}) holds.
\endofproof

\bs

Assume the {\it integer $\tm_0 \geq 1$ is chosen so large that} for any $z\in R$ and any
unstable cylinder $\cc$ of length $\geq \tm_0$ in $R$ we have $\diam(\tPsi(\cc)) \leq r_0$
and $\diam(\T_z(\cc)) \leq r_0$ for any $z\in \cc$. We will use the constant $\delta_0$ defined 
by (\ref{eq:delta0}).

{\bf Fix a constant $d > 1$} such that
\be\label{eq:d}
T = \frac{1}{2} \left( \frac{2d \beta}{\tau_0} - 1\right) > \frac{1}{\tbeta}
\ee
is a very large constant, to be specified later,
where $\beta > 0$ is the constant from Sect. 5.1 and $\tbeta \in (0, \beta)$ is the constant from Lemma 5.3.
As before we assume that $\beta > 0$ is sufficiently small so
that it satisfies the requirements of Lemma 5.1 and also those in Sect. 3. In what follows 
we will use Lemma 4.4 with the choice of the constant $d$ in (\ref{eq:d}).

{\bf Let the constants $\tm_0 \geq 1$ and $q_1 \geq q_0$ be as in Lemma 4.5.}
{\bf As before, we choose $N_0$ so that $1/\gamma^{N_0} < \delta' .$} Recall the constant $\alpha_2 \in (0,1)$
from Sect. 4.2 and {\bf fix a constant $\alpha_3 \in (0,\alpha_2)$}.

\bs

\noindent
{\bf Lemma 5.5.} 
{\it  
There exist global constants 
$d_3 > 1$, $d_4 \in (0,1)$, $T > 0$ and $0 < \hep_{13} \leq \con \hep$  which can be made arbitrarily small with $\hep$  such that for 
any point $z_0 \in \tP_0$ and any cylinder $\tcc$ in $\tR$ of length $m \geq \tm_0$ containing  $z_0$ there exist 
sub-cylinders $\tGamma$ and $\tLambda$ 
of $\tcc$ having co-lengths $q_0$ and  $q_1$ in $\tcc$, respectively, so that 
the following hold:}

\ms

(a) {\it For every $x \in \tGamma$, for $u = (\Phi^u_{z_0})^{-1}(\T_{z_0}(x)) \in E^u(z_0)$ we have
\be\label{eq:kappa}
 \|\uo\| \geq   \kappa\, \diam(\tcc) ,
\ee
where}
 $\kappa = d_2 e^{- m \hep_{13} }$.

\ms

(b) {\it For every constant $S = S(m) \geq 1$, possibly depending on the cylinder $\tcc$,
there exist an integer $j_0 = j_0(m) \geq 1$ and finite families $\{\tGamma_j\}_{j=1}^{j_0}$ and 
$\{\tLambda_j\}_{j=1}^{j_0}$ of sub-cylinders of $\tcc$  such that 
\be\label{eq:covering}
\cup_{j=1}^{j_0} \tGamma_j \subset \tGamma \quad , \quad \cup_{j=1}^{j_0} \tLambda_j = \tLambda ,
\ee
$$\di \diam(\tGamma_j)  \leq \frac{1}{S} \diam(\tcc) \quad , \quad \diam (\tLambda_j)  \leq \frac{1}{S} \diam(\tcc) ,$$
and
\be\label{eq:diamte}
\diamte(\tGamma_j)  \leq \frac{1}{\tS} ((\diamte(\tcc))^{1/\alpha_3}
\quad , \quad \diamte (\tLambda_j)  \leq \frac{1}{\tS} (\diamte(\tcc))^{1/\alpha_3} ,
\ee
for all $j = 1, \ldots, j_0$, where 
$$\tS = \frac{C_1^{1+1/\alpha_3}}{S^{\alpha_2}} .$$
We can choose all sub-cylinders $\tGamma_j$ and $\tLambda_j$ so that their lengths do not exceed $T \, m$ and
\be
\frac{1}{d_3} \leq \frac{\nu(\tGamma'_j)}{\nu(\tLambda'_j)} \leq d_3
\ee
for all $j = 1, \ldots, j_0$, where $\tGamma'_j = \piU(\tGamma_j)$ and $\tLambda'_j = \piU(\tLambda_j)$.}

\ms

(c)  {\it Moreover, for the sub-cylinders of $\tcc$ from parts (a) and (b),
any integer $N \geq N_0$, any $j = 1, \ldots, j_0$ and any  $i = 1,2$ there exist  a  (H\"older) continuous map 
$$B^u(z_0,\ep'') \ni x \mapsto v_{i,j} (z_0, x) \in U ,$$
such that $\sigma^N(v_{i,j}(z_0, x)) = x$  for all $x\in B^u(z_0,\ep'')$
and the following property holds: 
\be\label{eq:IN}
I_{N}(x',z') =  \left|  \psi(z_0,x') - \psi(z_0,z') \right|   \geq 6 c_0 \delta_0 \, \kappa\, \diam(\tcc) 
\ee
for all $z\in \tLambda_j$ and  $x \in \tGamma_j$, where $z' = \piU(z)$, $x' = \piU(x) \in U$, and}
$$\psi(z_0, x) = \tau_{N}(v_{1,j}(z_0,x)) - \tau_{N}(v_{2,j}(z_0 , x)) .$$

\bs

\noindent
{\it Proof.}  
We will now use the results of Lemma 4.5 with the constant $d > 1$ as in (\ref{eq:d}) and some of the set-up in its proof.

Define the global constant $\delta' \in (0,1)$ by (\ref{eq:delta'}) and choose the global constant  $N_0 \geq 1$ as in the proof of Lemma 5.4(b).

Take an arbitrary constant $C' > 8$, 
and set 
$\di C = \frac{3C'}{c_0 \delta_0 d_2} > 2 C',$
where $\delta_0 > 0$ is as in (\ref{eq:delta0}) and $c_0 > 0$ is as in Sect. 5.1. {\bf We will use Lemma 4.5 with this particular $C$.}

Fix for a moment  $z_0 \in \tP_0 $.   Assume $z_0\in \tR_{i_0}$. 
Let $\tcc$ be a cylinder of length $m \geq m_0$ with $z_0\in \tcc$, $Z_m =  \tpp^m(z_0) \in  \tR_{j_0}$ for some $j_0$.
Then choose $k$ with $0 < k \leq m$ with (4.28), set $p = [\ttau_{m-k}(z_0)]$
as in the proof of Lemma 4.5, and set $q = [\ttau_m(z_0)]$, $z_j = \varphi^j(z)$, $Z_j = \tpp^j(z_0)$,  etc., as in the proof of Lemma 4.5. 
As in the beginning of Sect. 4.4, consider the  locally defined near $z_0$ map
$$(\F0)^{-1}\circ \hvarphi_{z_0}^{-p} \circ \Fp : \tcc \longrightarrow \tcc ,$$
where
$$\Fp = (\Phi^u_{z_p})^{-1} \circ \T_{z_p} \circ \varphi^p : \tcc \longrightarrow E^u (z_p) \quad, \quad
\F0 = (\Phi^u_{z_0})^{-1} \circ \T_{z_0} : \tcc \longrightarrow E^u (z_0) .$$
Define the sets $H_1$, $H_2$, $\tH_1$ and $\tH_2$ as in the proof of Lemma 4.5.
It follows from the latter that there exists sub-cylinders $\tGamma_0$, $\tTheta$ and $\tLambda$ of $\tcc$
with co-lengths $q_0$, $q_2$ and $q_1$, respectively, in $\tcc$ with the properties
(\ref{eq:in-cond}) and (\ref{eq:out-cond}). Then $\tTheta \subset \tLambda \subset \tGamma_0$.

In what follows we will use arguments from the proofs of Lemmas 5.3 and 5.4.
It follows from the assumptions about $\tGamma$ in Lemma 4.5 that
(\ref{eq:kappa}) holds for  $u = (\Phi^u_{z_0})^{-1}(\T_{z_0}(x))$ for every $x \in \tGamma$.

Next, consider an arbitrary $z \in \tLambda$ so that it belongs to $\varphi^{-p}((\T_{z_p})^{-1}(\Phi^u_{z_p}(H'_1\setminus \tH_2)))$
and the corresponding 
$w_0 = (\Phi^u_{z_0})^{-1}(\T_{z_0}(z)) .$
Set
$$w_j = \hvarphi^j_{z_0}(w_0) \in E^u(z_j) \quad, \quad \hw_j = \hvarphi^j_{z_0}(0) \cdot w \in E^u(z_j)$$
for all $j = 0,1, \ldots, p$. 
It follows from our choice of $z$ that $w_p \in H'_1\setminus \tH_2$, so
\be
\|\wo_p\| \geq \frac{d_0 \ttc_1}{4C \gamma_1} e^{- 4d (m-k) \hep} .
\ee
It follows from (\ref{eq:in-cond}) that
there exists $x \in \tGamma$ so that 
$$C w_p = u_p = \hvarphi^p_{z_0}(u_0) ,$$
where 
$u_0 = (\Phi^u_{z_0})^{-1}(\T_{z_0}(x))$ satisfies (\ref{eq:kappa}), i.e.
$\|\uo_0\| \geq \kappa\, \diam(\tcc)$. Now, as in the proof of Lemma 5.4,  choose $\tv \in E^s_1(z_0)$ 
with $\|\tv\| = 1$ and $d\omega_{z_0}(\uo_0, \tv) \geq 2 c_0\, \|\uo_0\|$,
and then define $v_0 \in E^s_1(z_0)$ and $y_0 \in W^s_{\ep_1}(z_0)$ as in (\ref{eq:vo}). 

In most of what follows we will deal with the {\bf fixed $z\in \tLambda$ and the corresponding $x = x(z) \in \tGamma$}
which is related to $z$ as above; namely we determined successively: 
$$w_0 = (\Phi^u_{z_0})^{-1}(\T_{z_0}(z)) \quad, \quad
w_p = \hvarphi^p_{z_0}(w_0) \in E^u(z_p) \quad, \quad u_p = Cw_p = \hvarphi^p_{z_0}(u_0)$$
for some $u_0 \in E^u(z_0)$ and $u_0 = (\Phi^u_{z_0})^{-1}(\T_{z_0}(x))$ which determines $x$.
Later on we will consider neighbourhoods of $z$ and $x$ which will be corresponding little sub-cylinders of $\tLambda$ and $\tGamma$.

As in the proof of Lemma 5.3 we will use the notation
$$x_j = \varphi^j(x_0) \quad, \quad u_j = \hvarphi^j_{z_0}(u_0) \in E^u(z_j) \quad, 
\quad \hu_j = \hvarphi^j_{z_0}(0) \cdot u_0 \in E^u(z_j) 
\quad , \quad v_j = \hvarphi^j_{z_0}(v_0) \in E^s(z_j) .$$

{\bf Note.} Notice that $z_j = \varphi^j(z_0)$ are the iterates of the fixed point $z_0 \in \tP_0$, not those of the arbitrary point $z\in \Lambda$.

\ms

It follows from 
$u_p = C w_p \in H\setminus H_2$, the definition of $H_2$ (see the proof of Lemma 4.5) 
and $p = [\ttau_{m-k}(z_0)]  \in [ (m-k) \ttau_0, (m-k) \tau_0 + 1)$  that
\be
\|u_p\| \leq  e^{- 3d (m-k) \hep}  \leq  e^{- 3 d (p-1) \hep/\tau_0}  \leq  e^{- 2 d p \hep/\tau_0} \quad, \quad  
\|\uo_p\| \geq \frac{d_0 \ttc_1}{4 \gamma_1} e^{- 4d (m-k) \hep} \geq   \frac{d_0 \ttc_1}{4 \gamma_1} e^{- 4 d p \hep/\ttau_0} ,
\ee
assuming $p$ is sufficiently large.
Similarly for $w$ and $w_p$ we have
\be
\|w_p\| \leq \frac{1}{C}e^{- 2 d p \hep /\tau_0} \quad, \quad  \|\wo_p\| \geq \frac{d_0 \ttc_1}{4C \gamma_1} e^{-4 d p\hep/ \ttau_0} .
\ee

Next, it follows from Lemma 3.1  that
$\|\uo_p - \huo_p\| \leq L(z_p) \|u_p\|^{1+ \beta} \leq L_0 e^{p\hep} \|u_p\|^{1+\beta} ,$
so  $\|\huo_p\| \leq \|u_p\| (1+ L_0 e^{p\hep} \|u_p\|^\beta)$. Now (5.50) and (\ref{eq:d}) yield
$$L_0 e^{p \hep} \|u_p\|^\beta \leq L_0 e^{p\hep} e^{-2 d \beta p \hep/\tau_0}
\leq L_0 e^{- p \hep (2d\beta /\tau_0 - 1)} < L_0 e^{-p\hep}  < 1 , $$
since $2d \beta/\tau_0 -1 > 1$.
Thus,
$\|\huo_p\| \leq 2 \|u_p\| \leq 2 e^{- 2 d p\hep/\tau_0}$, and
\be
\|\uo_p - \huo_p\| \leq L_0 e^{- p \hep (2d \beta/\tau_0 - 1)} = L_0 e ^{- 2T p \hep} .
\ee

Similarly, again by Lemma 3.1,
$$\|\wo_p - \hwo_p\| \leq L(z_p) \|w_p\|^{1+ \beta} \leq L_0 e^{p\hep} \|w_p\|^{1+\beta} \leq  \frac{L_0}{C} e^{- 2T p \hep }  ,$$
and as above we derive
$$\|\hwo_p\| \leq 2 \|\wo_p\| \leq \frac{2}{C} e^{- 2 d p\hep/\tau_0} .$$
Since $u_p = C w_p$, it follows that $\uo_p = C \wo_p$, therefore
$$\|\huo_p - C \hwo_p\| \leq \|\huo_p - \uo_p\| + \|\uo_p - C \wo_p\| + \|C \wo_p - \hwo_p\|
= \|\huo_p - \uo_p\| + \|C \wo_p - C \hwo_p\| ,$$
and therefore
$$\|\huo_p - C \hwo_p\| \leq 2 L_0 e^{-2T p \hep} .$$
Notice that the above and (3.14) imply
\be\label{eq:wo}
\|\uo_0 - C \wo_0\| \leq \Gamma_0 e^{p \hep} \frac{2 L_0 e^{-2Tp \hep}}{\mu_1^p} \leq \frac{e^{-T p \hep}}{\lambda_1^p} ,
\ee
assuming that $T$ is sufficiently large.

We will now use (\ref{eq:wo}) to get an estimate for $|\Delta(z, y_0)|$ by means of $\|\uo_0\|$.
It follows from (5.53), (5.7) with $x_0$ replaced by $z$ and $\uo_0$ replaced by $\wo_0$ that
\begin{eqnarray}
|\Delta(z ,y_0)|
& \leq &  |d\omega_{z_0} (\wo_0, \vo_0 )| + C_0 L_0 \|\wo_0\| \, \|v_0\|^{1+\beta} + C_4 \, (\diam(\tcc))^{1+ \tbeta}\nonumber\\
& \leq & \frac{1}{C} |d\omega_{z_0} (\uo_0, \vo_0 )| + \frac{ C_0 L_0}{C} \|\uo_0\| \, \|v_0\|^{1+\beta} 
+ C_4 \, (\diam(\tcc))^{1+ \tbeta}\nonumber\\
&        & + \frac{1}{C} \frac{e^{-T p \hep}}{\lambda_1^p} + \frac{C_0 L_0}{C} \frac{e^{-T p \hep}}{\lambda_1^p} .
 \end{eqnarray}
Since $u$ satisfies (\ref{eq:kappa}), as in the proof of Lemma 5.4 we have (\ref{eq:C4-cond}), and therefore 
\be
C_4 \, (\diam(\tcc))^{1+ \tbeta} < \frac{c_0 \delta_0}{8} \, \kappa \, \diam(\tcc) 
\leq \frac{c_0 \delta_0}{8} \, \|\uo_0\| .
\ee
For the  second term in the right-hand-side of (5.54), it follows from (\ref{eq:d}), (\ref{eq:vo}) and (5.30) that
\be
C_0 L_0 \|\uo_0\| \|v_0\|^{1+\beta} \leq C_0L_0 \|\uo_0\| 10 \delta_0 \left(\frac{c_0}{32 C_0 L_0 \Gamma_0}\right)
< c_0 \delta_0 \, \|\uo_0\| \leq c_0 \delta_0 \,  \kappa\, \diam(\tcc) .
\ee

Next, notice that if we choose $T > 1$ sufficiently large, then
\be
4 C_0 L_0 \Gamma_0 \frac{e^{- T p \hep}}{\lambda_1^{p}} < \frac{c_0 \delta_0}{8} \,\kappa\, \diam(\tcc) \leq \frac{c_0 \delta_0}{8} \, \|\uo_0\| .
\ee
Indeed, using $q (1-\hep'_6) \leq p \leq q$ by (4.30) and $m \leq q \tau_0+1 < 2p\tau_0$, we get 
$$\kappa = d_2 e^{-m\hep_{13}} \geq d_2 e^{-2p\tau_0 \hep_{13}} .$$
Now it follows from (\ref{eq:diamC}) that 
\begin{eqnarray}
\frac{c_0 \delta_0}{8} \, \kappa\, \diam (\tcc) 
& \geq  &  c_0 \delta_0 \, d_2 e^{- 2p \tau_0 \hep_{13}} \frac{e^{-q \hep_7}}{8 C_3 \lambda_1^q} 
 >  c_0 \delta_0 \,  d_2 e^{- 2 p \tau_0 \hep_{13} } \frac{e^{-2p \hep_{7}} }{8 C_3 \lambda_1^{p/(1- \hep'_6)}} \nonumber\\
& =  &  c_0 \delta_0 \, d_2  \frac{e^{- 2 p (\tau_0 \hep_{13} + \hep_7)} \lambda_1^{- p\hep'_6 /(1-\hep'_6)}}{8 C_3\lambda_1^{p} }\nonumber\\
& =  & c_0 \delta_0 \,  d_2  \frac{e^{- p ( 2\tau_0 \hep_{13} + 2\hep_{7} + \hep'_6 (\log \lambda_1)/(1-\hep'_6) )} }{8 C_3\lambda_1^{p}} 
> 4 C_0L_0 \Gamma_0  \frac{e^{- T p \hep}}{\lambda_1^{p}} ,
\end{eqnarray}
assuming e.g. that $\frac{32 C_0L_0\Gamma_0 C_3}{c_0 \delta_0 \, d_2} < e^{T p \hep/2}$ and
$T  \hep/2 > 2\tau_0 \hep_{13} + 2 \hep_{7} + \hep'_6 (\log \lambda_1)/(1-\hep'_6)$, which will be satisfied if we take $T$ sufficiently large.
Thus, (5.57) holds.

Combining (5.54), (5.55), (5.56) and (5.57) implies
\begin{eqnarray}
|\Delta(z ,y_0)| 
& \leq & \frac{1}{C} |d\omega_{z_0} (\uo_0, \vo_0 )| + \frac{c_0 \delta_0}{C}  \,  \|\uo_0\|
+  \frac{c_0\delta_0}{8} \,  \|\uo_0\| +  2 \frac{c_0\delta_0}{8} \,  \|\uo_0\|\nonumber\\
& \leq & \frac{10C_0 \delta_0}{C}  \,  \|\uo_0\| + \frac{c_0 \delta_0}{C}  \|\uo_0\| +  3\frac{c_0 \delta_0}{8} \, \|\uo_0\|
 \leq  c_0 \delta_0 \, \|\uo_0\| ,
\end{eqnarray}
assuming e.g. that $C \geq \frac{20}{c_0 C_0}$. 

Before we continue we will derive estimates similar to (5.40) and (5.42) replacing $x_0$ by $z$. Assume as before that
$b_1, b_2 \in W^s_R(z_0)$ are such that $d(z_0, b_1) < \delta'$ and $d(y_0, d_2) < \delta'$, where $\delta'$ is defined by (\ref{eq:delta'}). 
Using (5.6) with $x_0$ replaced by
$z$ and $\uo_0$ by $\wo_0$, and taking into account (\ref{eq:wo}) and (5.55) and using the argument in the estimate of 
$|\Delta (x_0,y_0) - \Delta(x_0, b_2)|$ in the proof of Lemma 5.4, we get
\begin{eqnarray*}
|\Delta(z , y_0) - \Delta(z, b_2) | 
& \leq &\|\wo_0\| R_0 \delta' (1+ 2C_0 L_0) + C_4\, (\diam(\tcc))^{1+\tbeta}\\
& \leq & \frac{1}{C}\left(\|\uo_0\| + \frac{e^{-Tp\hep}}{\lambda_1^p}\right) \delta' 3 R_0 C_0 L_0 + C_4\, (\diam(\tcc))^{1+\tbeta}\\
& \leq & \frac{c_0 \delta_0\, \|\uo_0\|}{C}  + \frac{e^{-Tp \hep}}{C \lambda_1^p} + \frac{c_0 \delta_0}{8} \|\uo_0\| .
\end{eqnarray*}
It follows from (5.57), assuming that the constant $T > 0$ is chosen sufficiently large, that
$$\frac{c_0 \delta_0}{4 C} \|\uo_0\| \geq \frac{c_0 \delta_0\,  \kappa}{4 C} \diam(\tcc) 
\geq \frac{4 C_0 L_0 \Gamma_0}{4 C} \frac{e^{-Tp \hep}}{\lambda_1^p}  > \frac{1}{C} \frac{e^{-Tp \hep}}{\lambda_1^p}.$$
Therefore
\be
|\Delta(z , y_0) - \Delta(z, b_2) | \leq \frac{c_0 \delta_0 \,  \|\uo_0\|}{C} + \frac{c_0 \delta_0}{4C} \|\uo_0\| + \frac{c_0 \delta_0}{8} \|\uo_0\|
< \frac{c_0\delta_0}{2} \|\uo_0\| .
\ee

Similarly, using the estimates in the proof of (5.40), replacing $x_0$ by $z$ and $\uo$ by $\wo$, we get
\begin{eqnarray*}
|\Delta(z,b_1)| 
 \leq  2 C_0 R_0 \delta' \|\wo_0\| + \frac{1}{4} c_0 \|v_0\| \|\wo_0\|  \leq  \frac{1}{2} \delta_0\, \|\wo_0\| .
\end{eqnarray*}
On the other hand, (\ref{eq:wo}) and (5.57) yield
$$\|\wo_0\| \leq \frac{1}{C} \left(\|\uo_0\| + \frac{e^{-T p \hep}}{\lambda_1^p}  \right) 
\leq \frac{1}{C} \left(\|\uo_0\| + \frac{c_0\delta_0}{4} \, \|\uo_0\| \right) < \frac{2}{C} \|\uo_0\| . $$
Combining this with the previous estimates, we get
$$|\Delta(z,b_1)|  \leq \frac{\delta_0}{C} \|\uo_0\| \leq \frac{c_0 \delta_0}{2} \|\uo_0\| ,$$
assuming $C$ is sufficiently large, as before.
This, (5.59) and (5.60) now imply
\begin{eqnarray}\label{eq:twob}
|\Delta(z , b_1) - \Delta(z, b_2) | 
& \leq & |\Delta(z , b_1)|  + |\Delta(z , y_0)| +  |\Delta(z , y_0) - \Delta(z, b_2) | \nonumber\\
& \leq & \frac{c_0 \delta_0}{2} \|\uo_0\|  + c_0 \delta_0 \|\uo_0\| + \frac{c_0 \delta_0}{2} \|\uo_0\| 
\leq  2 c_0 \delta_0 \,  \|\uo_0\| .
\end{eqnarray}

Next, we will get estimates similar to (5.42) and (\ref{eq:twob}) replacing\footnote{We avoid here the notation $x'$ and $z'$ since 
these mean something specifically related to projections to $U$.}
 $x$ by $x'' \in \tcc$ close to $x$ and $z$ by $z'' \in \tcc$ close to $z$.
Recall that the function $\Delta(x,y_0)$ is uniformly H\"older, so there exist constants $D'_3 > 0$ and $\alpha \in (0,1]$ such that
$$|\Delta(x,y_0) - \Delta(x'',y_0)| \leq D'_3 (d(x,x''))^\alpha 
\quad, \quad |\Delta(z,y_0) - \Delta(z'',y_0)| \leq D'_3 (d(z,z''))^\alpha .$$
We {\bf fix an $\alpha$ with this property} so that $\alpha \in (0, \alpha_2 \alpha_3]$, where $0 < \alpha_3 < \alpha_2 < 1$ are as chosen just
before Lemma 5.5.
Fix for a moment a small constant $\delta > 0$ (we will determine later how small) and assume
$$d(x,x'') \leq \delta \quad, \quad d(z,z'') \leq \delta$$
for some $x'', z''\in \tcc$. 
It follows from (\ref{eq:diamC}) that we have the following estimate (in terms of $m$ now):
\begin{eqnarray*}
\kappa \, \diam (\tcc) 
 \geq   d_2 e^{-m \hep_{13}}  \frac{e^{- m\hep_7/\ttau_0}}{C_3 \lambda_1^{m\tau_0}}
= \frac{d_2}{C_3} \, e^{-m (\hep_{13}  + \hep_7/\ttau_0 + \tau_0 \log \lambda_1)}
 \geq   \frac{d_2}{C_3} \, e^{-m (1 + 1/\ttau_0 + \log \lambda_1)} ,
\end{eqnarray*}
where we used the fact that the constant $\hep_{13}$, $\hep_7$ and $\tau_0$ are all in $(0,1)$.

{\bf Assume now that $S = S(m) \geq 1$ is a given constant.}
We will determine $\delta$ so that
$$\di \frac{c_0 \delta_0\,  \kappa}{S} \,  \diam(\tcc) \geq 8 D'_3 \delta^\alpha ,$$
that is
\be\label{eq:delta}
\delta \leq \left(\frac{c_0 \delta_0}{8 \, S\, D'_3} \kappa \, \diam(\tcc) \right)^{1/\alpha} ,
\ee
for which it is enough to have
$\di \delta \leq D''_3 \, e^{-m (1 + 1/\ttau_0 + \log \lambda_1)/ \alpha}$,
where $D''_3 = \left(\frac{c_0 \delta_0 d_2}{8 S \,D'_3\, C_3} \right)^{1/\alpha}  > 0 $ is a global constant.
It follows from Lemma 4.1(c) that the sub-cylinders $\tGamma(x)$ and $\tLambda(z)$
of $\tcc$, containing  $x$ and $z$, respectively, and having lengths not less than
$$\tn = \frac{1}{|\log\theta|} ( - \log \delta + \log C_2 )$$
have diameters not exceeding $\delta$. Indeed, if $\tX$ is a cylinder in $\tcc$ of length $n \geq \tn$, then by
Lemma 4.2(c), 
$$\diam(\tX) \leq C_2 \, \diamte(\tX) = C_2 \theta^n = C_2 e^{-n |\log\theta|} \leq C_2 e^{\log \delta - \log C_2} = \delta .$$
Since
$$- \log \delta \geq m (1 + 1/\ttau_0 + \log \lambda_1)/ \alpha - \log D''_3 ,$$
it is enough to have $\tn \geq r'_m = T'_1 m + T'_2 ,$
where $T'_1 = \frac{1}{\alpha |\log\theta|} (1 + 1/\ttau_0 + \log \lambda_1) > 0$ and  $T'_2 = \frac{\log C_2}{|\log\theta|}$.

In this way we have shown that if  the sub-cylinders $\tGamma(x)$ and $\tLambda(z)$
of $\tcc$, containing  $x$ and $z$, respectively, having lengths $\geq r'_m$, then they have diameters not exceeding 
the number $\delta$ satisfying (\ref{eq:delta}). Since $0 < \alpha \leq \alpha_2\alpha_3 < 1$,  and $\kappa \leq 1$ for sufficiently large $m$,
(\ref{eq:delta}) implies
$$\diam(\tGamma(x)) 
\leq \frac{1}{S} (\diam (\tcc))^{1/(\alpha_2 \alpha_3)} \leq \frac{1}{S} \diam (\tcc) ,$$
and similarly $\diam(\tLambda(z))  \leq \frac{1}{S} (\diam (\tcc))^{1/(\alpha_2 \alpha_3)} \leq \frac{1}{S} \diam (\tcc)$.
Moreover,
it follows from Lemma 4.1 and the above that
$$\diamte(\tGamma(x)) \leq C_1\, (\diam(\tGamma(x))^{\alpha_2}
\leq  \frac{C_1}{S^{\alpha_2}} (\diam (\tcc))^{1/\alpha_3}
\leq  \frac{C_1}{S^{\alpha_2}} (C_1\, \diamte (\tcc))^{1/\alpha_3} 
\leq  \frac{1}{\tS} (\diamte (\tcc))^{1/\alpha_3} ,$$
where $\tS = \frac{C_1^{1+1/\alpha_3}}{S^{\alpha_2}} $. A similar estimate holds for $\tLambda(z)$, so
\be\label{eq:diamte}
\diamte(\tGamma(x)) 
\leq \frac{1}{\tS} (\diamte (\tcc))^{1/\alpha_3} \quad, \quad \diamte(\tLambda(z)) \leq  \frac{1}{\tS} (\diamte (\tcc) )^{1/\alpha_3} .
\ee

With such a choice of $\tGamma(x)$ and $\tLambda(z)$  it follows from (\ref{eq:twob}) that
for any $b_1, b_2 \in W^s_R(z_0)$ with $d(z_0, b_1) < \delta'$ and $d(y_0, b_2) < \delta'$ we have
\be
|\Delta(z'' , b_1) - \Delta(z'', b_2) | \leq 3 c_0 \delta_0\,  \|\uo_0\| 
\ee
for all $z'' \in \Lambda(z)$. Indeed, by the choice of $\delta$, for every $z''\in \tLambda(z)$ we have
$$|\Delta(z'',b_1) - \Delta(z,b_1)| \leq D'_3 \delta^\alpha
\quad, \quad |\Delta(z'', b_2) - \Delta(z, b_2)| \leq D'_3 \delta^\alpha ,$$
and now (\ref{eq:twob}) implies
$$|\Delta(z'' , b_1) - \Delta(z'', b_2) | \leq |\Delta(z , b_1) - \Delta(z, b_2) | + 2 D'_3 \delta^\alpha 
\leq 2 c_0 \delta_0 \|\uo_0\| + \frac{1}{4} c_0 \delta_0\, \kappa\, \diam (\tcc) \leq 3 c_0 \delta_0 \|\uo_0\| .$$

In a similar way, for $b_1, b_2$ as above, we can use (5.42) to derive
\be
|\Delta(x'',b_1) - \Delta(x'' , b_2)| \geq 12 c_0 \delta_0\,  \|\uo_0\|
\ee
for all $x'' \in \Gamma(x)$. Indeed, using (5.42) with $x_0$ replaced by $x$ and
$$|\Delta(x'',b_1) - \Delta(x,b_1)| \leq D'_3 \delta^\alpha
\quad, \quad |\Delta(x'', b_2) - \Delta(x, b_2)| \leq D'_3 \delta^\alpha ,$$
we get
$$|\Delta(x'' , b_1) - \Delta(x'', b_2) | \geq |\Delta(x, b_1) - \Delta(x, b_2) | - 2 D'_3 \delta^\alpha 
\geq 13 c_0 \delta_0 \,  \|\uo_0\| - \frac{1}{4} c_0 \delta_0 \,  \kappa\, \diam (\tcc)  \geq 12 c_0 \delta_0 \, \|\uo_0\| .$$

We will now use the above as follows.  Take an arbitrary {\bf $\hz\in \tLambda$ and fix it.}
Using the above with $z = \hz$, we construct a corresponding $x_0 = x_0(\hz) \in \tGamma$ as follows: define 
$w_0 = (\Phi^u_{z_0})^{-1}(\T_{z_0}(\hz))$, then $w_p = \hvarphi^p_{z_0}(w_0) \in E^u(z_p)$. This defines
$u_p = Cw_p = \hvarphi^p_{z_0}(u_0)$ for some $u_0 \in E^u(z_0)$. Then define $x_0 = x_0(\hz) = (\T_{z_0})^{-1} (\Phi^u_{z_0}(u_0))$
so that $u_0 = (\Phi^u_{z_0})^{-1}(\T_{z_0}(x_0))$. 

Next, we will use a construction from the proof of Lemma 5.4. 
Take a unit  vector $\tv \in E^s_1(z_0)$ with $d\omega_{z_0}(\uo_0, \tv) \geq 2 c_0 \|\uo_0\|$, and then define $v_0 \in E^s_1(z_0)$ and
$y_0 = \Phi^s_{z_0}(v_0) \in W^s_{\ep_1}(z_0)$ as in (\ref{eq:vo}). 
Then (\ref{eq:Delta-cond}) holds, that is there exist  a point $y_0 = y_0(z_0, x_0)\in  B^s (z_0,\ep_1)$ such that we have
\be\label{eq:DeltaAgain-cond}
 |\Delta(x_0, b_1)  - \Delta(x_0, b_2)| \geq 12 c_0 \delta_0 \, \|\uo_0\| \geq 12 c_0 \delta_0 \,\kappa\,   \diam (\tcc) 
\ee
for any $b_1,b_2 \in W^s_R(z_0)$ with $d(z_0,b_1) < 2\delta'$ and $d(y_0, b_2) < 2\delta'$.

Similarly, (\ref{eq:twob}) holds with $z$ replaced by $\hz$ for such $b_1,b_2$.

Next, fix a large integer $\tq \geq 1$.
We now choose a sub-cylinder $\tGamma(x_0)$ of $\tcc$ of co-length $\leq \tq$ with $x_0 \in \tGamma(x_0)$ and a 
sub-cylinder $\tLambda(\hz)$ of $\tcc$ of co-length $\leq \tq$ with $\hz \in \tLambda(\hz)$ 
both of length $r'_m$ so that for all $b_1,b_2 \in W^s_R(z_0)$ with $d(z_0,b_1) < \delta'$ and $d(y_0, b_2) < \delta'$,
(5.64) holds with $z''$ replaced by any $z \in \tLambda(\hz)$ and 
(5.65)  holds with $x''$ replaced by any $x\in \tGamma(x_0)$. It then follows that, for such $b_1,b_2$, we have
\begin{eqnarray}
&         &      \big |\,     \left|\Delta(x , b_1) - \Delta(x , b_2)\right| - \left|\Delta(z  , b_1) - \Delta(z, b_2)\right| \, \big|\nonumber\\
& \geq &  12 c_0 \delta_0 \, \|\uo_0\|  -  3 c_0 \delta_0 \,  \kappa \, \diam(\tcc)  \geq 9 c_0 \delta_0 \, \|\uo_0\|
\geq 9 c_0 \delta_0 \, \kappa \, \diam(\tcc)
\end{eqnarray}
for all $x \in \tGamma(x_0)$ and all $z \in \tLambda(\hz)$. {\bf Notice that we now no longer require that the sub-cylinders 
$\tGamma(x_0)$ and $\tLambda(\hz)$ are related by $C \hvarphi^p_{z_0}(\F0(\tLambda(\hz))) \subset \Fp(\tGamma(x_0)$;}
only the "central points $x_0$ and $\hz$ satisfy $C(\F0(\hz)) = x_0$. On the rest of the sub-cylinders $\tGamma(x_0)$
and $\tLambda(\hz)$ we just use the fact that each of them has diameter $< \delta$. 

We can now apply the above construction to the whole set $\tLambda$ obtained from Lemma 4.5, 
which itself is a union of sub-cylinders of $\tcc$. Given an arbitrary $\hz \in \tLambda$ take $x = x(\hz) \in \tGamma$
so that $C(\F0(\hz)) = x(\hz)$, and then consider sub-cylinders  $\tLambda(\hz)$ 
of $\tcc$ of co-length $\leq \tq$ and total length $\geq r'_m$ and a corresponding sub-cylinder $\tGamma(x(\hz))$ 
of $\tcc$ of co-length $\leq \tq$ and total length $\geq r'_m$ contained in  $\tGamma$ such that
$$C  \hvarphi^p_{z_0}(\F0 (\tLambda(\hz))) \subset \Fp (\tGamma(x (\hz)))$$
and so that the following condition is satisfied:

\ms

{\bf Binding condition:} {\it there exist a point $y = y(\hz) \in  B^s (z_0,\ep_1)$ such that
for all $b_1,b_2 \in W^s_R(z_0)$ with $d(z_0,b_1) < \delta'$ and $d(y, b_2) < \delta'$ we have
\begin{eqnarray}\label{eq:bin-cond}
  \big| \,  \left|\Delta(x , b_1) - \Delta(x , b_2)\right| - \left|\Delta(z  , b_1) - \Delta(z, b_2)\right| \,\big|
\geq 9 c_0 \delta_0 \, \kappa \, \diam(\tcc)
\end{eqnarray}
for all $x \in \tGamma(x(\hz))$ and all $z \in \tLambda(\hz)$. }

\ms

Let us make the obvious remark that if $\tGamma(x(\hz))$ and $\tLambda(\hz)$ are as above, then 
for any sub-cylinder $\tGamma_1$ of $\tGamma(x(\hz))$ and any sub-cylinder $\tLambda_1$ of $\tLambda(\hz)$,
(\ref{eq:bin-cond}) holds for all $z\in \tLambda_1$ and all $x\in \tGamma_1$.

\ms

In general $\nu(\hLambda(\hz))$ and  $\nu(\hGamma(x(z)))$ will not be the same. We apply the  following

\ms

{\bf Procedure:}

There are two possible cases to consider.

\ms

\noindent
{\bf Case 1.} $\nu(\tLambda'(\hz)) \leq \nu(\tGamma'(x(z)))$, where $\tLambda'(\hz) = \piU(\tLambda(\hz))$ and
$\tGamma'(x(z)) = \piU(\tGamma(x(z)))$.
Denote by $\tGamma_1(x(z))$ the sub-cylinder of  $\tGamma(x(z))$ containing $x(z)$ and
having maximal possible length so that $\nu(\tLambda'(\hz)) \leq \nu(\tGamma'_1(x(z)))$. 
If $t$ is the length of $\tGamma_1(x(z))$, then its sub-cylinder
$\tGamma_2(x(z))$ of length $t+1$ containing $x(z)$ satisfies $\nu(\tGamma'_2(x(z))) < \nu(\tLambda'(\hz))$, 
so using the projection of these cylinders to $U$ via $\piU$ and the point $x' = \piU(x(z))$, (4.2) implies
$$1 \leq  \frac{\nu(\tGamma'_1(x(z)))}{\nu(\tGamma'_2(x(z)))} \leq \frac{c_2}{c_1} \frac{e^{g_t(x')}}{e^{g_{t+1}(x')}}
= \frac{c_2}{c_1} \frac{1}{e^{g(\sigma(x'))}} \leq d_3$$
for some global constant $d_3 > 1$. 
Moreover it is easy to see, using (4.2), the fact that $\tGamma'_1(x(z))$ has maximal
possible length with $\nu(\tLambda'(\hz)) \leq \nu(\tGamma'_1(x(z)))$, and the fact that
the length of $\tLambda(\hz)$ is $r'_m = T'_1 m + T'_2$, that the length $t$
of $\tGamma_1(x(z))$ is still bounded above by  $r_m = T_1 m + T_2$ for some global constants $T_1, T_2 > 0$.
We now replace $\tGamma(x(z))$ by $ \tGamma_1(x(z))$. 
Then the pair of cylinders $(\tLambda(\hz), \tGamma_1 (x(z)))$ satisfies
$$\frac{1}{d_3} < 1 \leq  \frac{\nu(\tGamma'_1(x(z)))}{\nu(\Lambda'(\hz))} \leq d_3 .$$
As mentioned above, (\ref{eq:bin-cond}) clearly holds for all $z\in \tLambda (\hz)$ and all $x\in \tGamma_1(x(z))$.

\ms

\noindent
{\bf Case 2.} $\nu(\tLambda'(\hz)) > \nu(\tGamma'(x(z)))$. Repeating the argument in Case 1 above, changing the roles of $\tLambda(\hz)$
and $\tGamma(x(z))$, we construct a sub-cylinder $\tLambda_1(\hz)$ of $\tLambda(\hz)$ containing $\hz$ of length not exceeding $r_m$
and such that
$$\frac{1}{d_3} < 1 \leq  \frac{\nu(\tLambda'_1(\hz))}{\nu(\tGamma'(x(z)))} \leq d_3 .$$
We then replace $\tLambda(\hz)$ by $\tLambda_1(\hz)$, and consider the pair $\tLambda_1(\hz)$, $\tGamma(x(z))$.
Again, (\ref{eq:bin-cond}) holds for all $z\in \tLambda_1 (\hz)$ and all $x\in \tGamma(x(z))$.

\ms

Thus, reducing either $\tLambda(\hz)$ or $\tGamma(x(\hz)$ we get sub-cylinders of $\tcc$ of
of co-length $\leq \tq$ and total length $\geq r'_m$, containing $\hz$ and $x(\hz)$, respectively,
have total length $\geq r'_m$, satisfy the Binding condition and (5.47) as well for some global constant $d_3 > 1$.

\ms

Applying the above Procedure to each of the initial pairs $(\tLambda(\hz), \tGamma (x(z)))$ of cylinders in $\tLambda$ and $\tGamma$,
we construct families $\tLambda_j = \tLambda_j(z_0^{(j)})$, $\tGamma_j = \tGamma_j(x_0^{(j)})$ ($1 \leq j \leq j_0$)
of sub-cylinders of $\tcc_m$ of lengths $\leq r_m$, so that
$$\cup_{j=1}^{j_0}  \tLambda_j(z_0^{(j)}) = \tLambda \quad, \quad
\cup_{j =1}^{j_0} \tGamma_j (x_0^{(j)}) \subset \tGamma ,$$
where $\tLambda$ is exactly the initially constructed $\tLambda$. In general the construction will only produce
sub-cylinders $\tGamma_j(x_0^{(j)})$ whose union is contained in
the initially constructed $\tGamma$.  Importantly, for all $j = 1, \ldots, j_0$
there exists $y_0^{(j)} = y_0^{(j)} (z_0^{(j)}, x_0^{(j)}) \in B^s(z_0^{(j)}, \ep_1)$ such that
the relation (5.67) holds for all $x \in \tGamma_j(x_0^{(0)})$, all $z \in \tLambda_j (z_0^{(j)})$, and all 
$b_1,b_2 \in W^s_R(z_0^{(j)})$ with $d(z_0^{(j)},b_1) < \delta'$ and $d(y_0^{(j)}, b_2) < \delta'$.
Moreover,
\be
\di \frac{1}{d_3}  \leq  \frac{\nu(\tGamma'_j)}{\nu(\tLambda'_j)} \leq d_3 
\ee
for all $j= 1, \ldots, j_0$. 

Since $\cup_{j=1}^{j_0}  \tLambda_j = \tLambda$ and $\tLambda$ is a sub-cylinder of $\tcc_m$ of co-length $q_1$,
there exists a global constant $d_1 \in (0,1)$ so that
$\nu(\tLambda') \geq d_1 \nu(\tcc') .$
(See e.g. Remark 2(a) in Sect. 6.3 below.) It follows now from this and (5.69) that we have
$\di \nu(\cup_{j=1}^{j_0} \tGamma'_j) \geq \frac{d_1}{d_3} \nu(\tcc')$. 
Thus, replacing the initial global constant $d_1 > 0$ by $d_4 = d_1/d_3$, we get 
$$\nu(\cup_{j=1}^{j_0}\tGamma'_j) \geq d_4 \, \nu(\tcc') \quad , \quad \nu (\cup_{j=1}^{j_0} \tLambda_j') \geq d_4\, \nu (\tcc') .$$
In general, the number $j_0$ will depend on the cylinder $\tcc$. However $d_3 > 1$ and $d_4 \in (0,1)$ are {\bf global constants},
independent of $m$, $\tcc$ and $z_0 \in \tP_0$. 

\bs

Let $N \geq N_0$.
Consider now an arbitrary $j = 1,2, \ldots, j_0$. As we remarked earlier, (5.67) holds for every $x \in \tGamma_j$ and every $z \in \tLambda_j$
for some choice of the point $y_0 \in W^s_{\ep_1}(z_0)$ which we will now {\bf denote by $y_0^{(j)}$}.
As in Lemma 5.4, we construct corresponding points
$$\yj_1 = \yj_1(z_0,x_0) \in \tpp^N(B^u(z_0,\ep')) \cap B^s(z_0, \delta')  \:\: , \:\:
\yj_2 = \yj_2(z_0,x_0) \in \tpp^N(B^u(z_0,\ep')) \cap B^s(y_0, \delta') .$$
so that (5.33) holds 
for any $b_1,b_2 \in W^s_{\tR}(z_0)$ with $d(\yj_1,b_1) < \delta'$ and $d(\yj_2, b_2) < \delta'$, where as in Lemma 5.4, $\delta'$ is
the constant given by (\ref{eq:delta'}).

Given $i = 1,2$, there exists  a cylinder $\Lj_i = \Lj_i(z_0)$ of length $N$ in $W^u_{\tR_{i_0}}(z_0)$  so that 
$$\tpp^{N} : \Lj_i \longrightarrow W^u_{\tR_{i_0}}(\yj_i)$$
is a bijection; then it is a bi-H\"older homeomorphism. Consider its inverse and its H\"older continuous
extension  
$$\tpp^{-N} : W^u_{\tR_{i_0}}(\yj_i)   \longrightarrow \Lj_i$$ 
and the cylinder
$$\Mj_i = \Mj_i(z_0) = \piU(\Lj_i(z_0)) \subset U $$ 
of length $N$ in $U_{i_0}$. Define the maps
$$\tv_{i,j} (z_0, \cdot) : U_{i_0} \longrightarrow \Lj_i \subset B^u(z_0,\ep'')  
\quad , \quad v_{i,j} (z_0, \cdot)   : U_{i_0} \longrightarrow \Mj_i  \subset U$$
by 
$$\tv_{i,j}(z_0, y) = \tpp^{- N} (\phi_{\Delta(z_0,y)}(\pi_{\yj_i}(y)) \quad ,\quad  v_{i,j}(z_0, y) = \piU(\tv_{i,j}(z_0, y)) .$$
Then  
$$\tpp^{N}(\tv_{i,j}(z_0, y)) = \phi_{\Delta(z_0,y)}(\pi_{\yj_i}(y)) = \phi_{[-\ep_0,\ep_0]}(W^s_{\ep_0}(y)) \cap W^u_{\tR_{i_0}}(\yj_i) ,$$
and
\be
\tpp^{N}(v_{i,j}(z_0, y)) = \phi_{[-\ep_0,\ep_0]}(W^s_{\ep_0}(y)) \cap \tpp^{N}(\Mj_i) \in W^u_{\tR_{i_0}} (\bj_i) ,
\ee
where $\bj_i = \bj_i(z_0) \in W^s_{R}(z_0)$ is such that $\tpp^{N}(\Mj_i) = W^u_{\tR_{i_0}}(\bj_i)$.
Thus, $\sigma^N((v_{i,j}(z_0, y)) = y$.
Next, there exist $x' \in \Mj_i$ and $y'\in \Lj_i$ with $\piU(y') = x'$, $\tpp^{N} (x') = \bj_i$ and
$\tpp^{N}(y') = \yj_i$. Since stable leaves shrink exponentially fast, using (2.1) we get
$$\di d(\bj_i , \yj_i) \leq \frac{1}{c_0 \gamma^{N}} d(x',y') \leq  \frac{1}{\gamma^{N}} < \delta' .$$
As in the proof of Lemma 5.4(b) we derive that
 (5.67) holds for $x \in \tGamma_j$ with $b_i$ replaced by $\bj_i$ for $i = 1,2$, that is
\be 
9 c_0 \delta_0 \|\uo_0\| \leq |\Delta(x, \bj_1)  - \Delta(x, \bj_2)| 
\ee
for any $x \in \tGamma_j$.

Set $z'_0 = \piU(z_0) \in U_{i_0}$. 
If $x,z \in \tcc$, and $x' = \piU(x)$, $z' = \piU(z')$, then
\begin{eqnarray*}
 I_{N}(x',z') 
 & = & |[\tau_N(v_{1,j}(z_0,x')) - \tau_N(v_{2,j}(z_0,x'))] - [\tau_N(v_{1,j}(z_0,z')) - \tau_N(v_{2,j}(z_0,z'))]| \nonumber\\
& = & \left| \, \Delta(\pp^{N} (v_{1,j}(z_0, x')), \pp^{N}(v_{1,j}(z_0, z')))  
-  \Delta(\pp^{N} (v_{2,j}(z_0, x')), \pp^{N}(v_{2,j}(z_0, z'))) \,\right|\nonumber\\
& = & \left| \, \Delta(\pi_{\bj_1}(x'), \pi_{\bj_1}(z') ) -  \Delta(\pi_{\bj_2}(x') , \pi_{\bj_2} (z') ) \,\right| \nonumber\\
& = & \left| \, \Delta(x', \pi_{\bj_1}(z') ) -  \Delta(x' , \pi_{\bj_2} (z') ) \,\right| .
\end{eqnarray*}
In particular when $z = z_0$ and $z' = z'_0$ the latter gives
$$I_N(x', z'_0) =  \left| \, \Delta( x' , \bj_1 ) -  \Delta( x', \bj_2) \,\right|  = \left| \, \Delta( x , \bj_1) -  \Delta( x  , \bj_2 ) \,\right|.$$

Since  $\Delta(x,\pi_{y}(z)) = \Delta(x, y) - \Delta(z, y) $ for any $y \in W^s_\ep(z_0)$, it follows from (5.71) 
and (5.67) 
that for any $x \in \tGamma_j$ and any $z \in \tLambda_j$ we have
\begin{eqnarray*}
I_{N}(x' ,z') 
& \geq   &        
\left| \, \left|\Delta(x , \bj_1) - \Delta(x, \bj_2)\right| - \left|\Delta(z , \bj_1) - \Delta(z, \bj_2)\right| \, \right|\\
& \geq & 9 c_0 \delta_0 \, \|\uo_0\| - 3c_0 \delta_0 \|\uo_0\| \geq  6 c_0 \delta_0 \|\uo_0\| 
\geq 6 c_0 \delta_0 \kappa\, \diam(\tcc) .
\end{eqnarray*}
This proves (5.48) and thus completes the proof of the lemma.
\endofproof


\section{Contraction operators and their main properties}
\setcounter{equation}{0}

As in Sect. 5 here we assume that $M$ is a $C^2$ compact  Riemannian manifold and $\phi_t$ is 
a  $C^2$ contact Anosov flow on $M$.

We will use the notation in Sects. 2, 3 and 4. In particular, $\rr = \{ R_i\}_{i=1}^{k_0}$ will be a fixed
pseudo-Markov family for the flow and $\tRR$ will be the related Markov family as in Sect. 2. 

As in Sect. 4, 
we assume that  $\hep > 0$ is a  small constant as in  Sect. 3 (that can be taken smaller if necessary), 
and we will again assume that $R(x)$, $\Gamma(x)$, $D(x)$ and $L(x)$ are Lyapunov $\hep$-regularity functions, 
while $\hr(x)$ and $r(x)$ are $\hep$-slowly varying radius function so that they satisfy (3.7) -- (3.14)  and the 
conclusions of  Lemma 3.1 with $\hep'$  replaced by $\hep$. 
Replacing $r(x)$ with the smaller $\hep$-regularity function $\hr(x)$, 
without loss of generality we will assume that the conclusions of Lemmas 3.1, 3.2 and 3.3
hold with $\hr(x)$ replaced by $r(x)$.

\def\hm{\hat{m}}
\def\Gammam{\Gamma^{(m)}}
\def\Gammammm{\Gamma^{(m'')}}
\def\Lambdam{\Lambda^{(m)}}
\def\Gammaj{\Gamma^{(j)}}
\def\Lambdaj{\Lambda^{(j)}}
\def\Xj{X^{(j)}}
\def\omj{\omega^{(j)}}
\def\omt{\omega^{(t)}}
\def\ommt{\omega^{(m,t)}}
\def\ommo{\omega^{(m,1)}}
\def\ommtwo{\omega^{(m,2)}}

\def\Xm{X^{(m)}}
\def\omm{\omega^{(m)}}
\def\vm{v^{(m)}}
\def\dteo{D_{\theta_1}}
\def\diamteo{\diam_{\theta_1}}

\def\EE{{\mathcal E}}
\def\tEE{\widetilde{\EE}}
\def\diamteo{\diam_{\theta_0}}
\def\dteo{D_{\theta_0}}

\def\thetaL{\theta^{(\ell)}}
\def\lambdaL{\lambda^{(\ell)}}
\def\thetam{\theta^{(m)}}
\def\um{u^{(m)}}
\def\lambdam{\lambda^{(m)}}
\def\tlambdam{{\tilde{\lambda}^{(m)}}}
\def\tlambdaL{\tilde{\lambda}^{(\ell)}}
\def\tlambda{\tilde{\lambda}}
\def\hGammaj{\widehat{\Gamma}^{(j)}}
\def\hGammamq{\widehat{\Gamma}^{(m_q)}}
\def\hGammamp{\widehat{\Gamma}^{(m_p)}}
\def\tGammam{\widetilde{\Gamma}^{(m)}}
\def\tLambdam{\widetilde{\Lambda}^{(m)}}

\def\hLambdam{\widehat{\Lambda}^{(m)}}
\def\hLambdaj{\widehat{\Lambda}^{(j)}}

\def\tGammamt{\widetilde{\Gamma}^{(m,t)}}
\def\tGammamo{\widetilde{\Gamma}^{(m,1)}}
\def\tGammamtwo{\widetilde{\Gamma}^{(m,2)}}

\def\tddmt{\widetilde{\dd}^{(m,t)}}
\def\tddmo{\widetilde{\dd}^{(m,1)}}
\def\tddmtwo{\widetilde{\dd}^{(m,2)}}

\def\hLambdamt{\widehat{\Lambda}^{(m,t)}}
\def\hddmt{\widehat{\dd}^{(m,t)}}
\def\hddmo{\widehat{\dd}^{(m,1)}}
\def\hddmtwo{\widehat{\dd}^{(m,2)}}

\def\thetaj{\theta^{(j)}}
\def\lambdaj{\lambda^{(j)}}
\def\Omegaj{\Omega^{(j)}}
\def\Omegam{\Omega^{(m)}}
\def\tpsi{\tilde{\psi}}

\def\tsigma{\tilde{\sigma}}
\def\bb{\mathcal B}
\def\sbb{\Sigma^+_{\bb}}
\def\tnu{\tilde{\nu}}
\def\tS{\widetilde{S}}
\def\Psib{\Psi^{(b)}}
\def\labb{L_{a b_0}}

\def\tbeta{\tilde{\beta}}

\def\tqq{\widetilde{\qq}}
\def\diamteo{\diam_{\theta_1}}

\subsection{Constructing families of cylinders covering the given Pesin set}

As in Sect. 4.2, we will assume that $F_0 : M \longrightarrow \R$ is a {\bf fixed H\"older continuous function} and 
{\bf $\m$ is the  Gibbs measure determined by $F_0$ on $M$}, while
{\bf $\mu_0$ is the related Gibbs measure on $R$} with respect to the Poincar'e map $\pp: R \longrightarrow R$.
We will identify $\mu_0$ with a measure on $\tR$ so that $\tPsi : R \longrightarrow \tR$ is an isomorphism.
As in Sect. 4.2, {\bf $\nu$ will be the Gibbs measure on $U$ determined by the function $F_0$.}

Let {\bf $P_0$ be the Pesin set in $\ll\cap R$  fixed in Sect. 4.2 and let  $K_0 = \piU(P_0)$}. As in Sect. 4.2, there exist constants 
$r_0 > 0$, $R_0 > 0$, $\Gamma_0 > 0$, $L_0 > 0$ with $r(x) \geq r_0$ and $R(x) \leq R_0$, $\Gamma(x) \leq \Gamma_0$,
$L(x) \leq L_0$   for all $x\in P_0$. Then $\tP_0 = \tPsi(P_0)$ is a Pesin set for $\tpp$ on $\tRR$.
Replacing $r_0$ with a smaller positive constant, 
we may assume that for every $x\in \tP_0 \cap \tR_j$ for any rectangle $\tR_j$ in $\tRR$ we have $B^u(x, r_0) \subset \Int(\tR_j)$. 

Recall from Sect. 2 the constants $\alpha_1 > 0$ and $\gamma_1 > \gamma > 1$. As in Sect. 4.2
we will assume that $\alpha_2 \in (0,1)$ and $\theta \in (0,1)$ are so that
\be\label{eq:theta-cond}
\frac{1}{\gamma^{\alpha_1 \beta}} \leq \theta \leq \frac{1}{\gamma_1^{\alpha_2}}  < 1 ,
\ee
where $\beta \in (0,1)$ is as in Lemma 5.3.
In particular $\alpha_2 \in (0,1)$ satisfies $\alpha_2 \, \log \gamma_1 \leq \alpha_1 \beta\, \log \gamma$, i.e.
$$0 < \alpha_2 \leq \frac{ \alpha_1 \beta\, \log \gamma}{\log \gamma_1} .$$
We will also need the constant $\beta_0 > 0$ such that $\theta = e^{-\beta_0}$
and a constant $\alpha_3$ with $0 < \alpha_3 < \alpha_2$ so that
\be\label{eq:thetao}
\theta_1 = \theta^{1/\alpha_3} > \thetaoo ,
\ee
the initially chosen small constant.



In what follows {\bf $b \in \R$ will be such that $|b| > b_0$,} where $b_0 > e$ is a fixed large constant.
We will impose some conditions on $b_0$ later.

Recall the constants $C_3 > 0$ and $\hep_7 > 0$ from Lemma 4.3 and the constant $\hep_{13}$ from Lemma 5.5. Set
$\hep_{14} = \hep_7 + \hep_{13}/\ttau_0 .$
This is still a constant with $\hep_{14} \leq \con \hep$ that can be made arbitrarily small with $\hep$.
Assuming $\ttau_0 < 1/2$, we have $\hep_{14} > 2 \hep_{13}$.
{\bf Let $\hq\geq 1$  be the smallest integer} so that
$\di \frac{C_3 e^{-2 \hq \hep_{14}}}{\lambda_1^{\hq}}  \leq  \frac{1}{|b|} .$
Then
$\di \frac{1}{ |b|} <  \frac{C_3 e^{-2 (\hq-1) \hep_{14}}}{\lambda_1^{\hq-1}} , $
so
$\di \frac{\ep'_2}{ |b|} <  \frac{e^{-2 \hq \hep_{14}}}{\lambda_1^{\hq}} , $
for some constant $\ep'_2 = \frac{e^{-2\hep_{14}}}{C_3 \lambda_1} > 0$. Thus, our choice of $\hq$ is so that
\be\label{eq:ep2-cond}
\frac{\ep'_2}{ |b|} <  \frac{e^{-2\hq \hep_{14}}}{\lambda_1^{\hq}}  \leq  \frac{\ep'_1}{ |b|} 
\ee
for some constants $ \ep'_1 = \frac{1}{C_3} > \ep'_2 > 0$.
Notice that $\hq$ depends on $b$.

Next, for any $z \in \tP_0$ denote by $m = m(z)$ the length of the cylinder $\tqq(z)$ in $\tR$ containing $z$ so that $\ttau_m(z) \leq \hq$
and $m$ is maximal with this property, i.e. $\hq < \ttau_{m+1}(z)$, so $\hq < \ttau_{m+1}(z) \leq \ttau_m(z) + \tau_0$. 
Thus, $\hq < [\ttau_m(z)]+1 + \tau_0$, so
$q(z) = [\ttau_m(z)] \leq \hq < q(z) +  2 .$
Now $\ttau_{m+1}(z) > \hq$ implies $(m+1) \tau_0 > \hq$.
Then $m \ttau_0 \leq \hq \leq (m+1) \tau_0$, so 
$\hm_1 = [\hq/\tau_0 - 1] \leq \hq/\tau_0 - 1   \leq m(z) \leq  \hm_2  = [\hq/\ttau_0] .$
In what follows we will be considering cylinders in $\tR$ intersecting $\tP_0$ that have lengths in the interval $[\hm_1, \hm_2]$.
Notice that $\hm_1$ and $\hm_2$ depend on $b$.

Given $z\in \tP_0$,  let again $\tqq(z)$ be the cylinder of length $m(z)$ in $\tR$ containing $z$.
Then (\ref{eq:ep2-cond}), (\ref{eq:diamC}) and $\hep_7 < \hep_{14}$ imply 
$$\diam(\tqq(z)) \geq \frac{e^{-q(z) \hep_7}}{\lambda_1^{q(z)}} \geq  \frac{e^{-\hq \hep_{14}}}{\lambda_1^{\hq}}
= e^{\hq \hep_{14}}  \frac{e^{-2\hq \hep_{14}}}{\lambda_1^{\hq}} 
\geq e^{\hq \hep_{14}}\, \frac{\ep'_2}{ |b|} \geq e^{q(z)\hep_{14}}\, \frac{\ep'_2}{|b|} .$$
Similarly, using (\ref{eq:ep2-cond}) and (\ref{eq:diamC}) again, we get
$$\diam(\tqq(z)) \leq \frac{C_3 e^{q(z) \hep_7}}{\lambda_1^{q(z)}} \leq  \frac{C_3 e^{\hq  \hep_{14}}}{\lambda_1^{\hq - 2}}
\leq  C_3 \lambda_1^{2}\, \frac{e^{-2\hq \hep_{14}}}{\lambda_1^{\hq}} e^{3\hq \hep_{14}} 
 \leq  C_3 \lambda^{2}_1 \frac{\ep'_1}{ |b|} e^{3(q(z)+ 2) \hep_{14}} \leq \frac{C'_5 e^{3q(z) \hep_{14}}}{|b|},$$
where
$C'_5 = C_3 \lambda_1^{2} e^{6} > 0 $
is a global constant.
 
Let $\qq'(z) = \piU(\qq(z))$ be the corresponding projection of $\qq(z)$ in $U$. Here $\qq(z)$ is the cylinder in $R$ corresponding
to $\tqq(z)$, i.e. $\tPsi(\qq(z)) = \tqq(z)$.
Choose a finite set of points $Z_1, Z_2, \ldots, Z_{k_0}$ in $\tP_0$ such that the projections $\qq'(Z_k)$ cover completely $K_0$.
If $\qq'(Z_k)$ and $\qq'(Z_{k'})$  have common interior points (in the topology of $U$), then one of these 
cylinders contains the other. So, omitting some of the cylinders, we may assume the points $Z_1, Z_2, \ldots, Z_{k_0}$ in $\tP_0$
are chosen so that $\qq'(Z_k) \cap \qq'(Z_{k'}) \cap \hU = \e$ for all $k \neq k'$, and we still have
$K_0 \subset \cup_{k=1}^{k_0} \qq'(Z_k)$.
Then for each $k = 1,2, \ldots, k_0$, 
$\tqq_k = \tqq(Z_k)$ is a cylinder in $\tR$ containing $Z_k \in \tP_0$ with  ${\bf \piU(\qq_k) = \qq'_k} .$ 
Denote by {\it $t_k$ the length of the cylinder $\qq_k$}.
It follows from the above discussion, the construction of the cylinders $\qq_k$  and $t_k \ttau_0 \leq q(Z_k) \leq t_k\, \tau_0 + 1$ that 
there exist global constants $C_5 > 0$ and $\ep_2 > 0$ such that
\be\label{eq:Q-cond}
\frac{\ep_2 }{|b|} \, e^{t_k \hep_{15}} \leq \diam(\tqq_k) \leq \frac{C_5 e^{t_k \hep_{16}}}{|b|} 
\ee
for all $k = 1, 2, \ldots, k_0$, where $\hep_{15} = \hep_{14} \ttau_0$ and $\hep_{16} = 3 \hep_{14} > 3 \hep_{15}$. 
It follows from the definition of $\hep_{14}$ that $\hep_{15}  > \hep_{13}$, the constant  from Lemma 5.5 and 
$t_k \geq \tm_0$, the integer from Lemma 5.5.

Let us mention that (\ref{eq:diamC})  and the above imply estimates for $t_k$ by means of  $|b|$.
Indeed, by (\ref{eq:Q-cond}) and (\ref{eq:diamC}) we get 
$-  \log |b| + \log \ep_2 + t_k\hep_{15} \leq \log C_3 + t_k \tau_0 \hep_7 - t_k\ttau_0 \log \lambda_1 ,$
therefore 
$$t_k (\ttau_0 \log \lambda_1 + \hep_{15} - \tau_0 \hep_7) \leq \log |b| - \log \ep_2  + \log C_3 ,$$
so
$$\di t_k  = t_k(b) \leq \frac{1}{ \ttau_0 \log \lambda_1  - \tau_0 \hep_7} (\log |b| - \log \ep_2  + \log C_3 ) \leq D_1 \log |b|$$ 
for some global constant $D_1 > 0$, assuming $|b|$ is sufficiently large. 
In a similar way from the other sides of (\ref{eq:Q-cond}) and (\ref{eq:diamC}) we get
$t_k (\tau_0 \log \lambda_1 + \hep_7/\ttau_0 + \hep_{16}) \geq \log |b| - \log C_3 - \log C_5$, that is there exists a global constant 
$D_2 > 0$ so that (for sufficiently large $|b|$)
we have $t_k \geq D_2 \log |b|$. Thus,
\be\label{eq:D1-cond}
D_2 \log |b| \leq t_k (b) \leq D_1 \log |b| .
\ee

As a result of the construction in this section, we obtain a family $\{\qq_k\}_{k=1}^{k_0}$ of cylinders $\qq_k = \qq_k (b)$, 
where $k_0 = k_0(b)$ also depends on $b$, such that the length $t_k = t_k (b)$ of $\qq_k$ satisfies (6.4) and (6.5) and 
$${\bf K_0 = \piU(P_0) \subset \cup_{k=1}^{k_0} \qq'_k \subset U .}$$
Notice that {\bf $\di \gamma_2 = \frac{1}{2} \nu(K_0) > 0$ is a constant  independent of $b$.}

\ms

If $\cc$ is a cylinder in $U$, by Lemma 4.1(c) and our choice of $\alpha_2$ we have
$\diam_{\theta}(\tcc)  \leq  C_1 (\diam(\tcc))^{\alpha_2}  ,$
where $\tcc = \tPsi(\cc)\ \subset \tR$. Consequently, if $r$ is the length of $\cc$, then
$$\diam_{\theta_1}(\tcc) 
 =      \theta_1^{r} < \theta^r = \diamte(\tcc)  \leq C_1\, (\diam(\tcc))^{\alpha_2} .$$
A more precise estimate (which will be used below) is:
$$\diam_{\theta_1}(\tcc) = \theta_1^{r} = \theta^{r/\alpha_3 } = (\diamte(\tcc))^{1/\alpha_3} 
 \leq  (C_1 \diam(\tcc))^{\alpha_2/\alpha_3}  .$$
 
For the cylinders $\qq_k$  it follows from the above and (\ref{eq:D1-cond}) that
\begin{eqnarray*}
\diam_{\theta}(\tqq_k)   
& \leq &  
C_1 \left(\frac{C_5 e^{t_k \hep_{16}}}{|b|}\right)^{\alpha_2}
 \leq  C_1 \left(\frac{C_5 e^{\hep_{16} D_1 \log |b|}}{|b|}\right)^{\alpha_2}
 \leq  \frac{C_1 (C_5)^{\alpha_2}}{|b|^{\alpha_2(1 - D_1 \hep_{16})}}  \leq \frac{C_6}{|b|^{\alpha_3}} ,
\end{eqnarray*}
where we choose a global constant $C_6$ with  $C_6 \geq  C_1 (C_5)^{\alpha_2}$ and $0 < \alpha_3 < \alpha_2$. 
Thus,
$$\diamte(\qq_k) \leq \frac{C_6}{|b|^{\alpha_3}}  $$
for all $k = 1, \ldots, k_0(b)$ and all $|b| \geq b_0$. 
Apart from that, the above and (\ref{eq:D1-cond}) yield
\begin{eqnarray*}
\diamteo(\qq_k) 
& \leq & C_1 \left(\frac{C_5 e^{t_k \hep_{16}}}{|b|}\right)^{\alpha_2/\alpha_3}
 \leq  C_1 \left(\frac{C_5 e^{\hep_{16} D_1 \log |b|}}{|b|}\right)^{\alpha_2/\alpha_3}
  \leq  \frac{C_1 (C_5)^{\alpha_2/\alpha_3}}{|b|^{(\alpha_2/\alpha_3)(1 - D_1 \hep_{16})}}  \leq \frac{C_6}{|b|} ,
\end{eqnarray*}
assuming that $C_1 (C_5)^{\alpha_2/\alpha_3} \leq C_6$ and $D_1 \hep_{16}$ is sufficiently small so that
$(\alpha_2/\alpha_3)(1 - D_1 \hep_{16}) > 1$. 

We will now slightly enlarge our cylinders $\tqq_k$. For every $k$, consider a cylinder $\tcc_k$ 
such that $\tqq_k \subset \tcc_k$ and
$$\frac{\ep_4}{|b|} \leq \diamteo(\tcc_k)  \leq \frac{C_6}{|b|} $$
for some global constant $\ep_4 > 0$ (we can take e.g. $\ep_4 = \theta_1$).
Clearly such a cylinder $\tcc_k$ exists. In general the projections $\tcc'_k = \piU(\tcc_k)$ are not necessarily disjoint.
We now choose a disjoint sub-family of $\{ \tcc_k\}_{k=1}^{k_0}$ that covers the set $K_0$, and re-numbering its
elements, we get a family $\{\tcc_m\}_{m=1}^{m_0}$ of cylinders in  $\tR$ for some $m_0 = m_0(b) \leq k_0$
such that their projections $\tcc'_m$ in $\hU$ are disjoint and cover the compact set $K_0$. Then for each $m = 1, \ldots, m_0$
there exists $k = 1, \ldots, k_0$ with $\tqq_k \subset \tcc_m$, therefore the length $s_m$ of $\tcc_m$
satisfies $s_m \leq t_k$, and it follows from (6.4) that
$$\frac{\ep_2 e^{s_m \hep_{15}}}{|b|} \leq \frac{\ep_2 e^{t_k \hep_{15}}}{|b|} 
\leq \diam(\tqq_k) \leq \diam(\tcc_m) .$$

We have proved the following lemma.

\ms

\noindent
{\bf Lemma 6.1.} {\it There exists cylinders $\tcc_1, \ldots, \tcc_{m_0}$ in $\tR$ for some $m_0 = m_0(b) \geq 1$ such that
each of them has a common point with $\tP_0$, $\tcc'_m \cap  \tcc'_{m'}\cap \hU = \e$ for any $m \neq m'$, 
\be\label{eq:Vb-cond}
{\bf K_0 = \piU(P_0)  \subset V_b = \cup_{m=1}^{m_0} \cc'_m } ,
\ee
and for all $m = 1, \ldots, m_0$ we have
\be\label{eq:ccmo-cond}
\frac{\ep_4}{|b|} \leq \diamteo(\tcc_m)  \leq \frac{C_6}{|b|} 
\ee
and
\be\label{eq:C-cond}
\frac{\ep_2 }{|b|} \, e^{s_m \hep_{15}} \leq \diam(\tcc_m) ,
\ee
where $s_m$ is the length of the cylinder $\tcc_m$.}

\ms

We should stress that $\ep_2 > 0$, $\ep_4 > 0$ and $C_6 > 0$ are global constants, independent of $b$.

Since $\theta = \theta_1^{\alpha_3}$, it follows from (\ref{eq:ccmo-cond}) that
\be\label{eq:ccm-cond}
\frac{\ep_4}{|b|^{\alpha_3}} \leq \diamte (\tcc_m)  \leq \frac{C_6}{|b|^{\alpha_3}} 
\ee
for all $m = 1, \ldots, m_0$.

\def\qb{q^{(b)}}
\def\cb{c^{(b)}}
\def\fb{f^{(b)}}
\def\varphiso{\varphi^{(s_0)}}
\def\Ab{A^{(b)}}
\def\Bb{B^{(b)}}
\def\tSb{\tS^{(b)}}
\def\varphibobone{\varphi^{(b_1)}_{b_0}}
\def\varphibobtwo{\varphi^{(b_2)}_{b_0}}
\def\varphibonebtwo{\varphi^{(b_2)}_{b_1}}
\def\varphibob{\varphi^{(b)}_{b_0}}
\def\varphibobn{\varphi^{(b_n)}_{b_0}}
\def\Psibo{\Psi^{(b_0)}}
\def\Prf{\mbox{\rm\footnotesize Pr}}
\def\hdelta{\hat{\delta}}
\def\tG{\widetilde{G}}
\def\tH{\widetilde{H}}
\def\Psibb{\Psi^{(b')}}
\def\Psibn{\Psi^{(b_n)}}
\def\PsibN{\Psi^{(b_N)}}


\subsection{A consequence of Chebyshev's inequlaity}

Let $P_0$ be the fixed Pesin set in $R\cap \ll$ from Sect. 4.2 and let  $K_0 = \piU(P_0) \subset U$. 

Fix for a moment a large constant $b_0 > 1$ and assume that $|b| \geq b_0$. Recalling the construction of 
the sets $V_b$ in Sect. 6.1, we have $K_0 \subset V_b$, and $\cap_{|b| \geq b_0} V_b = K_0$.
Hence $\nu(V_b) \geq \nu(K_0) > 0$. Here, as in Sect. 4.2, $\nu$ is the Gibbs measure defined
by the function $g = f - P_f\, \tau$ on $U$ with $\Pr_{\sigma}(g) = 0$.


Consider the functions
$\di \Psib = \chi_{V_{b}}$ and $\Psi_0 = \chi_{K_0}$.

\ms

\noindent
{\bf Remark 1.} Let $\ep > 0$ be an expansivity constant (see e.g. \cite{Ba}) for $\sigma_A : \saa \longrightarrow \saa$ 
with respect to the metric $\dte$; one can take e.g. $\ep = \theta/2$. If $n \geq 1$ is an integer and $x = (x_0, x_1, \ldots)$,
$y = (y_0, y_1, \ldots) \in \saa$ are such that $\dte(\sigma^i(x), \sigma^i(y)) < \ep$ for all $ i = 0, 1, \ldots, n-1$,
then $x_i = y_i$ for all $i < n$, so $\dte(x,y) \leq \theta^n$. Hence
$$|\psi_n(x) - \psi_n(y)| \leq \sum_{i=0}^{n-1} |\psi(\sigma^i(x)) - \psi (\sigma^i(y))|
\leq |\psi|_\theta \sum_{i=0}^{n-1} \theta^{n-i} < \frac{|\psi|_\theta}{1-\theta} .$$
Notice that $|\Psib|_\theta$ is rather large. Indeed, by Lemma 4.1(a)  the natural estimate for $|\Psib|_\theta$ is 
$\di |\Psib|_\theta \leq \frac{1}{\theta^{s_m(b)}} = 1/\diamte (\tcc_m) $ which is of size $\Con |b|^{\alpha_3}$ by (\ref{eq:ccm-cond}).
The latter is much larger than $\Con \log |b|$. In what follows we will avoid using this norm, employing some simple 
estimates instead -- see e.g. (\ref{eq:K-cond}) below.

\bs


In what follows we will use sometimes the shift space $\saa$ and the shift $\sigma = \sigma_A : \saa \longrightarrow \saa$. 
We will identify $\nu$ with
a Gibbs measure on $\saa$ using the natural projection $\pi^+ : \saa \longrightarrow U$ which conjugates $\sigma_A$ and $\sigma$
(see Sects. 2 and 4.1). We will also identify the set $V_{b}$ with a subset of $\saa$ via this isomorphism, using the same notation for it.
We will use the notation $\Pr = \Pr_{\sigma}$, i.e. {\bf we will always consider topological pressure with respect to the map
$\sigma  : \saa \longrightarrow \saa$ or the related map $\sigma : U \longrightarrow U$.}

Recall that $V_b$ is a union of cylinders $\cc'_m$ of length $s_m(b) \leq t_k(b) \leq L_b$, where by (\ref{eq:D1-cond}) we can take
\be\label{eq:Lb-cond}
L_b = D_1\, \log |b| .
\ee
For two cylinders in $\saa$, either
one of them contains the other, or they are disjoint. It follows from this that if $X$ is a cylinder in $\saa$ of length $n > L_b$,
then either $X \subset V_b$ or $X \cap V_b = \e$. Hence either $\Psib(x) = 1$ for all $x\in X$ or $\Psib(x) = 0$ for all $x \in X$.
This implies a certain improvement of the estimate in Remark 1 in the special case $\psi = \Psib$. 
Assume $n > L_b$ and $D_\theta(x, y) \leq \theta^n$ for some $x,y\in \saa$. Then for any integer $i < n- L_b$ we have
$D_\theta(\sigma^ix, \sigma^iy) \leq \theta^{L_b}$, so $\sigma^iy$ belongs to the cylinder  of length $L_b$ determined by $\sigma^ix$.
Thus, either both $\sigma^ix$ and $\sigma^iy$ belong to $V_b$ or both do not belong to $V_b$. That is $\Psib(\sigma^i x) = \Psib(\sigma^iy)$.
This implies
\begin{eqnarray}\label{eq:K-cond}
|\Psib_n(x) - \Psib_n(y)| 
& \leq & \sum_{i=0}^{n-1} |\Psib(\sigma^i(x)) - \Psib (\sigma^i(y))|\nonumber\\
& = & \sum_{i = n-L_b}^{n-1} |\Psib(\sigma^i(x)) - \Psib (\sigma^i(y))| \leq L_b 
\end{eqnarray}
for all $n > L_b$ and all $x,y\in \saa$ with $\dte (x,y) \leq \theta^n$. 
Notice that (\ref{eq:K-cond}) also holds for all $n = 1,\ldots, L_b$ and all $x,y\in \saa$.

\bs

For any integer $n \geq 1$, any $|b| \geq b_0$ and any $q\in \R$ set
$$\fb_n (q)= \frac{1}{n} \log \int_{U} e^{q \Psib_n(x)} \, d\nu(x) .$$

The following lemma can be derived from the proof of Lemma 4.1 in \cite{MV}.

\bs

\noindent
{\bf Lemma 6.2.} {\it
There exists a global constant $C_7 = C(\ep) \geq 2$  with the following property: if
$\Psi \in \ff_\theta(\saa)$ and the constant $K > 0$ 
are so that
$|\Psi_n(x) - \Psi_n(y)| \leq K$ for all $n \geq 1$ and all $x,y\in \saa$ with $\dte(x,y) \leq \theta^n$,
then the function
$$f_n(q) = \frac{1}{n} \log \int_U e^{q \Psi_n(x)}\, d\nu(x) \quad, \quad q\in \R ,$$
satisfies
\be\label{eq:n-cond}
f_n(q) \leq \Pr (g + q \Psi) +  \frac{2|q| K}{n}  + \frac{1}{n} \log C_7
\ee
for all $q\in \R$ and all $n \geq 1$.}

\bs

It follows from Lemma 4.1 in \cite{MV} (see also Proposition 3.1 in \cite{Den}) that for the function $f_n(q)$ in the above lemma there exists 
$\lim_{n\to\infty} f_n(q) = \Pr (g+ q \Psi)$
for every $q \in \R$.  Here we only need one half of that argument.
For completeness  we prove Lemma 6.2 in Appendix II using part of the argument in \cite{MV} and some lemmas from \cite{B3}.

In what follows it is more convenient to work with $U$ and the corresponding first-return map $\sigma : U \longrightarrow U$ rather
than the shift space $\saa$ and the shift $\sigma_A : \saa \longrightarrow \saa$. 

We always consider $b\in \R$ with $|b| \geq b_0$. 
Set
$$
\tG =  \lim_{|b| \nearrow \infty} \Pr (g - \Psib) \leq 0 .$$
Notice that $\Psib \geq \Psibb$ when $|b| \leq |b'|$, so the function 
$\Pr (g - \Psib) \leq 0$ is increasing. 

Given $\alpha \in [0,1]$, the convexity of $\Pr$ and $\Pr(g) = 0$ imply
$$\Pr (g - \alpha\, \Psib) = \Pr (\alpha (g - \Psib) + (1-\alpha) g) 
\leq \alpha\, \Pr(g - \Psib) + (1-\alpha)\, \Pr (g) = \alpha\, \Pr(g - \Psib) $$
Thus,
\be\label{eq:alpha-cond}
\Pr (g- \alpha\, \Psib) \leq \alpha\, \Pr (g - \Psib) 
\ee
for all $\alpha \in [0,1]$ and all $|b| \geq b_0$.

The following lemma is an easy consequence of Lemma 4.3 and the choice of the Pesin set $P_0$
and its projection $K_0$ in $U$ in Sect. 4.2.

\bs

\noindent
{\bf Lemma 6.3.} {\it We have $\tG < 0$.}

\bs

\noindent
{\it Proof.} By the choice of the Pesin set $P_0$ and its projection $K_0$ we have
$$\lim_{n\to\infty} \Pr (g - \chi_{K_0^{(n)}}) = - c < 0$$
for some constant $c > 0$. Then $\Pr (g - \chi_{K_0^{(n)}}) \leq - c$ for all $n \geq 1$, since
$ \Pr (g - \chi_{K_0^{(n)}})$ is increasing with respect to $n \geq 1$.
Given $b\in \R$ with $|b| \geq b_0$ we can choose an integer $n \geq 1$ so large that
$K_0^{(n)} \subset V_b$.
Then $\chi_{K_0^{(n)}} \leq \chi_{V_b}$, so $g - \chi_{K_0^{(n)}} \geq g - \chi_{V_b}$ and therefore
$$\Pr (g - \Psib)  = \Pr (g - \chi_{V_b}) \leq \Pr (g - \chi_{K_0^{(n)}}) \leq - c .$$
Hence $\tG = \lim_{|b|| \to \infty} \Pr(g- \Psib) \leq -c < 0$.
\endofproof

\bs

Set
\be
\ep_0 = \frac{1}{4} \, |\Pr(g - \Psi_0)| >   0 .
\ee
Then (\ref{eq:alpha-cond}) implies
$$\Pr (g- \alpha\, \Psib) \leq \alpha\, \Pr (g - \Psib) \leq - 4 \alpha\, \ep_0$$
for all $\alpha \in [0,1]$ and all $|b| \geq b_0$. That is, for any $-1 < q < 0$ we have
\be
\Pr (g + q \, \Psib) \leq 4 q \, \ep_0 .
\ee

{\bf Fix a constant}
$$- 1 < \tq_0  < 0 .$$
Let $\ep_0 > 0$ be as in (6.14), and let $|b|\geq b_0$ and $n \geq 1$.
We will now use Lemma 6.2 with $\Psi = \Psib$, $q = \tq_0 < 0$ and $K = L_b$.
It follows from (\ref{eq:n-cond}) and (6.15)
that
$$\fb_n(\tq_0) \leq \Pr_\sigma(g + \tq_0 \,\Psib) +  \frac{2\, |\tq_0|}{n} L_b + \frac{1}{n} \log C_7
\leq 4\, \tq_0 \, \ep_0  + \frac{2 \, |\tq_0| \, L_b}{n} + \frac{1}{n} \log C_7$$
for all $n \geq 1$. 
Using (\ref{eq:Lb-cond}), we have $L_b = D_1 \log |b|$.  Define an integer $n_0(b)$ by
\be\label{eq:nb-cond}
n_0(b) = D_3\, \log |b| ,
\ee
where $D_3 > 1$ is a {\bf fixed constant} with
$$\di D_3 \geq  \max \left\{ \frac{2 D_1}{\ep_0}\, , \,  \frac{\log C_7}{|\tq_0|\, \ep_0 } \right\} . $$
Then for every $|b| \geq b_0$ and every integer $n \geq n_0(b)$ we have
\be\label{eq:fbfinal-cond}
\fb_n(\tq_0) \leq 4 \tq_0 \ep_0 + 2 |\tq_0|\,  \ep_0 = 2 \tq_0 \, \ep_0 < 0 . 
\ee

Set 
$$\Ab(n) = \{ x\in U : \Psib_n(x) \leq n \ep_0 \} .$$
Since $\tq_0 < 0$, we have
$\di \Ab (n) = \{ x \in U : e^{\tq_0\, \Psib_n(x)} \geq e^{n \tq_0 \, \ep_0} \} .$
Using Chebyshev's inequality\footnote{Estimates of this kind appear e.g. in  \cite{BLM} or Sect. 2.4 in \cite{Var}.}
(sometimes called Markov's inequality; see e.g. Ch. 3 in \cite{G}), it follows that
$$\nu(\Ab (n)) \leq \frac{\int_{U} e^{\tq_0 \Psib_n(x)}\, d\nu(x)}{e^{n \tq_0 \ep_0} }
=  e^{-n \tq_0 \ep_0} \,  \int_{U} e^{\tq_0 \Psib_n(x)}\, d\nu(x) .$$
Assuming $n \geq n_0(b)$ and combining this with (\ref{eq:fbfinal-cond}) gives
$$\frac{1}{n} \log \nu(\Ab(n))
\leq - \tq_0 \ep_0 + \frac{1}{n} \log \int_{U} e^{\tq_0 \Psib_n(x)} \, d\nu(x)
= - \tq_0 \ep_0 + \fb_n(\tq_0)  \leq - \tq_0 \ep_0 +  2 \tq_0 \ep_0 = \tq_0 \ep_0  $$
for all $|b| \geq b_0$, $n \geq n_0(b)$.
Thus, for such $b$ and $n$ we have $\nu(\Ab(n)) \leq e^{n \tq_0 \ep_0}$.

{\bf Set}
$c =  |\tq_0| \ep_0  > 0 .$
It follows from our discussion that for the set
$$U_b(n)  = \left\{ x\in U : \frac{\Psib_n(x)}{n} \leq \ep_0 \right\} ,$$
we have
$\di \nu \left( U_b(n) \right) \leq e^{- n c} $
for all $|b| \geq b_0$ and all integers $n \geq n_0(b)$.

We will use the above for integers of the form $n = M N$, where $M \geq n_0(b)$ and $N \geq 1$.
For any $x \in U$ set
$$\tSb_{M,N}(x) = \sharp \{ j : 1 \leq j \leq  M\;, \; \sigma^{j N}(x) \in \cup_{r=0}^{N-1} \sigma^{-r}(V_{b})\}  ,$$
and
$$\tU_b(M,N) = \left\{ x\in U : \frac{\tSb_{M,N}(x)}{M} \leq \ep_0 \right\} .$$

As a consequence of the above we obtain

\bs

\noindent
{\bf Lemma 6.4.} {\it For any $b\in \R$ with $|b| \geq b_0$ we have
\be\label{eq:Ub-cond}
\nu(\tU_b(M,N)) \leq  e^{- M N c}
\ee
for all integers $M \geq n_0(b)$ and $N \geq 1$, where $n_0(b)$ is defined by} (\ref{eq:nb-cond}).

\bs

\noindent
{\it Proof.} Let $x \in \tU_b(M,N)$. Then $x \in U_b(MN)$. Indeed, if $\Psib_{MN}(x) > \ep_0 \, MN$,
then $\sigma^k(x) \in V_{b}$ for more than $\ep_0 MN$ values of $k = 0,1, \ldots, MN-1$. 
So, for more than $\ep_0 M$  values of 
$j = 1, \ldots, M$ there exists $r= 0,1, \ldots, N-1$ such that $\sigma^{j N + r}(x) \in V_{b}$, that is $\tSb_{M,N}(x) > \ep_0 M$.
The latter means $x \notin \tU_b(M,N)$, a contradiction. Thus, $\tU_b(M,N) \subset U_b(MN)$, so (\ref{eq:Ub-cond}) holds.
\endofproof

\bs

The above will play an important role in Sect. 8.


\def\ta{\tilde{a}}
\def\tS{\widetilde{S}}
\def\tE{\widetilde{E}}
\def\tWm{\tW^{(m)}}
\def\tDm{\widetilde{D}^{(m)}}
\def\tdd{\widetilde{\dd}}
\def\tddm{\tdd^{(m)}}
\def\hddm{\widehat{\dd}^{(m)}}
\def\Xm{X^{(m)}}
\def\Xmt{X^{(m,t)}}
\def\Ym{Y^{(m)}}
\def\xim{\xi^{(m)}}
\def\etam{\eta^{(m)}}


\subsection{The contraction operators}

Let $|b| \geq b_0$ and let  {\bf $c =  |\tq_0| \ep_0 > 0$
be the constant  from Lemma 6.4} so that (\ref{eq:Ub-cond}) holds for all integers $M \geq n_0(b)$ and $N \geq 1$. 
In this section we will fix several more objects that will stay the same throughout the rest of the paper.

In what follows we will use the constants $\theta$ with (6.1) and $\theta_1$ defined by (\ref{eq:thetao})
for some $0 < \alpha_3 < \alpha_2 < 1$.
As in Sect. 6.1 we will also need the constant $\beta_0 > 0$ such that $\theta = e^{-\beta_0}$.
{\bf The constants $\theta_1 $ and $\theta$ will stay fixed throughout the rest of the paper.}

\bs

\noindent
{\bf Remarks 2:} (a)
It follows from (4.2) that 
if $X$ and $Y$ are cylinders in $U$ of lengths $\ell_X$ and $\ell_Y$, $X \subset Y$
and $\nu(Y) \leq d\, \nu(X)$ for some constant $d > 0$, then $\ell_X \leq \ell_Y + d'$,
where we can take $d' = \frac{1}{|g_0|} \log (d c_2/c_1)$.
Here $g = f - P_f\tau$ is the function from Sect. 4.2 with $g_0 = \max g < 0$, while $c_1$ and $c_2$
are the constants from (4.2).
Indeed, take any $x\in X$ and use (4.2) to get
$$d \geq \frac{\nu(Y)}{\nu(X)} \geq \frac{c_1 \, e^{g_{\ell_Y(x)}}}{c_2 \, e^{g_{\ell_X(x)}}}
= \frac{c_1}{c_2 \, e^{g_{(\ell_X - \ell_Y)(x)}}} \geq \frac{c_1}{c_2 \, e^{(\ell_X - \ell_Y)\,g_0}} . $$
Therefore $e^{(\ell_X - \ell_Y)\,g_0}\geq \frac{c_1}{d c_2}$, so $(\ell_X - \ell_Y)\,g_0\geq \log \frac{c_1}{d c_2}$.
Thus, $(\ell_X - \ell_Y) |g_0| \leq \log \frac{d c_2}{c_1}$, and therefore $\ell_X - \ell_Y \leq d'$.

\ms

(b) Notice also that if $X$ is a cylinder in $\hU$ and $d \in (0,1)$ is a given constant, then there exists an integer
$\ell = \ell(d) \geq 1$, uniquely determined by $d$, such that if $X_1, \ldots, X_k$ are any disjoint sub-cylinders
of $X$ with $\sum_{i=1}^k \nu(X_i) \geq d\, \nu(X)$, then the smallest sub-cylinder $Y$ of $X$ containing
$\cup_{i=1}^k X_i$ has co-length at most $\ell$ in $X$. Indeed, for such $Y$ we must have $\nu(Y) \geq d\, \nu(X)$,
and as in part (a), this gives an upper bound for $\ell$. 

\ms

(c) Using the last Remark, there exists a global constant $\ell_0$ such that if $\tX_1, \ldots, \tX_k$ are disjoint sub-cylinders of 
$\tcc_m$ for some $m$, and $\sum_{i=1}^k \nu(\tX_i) \geq \frac{d_4}{8d_3} \nu (\cc'_m)$, where $d_3 > 1$
and $d_4 \in (0,1)$ are as in Lemma 5.5, then the co-length 
of the smallest sub-cylinder $\tY$ of $\tcc_m$ containing $\cup_{i=1}^k \tX_i$ is $\leq \ell_0$. 
{\bf Fix $\ell_0 \geq 1$ with this property.}

\bs

Next, let $M_1 \geq \tp_0 + 1$, where $\tp_0 > 0$ is as in  Sect. 2. We will also assume that
$M_1 \geq 2 |P_f| + 1$  (see Sect. 4.2 for $P_f$).  Let  $T_0 $ be as in (4.4), and let 
\be\label{eq:C8-cond}
\di C_{8} = \frac{T_0 C_{6} C_4}{1-\theta_1} + \frac{M_1 C_4}{d_0^2} ,
\ee
where $d_0 \in (0,1]$ is the constant fixed at the end of Sect. 2 and $C_4 \geq 1$ is the constant from Lemma 5.3.
Next, {\bf fix an arbitrary constant $\mu_0$ with}
\be\label{eq:muo-cond}
0 < \mu_0 \leq \min \left\{ \frac{1}{4} ,  \frac{1- \cos \ep_3}{20}   \right\} ,
\ee
where 
$$\ep_3 = \frac{1}{2} \min \left\{ \frac{\delta_1 \ep_2}{16}\, , \, \frac{\pi}{32} \,  , \, \ln \frac{19}{16} \right\} $$
with $\delta_1 = 6 c_0  \delta_0 d_2 > 0$ and $\ep_2 > 0$ as in (\ref{eq:C-cond}).
Then {\bf fix a large constant $E > 0$} 
so that
\be\label{eq:E-cond}
E \geq \max\left\{ 240 e^{C_{8}}\, C_{8}\, ,\, \frac{1}{16\theta_1^{\ell_0}}\, , 
\, \frac{9 C_1 T_0 e^{T_0/(1-\theta_1)}}{(1- \theta_1)\theta_1^{\ell_0} \, \ep_4} \, , \, \frac{2 C_4}{d_0^2} \, , \,
12 T_0 e^{T_0/(1-\theta)}  \right\} ,
\ee
where $\ell_0 \geq 1$ is the global constant from Remark 2(c)
and $\ep_4 > 0$ is the constant from (\ref{eq:ccmo-cond}).
Let $N_0 \geq 1$ be the global constant from Lemma 5.5. We will assume $N_0$ is chosen so that
\be\label{eq:No-cond}
\theta^{N_0/2} < \min \left\{ \frac{1}{100}\, , \, \frac{\theta^{\ell_0}}{2 E\, e^{C_{8}}} \, , \, 
\frac{1}{6 e^{T_0/(1-\theta)}} \, , \, \frac{d_6\, \ln 2}{E \, C_6 } 
\right\} ,
\ee
where $d_6 = \theta^{\frac{1}{|g_0|} \ln (\frac{c_2 d_3}{c_1})}$.
{\bf We will assume from now on that $N \geq N_0$}.  

\ms

Recall from Sect. 6.1 the family $\{\tcc_m\}_{m=1}^{m_0}$ of cylinders $\tcc_m = \tcc_m(b)$  in $\tR$ with (\ref{eq:ccmo-cond}) 
and (\ref{eq:ccm-cond}) having lengths $s_m$ satisfying (\ref{eq:C-cond}). It follows from Lemma 5.5 with 
$$\tS = \tS(m) = \frac{32\, E \, (C_6)^{1/\alpha_3}}{\ep_3}$$
that there exist finite families $\{\tGammam_j\}_{j=1}^{j_m} = \{\tGammam_j(b)\}_{j=1}^{j_m}$ and 
$\{\tLambdam_j\}_{j=1}^{j_m} = \{\tLambdam_j(b)\}_{j=1}^{j_m}$ of sub-cylinders of $\tcc_m = \tcc_m(b)$ 
for some integer  $j_m$ (depending on $m$ and $b$) such that 
$$\cup_{j=1}^{j_m} \tGammam_j \subset \tGammam \quad, \quad \cup_{j=1}^{j_m} \Lambdam_j = \Lambdam ,$$
where {\bf $\tGammam = \tGammam(b)$ and $\tLambdam = \tLambdam(b)$ are disjoint sub-cylinders  of $\tcc_m = \tcc_m(b)$} 
of co-lengths $\tq_0$ and $\tq_1$, respectively.
With the choice of $\tS$ we made, using appropriately (\ref{eq:diamte}) in 
Lemma 5.5, we get\footnote{In a similar way, it follows from the above and  (\ref{eq:C-cond}) that
$\diam(\tGammam_j) \leq \frac{\ep_3}{32 E |b|}$ and $\diam(\tLambdam_j) \leq \frac{\ep_3}{32 E |b|} $,
for all $j = 1, \ldots, j_m$, but we are not going to use these in what follows.}
$$\diamte(\tGammam_j) \leq \frac{1}{\tS} (\diamte (\tcc_m))^{1/\alpha_3} \quad , \quad
\diamte(\tLambdam_j) \leq \frac{1}{\tS} (\diamte (\tcc_m))^{1/\alpha_3} $$
for all $j = 1, \ldots, j_m$. 
Moreover, it follows from (5.47) in Lemma 5.5 that
\be
\frac{1}{d_3} \leq \frac{\nu(\hGammam_j)}{\nu(\hLambdam_j)} \leq d_3
\ee
for all $j = 1, \ldots, j_0$, where $\hGammam_j = \piU(\tGammam_j)$ and $\hLambdam_j = \piU(\tLambdam_j)$
and $d_3 > 1$ is the global constant from Lemma 5.5.

For later convenience we will now slightly {\bf change the notation involving the sub-cylinders $\tGammam_j$ and $\tLambdam_j$}.
Namely, set 
$$\tddmo = \tGammam \quad, \quad \tddmo_j = \tGammam_j \quad, \quad \tddmtwo = \tLambdam \quad, \quad
 \tddmtwo_j = \tLambdam _j . $$
We will also use the notation
$$\hddmt = \piU(\tddmt) \subset U \quad, \quad \hddmt_{j} = \piU(\tddmt_j) \subset U .$$
We should stress that $\tq_0$ and $\tq_1$ are global constants, so for the
global constant $d_4 \in (0,1)$ from Lemma 5.5 we have $\nu(\hddmt) \geq d_4 \nu(\tcc'_m) $
for all $t = 1,2$, all $m = 1, \ldots, m_0$ and all $|b|\geq b_0$. That is, we have
\be\label{eq:d4-cond}
\sum_{j=1}^{j_m} \nu(\hddmt_{j}) \geq d_4 \, \nu(\tcc'_m) 
\ee
for all $t = 1,2$ and all $m = 1, \ldots, m_0$, where $\tcc'_m = \piU(\tcc_m)$. With the new notation,
using (\ref{eq:ccm-cond}) as well,  the conditions
on the sub-cylinders $\tGammam_j$ and $\tLambdam_j$ translate as
\be\label{eq:dd-cond}
\diamte(\tddmt_j)  < \frac{\ep_3}{32 \, E\, |b|} 
\ee
for all $t = 1,2$, all $m = 1, \ldots, m_0$ and all $j = 1, \ldots, j_m$.
Apart from that, (6.23) gives
\be\label{eq:d3-cond}
\frac{1}{d_3} \leq \frac{\nu(\hddmo_j)}{\nu(\hddmtwo_j)} \leq d_3
\ee
for all $t = 1,2$, all $m = 1, \ldots, m_0$ and all $j = 1, \ldots, j_m$.



\ms




Next, {\bf fix a constant} $C_{9} > 0$ such that
$$\di C_{9} \geq \frac{32  d_3 E^2 C_6^2}{d_4}  .$$
Then {\bf fix a constant  $\beta_4  = \beta_4(N)$ so that}
$$0 < \beta_4 < \min\left\{ \frac{\mu_0\, \ep_0\, e^{-N T_0}}{10 C_{9}}\; , \; \frac{c}{8} \right\} ,$$
where $c > 0$ is the constant from Lemma 6.4 and {\bf a constant $a_0 = a_0(N) > 0$} so that
\be\label{eq:ao-cond}
0 < a_0 < \min\left\{ \frac{\mu_0 \ep_0\, e^{-NT_0}}{40 C_{9} D_1  N T_0} \; , \; \frac{\beta_4}{4 NT_0} 
\; , \; \frac{c}{32 NT_0} \right\}
\ee
where $D_1 > 1$ is one of the constants from (\ref{eq:D1-cond}).
Set
\be\label{eq:rho3}
\rho_3 = \frac{e^{a_0NT_0}}{1+ \frac{\mu_0 \,  e^{-NT_0}}{C_{9}}} < 1 \quad ,  \quad S_0 = e^{a_0 N T_0} > 1 .
\ee
Since $\log(1+y) > y/2$ for $y\in (0,1)$, it follows from (\ref{eq:ao-cond}) and $\ep_0 < 1$ that 
$$\log \left(1 + \frac{\mu_0 \, e^{-NT_0}}{C_{9}} \right) 
> \ep_0 \, \log \left(1 + \frac{\mu_0 \, e^{-NT_0}}{C_{9}} \right) 
> \frac{\mu_0 \, \ep_0\,  e^{-NT_0}}{2 C_{9}} > a_0 N T_0 ,$$ 
so $\rho_3  < 1$. 
Again by (\ref{eq:ao-cond}), $a_0NT_0 < c/2$, so 
$S_0 e^{-c} = e^{a_0 N T_0 - c} \leq e^{-c/2} < 1 .$
Also notice that the above estimate implies
$$\rho_3^{\ep_0} S_0^{1- \ep_0} 
= \frac{e^{a_0NT_0}}{\left(1+ \frac{\mu_0 \, e^{-NT_0}}{C_{9}}\right)^{\ep_0}} < 1 .$$
Thus, 
\begin{equation}\label{eq:rho-cond}
\rho_4 = \max \{ \rho_3^{\ep_0} S_0^{1- \ep_0}\; , \; S_0 e^{- c} \}   < 1 .
\end{equation}
Using again $\log(1+y) > y/2$ for $y\in (0,1)$ and the choice of $\beta_4$, we get
$$\ep_0 \, \log \left(1 + \frac{\mu_0 \,  e^{-NT_0}}{C_{9}} \right) > 5 \beta_4 .$$
Combining it with $a_0 NT_0 < \beta_4$, which follows from (\ref{eq:ao-cond}), the latter yields
$$e^{4 \beta_4}\, \rho_3^{\ep_0} S_0^{1- \ep_0} 
 = \frac{e^{4\beta_4}\, e^{a_0NT_0}}{\left(1+ \frac{\mu_0 \,  e^{-NT_0}}{C_{9}}\right)^{\ep_0}} 
< \frac{e^{5\beta_4}}{\left(1+ \frac{\mu_0 \,  e^{-NT_0}}{C_{9}}\right)^{\ep_0}} < 1. $$
Hence $\rho_3^{\ep_0} S_0^{1- \ep_0}  < e^{-4\beta_4}$. 
As we observed earlier, $S_0 e^{-c} \leq e^{-c/2}$, so the choice of $\beta_4$
now yields $S_0 e^{-c} < e^{-4 \beta_4}$. Therefore (\ref{eq:rho-cond}) implies
\be\label{eq:rho4}
\rho_4 <  e^{- 4 \beta_4} .
\ee

{\bf We now apply the construction at the end of the proof of Lemma 5.5.}
Given $|b| \geq b_0$ and $m = 1, \ldots, m_0(b)$,  for  any $j = 1, \ldots, j_m$ and any 
$i = 1,2$, using the fixed point  $Z_m \in \tP_0 \cap \tcc_m$, there exists  a  (H\"older) continuous map 
$$B^u(Z_m,\ep'') \ni x \mapsto \vm_{i,j}(x) = \vm_{i,j} (Z_m, x) \in U ,$$
such that $\sigma^N(\vm_{i,j}(x)) = x$  for all $x\in B^u(Z_m,\ep'')$
and the following property holds: 
\be\label{eq:IN-cond}
 I_{N}(x',z') =  \left|  \psi_m(x') - \psi_m (z') \right|   \geq \delta_1\, e^{-s_m \hep_{13}} \diam(\tcc_m)
\ee
for all $z\in \tddmtwo_{j}$ and  $x \in \tddmo_{j}$, where $z' = \piU(z)$, $x' = \piU(x) \in U$, and
$$\psi_m ( x) = \tau_{N}(\vm_{1,j}(x)) - \tau_{N}(\vm_{2,j}(x)) .$$
In the above we use the constant
$\delta_1 = 6 c_0  \delta_0 d_2 > 0$,
where $c_0$, $\delta_0$ and $d_2$ are as in Lemma 5.5.

For every $m = 1, \ldots, m_0$, using the point $Z'_m \in U$, {\bf fix maps $\vm_{i,j}(Z'_m, \cdot)$}
with the properties described above.
For any $i = 1,2$, $m = 1, \ldots, m_0$, and $j = 1, \ldots, j_m$ set 
$$\vm_{i,j} = \vm_{i, j}(Z'_m, \cdot) \quad , \quad \Xmt_{i,j} = \vm_{i,j}(\tddmt_{j}) \subset U .$$
By construction  the sets $\Xmt_{i,j}$ are disjoint and satisfy $\sigma^N(\Xmt_{i,j}) = \tddmt_{j}$.
By Lemma 4.1(a), the {\it characteristic functions} 
$$\di \ommt_{i,j} = \chi_{\Xmt_{i,j}} : \hU \longrightarrow [0,1] $$
of $\Xmt_{i,j}$ belongs to $\ff_\theta(\hU)$ and  
$\di \Lip_\theta(\ommt_{i,j}) \leq 1/\diamte(\Xmt_{i,j}) .$

The so called contraction operators are now defined similarly to what was done in \cite{St4}.
A subset $J$ of the set
$$\Pi(b)  = \{\; (t, m, i, j ) \; :  \;  t = 1,2 \; ,\;  
1 \leq m \leq m_0\;, \; i = 1,2\; , \;1\leq j\leq j_m \; \;\} $$
will be called {\it representative}  if either all elements of $J$ have the form $(1,m, i, j)$ and
\be\label{eq:1-cond}
\sum_{(1,m, i, j) \in J} \nu(\hddmo_{j}) \geq \frac{d_4}{8 d_3} \nu(\tcc'_m) ,
\ee
or  all elements of $J$ have the form $(2,m,i, j)$ and
\be\label{eq:2-cond}
\sum_{(2,m,i, j) \in J} \nu(\hddmtwo_{j}) \geq \frac{d_4}{8d_3} \nu(\tcc'_m) .
\ee
It follows from Remark 2(c) above and (\ref{eq:1-cond}) that 
whenever $J$ is representative and its elements are of type $(1,m,i,j)$, then 
the smallest sub-cylinder $X$ of $\hddmo$ that contains all $\hddmo_j$ for all
$(1,m,i,j) \in J$ has co-length $\leq \ell_0$ in $\cc'_m$. Similarly, 
the smallest sub-cylinder $Y$ of $\hddmtwo$ that contains all $\hddmtwo_j$ for all
$(2,m,i,j) \in J$ has co-length $\leq \ell_0$ in $\cc'_m$.

Let $\jj(b)$ be the {\it family of all representative subsets} $J$  of $\Pi(b)$.
Given $J  \in \jj(b)$, define the function    $\omega_{J}  : \hU \longrightarrow [0,1]$ by
$$\di \omega_J  = 1- \mu_0 \,\sum_{(t,m,i,j) \in J} \ommt_{i,j} .$$
Then $\omega_J \in \ff_\theta(\hU)$ and 
$\frac{3}{4} \leq 1-\mu_0 \leq \omega_J (u) \leq 1$ for all  $u \in \hU .$
Define the {\it contraction operator} 
$$\nn = \nn_J(a,b) : \ff_\theta (\hU) \longrightarrow \ff_\theta (\hU) \quad  \mbox{\rm by } \quad \nn h = \ma^{N} (\omega_J \cdot h) ,$$
where $\ma = L_{\fa}$ is defined in Sect. 4.2. 

\bs

\noindent
{\bf Remark 3.} The contraction operators $\nn_J$  look similar to the operators 
defined and studied by Dolgopyat in Sect. 7 in \cite{D1},
although here the construction is necessarily much more complicated. In fact it is significantly more complicated than the constructions
used in \cite{St2} and \cite{St4} as well.

\bs

We will now prove some basic properties of the contraction operators.
Our exposition here is similar to that in Sect. 6.2 in \cite{St4}, however there are some substantial differences.


\medskip


\subsection{Estimates for  eigenfunctions and Lasota-York type inequality}

For any $u, u'\in \hU$, we denote by $\ell(u,u') \geq 0$  {\it the length of the smallest cylinder $Y(u,u')$  
 in $\hU$ containing} $u$ and $u'$. For later use,
notice that for every $p \geq 1$, $\sigma^p(Y(u,u'))$ is the smallest cylinder\footnote{Indeed, let $X$ be a cylinder in 
$U$ containing both  $\sigma^p(u)$ and $\sigma^p(u')$ and let $X \subset \sigma^p(Y(u,u'))$. Then every $x \in X$ 
has the form $x = \sigma^p(y)$ for some $y \in Y(u,u')$.  Let 
$Y' = \{y \in Y(u,u') : \sigma^p(y) \in X\} .$
Since $X$ is a cylinder in $U$, $Y'$ is a cylinder in $U$, too. Now $u, u' \in Y'$ imply
$Y(u,u') \subset Y'$, therefore $\sigma^p(Y(u,u')) \subset X$. This proves that 
$\sigma^p(Y(u,u')) = Y(\sigma^p(u), \sigma^p(u'))$.}
 $Y(\sigma^p(u), \sigma^p(u'))$ in $U$ containing both $\sigma^p(u)$ and $\sigma^p(u')$.

Consider the following 

\ms

{\bf Assumption: } for points $u,u' \in \cc'_m$ for some $m = 1, \ldots, m_0(b)$, 
$i= 1, 2$, $j = 1, \ldots, j_m(b)$, an integer $p \geq 0$ and points $v,v'\in U$ we have:
\be\label{eq:u-cond}
\sigma^p(v) = \vm_{i,j}(u)\;, \; \sigma^p(v')  = \vm_{i,j}(u')\;,\; \ell (v,v') \geq p  .
\ee


Notice that the latter implies $\ell (v,v') \geq N + p$, $\sigma^{N+p}(v) = u$ and $\sigma^{N+p}(v') = u'$.

The  following estimate is similar to the one in Lemma 6.6 in \cite{St4}.

\bs

\noindent
{\bf Lemma 6.5.} {\it 
If the points  $u,u' \in U$, the  cylinder $\cc'_m$, the integer 
$p \geq 0$ and the points $v,v'\in U$ satisfy {\rm (\ref{eq:u-cond})}  for some $j = 1, \ldots, j_m$ and
$i = 1,2$, and $w,w'\in U$ are such that $\sigma^N w = v$, 
$\sigma^N w' = v'$ and $\ell (w,w') \geq N$, then
$$|\tau_N(w) - \tau_N(w')| \leq \frac{C_4}{d_0^2} \, \theta^{p+N} \, \diam(\tcc_m) 
 \leq \frac{C_1 C_4}{d_0^2} \, \theta^{p+N} \, \diamte(\tcc_m),$$
where $d_0 \in (0,1]$ is the constant fixed at the end of Sect. 2 and $C_4 \geq 1$ is the constant from
Lemma 5.3.
If we assume in addition that $u,u' \in \hddmt_j$ for some $t = 1,2$ and some $j = 1, \ldots, j_m$, then }
$$|\tau_N(w) - \tau_N(w')| \leq \frac{C_4}{d_0^2} \, \theta^{p+ N} \, \diam(\tddmt_j)
\leq \frac{C_1 C_4}{d_0^2} \, \theta^{p+ N} \, \diamte(\tddmt_j)  .$$

\bs

\noindent
{\it Proof.} Assume that the points $u,u',v,v',w,w'$ and the cylinder $\cc'_m$
satisfy the assumptions in the lemma. Clearly, $\ell(w,w') \geq p+2N$ and
\be
\tau_N(w) - \tau_N(w') = [\tau_{p+2N}(w) - \tau_{p+2N}(w')] - [\tau_{p+N}(v) - \tau_{p+N}(v')] .
\ee

Consider now some  fixed $j = 1, \ldots, j_m$ and $i = 1,2$.
Recall the construction of the map $\vm_{i,j}$ from the proof of Lemma 5.5. In particular by (5.70),
$$\pp^N(\vm_{i,j}(u)) = \phi_{[-\ep_0,\ep_0]}(W^s_{\ep_0} (u)) \cap W^u_{R}(\bj_i) ,$$
where we set $\bj_i = \bj_i (Z_m) \in W^s_{R}(Z_m)$ for brevity. 
Since $\sigma^{p}(v) = \vm_{i,j}(u)$ and $\sigma^p(v') = \vm_{i,j}(u')$, we have $\sigma^{p+2N}(w) = \sigma^{p+N}(v) = u$ and 
$\sigma^{p+2N}(w') = \sigma^{p+N}(v') = u'$, so both
$x' = \pp^{p+N}(v)$ and $ z' = \pp^{p+N}(v')$ belong to  $W^u_{R}(b')$
for some $b' \in W^s_{R}(Z_m)$. Then $\piU(x') = u$ and $\piU(z') = u'$.
Moreover, $\pp^p(v) \in W^s_{R}(\vm_{i,j}(u))$ and the choice of $N$ imply (as in the proof of Lemma 5.5) that 
$d(\bj_i, b') < \delta'$, the constant with (\ref{eq:delta'}) from Lemmas 5.4 and 5.5.
Similarly, 
$x'' = \pp^{p+2N}(w)$ and $z'' = \pp^{p+2N}(w')$ belong to $ W^u_{R}(b'')$ 
for some $b'' \in W^s_{R}(Z_m)$ with $d(\bj_i, b'') < \delta'$, and $\piU(x'') = u$, $\piU(z'') = u'$. 
Thus, $x', x'' \in W^s_R(u)$ and $z', z'' \in W^s_R(u')$.
Moreover, since the local stable/unstable holonomy maps are uniformly 
$\alpha_1$-H\"older, by the choice of $d_0 \in (0,1]$ at the end of Sect. 2,
$$d(b', b'') \leq \frac{1}{d_0} (d(\pp^{p+N}(v), \pp^{p+2N}(w)))^{\alpha_1} .$$
Using this and (2.1) for points on local stable manifolds, i.e. going backwards along the flow, we get
\be
d(b', b'') \leq \frac{1}{d_0} (d(\pp^{p+N}(v), \pp^{p+2N}(w)))^{\alpha_1} 
\leq \frac{1}{d_0}\, \left( \frac{d(v, \pp^N(w))}{d_0 \gamma^{p+N}} \right)^{\alpha_1} 
\leq \frac{1}{d_0^{2} \gamma^{\alpha_1 (p+N)}} .
\ee
By the choice of $\theta$ (see Sect. 6.1) we have $\di \frac{1}{\gamma^{\alpha_1 \beta}} \leq \theta$, hence
$$(d(b', b''))^{\beta} \leq  (1/d_0^{2})^{\beta} (1/\gamma^{p+N})^{\alpha_1 \beta} 
\leq  \frac{\theta^{p+N}}{d_0^2} .$$

We are preparing to use Lemma 5.3. Let $\tu \in \tR$ and $\tu' \in \tR$ be the shifts along the flow of 
the points $\pi_{Z_m}(u)$ and $\pi_{Z_m}(u')$.
Then we have $\tu = \phi_{t(u)} (\pi_{Z_m}(u))$ and  $\tu' = \phi_{t(u')} (\pi_{Z_m}(u'))$ 
for some small $t(u), t(u') \in \R$. So
\begin{eqnarray*}
\tau_{p+N}(v) - \tau_{p+N}(v') 
& = & \Delta(\pp^{p+N}(v), \pp^{p+N}(v')) = \Delta(x',z') = \Delta(u, \pi_{b'}(u'))\\
& = & \Delta(\pi_{Z_m}(u) ,\pi_{b'}( \pi_{Z_m}(u')) ) =  \Delta(\tu, \pi_{b'}(\tu')) + t(u) - t(u') ,
\end{eqnarray*}
and similarly
\begin{eqnarray*}
\tau_{p+2N}(w) - \tau_{p+2N}(w') 
& = & \Delta(\pp^{p+2N}(w), \pp^{p+2N}(w')) 
=  \Delta(\tu, \pi_{b''}(\tu')) + t(u) - t(u') .
\end{eqnarray*}
This, the above estimate  and Lemma 5.3(a)  yield
\begin{eqnarray*}
|\tau_N(w) - \tau_N(w')|
 =      |  \Delta(\tu, \pi_{b'}(\tu')) -  \Delta(\tu, \pi_{b''}(\tu'))| 
\leq C_4 \diam (\tcc_m) \, (d(b',b''))^{\beta}
 \leq  \frac{C_4}{d_0^2} \theta^{p+N}  \diam(\tcc_m).
\end{eqnarray*}
Finally, by Lemma 4.1 we have $ \diam(\tcc_m) \leq C_1  \diamte(\tcc_m)$.

Next, assume that $u,u' \in \hddmt$ for some $t$. We can then apply the above argument
replacing $\cc'_m$ by $\hddmt_{j}$ to get
$\di |\tau_N(w) - \tau_N(w')| \leq \frac{C_4}{d_0^2} \, \theta^{p+N} \, \diam(\tddmt_j) .$
This proves the lemma.
\endofproof

\bs


Set  
$\di E_1 =   \frac{E}{120} ,$
where $E$ is as in (\ref{eq:E-cond}). Then $E_1 \geq C_{8} e^{C_{8}}$, where
$C_{8} = \frac{T_0 C_{6} C_4}{1-\theta_1} + \frac{M_1 C_4}{d_0^2}$, as defined  in Sect. 6.3.
Here  $M_1$ is a constant so that $M_1 \geq 2 |P_f| + 1$  (see  Sect. 4.2 for $P_f$). 

\ms

Denote by $\kk_0$
{\it the set of all {\bf $h \in \ff_{\theta}(U)$} such that $h \geq 0$ on $U$ and for any 
$u,u' \in U$ contained in some cylinder $\tddmt_j$ for some $t = 1,2$, some $1\leq m \leq m_0$
and some $j = 1, \ldots, j_m$, 
any integer $p \geq 0$ and any points $v,v'\in U$ satisfying {\rm (\ref{eq:u-cond})} we have}
\be
|h(v) - h(v') | \leq E_1 \, \theta^{p+N} \, h(v')\, \diamte (\tddmt_j) .
\ee

\ms

It turns out that the eigenfunctions $h_a \in \kk_0$ for $|a| \leq a_0$ (see Sect. 4.2). 
This follows from the following lemma whose proof  is almost identical, modulo the different assumption 
in (\ref{eq:u-cond}), with that of Lemma 6.7 in \cite{St4}. 
See Appendix I for a sketch of the proof.

\bs

\noindent
{\bf Lemma 6.6.} {\it For any real constant $s$ with $|s| \leq M_1$ we have 
$L^{m N}_{f - s \tau}(\kk_0) \subset \kk_0$ for all $m \geq 1$.} 

\bs

\noindent
{\bf Corollary 6.7.} {\it For any real constant $a$ with $|a| \leq a_0$ we have $h_a \in \kk_0$.} 

\bs

\noindent
{\it Proof.} Let $|a| \leq a_0$. Since the constant  function $h = 1 \in \kk_0$, it follows from 
Lemma 6.6 that $L^{m N}_{f-(P+a)\tau} 1 \in \kk_0$ for all $m \geq 1$. Now the Ruelle-Perron-Frobenius 
Theorem (see e.g. \cite{PP}) and the fact that $\kk_0$ is closed in $\ff_{\theta}(\hU)$ imply  $h_a\in \kk_0$. 
\endofproof

\bs



\def\Gammamm{\Gamma^{(m')}}
\def\ttheta{\tilde{\theta}}
\def\tdet{D_{\ttheta}}
\def\tdiamte{\diam_{\ttheta}}
\def\dteo{D_{\theta_1}}

\def\ttheta{\tilde{\theta}}
\def\dtte{D_{\ttheta}}
\def\diamtte{\diam_{\ttheta}}
\def\tff{\widetilde{\ff}}


The following lemma  is similar to  Lemma 6.4 in \cite{St4}.

\bs


\noindent
{\bf Lemma 6.8.} {\it Assume that $u,u' \in \hU$, $u \neq u'$, and $\sigma^N(v) = u$, $\sigma^N(v') = u'$
for some $v,v'\in \hU$ with $\ell(v,v') \geq N$.
Let $\omega_J(v) < 1$ and $\omega_J(v') = 1$ for some $J \in \jj(b)$.  Then 
$$\di |\omega_J (v) - \omega_J (v')| \leq \frac{\mu_0 \theta_1^{-\ell_0}}{\ep_4} \, |b|\,  \dteo (u,u') ,$$
where $\ell_0 \geq 1$ is the constant from Remark} 2(c).

\bs

\noindent
{\it Proof.}  Under the given assumption we have $|\omega_J (v) - \omega_J (v')|  = \mu_0$.
Assume e.g. that all elements of $J$ have the form $(1,m,i, j)$ and (\ref{eq:1-cond}) holds.
Let $\omega_J(v) < 1$; then $v\in \Xmt_{1,j}$ for some $(1,m,i, j) \in J$, so  $u = \sigma^N(v) \in \tddmt_{j}$.  
However $u' = \sigma^N(v') \notin \tddmt_{j'}$ for any $j'$ with $(1,m,i,j') \in J$. 
Using (\ref{eq:1-cond}), it follows from Remark 2(c)  that 
the smallest sub-cylinder $X$ of $\hddmo$ that contains all $\hddmo_j$ for all
$(1,m,i,j) \in J$ has co-length $\leq \ell_0$ in $\cc'_m$. Since $u'\notin X$ and $u \in X$, it follows from (\ref{eq:ccmo-cond})  that
$$\dteo (u,u') \geq \diamteo (X) \geq \theta_1^{\ell_0} \diamteo(\tcc_m) \geq \frac{\theta_1^{\ell_0}\, \ep_4}{|b|} .$$ 
Hence
\begin{eqnarray*}
|\omega_J (v) - \omega_J (v')| 
 =  \mu_0 \, \frac{\dteo(u,u')}{\dteo(u,u')}   \leq \mu_0 \frac{\dteo(u,u')}{\theta_1^{\ell_0}\, \ep_4/|b| } 
 = \frac{\mu_0 \theta_1^{-\ell_0}}{\ep_4} |b|\,  \dteo(u,u') .
\end{eqnarray*}
This proves the lemma. 
\endofproof

\bs


Recall the large constant $E$ with (\ref{eq:E-cond}).

\bs

\noindent
{\bf Definition 6.9.}
Denote by $\kk_{E|b|}$  be {\it the set of all functions $H \in \ff_{\theta}(\hU)$  such that } $0 < H \leq 1$ on $\hU$ and
$$\frac{|H(u) - H(u')|}{H(u')} \leq E\,|b|  \dteo (u,u') $$ 
for all $u,u'\in \hU$.

\ms

We can now derive a Lasota-Yorke type inequality for functions in $\kk_{E|b|}$.
Its proof is similar to that of Proposition 6 in \cite{D1} and that of Lemma 5.6 in \cite{St2}. 
We prove it in  Appendix I for completeness.

\bs

\noindent
{\bf Lemma 6.10.} {\it Let $f \in \ff_{\theta_1}(\hU)$. Then for any $J \in \jj(b)$ we have  $\nn_J(\kk_{E|b|}) \subset \kk_{E|b|}$.} 





\section{Iteration procedure -- the role of the contact structure}
\setcounter{equation}{0}

We continue here with the notation and the assumptions in Sect. 6. 
Let $|b| \geq b_0$.

Denote by   {\it $\kk_b$ the set of all pairs $(h, H)$ such that $h \in \ff_{\theta}(\hU)$, 
$H \in \kk_{E|b|}$, and the following two conditions are  satisfied:}

\ms

$(T1)$\quad  $|h| \leq H \leq 1$ on $\hU$,

\ms

$(T2)$ \quad for any $u,u' \in \hddmt_j =  \hddmt_j(b)$ for some $m = 1, \ldots, m_0$, some $t = 1,2$ 
and some $j = 1, \ldots, j_m$, any integer  
$p \geq 0$ and any points $v,v'\in \hU$ satisfying (\ref{eq:u-cond})  for $ m$ we have
\be\label{eq:h-cond}
 |h(v) - h(v') | \leq E \, |b|\, \theta^{p+N} \, H(v') \, \diamte (\tddmt_j) .
\ee

Notice that (\ref{eq:dd-cond}) implies
\be\label{eq:Eb-cond}
E\, |b|  \, \diamte (\tddmt_j)  < \frac{\ep_3}{32} .
\ee

The following lemma is fundamental for the iteration procedure that will be used in Sect. 8 which will  show that
the so called contraction operators are "eventually" contracting. The idea behind all this is in Lemma $10''$ in \cite{D1}
(see Sects. 6, 7 and 8 in \cite{D1}), although here we have to proceed in a slightly more complicated way.

\bs

\noindent
{\bf Lemma 7.1.} {\it For any $|a| \leq a_0$, any $|b| \geq b_0$, and any $(h,H) \in \kk_b$ 
there exists $J \in \jj(b)$ such that $(\lab^{N} h, \nn_J H ) \in \kk_b$.}

\bs

To prove this  we need the following lemma, whose proof is very similar to that of 
Lemma 14 in \cite{D1} (and essentially the same as that of Lemma 6.10 in \cite{St4}). 
For completeness we prove it in Appendix I.

\bs

\noindent
{\bf Lemma 7.2.}  {\it Let $(h,H) \in \kk_b$. Then for any $m = 1, \ldots, m_0$,  any $j = 1, \ldots, j_m$, any $t = 1,2$
and any $i = 1,2$ we have:}

\ms

(a) {\it $\di\frac{1}{2} \leq \frac{H(\vm_{i,j} (u'))}{H(\vm_{i,j}(u''))} \leq 2$ for all} $u', u'' \in \hddmt_j$;

\ms

(b) {\it Either $H(\vm_{1,j}(u)) \geq H(\vm_{2,j}(u))/4$ for all $u \in \hddmt_j$ or $H(\vm_{2,j}(u)) \geq H(\vm_{1,j}(u))/4$ 
for all $u \in \hddmt_j$.}

\ms

(c) {\it Either for all $u\in \hddmt_j$ we have
$|h(\vm_{i,j} (u))|\leq \frac{3}{4}H(\vm_{i,j} (u))$, or $|h(\vm_{i,j} (u))|\geq \frac{1}{4}H(\vm_{i,j} (u))$ 
for all $u\in \hddmt_j$.} 

\bs

\noindent
{\it Proof of Lemma} 7.1. 
Let $|a| \leq a_0$, $|b| \geq b_0$ and let $(h, H) \in \kk_b$. We will construct a representative set 
$J \in \jj(b)$ such that  $(\lab^N h , \nn_J H) \in \kk_b$. 

Notice that, since $1/2 \leq \omega_J \leq 1$, we have
$\frac{1}{2} \ma^N(H) \leq \nn_J(H) = \ma^N(\omega_J\, H)$.

Consider for a moment an arbitrary (at this stage) representative set $J \in \jj(b)$. 
We will first show that  $(\lab^N h , \nn_J H) $ has property $(T2)$. This is done as in the proof
of Lemma 6.9 in \cite{St4}. 

Assume that the cylinder $\tddmt_j$ in $U$, the points $u,u' \in \tddmt_j$, the integer $p \geq 0$ and the points 
$v,v'\in \hU$  are satisfying (\ref{eq:u-cond}) and
(\ref{eq:h-cond})  for some $t = 1,2$, $m = 1, \ldots, m_0$ and $j = 1,\ldots, j_m$.

From the definition of $\fa$, for any $w, w'$ with $\sigma^N w = v$, $\sigma^N(w') = v'$  and $\ell (w,w') \geq N$ we have 
$$\fa_N(w) 
 =   f_N(w) - (P+a)\tau_N(w) + \ln h_a (w)  - \ln h_a (v) - N \lambda_a .$$
Since $h_a \in \kk_0$ by Corollary 6.7,
\begin{eqnarray*}
|\ln h_a(w) - \ln h_a(w') | 
& \leq &\frac{| h_a(w) - h_a(w')|}{\min \{ |h_a(w)|, |h_a(w')|\}} 
\leq E_1 \, \theta^{p+2N} \, \diamte (\tddmt_j) .
\end{eqnarray*}
and similarly, $|\ln h_a(v) - \ln h_a(v') | \leq E_1 \,\theta^{p+ N} \, \diamte (\tddmt_j)$.
Using   $|f|_{\theta} \leq T_0$, we get
\begin{eqnarray*}
         |f_N(w) - f_N(w')| 
& \leq  & \sum_{j=0}^{N-1} |f (\sigma^j(w)) - f (\sigma^j(w'))|
 \leq  \sum_{j=0}^{N-1} |f|_{\theta} \,\theta^{N-j}\, D_{\theta} (v,v')\\
& \leq &  \frac{T_0 }{1-\theta}\,\theta^{p+N} \,D_{\theta}(u,u') 
\leq \frac{T_0}{1-\theta} \theta^{p+ N}\, \diam_{\theta} (\tddmt_j)\\
& \leq & E_1 \, \theta^{p+N} \, \diamte (\tddmt_j) .
\end{eqnarray*}
Apart from that it follows from Lemma 6.5 and the choice of $E_1$ that
$$ |P+a|\, |\tau_N(w) - \tau_N(w') | 
\leq \frac{2|P+a|\, C_4}{d_0^2}\,  \theta^{p+N}\, \diamte (\tddmt_j)
\leq E_1\, \theta^{p+N}\, \diamte (\tddmt_j) .$$
Using the above we get
\begin{eqnarray}
 |\fa_N(w) - \fa_N(w')| 
\leq  4 E_1 \, \theta^{p+N}\, \diamte (\tddmt_j) < \frac{\ep_3}{32} ,
\end{eqnarray}
using (\ref{eq:Eb-cond}) for the last estimate. This implies 
$$ \left| (\fa_N - \i b \tau_N)(w) -  (\fa_N - \i b \tau_N) (w')\right|
\leq  \left(4E_1 + E_1\, |b| \right)\, \theta^{p+N}\, \diamte (\tddmt_j) < \frac{\ep_3}{16} < \ep_3 ,$$
using (\ref{eq:Eb-cond}) again for the latter.

Notice that (7.1) for the pair $w,w'$ gives 
$$ |h(w) - h(w') | \leq E \, |b|\, \theta^{p+2N} \, H(w') \, \diamte (\tddmt_j) .$$
Using  this and (7.3) we now derive (as we did in the proof of Lemma 6.9 in \cite{St4})  
\begin{eqnarray*}
 & &      |(\lab^N h)(v) - (\lab^N  h)(v')|
=       \left| \sum_{\sigma^N w = v} e^{(\fa_N - \i b \tau_N) (w)}\, h(w)   
-  \sum_{\sigma^N w = v} e^{(\fa_N - \i b \tau_N)(w'(w))}\, h(w'(w)) \right| \\
& \leq & \left| \sum_{\sigma^N w = v} e^{(\fa_N - \i b \tau_N)(w)}\, [h(w) -  h(w')]\right| 
        + \sum_{\sigma^N w = v}  \left|e^{(\fa_N - \i b \tau_N)(w)} -  e^{(\fa_N - \i b \tau_N) (w')}\right| \, |h(w')| \\
& \leq & \sum_{\sigma^N w = v} e^{(\fa_N (w) - \fa_N(w'))} e^{\fa_N (w')}\, E |b| \, \theta^{p+2 N} \,H(w') \, \diamte (\tddmt_j)  \\
&   &     + \sum_{\sigma^N w = v}  \left| e^{(\fa_N - \i b \tau_N)(w) -  (\fa_N - \i b  \tau_N)(w')} - 1\right| \,  e^{\fa_N(w')} H(w') \\
& \leq & e^{\ep_3}\, E |b|\, \theta^{p+2N}\,  \diamte (\tddmt_j)\,  (\ma^N H)(v')        
+  e^{\ep_3} \, (4E_1 + E_1 |b| ) \theta^{p+N} \, \diamte  (\tddmt_j)\,  (\ma^N H)(v') .
\end{eqnarray*}
As remarked earlier, $\ma^N H \leq 2 \ma^N (\omega_J H) =  2 \nn_J H$. 
This and the above yield
\begin{eqnarray*}
 |(\lab^N h)(v) - (\lab^N  h)(v')|
 & \leq &  3 E |b|\, \theta^{p+2N}\, \diamte (\tddmt_j)\, 2 (\nn_J H)(v') \\
& &       +  30 E_1\, |b| \, \theta^{p+N} \, \diamte (\tddmt_j)\,  2 (\nn_J H)(v') \\
& \leq & (6 \, \theta^{N} + 60 E_1/E)  \, E \, |b|\, \theta^{p+N}\,  \diamte (\tddmt_j)\,  (\nn_J H)(v')\\
& \leq &  E \,|b|\, \theta^{p+N}\,   \diamte (\tddmt_j)\,  (\nn_J H)(v') ,
\end{eqnarray*}
since $6 \, \theta^{N} + 60 E_1/E \leq 1$ by (\ref{eq:No-cond}) and the choice of $E_1$.
Thus, $(\lab^N h , \nn_J H)$ has property $(T2)$.

\ms

\def\ffo{\ff^{(1)}}
\def\fftwo{\ff^{(2)}}
\def\tJ{\widetilde{J}}

So far the choice of $J$ was not important. We will now construct a representative set $J = \jj(b)$ so that 
$(\lab^N h , \nn_J H)$ has property $(T1)$, namely 
\be\label{eq:Lab}
|\lab^N h|(u) \leq (\nn_J H) (u)
\ee
for all $u \in \hU$.

Notice that (\ref{eq:Lab}) is trivially satisfied for $u \notin V_b$ for any choice of $J \in \jj(b)$.
So we need to deal with those $u$ that belong to $\cc'_m = \cc'_m(b)$ for some $m = 1, \ldots,  m_0 = m_0(b)$.

{\bf Fix an arbitrary  $m = 1, \ldots, m_0$.}
We will construct  a family of triples $(t,i,j)$ with $t = 1,2$, $i = 1,2$ and $j = 1, \ldots, j_m$
so that $(t,m,i,j)$ will be included in $J$, namely a family which satisfies (\ref{eq:1-cond})
or (\ref{eq:2-cond}).
That is, for the given $m$, we need to construct either a family
$$\ffo_m \subset  \{ (1,i,j) : i = 1,2 \; ;\;  1 \leq j \leq j_m\} ,$$
such that 
\be\label{eq:Fmo-cond}
\sum_{(1,m,i,j) \in \ffo_m} \nu(\hddmo_j)  \geq \frac{d_4}{8 d_3} \nu(\tcc'_m) ,
\ee
and  (\ref{eq:Lab}) holds for all $u \in \hddmo_j$ whenever $(1,i,j) \in \ffo_m$,
or a family
$$\fftwo_m \subset  \{ (2,i,j) : i = 1,2 \; ;\;  1 \leq j \leq j_m\} ,$$
such that 
\be\label{eq:Fmtwo-cond}
\sum_{(2,m,i,j) \in \fftwo_m} \nu(\hddmtwo_j)  \geq \frac{d_4}{8 d_3} \nu(\tcc'_m) ,
\ee
and  (\ref{eq:Lab}) holds for all $u \in \hddmo_j$ whenever $(2,i,j) \in \fftwo_m$.

Once $\ffo_m$ or $\fftwo_m$ is constructed, we will set $J = \ffo_m$ or $J = \fftwo_m$, respectively.
Set
$$\tff_m = \{ (t,i,j) : t= 1,2\; , \; i = 1,2 \; ;\;  1 \leq j \leq j_m\} .$$

Define the functions $\tpsi_m, \gao_m, \gat_m: \hU  \longrightarrow \C$  by
$$\di \tpsi_m(u) = e^{(\fa_{N}+\i b\tau_{N})(\vm_{1,j}(u))} h(\vm_{1,j}(u)) 
+ e^{(\fa_{N}+\i b\tau_{N})(\vm_{2,j}(u))} h(\vm_{2,j}(u)) ,$$
$$\di \gao_m (u) = (1-\mu_0)\, e^{\fa_{N} (\vm_{1,j}(u))} H(\vm_{1,j}(u)) + e^{\fa_{N}(\vm_{2,j}(u))} H(\vm_{2,j}(u)) ,$$
while $\gat_m (u)$ is defined similarly with a coefficient $(1-\mu_0)$ 
in front of the second term. 
Recall the functions 
$$\di \psi_m (u) = \tau_{N}(\vm_{1,j}(u)) - \tau_{N}(\vm_{2,j}(u)) \quad, \quad u \in U ,$$
that appear in (\ref{eq:IN-cond}).

Next, denote by $\ff'_m$ the set of those $(t,i,j) \in \tff_m$ so that 
the  first alternative in Lemma 7.2(c) holds for $t, i, j$, and by $\ff''_m$  the set of those $(t, i,j) \in \tff_m\setminus \ff'_m$ so that 
the second alternative in Lemma 7.2(c) holds for $t, i, j$. 

\ms

{\bf Case 1.} Assume that
$\di \sum_{(1,1, j)\in \ff'_m} \nu(\hddmo_j) 
\geq \frac{d_4}{8 d_3} \nu(\tcc'_m)$.
Consider an arbitrary $(1,1, j) \in \ff'_m$. Given $u \in \hddmo_j$, $1- \mu_0 \geq 3/4$ and $|h(\vm_{1,j} (u))|\leq \frac{3}{4}H(\vm_{1,j} (u))$
imply $|\tpsi_m(u)| \leq \gao_m(u)$. So if $(1,m,1, j) \in J$ with $(1,1,j) \in \ff'_m$, then
\begin{eqnarray*}
&        & \left| (\lab^N h)(u)\right| 
 \leq  \left| \sum_{\sigma^N v = u, \;v\neq \vm_{1,j}(u), \vm_{2,j}(u)} e^{(\fa_N+\i b\tau_N)(v)} h(v) \right|  + |\tpsi_m(u)| \nonumber\\
& \leq & \sum_{\sigma^N v = u, \;v\neq \vm_{1,j}(u), \vm_{2,j}(u)} e^{\fa_N(v)} |h(v)|  + \gao_m (u)\nonumber\\
& \leq & \sum_{\sigma^N v = u, \;v\neq \vm_{1,j}(u),\vm_{2,j}(u)} e^{\fa_N(v)} \omega_J(v) H(v)\nonumber\\
&      & + \left[e^{\fa_N(v_1(u))} \omega_J(\vm_{1,j}(u)) H(\vm_{1,j}(u)) 
+  e^{\fa_N(\vm_{2,j}(u))} \omega_J(\vm_{2,j}(u)) H(\vm_{2,j}(u))\right]  
 \leq   (\nn_J H) (u) . 
\end{eqnarray*}
Thus, in this case we can simply take $\ffo_m = \{(1,1, j) : (1,1,j) \in \ff'_m\}$ and then (\ref{eq:Lab}) and (\ref{eq:Fmo-cond}) will be satisfied.

\ms

{\bf Case 2.} Assume that 
$\di \sum_{(1,1,j)\in \ff'_m} \nu(\hddmo_j) < \frac{d_4}{8 d_3} \nu(\tcc'_m)$.
Since $\sum_{j=1}^{j_m} \nu(\hddmo_j) \geq d_4 \nu(\tcc'_m)$ by (\ref{eq:d4-cond}), $d_3 > 1$ and $(1,i,j) \notin \ff'_m$ 
implies $(1,i,j) \in \ff''_m$, it follows that 
\be
\sum_{(1,1,j) \in \ff''_m}   \nu(\hddmo_j)  \geq \frac{7d_4}{8 d_3} \nu(\tcc'_m) .
\ee

\ms

{\bf Sub-case 2.1.} Assume that
\be
\sum_{(1,2,j)\in \ff'_m } \nu(\hddmo_j) 
\geq \frac{d_4}{8 d_3} \nu(\tcc'_m) .
\ee
As in Case 1 one shows that for all $j$ with $(1,2,j) \in \ff'_m$, including $(1,2,m,j)$ in $J$, we have
$\left| (\lab^N h)(u)\right| \leq (\nn_J H) (u)$ for all $u \in \hddmo_j$.
Thus, setting 
$$\ffo_m = \{(1,2,j) : (1,2,j) \in \ff'_m \} ,$$ 
(\ref{eq:Lab}) and (\ref{eq:Fmo-cond}) will be satisfied.

\ms

{\bf Sub-case 2.2.} Assume that (7.8) does not hold, that is 
$\di \sum_{(1,2,j)\in \ff'_m } \nu(\hddmo_j) < \frac{d_4}{8 d_3} \nu(\tcc'_m)$.
It follows now that
$$\sum_{(1,1,j)\in \ff'_m } \nu(\hddmo_j) + \sum_{(1,2,j)\in \ff'_m } \nu(\hddmo_j) < \frac{d_4}{4 d_3} \nu(\tcc'_m) .$$
Using (\ref{eq:d3-cond}) we can write 
\be
\sum_{(1,1,j)\in \ff'_m } \nu(\hddmtwo_j) + \sum_{(1,2,j)\in \ff'_m } \nu(\hddmtwo_j) < \frac{d_4}{4} \nu(\tcc'_m) .
\ee
Denote by {\it $J'_m$ the set of those $j \in \tff_m$ such that both $(1,1,j) \in \ff''_m$ and $(1,2,j) \in \ff''_m$}.
Then the above and (\ref{eq:d4-cond}) with $t = 2$ imply
\be
\sum_{j\in J'_m } \nu(\hddmtwo_j) > \frac{3d_4}{4} \nu(\tcc'_m) .
\ee

In what follows we assume (7.10).

\ms

{\bf Sub-case 2.2.1.}  Assume now that
$$\di \sum_{(2,1, j)\in \ff'_m} \nu(\hddmtwo_j)  \geq \frac{d_4}{8 d_3} \nu(\tcc'_m) .$$
Consider an arbitrary $(2,1, j) \in \ff'_m$. Given $u \in \hddmtwo_j$, $1- \mu_0 \geq 3/4$ and $|h(\vm_{1,j} (u))|\leq \frac{3}{4}H(\vm_{1,j} (u))$
imply $|\tpsi_m(u)| \leq \gao_m(u)$. So if $(2,m,1, j) \in J$ with $(2,1,j) \in \ff'_m$, then as in Case 1 we derive
$\left| (\lab^N h)(u)\right|  \leq   (\nn_J H) (u)$.
Thus, in this case we can  take $\fftwo_m = \{(2,1,j) : (2,1,j) \in \ff'_m\}$ and then (\ref{eq:Lab}) and (\ref{eq:Fmtwo-cond}) will be satisfied.

\ms

{\bf Sub-case 2.2.2.}  Assume that
\be
\di \sum_{(2,1, j)\in \ff'_m} \nu(\hddmtwo_j)  < \frac{d_4}{8 d_3} \nu(\tcc'_m) .
\ee
Next, consider the case when
\be
\di \sum_{(2,2, j)\in \ff'_m} \nu(\hddmtwo_j)  \geq \frac{d_4}{8d_3} \nu(\tcc'_m) .
\ee
As in Sub-case 2.2.1 we observe that if we include in $J$ those $(2,m,2, j) \in J$ with $(2,2,j) \in \ff'_m$, then  we derive
$\left| (\lab^N h)(u)\right|  \leq   (\nn_J H) (u)$.
Thus, we can take $\fftwo_m = \{(2,2,j) : (2,2,j) \in \ff'_m\}$ and then (\ref{eq:Lab}) and (\ref{eq:Fmtwo-cond}) will be satisfied.

\ms

{\bf Sub-case 2.2.3.}  Assume next that
$$\di \sum_{(2,2, j)\in \ff'_m} \nu(\hddmtwo_j)  < \frac{d_4}{8d_3} \nu(\tcc'_m) .$$
This and (7.11) yield
$$\sum_{(2,1,j)\in \ff'_m} \nu(\hddmtwo_j) + \sum_{(2,2,j)\in \ff'_m} \nu(\hddmtwo_j) < \frac{d_4}{4 d_3} \nu(\tcc'_m)
< \frac{d_4}{4} \nu(\tcc'_m) .$$
Denote by {\it $J''_m$ the set of those $j \in J_m$ such that both $(2,1,j) \in \ff''_m$ and $(2,2,j) \in \ff''_m$}.
Then the above and (\ref{eq:d4-cond}) with $t = 2$ imply
\be
\sum_{j\in J''_m } \nu(\hddmtwo_j) > \frac{3d_4}{4} \nu(\tcc'_m) .
\ee
Finally, denote by {\it $J'''_m$ the set of those $j \in J_m$ such that both $(1,1,j) \in \ff''_m$ and $(1,2,j) \in \ff''_m$
and both $(2,1,j) \in \ff''_m$ and $(2,2,j) \in \ff''_m$,} i.e. $J'''_m = J'_m \cap J''_m$.  It now follows from 
(7.10), (7.13) and (\ref{eq:d4-cond}) with $t = 2$ that
\be
\sum_{j\in J'''_m } \nu(\hddmtwo_j) > \frac{d_4}{2} \nu(\tcc'_m) .
\ee

{\bf Fix now an arbitrary $j \in J'''_m$.} 

Then $(1,1,j) \in \ff''_m$, $(1,2,j) \in \ff''_m$, $(2,1,j) \in \ff''_m$
and $(2,2,j) \in \ff''_m$, so for any $t = 1,2$ the second alternative in Lemma 7.2(c) holds for $(t,1,j)$ and $(t,2,j)$, that is
\be
|h(\vm_{i,j}(u))|\geq \frac{1}{4}\, H(\vm_{i,j}(u)) > 0 
\ee 
for all $u \in \hddmt_j$ for both $i = 1$ and $i = 2$.

We will prove the following

\ms

{\bf Claim:} Either (\ref{eq:Lab}) holds for all $u \in \hddmo_{j}$ or (\ref{eq:Lab}) holds for all $u \in \hddmtwo_{j}$.

\ms

{\it Proof of Claim.} 
Let $u,u' \in \hddmt_{j}$ for some $t = 1,2$.
Using the assumption  $(h,H) \in \kk_b$, and in particular property 
$(T2)$ with $p = 0$, $v = \vm_{i,j}(u)$ and $v' = \vm_{i,j}(u')$, and assuming e.g.
$$\min\{ |h(\vm_{i,j}(u))| , |h(\vm_{i,j}(u'))|\}  = |h(\vm_{i,j}(u'))| ,$$
it follows from (\ref{eq:h-cond}), (7.2) and (7.15) that
\begin{eqnarray*}
\frac{|h(\vm_{i,j}(u)) - h(\vm_{i,j}(u'))|}{\min\{ |h(\vm_{i,j}(u))| , |h(\vm_{i,j}(u'))| \}}
& \leq &   \frac{E|b|\,\theta^N H(\vm_{i,j}(u'))}{|h(\vm_{i,j}(u'))| } \diamte (\tddmt_{j})\\
& \leq & 4 E |b|\,   \theta^{N} \diamte (\tddmt_{j}) < \frac{\ep_3}{8} .
\end{eqnarray*}
So, the difference between the arguments of  the complex numbers 
$h(\vm_{i,j}(u))$ and $h(\vm_{i,j}(u'))$ (regarded as vectors in $\R^2$)  is   $< \frac{\ep_3}{8} < \frac{\pi}{8}$.  
In particular,  for any $i = 1,2$ we can choose a real continuous function
$\thetam_{i,j}(u)$, $u \in  \hddmt_{j}$, with values in $[0, \ep_3/8]$  and $\lambdam_{i,j} \in [0,2\pi)$ such that
\be
\di h(\vm_{i,j}(u)) = e^{\i(\lambdam_{i,j} + \thetam_{i,j}(u))}|h(\vm_{i,j}(u))| \quad , \quad u\in \hddmt_{j}  .
\ee
The above, 
yields
\be
|\thetam_{i,j}(u) - \thetam_{i,j}(u')|\leq  \frac{\ep_3}{4}  < \frac{\pi}{4} 
\ee
for all $u,u' \in \hddmt_{j}$.


Choose an arbitrary point $\um_{i,j} \in \hddmt_{j}$, and set
$\tlambdam_{i,j} = |b| \psi_m(\um_{i,j}) + 2 k\pi ,$
where we choose $k = k(i,m,j)\in \Z$ so that
\be
|\lambdam_{2,j} - \lambdam_{1,j} + \tlambdam_{i,j} | \leq \pi .
\ee

By (7.16), the difference between the arguments of the complex numbers
$e^{\i \,b\,\tau_N(\vm_{1,j}(u))} h(\vm_{1,j}(u))$ and $e^{\i \,b\, \tau_N(\vm_{2,j}(u))} h(\vm_{2,j}(u))$
is given by the function
\begin{eqnarray*}
\Omegam_{j}(u) 
& = & [b\,\tau_N(\vm_{2,j}(u)) + \thetam_{2,j}(u) + \lambdam_{2,j}] -  
[b\, \tau_N(\vm_{1,j}(u)) + \thetam_{1,j}(u) + \lambdam_{1,j}]\\
& = & (\lambdam_{2,j}-\lambdam_{1,j}) - b \psi_{m}(u) + (\thetam_{2,j}(u) - \thetam_{1,j}(u)) .
\end{eqnarray*}
It follows from Lemma 6.5 (or rather its proof) and (7.2), using (7.3) as well, that for all $u, u'  \in \hddmt_{j}$ we have
\begin{eqnarray*}
        |b|\, |\psi_{m}(u) - \psi_m(u')| 
& =    & |b|\, \left| [\tau_{N}(\vm_{1,j}(u)) - \tau_{N}(\vm_{2,j}(u))] 
- [\tau_{N}(\vm_{1,j}(u')) - \tau_{N}(\vm_{2,j}(u'))] \right|\nonumber\\
& \leq & |b|\, \frac{2 C_4}{d_0^2} \, \theta^{N} \, \diamte (\tddmt_{j}) \leq \frac{2 C_4}{d_0^2} \,  \frac{\ep_3}{32 E} < \frac{\ep_3}{16} .
\end{eqnarray*}
The latter, (7.17), (7.18) and the choice of $\ep_3$ in Sect. 6.3 imply that for any $u \in \hddmt_{j}$ we have
\begin{eqnarray}
|\Omegam_j(u)| 
& \leq & |\lambdam_{2,j}-\lambdam_{1,j} +\tlambdam_{1,j} | 
+  |b|\, |\psi_m (u) -\psi_m (\um_{i,j})| + |\thetam_{2,j}(u) - \thetam_{1,j}(u)| \nonumber\\
& \leq & \pi + \frac{\ep_3}{16} + \frac{\ep_3}{4} < \frac{3\pi}{2} .
\end{eqnarray}
We will show now that either for $\hddmo_{j}$ or for $\hddmtwo_{j}$ we can bound
$|\Omegam_j(u)|$ from below by a positive constant.

It follows from the properties of the cylinders $\hddmo_{j}$ and $\hddmtwo_{j}$ in (\ref{eq:IN-cond}) that for
$u \in \hddmo_{j}$  and $u'\in \hddmtwo_{j}$  we have 
$$|\psi_{m}(u) - \psi_{m}(u')| \geq  \delta_1  e^{- s_m \hep_{13}}\, \diam(\tcc_m) .$$   
This, (\ref{eq:C-cond}) and $\hep_{15} \geq \hep_{13}$ imply
$$|b| |\psi_{m}(u) - \psi_{m}(u')| \geq |b| \delta_1  e^{- s_m \hep_{13}}\, \frac{\ep_2e^{s_m \hep_{15}}}{ |b|}
\geq \delta_1 \, \ep_2$$
for all $u \in \hddmo_{j}$ and $u'\in \hddmtwo_{j}$.
Then for such $u$ and $u'$ we have
\begin{eqnarray*}
|\Omegam_{j} (u)- \Omegam_{j}(u')|
& \geq & |b|\, |\psi_{m}(u) - \psi_{m}(u')| - |\thetam_{1,j}(u) - \thetam_{1,j}(u')| - |\thetam_{2,j}(u) - \thetam_{2,j}(u')|\\
& \geq & \delta_1 \, \epsilon_{2} - \ep_3 >  2\ep_3 ,
\end{eqnarray*}
by the choice of $\ep_3$.

Thus,  $|\Omegam_j (u)- \Omegam_j (u')|\geq 2\epsilon_{3}$ for all  $u\in \hddmo_{j}$ and all $u'\in \hddmtwo_{j}$. Hence we either 
have $|\Omegam_j (u)| \geq \ep_3$ for all $u\in \hddmo_{j}$ or $|\Omegam_j(u')| \geq \ep_3$ for all $u'\in \hddmtwo_{j}$. Indeed, if
$|\Omegam_j (u')| < \ep_3$ for some $u'\in \hddmo_{j}$, then for every $u\in \hddmtwo_{j}$ we get
$$|\Omegam_j (u)| = |(\Omegam_j (u)- \Omegam_j (u')) + \Omegam_j (u')| 
\geq |\Omegam_j (u)- \Omegam_j (u')| - |\Omegam_j (u')| > \ep_3 . $$
Similarly, if $|\Omegam_j (u)| < \ep_3$ for some $u\in \hddmo_{j}$, then $|\Omegam_j(u')| \geq \ep_3$ for every $u'\in \hddmtwo_{j}$.

\bs

Hence we either have

\ms

{\bf A:} $|\Omegam_j (u)| \geq \ep_3$ for all $u\in \hddmo_{j}$,

or 

{\bf B:} $|\Omegam_j (u)| \geq \ep_3$ for all $u\in \hddmtwo_{j}$.

\ms

Assume for example that we have {\bf A:} $|\Omegam_j(u)| \geq \ep_3$ for all $u\in \hddmo_{j}$. 
It follows from this and (7.19) that $\ep_3 \leq |\Omegam_j (u)| <  \frac{3\pi}{2}$ for all $u \in \hddmo_{j}$.
Hence, we see that for $u\in \hddmo_{j}$ the difference $\Omegam_j(u)$  between the  arguments of the complex numbers
$e^{\i \,b\,\tau_N(\vm_{1,j}(u))} h(\vm_{1,j}(u))$ and $e^{\i \,b\, \tau_N(\vm_{2,j}(u))} h(\vm_{2,j}(u))$,
defined as a number in the interval $[0, 2\pi)$, satisfies
$\ep_3 \leq \Omegam_j (u) < 3\pi/2$ for all $u\in \hddmo_{j}$.

\def\halpha{\hat{\alpha}}
\def\hgamma{\hat{\gamma}}

It follows from Lemma 7.2(b) that either $H(\vm_{1,j}(u)) \geq H(\vm_{2,j}(u))/4$ for all
$u \in \hddmo_{j}$ or $H(\vm_{2,j}(u)) \geq H(\vm_{1,j}(u))/4$ for all $u \in \hddmo_{j}$. 
Assume e.g. that 
\be
H(\vm_{1,j}(u))/4 \leq H(\vm_{2,j}(u)) \quad , \quad u \in \hddmo_{j} .
\ee
As in \cite{D1} (see also \cite{St4}) we will show that 
$|\tpsi_m (u)| \leq \gamma^{(1)}_m (u)$ for all $u \in \hddmo_{j}$.
Given such $u$, consider the points 
$$z_1 = e^{(\fa_{N}+\i b\tau_{N})(\vm_{1,j}(u))} h(\vm_{1,j}(u)) \:\:\: ,  \:\:\: z_2 =  e^{(\fa_{N}+\i b\tau_{N})(\vm_{2,j}(u))} h(\vm_{2,j}(u))$$
in the complex plane $\C$, and let $\varphi$ be the smaller angle between the arguments of $z_1$ and $z_2$.
It then follows from the above estimate for $\Omegam_j(u)$ that $\epsilon_3 \leq \varphi \leq 3\pi/2$. 
Moreover,  (7.3), $|h| \leq H$, (7.15)  and our assumption (7.20) imply
\begin{eqnarray*}
\frac{|z_1|}{|z_2|} 
=     e^{\fa_N(\vm_{1,j}(u)) - \fa_N(\vm_{2,j}(u))}   \frac{|h(\vm_{1,j}(u))|}{|h(\vm_{2,j}(u))|}
  \leq    e^{\ep_3} \, \frac{H(\vm_{1,j}(u))}{H(\vm_{2,j}(u))/4} \leq 16 e^{\ep_3} < 19 ,
 \end{eqnarray*}
by the choice of $\ep_3$. This yields
\begin{equation}
|z_1 + z_2| \leq (1- t)  |z_1| + |z_2| ,
\end{equation}
where we can take e.g. 
$\di t = \frac{1 - \cos (\epsilon_3)}{20} .   $
Indeed, we have
$$|z_1 + z_2|^2 = |z_1|^2 + |z_2|^2 + 2 \langle z_1, z_2\rangle \leq  |z_1|^2 + |z_2|^2 + 2 |z_1|\,  |z_2| (1- s) ,$$
where $s = 1- \cos \ep_3$. Thus, (7.21) will hold if
$$|z_1|^2 + |z_2|^2 + 2 |z_1|\,  |z_2| (1- s)  \leq  (1-t)^2 |z_1|^2 + |z_2|^2 + 2 (1-t) |z_1|\,  |z_2| ,$$
that is if 
$$(1 - (1-t)^2) |z_1| + 2   |z_2| (1- s)  \leq  2 (1-t)   |z_2| ,$$
which equivalent to
$\di |z_1| \leq 2 \frac{s - t}{t(2-t)} \, |z_2| .$
Since $t = s/20$, it follows that
$$\frac{|z_1|}{|z_2|} < 19 <  2 \frac{s - t}{t(2-t)} = \frac{38}{2- s/20} ,$$ 
so the above inequality holds.
This proves (7.21) with the given choice of $t$.

Since $\mu_0 \leq t$ by (\ref{eq:muo-cond}), it now follows from (7.21) that $|z_1+ z_2| \leq (1-\mu_0) |z_1| + |z_2|$.
Therefore  $|\tpsi_m(u)| \leq \gamma^{(1)}_m(u)$   for all $u \in \hddmo_{j}$.  
Now the argument from Case 1 proves that (\ref{eq:Lab}) holds for all $u \in \hddmo_{j}$.

Thus, in the case {\bf A}, (\ref{eq:Lab}) holds for all $u \in \hddmo_{j}$. 
In a similar way we prove that in the case {\bf B}, (\ref{eq:Lab}) holds for all $u \in \hddmtwo_{j}$.
This proves the Claim.

\ms

We will now define the set $\ff_m$ in the Sub-case 2.2.3. Consider the set $\tJ_m$ of all $j \in J''''_m$ so that
(\ref{eq:Lab}) holds for all $u \in \hddmtwo_{j}$ (i.e. case {\bf B}). If
$$\sum_{j\in \tJ_m } \nu(\hddmtwo_j) \geq \frac{d_4}{4} \nu(\tcc'_m) ,$$
then we just set $\fftwo_m = \{ (2,m,j) : j \in \tJ_m\}$, and then (7.6) will be satisfied. We then set $J = \fftwo_m$.

Assume that
$$\sum_{j\in \tJ_m } \nu(\hddmtwo_j) < \frac{d_4}{4} \nu(\tcc'_m) .$$
This and (7.14) imply
$$\sum_{j\in J'''_m \setminus \tJ_m } \nu(\hddmtwo_j) \geq \frac{d_4}{4} \nu(\tcc'_m) .$$
For $j \in J'''_m\setminus \tJ_m$ we have the case {\bf A}, that is (\ref{eq:Lab}) holds for all $u \in \hddmo_{j}$.
The above and  (\ref{eq:d3-cond}) yield
$$\sum_{j\in J'''_m \setminus \tJ_m } \nu(\hddmo_j) \geq \frac{d_4}{4 d_3} \nu(\tcc'_m) .$$
Now we set $\ffo_m = \{(1,m,j) : j \in  J'''_m\setminus \tJ_m\}$ and then (7.5) will be satisfied.
So, we just define $J = \ffo_m$ in this case.

This completes the construction of the set $J$ in all possible cases.
Clearly $J   \in \jj(b)$ and
(\ref{eq:Lab}) holds for all $u \in V_b$. As we mentioned in the beginning of the proof, (\ref{eq:Lab}) 
always holds for $u \in \hU\setminus V_b$.
\endofproof

\newpage


\def\hddmpt{\widehat{\dd}^{(m_p,t)}}
\def\hddmpo{\widehat{\dd}^{(m_p,1)}}
\def\hddmptwo{\widehat{\dd}^{(m_p,2)}}

\section{$L^1$ contraction estimates}
\setcounter{equation}{0}

Here we obtain $L^1$-contraction estimates for large powers of the contraction operators $\nn_J$.
We continue to use the notation from Sections 5 and 6. Here 
{\bf we assume that $f \in \ff_{\theta_1}(\hU)$ and  $N \geq N_0$ is a fixed integer} as in (\ref{eq:No-cond}).
Recall also the constants $b_0 > 1$, $D_3 > 1$ from (\ref{eq:nb-cond}), the constant $c =  |q_0|\ep_0$ 
that appears in (\ref{eq:Ub-cond}) in Sect. 6.2 and $C_{9} > 0$, $a_0 = a_0(N) > 0$, $\rho_3 = \rho_3(N) \in (0,1)$ and
$\rho_4 = \rho_4(N) \in (0,1)$ from Sect. 6.3.


Given $b$ with $|b| \geq b_0$ and a representative set $J \in \jj(b)$ set
$$W_J = \cup_{(t,m,i, j) \in J} \hddmt_j \subset V_{b} .$$
Recall from (\ref{eq:Vb-cond}) that $V_b$ is a union of cylinders $\cc'_m = \cc'_m(b)$ and 
for each $m = 1, \ldots, m_0= m_0(b)$, $\hddmt_{j} = \hddmt_{j} (b)$ are sub-cylinders of 
$\cc'_m$ for $t = 1,2$, $j=1,\ldots,j_m $ for some $j_m = j_m(b)$.

The following lemma is the analogue of Lemma 12 in \cite{D1}. It is similar to Lemma 7.1 in \cite{St4}.

\bs

\noindent
{\bf Lemma 8.1.}  
(a) {\it For any $H\in \kk_{E|b|}$, any $J\in \jj(b)$ and any integer $r$ with $0 \leq r \leq N-1$ we have}
\be
\int_{\cup_{r=0}^{N-1}\sigma^{-r}(V_{b})} H^2 \, d\nu 
\leq C_{9}\, \int_{\cup_{r=0}^{N-1} \sigma^{-r}(W_J)} H^2 \, d\nu .
\ee

\ms

(b) {\it For any $H\in \kk_{E|b|}$,  any $a\in \R$ with $|a| \leq a_0$ and any $J\in \jj(b)$ we have}
\be
\int_{\cup_{r=0}^{N-1} \sigma^{-r}(V_{b})} (\nn_J H)^2 \, d\nu 
\leq \rho_3\, \int_{\cup_{r=0}^{N-1} \sigma^{-r}(V_{b})} L^N_{\f0} (H^2)\, d\nu .
\ee

\ms

\noindent
{\it Proofs.} 
(a) 
Let $H \in \kk_{E|b|}$. Then $0 < H \leq 1$. Consider an arbitrary $J \in \Pi(b)$. 
Then either all elements of $J$ have the form $(1,m,i,j)$ or all of them 
have the form $(2,m,i,j)$. Assume e.g.  that we have the first case, and that (\ref{eq:1-cond}) holds.

Notice that the set $V = \cup_{r=0}^{N-1} \sigma^{-r}(V_{b})$ is a finite union of cylinders in $\hU$.
Any two cylinders in the latter are either disjoint or one of them contains the other.
Since $V_{b}$ is a disjoint union of cylinders $\cc'_m$ ($1 \leq m \leq m_0 $), $V$ is a union of cylinders in
$\sigma^{-r}(\cc'_m)$, $0 \leq r \leq N-1$, $m = 1, \ldots, m_0$. Taking maximal cylinders of this kind, we get a {\bf disjoint} family 
$$\{Y_p \, :\,  1 \leq p \leq k \}$$
of cylinders in $\sigma^{-r}(V_b)$ so that for every $p = 1, \ldots,k$ we have $\sigma^{r_p}(Y_p) = \cc'_{m_p}$
for some $r_p = 0,1, \ldots,N-1$ and some $m_p = 1, \ldots, m_0$. Naturally, $r_p$ may take the same
value $r$ for many different $p$, and similarly $m_p$ might take the same value $m$ for various different $p$.
Then we have 
$$V = \cup_{p=1}^k Y_p \quad , \quad V_b = \cup_{p=1}^k \sigma^{r_p}(Y_p).$$ 
For every $p = 1, \ldots,k$ and every $(1,m_p, i,j) \in J$
let $Z_{(1,m_p,j)}$ be the cylinder in $Y_p$ such that $\sigma^{r_p}(Z_{(1,m_p,j)}) = \hddmpo_{j}$.
Then
$$W_J = \cup_{p=1}^k \cup_{(1,m_p,i,j) \in J} \sigma^{r_p}(Z_{(1,m_p,j)}) .$$

Consider now a fixed $p$ and corresponding sets $Y_p$ and $Z_{(1,m_p,j)} \subset Y_p$.
By (\ref{eq:1-cond}),
\be
\di \sum_{(1,m_p,i, j)\in J} \nu(\hddmpo_{j}) \geq \frac{d_4}{8 d_3} \nu(\cc'_{m_p}) .
\ee
Given $u,u'\in \cc'_{m_p}$, it follows from $H \in \kk_{E|b|}$ and (\ref{eq:ccmo-cond}) that
$$\di \frac{|H(u) - H(u')|}{H(u')} \leq E \,|b|\,  \dteo(u,u') \leq E C_6 ,$$
therefore
$H(u)/H(u') \leq 1+ EC_6 \leq 2E C_6 .$
Thus, for
$L_1 = \max_{\cc'_{m_p}} H$ and  $ L_2 = \min_{\cc'_{m_p}}H$,  we have $1 \leq L_1/L_2 \leq 2E C_6$. 
Using (8.3) we derive
$$\int_{\cc'_{m_p}} H^2 \, d\nu \leq L_1^2 \nu(\cc'_{m_p}) \leq \frac{8 d_3 L_1^2}{d_4}\,  \sum_{(1,m_p,i,j)\in J} \nu(\hddmpo_{j})
\leq  \frac{8 d_3 L_1^2 }{d_4 L_2^2}\, \sum_{(1,m_p,i,j) \in J} \int_{\hddmpo_{j}} H^2\, d\nu . $$
On the other hand, using $\sigma^{r_p}(Y_p) = \cc'_{m_p}$ and the $\sigma$-invariance of $\nu$, we get
\begin{eqnarray}
\int_{Y_p} H^2\, d\nu 
& = & \int_{U} H^2(y)\, \chi_{Y_p}(y)\, d\nu(y) = \int_{U} H^2(\sigma^{r_p} y)\, \chi_{Y_p}(\sigma^{r_p}y)\, d\nu(y)\nonumber\\
& = & \int_{U} H^2(x)\, \chi_{\cc'_{m_p}}(x)\, d\nu(x) = \int_{\cc'_{m_p}} H^2 \, d\nu .
\end{eqnarray}
Similarly, 
$$ \int_{Z_{(1,m_p, j)}} H^2\, d\nu =  \int_{\hddmpo_{j}} H^2\, d\nu .$$
It now follows from (8.4) and $C_{9} \geq   \frac{32 d_3 E^2 C_6^2}{d_4 } \geq \frac{8 d_3 L_1^2 }{d_4 L_2^2}$ (see Sect. 6.3) that
\be\label{eq:Yp-cond}
\int_{Y_p} H^2 \, d\nu 
\leq  \frac{8 d_3 L_1^2 }{d_4 L_2^2}\, \sum_{(1,m_p,i, j) \in J} \int_{Z_{(1,m_p,j)}} H^2\, d\nu 
\leq  C_{9}\, \sum_{(1,m_p,i, j) \in J} \int_{Z_{(1,m_p,j)}} H^2\, d\nu .
\ee
This, $\cup_{p=1}^k Y_p = V = \cup_{r=0}^{N-1} \sigma^{-r}(V_b)$ and 
$\cup_{p=1}^k \cup_{(1,m_p,i, j) \in J} Z_{(1,m_p,j)} = \cup_{r=0}^{N-1} \sigma^{-r}(W_J)$ imply
$$\int_{\cup_{r=0}^{N-1}\sigma^{-r}(V_{b})} H^2 \, d\nu = \sum_{p=1}^k \int_{Y_p} H^2 \, d\nu \leq 
C_{9} \, \sum_{p=1}^k \sum_{(1,m_p,i, j) \in J} \int_{Z_{(1,m_p,j)}} H^2 \, d\nu = C_{9} \, \int_{\cup_{r=0}^{N-1} \sigma^{-r}(W_J)} H^2 \, d\nu ,$$
which proves (8.1).

\bs

(b) The proof of this part is very similar to the proof of Lemma 7.1(b) in \cite{St4}.
We provide some details since they will be used later. 

Let again $H \in \kk_{E|b|}$ and $J \in \jj(b)$, and let $0 \leq r \leq N-1$ be an integer. We continues to use the notation from
the proof of part (a) and the assumption about $J$.

We will use the representations 
$$\cup_{p=1}^k Y_p = V = \cup_{r=0}^{N-1} \sigma^{-r}(V_{b}) \quad, \quad 
W_J = \cup_{p=1}^k \cup_{(1,m_p,i,j) \in J} \sigma^{r_p}(Z_{(1,m_p,j)}) $$
from the proof of part (a). Recall that $\sigma^{r_p}(Y_p) = \cc'_{m_p}$ and 
$Z_{(1,m_p,j)}$ is the cylinder in $Y_p$ with $\sigma^{r_p}(Z_{(1,m_p,j)}) = \hddmpo_{j}$
for every $(1,m_p, i, j) \in J$.

 By Lemma 6.10, $\nn_J H \in \kk_{E|b|}$, while the Cauchy-Schwartz  inequality implies
\be
(\nn_J H)^2 = (\ma^N \omega_J H)^2 \leq (\ma^N \omega_J^2)\, (\ma^N H^2) \leq 
(\ma^N \omega_J)\, (\ma^N H^2) \leq \ma^N H^2 .
\ee
If $u \notin W_J$, then $\omega_{J}(u) = 1$. Let $u \in W_J$; then $u \in \hddmo_{j}$ 
for some (unique) $(1,m,i,j) \in J$. 
Assuming e.g. $i = 1$ so that $\sigma^N(v) = u$ for some $v = \vm_{1,j}(u)$ with $\omega_J(v) = 1-\mu_0$, we have
\begin{eqnarray*}
(\ma^N \omega_{J})(u)
& = &  \sum_{\sigma^N v = u, \;v\neq \vm_{1,j}(u)} e^{\fa_N(v)}   + e^{\fa_N(\vm_{1,j}(u))} \omega_{J} (\vm_{1,j}(u))  \\
& = & \sum_{\sigma^N v = u, \;v\neq \vm_{1,j}(u)} e^{\fa_N(v)}   + (1-\mu_0)  e^{\fa_N(\vm_{1,j}(u))} \\ 
&  =  & \sum_{\sigma^{N}v = u} e^{\fa_{N}(v)} - \mu_0 \,e^{\fa_{N}(\vm_{1,j}(u))}  
 \leq  (\ma^{N}\; 1)(u) - \mu_0 \, e^{-N T_0} = 1 - \mu_0 \, e^{-NT_0} .
\end{eqnarray*}
This holds for all $u \in W_J$, so by (8.6),
$$\di (\nn_J H)^2\circ \sigma^{r_p} \leq (1- \mu_0 e^{-NT_0})\, (\ma^N H^2)\circ \sigma^{r_p}$$
on $Z_{(1,m_p,j)}$. Therefore
\begin{eqnarray*}
\int_{Z(1, m_p,j)} (\nn_J H)^2\, d\nu
& = & \int_{U} (\nn_J H)^2(y) \chi_{Z(1,m_p,j)}(y)\, d\nu = \int_{U} (\nn_J H)^2(\sigma^{r_p} y) \chi_{Z(1,m_p,j)} (\sigma^{r_p} y)\, d\nu \\
& \leq & (1- \mu_0 e^{-NT_0}) \int_{U} (\ma^N H)^2(\sigma^{r_p} y) \chi_{Z(1,m_p,j)} (\sigma^{r_p} y)\, d\nu (y)\\
& = & (1- \mu_0 e^{-NT_0}) \int_{U} (\ma^N H)^2(x) \chi_{Z(1,m_p,j)} (x)\, d\nu (x)\\
& = &(1- \mu_0 e^{-NT_0}) \int_{Z(1,m_p,j)} (\ma^N H)^2 \, d\nu .
\end{eqnarray*}
Setting 
$$Z_p = \cup_{(1,m_p,j) \in J} Z_{(1,m_p,j)} \subset Y_p$$ 
for every $p = 1, \ldots, k$, the above implies
$$ \int_{Z_p} (\nn_J H)^2 \, d\nu \leq (1-\mu_0 e^{-NT_0})\,  \int_{Z_p} (\ma^N H)^2 \, d\nu ,$$
while (8.5) with $H$ replaced by $\ma^N H$ yields
$$\int_{Y_p} (\ma^N H)^2 \, d\nu  \leq  C_{9}\,  \int_{Z_p} (\ma^N H)^2\, d\nu .$$
As in (8.4) we have $\di \int_{Y_p} (\nn_J H)^2\, d\nu = \int_{\cc'_{m_p}} (\nn_J H)^2 \, d\nu$.
This, (8.6)   and the above yield
\begin{eqnarray*}
\int_{Y_p} (\nn_J H)^2\, d\nu
& =     & \int_{Y_p \setminus Z_p} (\nn_J H)^2 \, d\nu  +  \int_{Z_p} (\nn_J H)^2 \, d\nu\\
& \leq &  \int_{Y_p \setminus Z_p} (\ma^N H^2) \, d\nu  + (1-\mu_0 e^{-NT_0})\,  \int_{Z_p} (\ma^N H)^2 \, d\nu\\
& =     & \int_{Y_p} (\ma^N H)^2 \, d\nu  - \mu_0 e^{-NT_0}\,  \int_{Z_p} (\ma^N H)^2 \, d\nu\\
& \leq & \int_{Y_p} (\ma^N H)^2 \, d\nu  - \frac{\mu_0 e^{-NT_0}}{C_{9}} \,  \int_{Y_p} (\nn_J H)^2 \, d\nu .
\end{eqnarray*}
From the latter and 
\be
(\ma^N H)^2 \leq (\ma^N 1)^2 (\ma^NH^2) \leq \ma^N H^2 
= L^N_{\f0} (e^{\fa_N - \f0_N} H^2)\leq e^{a_0 NT_0} (L^N_{\f0} H^2) ,
\ee
we get
\begin{eqnarray*}
(1+\mu_0 \,  e^{-NT_0}/C_{9})\, \int_{Y_p} (\nn_J H)^2\, d\nu 
\leq  \int_{Y_p} (\ma^N H)^2 \, d\nu \leq e^{a_0NT_0}\, \int_{Y_p} L^N_{\f0} H^2 \, d\nu .
\end{eqnarray*}
Thus,
$$\di  \int_{Y_p} (\nn_J H)^2\, d\nu \leq \rho_3 \,  \int_{Y_p} L^N_{\f0} H^2 \, d\nu $$
for every $p = 1,\ldots, k$. This yields
\begin{eqnarray*}
\int_{\cup_{r=0}^{N-1} \sigma^{-r}(V_{b})} (\nn_J H)^2 \, d\nu 
 =  \sum_{p=1}^k \int_{Y_p} (\nn_J H)^2\, d\nu \leq \rho_3 \,  \sum_{p=1}^k \int_{Y_p} L^N_{\f0} H^2 \, d\nu 
 =  \rho_3\, \int_{\cup_{r=0}^{N-1} \sigma^{-r}(V_{b})} L^N_{\f0} (H^2)\, d\nu ,
\end{eqnarray*}
which proves (8.2). 
\endofproof

\bs

We can now prove that iterating sufficiently many contraction operators provides an $L^1$-contraction on $U$.

Let again $|b| \geq b_0$ and let $n_0(b) \geq 1$ be as in (\ref{eq:nb-cond}).
Define the function $\hath: \hU \longrightarrow [0,\infty)$ by
$$\hath = \rho_3\, \chi_{\cup_{i=0}^{N-1}\sigma^{-i}(V_{b})} + S_0\, \chi_{\hU\setminus \cup_{i=0}^{N-1} \sigma^{-i}(V_{b})} .$$
Then $\hath(x) = \rho_3$ whenever $\sigma^i(x) \in V_{b}$ for some $i = 0,1, \ldots, N-1$ and $\hath(x) = S_0$ otherwise.
Since $V_{b}$ and $\hU\setminus V_{b}$ are unions of finitely many cylinders, we have $\hath \in \ff_{\theta}(\hU)$. 

\ms

At some stage later on we will need the {\it Perron-Ruelle-Frobenius Theorem} (see e.g. \cite{PP}):
there exist global constants $C_{10} \geq 1$ and $\rho_0 \in (0,1)$ such that 
\be \label{eq:PRF}
\| L_{\f0}^n h - h_0 \, \int_{U} h \, d\nu \| \leq C_{10}\, \rho_0^n \, \|h\|_{\theta}
\ee
for all $h \in \ff_{\theta}(\hU)$ and all integers $n \geq 0$, where $h_0 > 0$ is the normalised eigenfunction of 
$L_{f-P_f \tau}$ in $\ff_{\theta}(\hU)$ (see Sect. 4.2). 
Fix such constants $C_{10} > 0$ and $\rho_0 \in (0,1)$. 

Here $\|h\|_\theta = \|h\|_0 + |h|_\theta$ as defined in Sect. 2.

The following  is the main result in this section. It is similar to Lemma 7.3  in \cite{St4} and its proof has
some similarity with that in \cite{St4}. 

\bs

\noindent
{\bf Theorem 8.2.}  {\it Let $|b| \geq b_0$.}

\ms

(a) {\it  For any sequence $J_1, J_2, \ldots, J_r \ldots $ of elements  of $\jj(b)$, setting 
$H^{(0)} = 1$ and\\ $H^{(r+1)} = \nn_{J_r} (H^{(r)})$ ($r \geq 0$) we have
\be
\di\int_{U} (H^{(M)})^2 \, d\nu \leq  2 \, \rho_4^M 
\ee
for all $M \geq n_0(b)$.}

\ms

(b) {\it Let $k \geq D_3$ be a constant. 
For all $|a| \leq a_0$ and $|b| \geq b_0$, and all $h \in \ff_{\theta}(U)$ for\\ $m =  (k/\beta_4) \, \log |b|$ we have}
\be
\int_{\hU} |\lab^{m N} h |\, d\nu \leq   \frac{2}{|b|^{2 k}} \, \|h\|_{\theta,b} .
\ee

(c) {\it Let $s \geq D_3$ be a constant. There exists a global constant  $C_{11} > 0$ such that for any integer $k \geq 2s$
and all $|a| \leq a_0$, $|b| \geq b_0$ and  $h \in \ff_{\theta}(U)$ for $m = (k/\beta_4) \, \log |b|$ we have}
\be
\|\lab^{2 m  N} h \|_0 \leq  \frac{C_{11}}{|b|^{k/2}} \, \|h\|_{\theta,b}
\leq  \frac{C_{11}}{|b|^{s}} \, \|h\|_{\theta,b} .
\ee

\ms

\noindent
{\it Proof of Theorem} 8.2. (a) Set  $\omega_r = \omega_{J_r}$ and $\nn_r = \nn_{J_r}$.
Since $H^{(0)} = 1 \in \kk_{E|b|}$, it follows from Lemma 6.10 that $H^{(r)} \in \kk_{E|b|}$ for all $r \geq 1$. 

Let $M \geq 1$ be an arbitrary integer. 
Using 
$L^N_{\f0}((\hath\circ \sigma^N) \, H) = \hath\, (L_{\f0}^N H)$, 
Lemma 8.1(b), (8.2) and (8.7),  we get
\begin{eqnarray*}
\int_{U} (H^{(M)})^2 \, d\nu
& =      & \int_{\cup_{r=0}^{N-1} \sigma^{-r}(V_b)} (H^{(M)})^2 \, d\nu + \int_{U\setminus \cup_{r=0}^{N-1} \sigma^{-r}(V_b)} (H^{(M)})^2 \, d\nu\\
& \leq  & \rho_3\, \int_{\cup_{r=0}^{N-1} \sigma^{-r}(V_b)} L^{N}_{\f0} (H^{(M-1)})^2\, d\nu 
+ e^{a_0 NT}\, \int_{U\setminus \cup_{r=0}^{N-1} \sigma^{-r}(V_b)} L^N_{\f0} (H^{(M-1)})^2\, d\nu\\
& =     & \int_{U} \hath\,  (L_{\f0}^N (H^{(M-1)})^2)\, d\nu 
= \int_{U}   L_{\f0}^N ( (\hath \circ \sigma^N)\, (H^{(M-1)})^2)\, d\nu\\
& =     &  \int_{U} (\hath\circ \sigma^N)\,  (H^{(M-1)})^2\, d\nu ,
\end{eqnarray*}
using here also one of the main properties of $L_{\f0}$, namely 
$\di \int_U L_{\f0} F\, d\nu = \int_U F\, d\nu $ for all $F \in C(U)$.
Continuing by induction and using $H^{(0)} = 1$, we get
\be\label{eq:HM-cond}
\int_{U} (H^{(M)})^2 \, d\nu \leq   \int_{U}  (\hath\circ \sigma^{MN})\, (\hath\circ \sigma^{(M-1)N}) 
 \ldots  (\hath\circ \sigma^{2N})\, (\hath\circ \sigma^N)\, d\nu .
\ee
 
By Lemma 6.4, the set $\tU_b = \tU_b(M,N)$ defined in Sect. 6.2 satisfies (\ref{eq:Ub-cond}) for any integers $M \geq n_0(b)$ and $N \geq 1$.
Using the rough estimate $\hath \leq S_0$ on $\tU_b$, it follows from (\ref{eq:HM-cond}), (\ref{eq:Ub-cond}) 
and (\ref{eq:rho-cond}) that for $M \geq n_0(b)$ we have
\be
\int_{\tU_b} (H^{(M)})^2 \, d\nu \leq S_0^M\, \nu(\tU_b) <  S_0^M e^{-cM N} \leq \rho_4^M .
\ee
On the other hand when $x \in U \setminus \tU_b$, the definition of $\tU_b$ implies
$\sigma^{jN}(x) \in \cup_{r=0}^{N-1} \sigma^{-r}(V_b)$ for at least $M \ep_0$ values of $j = 1,2,  \ldots, M$,
i.e. $\tSb_{M,N}(x) \geq M\ep_0$.
For such $j$ the definition of $\hath$ gives $\hath(\sigma^{jN}(x)) = \rho_3$. For all other $j$ we can
still use $\hath(\sigma^{jN}(x)) \leq S_0$. Thus, from (\ref{eq:HM-cond}) (or rather the analogous estimate we get for the integral over 
$U\setminus \tU_b$), using (\ref{eq:rho-cond}), we derive
\begin{eqnarray*}
\int_{U \setminus \tU_b} (H^{(M)})^2\, d\nu 
& \leq &  \int_{U\setminus \tU_b} \rho_3^{\tSb_{M,N}(x)}\, S_0^{M - \tSb_{M,N}(x)} \, d\nu (x)
 \leq  \rho_3^{M \ep_0}\, S_0^{M - M \ep_0} \nu(U\setminus \tU_b) \\
& \leq & (\rho_3^{\ep_0} S_0^{1- \ep_0})^M \leq \rho_4^M .
\end{eqnarray*}
This, (\ref{eq:HM-cond}) and (8.13) yield
\be
\int_{U} (H^{(M)})^2\, d\nu \leq \rho_4^M + \rho_4^M  = 2 \, \rho_4^M .
\ee
The latter holds for all integers $M \geq n_0(b)$.
This proves part (a).

\ms

(b) Let $h\in \ff_{\theta}(\hU)$ be such that 
$\| h\|_{\theta,b} \leq 1$.  Then  $|h(u)| \leq 1$ for all $u\in \hU$ and  $|h|_{\theta} \leq |b|$.  


Assume that the points $u,u'$, the cylinder $\hddmt_{j}$ in $U$, the integer $p \geq 0$ and the points $v,v'\in U$
satisfy (\ref{eq:u-cond})  for some $t = 1,2$, $j = 1, \ldots, j_m$. This implies
\begin{eqnarray*}
|h(v) - h(v')| 
& \leq &  |h|_\theta\, D_{\theta} (v,v') \leq  |b|\, D_{\theta} (v,v') =  \theta^{p+N}\, |b|\, D_{\theta} (u,u') 
\leq   \theta^{p+N}\, |b|\, \diam_{\theta}(\tddmt_{j}) \\
& \leq & E \, |b|\, \theta^{p+N} \, \diamte (\tddmt_{j}) .
\end{eqnarray*}
Thus, $(h, 1) \in \kk_{b}$ (see Sect. 7).

Set $(h^{(0)}, H^{(0)})= (h,1) \in \kk_{b}$ and $h^{(m)} = \lab^{m N}{h}$ for $m \geq 0$. 
Define the sequence of functions $\{ H^{(m)}\}$ recursively by
$H^{(0)} = 1$ and $H^{(m+1)} = \nn_{J_m} H^{(m)}$, where $J_m \in \jj(b)$ is chosen by induction as follows.
Since $(h^{(0)}, H^{(0)}) \in \kk_{b}$,  using Lemma 7.1  we find $J_0 \in \jj(b)$ such that
for $h^{(1)} = \lab^{N} h^{(0)}$ and $H^{(1)} = \nn_{J_0} H^{(0)}$ we have
$(h^{(1)} , H^{(1)}) \in \kk_{b}$. Continuing in this way we construct by induction
an infinite sequence of functions $\{ H^{(m)}\}$ with
$H^{(0)} = 1$, $H^{(m+1)} = \nn_{J_m} H^{(m)}$ for all $m \geq 0$, such that $(h^{(m)} , H^{(m)}) \in \kk_{b}$.
The latter implies $|h^{(m)}| \leq H^{(m)}$ for all $m \geq 0$.

It follows from part (a) that
$\di \int_U (H^{(m)})^2 d\nu \leq 2 \rho_4^m$ for all $m \geq n_0(b) = D_3\, \log |b|$.
Hence
$$\int_{U} |\lab^{m N} h|^2 \; d\nu  = \int_{U}  |h^{(m)}|^2\; d\nu \leq
\int_{U}  (H^{(m)})^2\; d\nu   \leq 2  \, \rho_4^m $$
for all $m \geq n_0(b)$.
The latter applies to any $h\in \ff_{\theta}(\hU)$  such that $\| h\|_{\theta,b} \leq 1$.
From this it follows that for any $h \in \ff_{\theta}(\hU)$ and any $m\geq n_0(b)$ we have 
$\di \int_{U} |\lab^{m N} h|^2 \; d\nu \leq  2 \, \|h\|^2_{\theta, b} \,  \rho_4^m$. Now H\"older's inequality yields
\be
\int_{U} |\lab^{m N} h| \; d\nu \leq  2 \, \rho_4^{m/2}  \,\|h\|_{\theta,b} 
\ee
for any $h \in \ff_{\theta}(\hU)$ and any $m\geq n_0(b)$.

Let $k \geq D_3$ be a fixed integer. Let
\be
m = \frac{k}{\beta_4}\, \log |b|  .
\ee
Recall from (\ref{eq:rho4}) that $\rho_4 \leq e^{-4\beta_4}$. Let $h \in \ff_\theta(\hU)$. Then (8.15) implies that
for any $m$ with (8.16) we have
$$\int_{U} |\lab^{m N} h| \; d\nu \leq  2  \, \rho_4^{m/2}  \,\|h\|_{\theta,b} 
= 2 e^{-4 \beta_4 \frac{k}{2 \beta_4}\, \log |b| }  \,\|h\|_{\theta,b} 
= 2 \frac{1}{|b|^{2k}}  \,\|h\|_{\theta,b} .$$
This proves (8.10).

\bs

(c)
We will now use a {\bf standard procedure} (see \cite{D1}) to derive an estimate of the form (8.11) from (8.10).
This partially repeats Sect. 7.4 in \cite{St4}.

The main tool to use here is the Perron-Ruelle-Frobenius Theorem -- see (8.8).

Fix an arbitrary constant $s \geq D_3$.  Recall from Sect. 6.3 that
$a_0 NT_0 \leq \beta_4/4$.
Let $m$ and $k$ satisfy (8.16), however this time we assume $k \geq 2s$.

Let $h \in \ff_{\theta}(\hU)$ be such that  $\|h\|_{\theta,b} \leq 1$. Then $\|h\|_0 \leq 1$ and 
$|h|_{\theta} \leq |b|$, so using Lemma 4.2 with $H = 1$ yields
\be
|\lab^{r} h|_{\theta} \leq A_0 [ |b| \theta^{r} + |b| ]\leq 2A_0 |b|
\ee
for any integer $r \geq 1$. It follows from H\"older's inequality that
\begin{eqnarray}
|\lab^{2mN} h|
& =    & |\lab^{mN} \left(|\lab^{mN} h |\right)| \leq \ma^{mN} |\lab^{mN} h| 
= L_{\f0}^{mN} \left(e^{\fa_{mN} - \f0_{mN}}\, |\lab^{mN} h|\right) \nonumber\\
& \leq & \left(L_{\f0}^{mN}\left(e^{\fa_{mN} -\f0_{mN}}\right)^{2}\right)^{1/2}\,
\left( L_{\f0}^{mN} |\lab^{mN} h|^{2} \right)^{1/2} .
\end{eqnarray}
For the first term in this product (\ref{eq:T0-cond}) and $a_0NT_0 \leq \beta_4/4$  imply
$$\di \left(L_{\f0}^{mN} \left(e^{\fa_{mN} - \f0_{mN}})\right)^{2}\right)^{1/2} \leq e^{a_0 NT_0 m} 
\leq e^{(k/4) \log |b|} = |b|^{k/4} .$$

For the second term in (8.18), using (\ref{eq:PRF}) replacing $h$ by  $|\lab^{mN} h|$ and $n$ by $mN$, we get
$$L_{\f0}^{mN} |\lab^{mN} h|^2 \leq L_{\f0}^{mN} |\lab^{mN} h| \leq \|h_0\|_0\, \int_U |\lab^{mN} h|\, d\nu 
+ C_{10}\, \rho_0^{mN} \, \|\lab^{mN} h\|_{\theta} .$$
By Lemma 4.2,  $\|\lab^{mN} h\|_{\theta} \leq  A_0 [ |b| \theta^{mN} + |b|] \leq 2 A_0 |b|$, so the above and  (8.10) imply
$$L_{\f0}^{mN} |\lab^{mN} h|^2 \leq  \frac{2}{|b|^{2k}} \|h_0\|_0  + 2 A_0 C_{10} |b| \, \rho_0^{mN} .$$
Assuming e.g. $\rho_0 \leq \rho_4 = e^{-4\beta_4}$ 
we get 
$$|b|\, \rho_0^{mN} \leq |b|\, e^{- 4 m N \beta_4} = |b|\, \frac{1}{|b|^{4N k}} < \frac{1}{|b|^{2k}} .$$
Thus,
$$L_{\f0}^{mN} |\lab^{mN} h|^2 \leq \|h_0\|_0\, \frac{2}{|b|^{2k}}  +
2 A_0 C_{10} |b|\, \rho_4^{mN} \leq \frac{C'_{11}}{|b|^{2k}} $$
for some constant $C'_{11} > 0$.
Therefeore
$$\left(L_{\f0}^{mN} |\lab^{mN} h |^2 \right)^{1/2} \leq \left(\frac{C'_{11}}{|b|^{2k}} \right)^{1/2}
\leq \frac{C_{11}}{|b|^{k}}  $$
for some constant $C_{11} > 0$.
Combining the estimates of the two terms in (8.18), it follows that
$$|\lab^{2mN} h | \leq   |b|^{k/4} \, \frac{C_{11}}{|b|^{k}} \leq \frac{C_{11}}{|b|^{k/2}} \leq \frac{C_{11}}{|b|^{s}} , $$
using $k \geq 2s$. This was proved assuming $\|h\|_{\theta, b} \leq 1$. For general
$h \in \ff_{\theta} (\hU)$ we now get  (8.11).
\endofproof

\bs

\noindent
{\it Proof of Theorem } 1.2. This is derived from Theorem 8.2(c) by a standard argument as in \cite{D1}.
We sketch the argument here for completeness. It is very similar to the one in Sect. 7 in \cite{St4}.

Recall that  $\theta = e^{-\beta_0}$ for some $\beta_0 > 0$.

Let again $s  \geq D_3$ and let $C_{11}$ be as in Lemma  8.2(c).  
We now assume 
$$m = \frac{k}{\beta_4 \beta_0} \, \log|b| $$ 
for some $k \geq 2s$. 
In particular we have $m \geq \frac{k}{\beta_4} \, \log|b|$ and therefore (8.11) holds for all $h \in \ff_{\theta}(\hU)$.

Consider an arbitrary  $h\in \ff_{\theta}(\hU)$ with
$\| h\|_{\theta,b} \leq 1$.  Then  $\|h\|_0 \leq 1$ and  $|h|_{\theta} \leq |b|$.

Since $|h|_{\theta} \leq |b|$, as before, using Lemma 4.2 with $H = 1$ yields (8.17)
for any integer $r \geq 0$ and any $h \in \ff_{\theta}(\hU)$.
Using again Lemma 4.2 and (8.11), we obtain
\begin{eqnarray*}
|\lab^{4 m N} h|_{\theta} 
 &   =  &  |\lab^{2 m N } ( \lab^{2 m N} h) |_{\theta} 
 \leq   A_0 \left[ 2 A_0 |b| \, \theta^{2 m N}  + |b| \, \| \lab^{2 m N} h\|_0 \right] \\
 & \leq & A_0\left[2 A_0  |b| \; e^{-2 m N \beta_0}  + |b| \; \frac{C_{11}}{|b|^{k/2}}\right]
\leq A_0 \, \left[2 A_0  |b| \; e^{-2 \frac{k N}{\beta_4} \log |b|}  + |b| \; \frac{C_{11}}{|b|^{k/2}}\right]\\
& \leq & A_0\left[2 A_0  |b| \; \frac{1}{|b|^{2 k N/\beta_4} } + |b| \; \frac{C_{11}}{|b|^{k/2}}\right]
\leq A_0\left[2 A_0  |b| \; \frac{1}{|b|^{2 k} } + |b| \; \frac{C_{11}}{|b|^{k/2}}\right]\leq \frac{C'_{12}}{|b|^{k/2 -1}}  ,
\end{eqnarray*}
for some constant $C'_{12} > 0$. This was done under the assumption $\| h\|_{\theta,b} \leq 1$.
Therefore
$$|\lab^{4 m N} h|_{\theta}  \leq \frac{C'_{12}}{|b|^{k/2 -1}} \, \|h\|_{\theta, b}$$
for all $h \in \ff_\theta(\hU)$.
Using (8.11) again, we get 
$\di \|\lab^{4 m N} h\|_{0} 
\leq  \frac{C'_{12}}{|b|^{k/2}} \,  \|h\|_{\theta,b} ,$
therefore there exists a global constant $C_{12} > 0$ such that
$\di \|\lab^{4 m N} h\|_{\theta,b} \leq \frac{C_{12}}{|b|^{k/2}}\, \|h\|_{\theta,b}$
for all $h \in \ff_\theta(\hU)$.
Assuming  $|b| \geq b_0$ and $b_0 > 1$ is sufficiently large, we get
$$\di \|\lab^{4 m N} h\|_{\theta,b} \leq \frac{C_{12}}{|b|^{k/2}} \, \|h\|_{\theta, b}  \leq \frac{1}{|b|^{k/4}}\, \|h\|_{\theta,b}$$
for all $h \in \ff_{\theta}(\hU)$.

Let $n \geq 4 m N $ be an arbitrary integer. Writing $n = r \, ( 4 m N)+ \ell$ for some
$\ell  = 0,1, \ldots, 4 m N -1$, and using the above $r$ times we get
$\|\lab^{r 4 m N} h\|_{\theta,b}  
\leq \frac{1}{|b|^{r k/4}}\, \|h\|_{\theta,b} .$
As before, using Lemma 4.2 with $H = 1$ and $B = |\lab^{ r 4 m N} h|_{\theta}$, implies
\begin{eqnarray*}
|\lab^{n} h|_{\theta} = |\lab^{\ell} ( \lab^{  r 4 m N} h) |_{\theta} 
 \leq   A_0 \left[ |\lab^{r 4 m  N} h|_{\theta} \, \theta^{\ell} + |b| \, \| \lab^{r 4 m N} h\|_0 \right] ,
\end{eqnarray*}
so  $\di \frac{1}{|b|} |\lab^{n} h|_{\theta}  \leq \frac{2A_0 }{|b|^{r k/4 }}\, \|h\|_{\theta,b}$.
This and  $\di \|\lab^{n} h\|_0 \leq \|\lab^{r 4 m N} h\|_{0} \leq \frac{1}{|b|^{r k/4}}\, \|h\|_{\theta,b}$ imply
$$\di \|\lab^{n} h\|_{\theta,b} 
 \leq  \frac{3A_0}{|b|^{r k/4}}\, \|h\|_{\theta,b} = 3A_0 \, e^{- (r k/4) \log |b|} \|h\|_{\theta,b} .$$
We have $r \geq (r+1)/2$ for all $r \geq 1$, so the above and 
$$(r+1) k \log |b| = \beta_4 \beta_0\, (r+1) \frac{k}{\beta_4 \beta_0}\, \log |b| 
= \frac{\beta_4 \beta_0}{4N}\, (r+1) 4 m N > \frac{\beta_4 \beta_0}{4N} \, n $$
yield 
\begin{eqnarray} 
\|\lab^{n} h\|_{\theta,b} 
& \leq & 3A_0 \, e^{- (r k/4) \log |b|} \|h\|_{\theta,b}
  \leq   3A_0  e^{-\frac{(r+1) k  \log |b|}{8}} \|h\|_{\theta,b}\nonumber\\
& \leq  & 3A_0  e^{-\frac{\beta_4 \beta_0 4 (r+1) m N}{32 N}} \|h\|_{\theta,b}
\leq 3A_0  e^{- \frac{\beta_4 \beta_0}{32 N}\, n } \|h\|_{\theta,b} 
 \leq  3A_0  \, \rho_6^n\, \|h\|_{\theta,b} ,
\end{eqnarray}
where $\rho_6 = \rho_6(N) = e^{-\frac{\beta_4 \beta_0}{32 N}} \in (0,1)$. 

Thus, (8.19) holds for all $h \in \ff_{\theta}(\hU)$ and all integers $n \geq 4 m N = 4 (k/\beta_4 \beta_0)\,   N \log |b|$. 
Finally, recall the eigenfunction $h_a \in \ff_{\theta}(\hU)$ for the operator $L_{f- (P_f+a)\tau}$
from Sect. 4.2. It is known that $\|h_a\|_{\theta} \leq \Con$ for bounded $a$, e.g. for $|a|\leq a_0$.
For $|a| \leq a_0$ and $a_0 > 0$ sufficiently small we have $\lambda_a \rho_6 \leq \rho$ for some
global constant $\rho \in (0,1)$. Now 
$$\di \lab^n(h/h_a) = \frac{1}{\lambda_a^n h_a}\, L^n_{f- (P+a+\i b) \tau} h$$
and the above estimate show that there exist constants $0 < \rho < 1$, $a_0 > 0$, 
and  $C > 0$ such that if $a,b\in \R$  satisfy $|a| \leq a_0$ and $|b| \geq b_0$,  then 
$\di \|L_{f -(P_f+a+ \i b)\tau}^n h \|_{\theta,b} \leq C \;\rho^n \; \| h\|_{\theta,b}$
for any integer $n \geq \frac{4 k N}{\beta_4 \beta_0}\,  \log |b|$ and any  $h\in \ff_{\theta} (\hU)$. So, we can just set 
$T = T(k,N) = \frac{4 k N}{\beta_4\beta_0}$, where we can take e.g. $k \geq 4 D_3$.

This completes the proof of Theorem 1.2.
\endofproof

\bs

\noindent
{\it Proof of Theorem} 1.1. This follows from the arguments in Sect. 4 and Appendix 1 in \cite{D1}.
Here we very briefly describe\footnote{Similarly to what we did in Sect. 8 in \cite{St1}, although here we provide 
a bit more details.}  the procedure from \cite{D1} that uses
Corollary 3.3 and leads to the proof of Theorem 1.2. It works in exactly the same way, 
since it is practically carried out on the symbol space $\Sigma_A$. 

Consider an arbitrary $F = F_M \in \ff_\alpha(M)$, and let $f\in \ff_\alpha (U)$ be defined by 
$f(x) = \int_0^{\tau(x)} F(\phi_s(x))\; ds$ (cf. Sect. 2). Given continuous functions
$A$ and $B$ on $M$, let
$$\rho_{A,B}(t) = \int_M A(x) B(\phi_t(x))\; d\nu_F(x)\; $$
be the corresponding {\it correlation function}, where as in Theorem 1.2, $\nu_F$ is the Gibbs
measure on $M$ determined by $F$. One also needs the {\it Laplace transform}
$\di \hat{\rho}_{A,B}(z) = \int_0^\infty \rho_{A,B}(t) e^{-zt}\; dt$
of $\rho_{A,B}$, where $z = a + \i b \in \C$. As before, we will denote by $\nu$ 
the Gibbs measure on $U$  determined by $f-P\tau$, where $U$ is identified with $\saa$ as before.
We will also use the notation $\rho_{A,B}(t)$ and $\hat{\rho}_{A,B}(z)$ for functions $A,B$ on
$\ff_\theta(\sa)$.

The first step comprises Propositions 3 and 4 in \cite{D1} and is based on previous arguments
of Pollicott \cite{Po} and Ruelle \cite{R3}. Proposition 3 deals with the representatives of $F_M$ and $f$
on the symbol spaces $\ff_{\theta^2}(\sa)$ and $\ff_{\theta}(\saa)$.

\medskip

{\bf Step 1.} {\it Assume that $\tau \in \ff_\theta(U)$ and $F \in \ff_{\theta^2}(\sa)$. 
There exist constants $\ep_4 > 0$, $\theta \in (0,1)$, $K > 0$ and $C_{13} > 0$
such that for every $z = a + \i b$ with $|a| \leq \ep_4$ and $1/|b| \leq \ep_4$ and every 
integer $n \geq 0$ there exist linear operators 
$Q_n(z), R_n(z) : \ff_\theta(\sa) \longrightarrow \ff_\theta (\saa)$ such that:}

\smallskip

(a) $\|Q_n(z) A\|_0 \leq C_{13} (1-\ep_4)^n |b|\; \|A\|_\theta $ {\it and
$\|R_n(z) A\|_0 \leq C (1-\ep_4)^n |b|\; \|A\|_\theta $ for all $A\in \ff_\theta(\sa)$};

\smallskip

(b) $|Q_n(z) A|_\theta \leq C_{13} \, K^n |b|^2\; \|A\|_\theta $ {\it and
$|R_n(z) A|_\theta \leq C \, K^n |b|^2\; \|A\|_\theta $ for all $A\in \ff_\theta(\sa)$};

\smallskip

(c) {\it The function
$$\frac{1}{\| A\|_0 \|B\|_0} \left[ \hat{\rho}_{A,B}(z) -  \sum_{i,j=0}^\infty \int_{\saa} \left( L^{i+j}_{f - (P+z)\tau}(1 - 
L_{f-(P+z)\tau})^{-1} [Q_i(z)A]\right)[R_j(z)B]\; d\nu\right]$$
is uniformly bounded}.

\medskip

From this Proposition 4 in \cite{D1} is derived:

\ms

{\bf Step 2.} {\it Assume that $F_M \in \ff_{\alpha}(M)$ with corresponding representative $F \in \ff_{\theta^2}(\sa)$.
With the notation and assumptions in Step 1, consider the operators $Q_n(z), R_n(z) $ as operators\\
$Q_n(z), R_n(z) : C^\alpha(M) \longrightarrow C^\alpha (U)$. Then:}

\smallskip

(a) $\|Q_n(z) A\|_0 \leq C_{13} (1-\ep_4)^n |b|\; \|A\|_\alpha $ {\it and
$\|R_n(z) A\|_0 \leq C (1-\ep_4)^n |b|\; \|A\|_\alpha $ for all $A\in C^\alpha(M)$};

\smallskip

(b) $|Q_n(z) A|_\alpha \leq C_{13} \, K^n |b|^2\; \|A\|_\alpha $ {\it and
$|R_n(z) A|_\alpha \leq C \, K^n |b|^2\; \|A\|_\alpha $ for all $A\in C^\alpha(M)$}.



\bs

The next step is given by Corollary 3 (b) in \cite{D1} (cf. also Lemma 18 there).

\medskip

{\bf Step 3.} {\it For every $\beta \in (0,\alpha)$ there exists a constant $C_{14} > 0$ such that
for $z = a + \i b$ with $|a| \leq \ep_4$ and $1/|b| \leq \ep_4$ we have
$\| (1 - L_{f-(P+z)\tau})^{-1} h\|_0  \leq C_{14} \ln |b| \; \| h\|_\beta$ for all $h \in \ff_\beta(U)$}. 

\medskip

Using this and applying  Corollary 4 in \cite{D1} one obtains

\medskip

{\bf Step 4.} {\it For all $\beta \in (0,\alpha)$ and $A,B\in \ff_\beta(M)$ there exists
$\ep_5 > 0$ such that $\hat{\rho}_{A,B}(z)$ is analytic in the region $|a| \leq \ep_5$,
$1/|b|\leq \ep_5$. Moreover, 
$|\hat{\rho}_{A,B}(z)| \leq C_{15} |b|^2 \ln |b|\; \|A\|_\beta\; \|B\|_\beta$ for
some constant} $C_{15} > 0$.

\medskip

As in \cite{D1}, this yields the exponential decay of the function
$\rho_{A,B}(t) - (\int_M A\; d\nu_F)(\int_M B\; d\nu_F)$.

\bs

\noindent
{\it Proof of Corollary 1.3.}  This is essentially the same as  the proof of Corollary 1.4 in Sect. 8 in \cite{St4}.
\endofproof

\ms

\section{Appendix I}
\setcounter{equation}{0}

\noindent
{\it Proof of Lemma} 6.6. We will use Lemma 6.5 and a standard argument.

Let $h \in \kk_0$ and let $|s| \leq M_1$. We will prove that $L^N_{f-s \tau} h \in \kk_0$.

Assume that $u,u' \in \tddmt_j$ for some $t = 1,2$, some $m = 1,\ldots, m_0$ and some $j = 1, \ldots, j_m$, and the integer $p \geq 0$ 
and the points $v,v'\in \hU$ satisfy (\ref{eq:u-cond}) for some $i = 1,2$, and 
$w,w'\in \hU$ are such that $\sigma^N w = v$, $\sigma^Nw' = v'$ and $\ell(w,w') \geq N$; then $w' = w'(w)$ is
uniquely determined by $w$. 

We have 
$D_{\theta} (\sigma^j(w),\sigma^j(w')) = \theta^{N-j}\, D_{\theta} (v,v')$
for all $j = 0,1, \ldots, N-1$, and
$D_{\theta}(v,v')= \theta^{N+p} \, D_{\theta} (u,u')$. Assuming $|f|_{\theta} \leq T_0$, we get
\begin{eqnarray*}
         |f_N(w) - f_N(w')| 
 \leq   \sum_{j=0}^{N-1} |f (\sigma^j(w)) - f (\sigma^j(w'))|
 \leq  \sum_{j=0}^{N-1} |f|_{\theta} \,\theta^{N-j}\, D_{\theta} (v,v')
 \leq   \frac{T_0  D_{\theta} (v,v')}{1-\theta} .
\end{eqnarray*}
It follows from this and (\ref{eq:u-cond}) that
\begin{eqnarray*}
|f_N (w) - f_N(w')| 
 \leq  \frac{T_0}{1-\theta} \theta^{p+N}\, D_{\theta}(u,u')
 \leq   \frac{T_0 }{1-\theta}  \theta^{p+N} \, \diam_{\theta}(\tddmt_j) .
\end{eqnarray*}
This and Lemma 6.5 imply
$$|(f-s\tau)_N(w) - (f- s\tau)_N(w')| 
\leq \left(\frac{T_0 C_{6} C_4}{1-\theta} + \frac{M_1 C_4}{d_0^2} \right)\, \theta^{p+ N} \, \diamte(\tddmt_j)
= C_{8}\, \theta^{p+ N} \, \diamte (\tddmt_j) $$
for all $s\in \R$ with $|s| \leq M_1$, where $C_{8} > 0$ is as in (\ref{eq:C8-cond}).

Since $h \in \kk_0$, using the assumptions about $w,w'$ we get
$$|h(w) - h(w')| \leq E_1\, \theta^{p+2N}\, h(w') \, \diamte (\tddmt_j).$$
Thus, given $s$ with $|s| \leq M_1$ we have:
\begin{eqnarray*}
&        & |(L^N_{f-s\tau}h)(v) - (L^N_{f- s\tau} h)(v')|
 =       \left| \sum_{\sigma^N w = v} e^{(f-s \tau)_N(w)}\, h(w)   -  \sum_{\sigma^N w = v} e^{(f - s \tau)_N(w'(w))}\, h(w'(w)) \right| \\
& \leq & \left| \sum_{\sigma^N w = v} e^{(f- s\tau)_N(w)}\, [h(w) -  h(w')]\right| 
        + \sum_{\sigma^N w = v}  \left|e^{(f- s\tau)_N(w)} -  e^{(f- s\tau)_N(w')}\right| \, h(w') = (I) + (II) .
        \end{eqnarray*} 
Estimating separately each of the sums $(I)$ and $(II)$, we get
\begin{eqnarray*}
(I)        
& \leq & \sum_{\sigma^N w = v} e^{(f- s\tau)_N(w) - (f- s\tau)_N(w')}  
e^{(f- s \tau)_N(w')}\, E_1 \theta^{p+2N} \, h(w') \, \diamte (\tddmt_j)  \\  
& \leq & E_1 \, e^{C_{8}}\, \theta^{p+2N} \diamte (\tddmt_j)\, (L^N_{f-s \tau} h)(v') ,
\end{eqnarray*}
and
\begin{eqnarray*}
(II)        
& \leq &  \sum_{\sigma^N w = v}  \left| e^{(f- s \tau)_N(w) - (f- s \tau)_N(w')} - 1\right|  \,  e^{(f- s \tau)_N(w')} h(w') \\
& \leq & e^{C_{8}} \,C_{8} \, \theta^{p+N} \, \diamte (\tddmt_j)\, (L^N_{f-s \tau} h)(v') .
\end{eqnarray*}
Therefore
$$\di (I) + (II) \leq E_1 \,\theta^{p+N}\,   \diamte (\tddmt_j) \, (L^N_{f-s \tau} h)(v') ,$$
using $e^{C_{8}} \,C_{8} \leq E_1/2$ and $\theta^N \, e^{C_{8}} < 1/2$ by $N \geq N_0$ and the choices of $N_0$ and $E$ in Sect. 6.3.

Hence $L^N_{f-s \tau} h \in \kk_0$.

A simple induction implies now that $L^{m N}_{f - s \tau} h \in \kk_0$ for all $h \in \kk_0$ and all integers $m \geq 1$.
\endofproof

\bs

\noindent
{\it Proof of Lemma} 6.10.  
Let $J \in \jj(b)$ and let $H \in \kk_{E|b|}$. Set $\nn = \nn_J$. We will show that $\nn H \in \kk_{E|b|}$. 
Since $0 < H \leq 1$, it follows from the definition of $\nn$ that $0 < \nn H \leq 1$.

Let $u,u' \in \hU$.
Given $v \in \hU$ with $\sigma^{N}(v) = u$, let $C[\ii] = C[i_0, \ldots,i_{N}]$ 
be the  cylinder of length $N$ containing $v$.  
Set $\hC[\ii] = C[\ii]\cap \hU$. Then $\sigma^{N}(\hC[\ii]) = \hU_{i_0}$. Moreover, 
$\sigma^N : \hC[\ii] \longrightarrow \hU_{i_0} $ is a homeomorphism, so
there  exists a unique $v' = v'(v)\in \hC[\ii]$ such that $\sigma^N(v') = u'$.  Then 
$$\dteo  (\sigma^j(v),\sigma^j(v'(v))) = \theta_1^{N-j}\, \dteo (u,u')$$
for all $j = 0,1, \ldots, N-1$. 
Also $\dteo (v,v'(v))= \theta_1^N \dteo  (u,u')$. Hence
\begin{eqnarray}
         |\fa_N(v) - \fa_N(v')| 
& \leq & \sum_{j=0}^{N-1} |\fa(\sigma^j(v)) - \fa (\sigma^j(v'))|
 \leq  \sum_{j=0}^{N-1} |\fa|_{\theta_1} \,\theta_1^{N-j}\, \dteo (u,u') \nonumber\\
& \leq &  \sum_{j=0}^{N-1} T_0 \,\theta_1^{N-j}\, \dteo (u,u') \leq \frac{T_0}{1-\theta_1}\, \dteo (u,u') .
\end{eqnarray}

The definition of $\omega$ implies that either  $\omega(v) = \omega(v')$ or at least one of these numbers is  $ < 1$. 
Using Lemma 6.8 we get 
$$|\omega(v) - \omega(v')| \leq \frac{\theta_1^{-\ell_0} }{\ep_4}\, |b|\, \dteo(u,u') .$$
Since $\sigma^{N}(v) = u$, $\sigma^{N}(v') = u'$ and $\ell(v,v') \geq N$, the assumption $H \in \kk_{E|b|}$ implies
$$|H(v) - H(v')| \leq E \, |b|\, H(v') \, \dteo(v,v') \leq E \, |b|\, H(v') \theta_1^{N} \dteo (u,u') .$$
It follows from (9.1) that
$$\left| e^{\fa_N(v) - \fa_N(v')} - 1\right| 
\leq e^{T_0/(1-\theta_1)} \frac{T_0}{1-\theta_1}\, \dteo (u,u') .$$
Also, $2 \omega_J (v') \geq 1 \geq \omega_J(v)$.

Using  the above  and the definition of $\nn = \nn_J$, setting  $v' = v'(v)$ and $\omega = \omega_J$ for brevity, we get
\begin{eqnarray*}
&        & \frac{|(\nn H)(u) - (\nn H)(u')|}{\nn H (u')} 
 =      \frac{\di \left| \sum_{\sigma^N v = u} e^{\fa_N(v)}\, \omega(v) H(v) -  
\sum_{\sigma^N v = u} e^{\fa_N(v'(v))}\, \omega(v'(v)) H(v'(v)) \right|}{\nn H (u')} \nonumber\\
& \leq & \frac{\di \sum_{\sigma^N v = u} e^{\fa_N(v) - \fa_N(v')}  e^{\fa_N(v')}\, |\omega(v) - \omega(v')| H(v')}{\nn H (u')}
            + \frac{\di \sum_{\sigma^N v = u} e^{\fa_N(v)}\, \omega(v) | H(v) -   H(v')|}{\nn H (u')} \nonumber \\
&       & + \frac{\di \sum_{\sigma^N v = u}  \left| e^{\fa_N(v) - \fa_N(v')} - 1\right| 
\,  e^{\fa_N(v')} \omega(v') H(v')}{\nn H (u')} .
\end{eqnarray*}
Estimating separately each of the three sums above, we get
\begin{eqnarray*}
&    &        \frac{|(\nn H)(u) - (\nn H)(u')|}{\nn H (u')} 
 \leq  e^{T_0/(1-\theta_1)} \frac{\theta_1^{-\ell_0} }{\ep_4}\, |b|\, 
  \frac{\di \sum_{\sigma^N v = u} e^{\fa_N(v')}\, H(v')\, \dteo(u,u')}{\nn H (u')} \\
&        &  + \frac{\di \sum_{\sigma^N v = u}  e^{\fa_N(v)- \fa_N(v')}\,e^{\fa_N(v')}\,  
2\omega(v')\, E \, |b| H(v') \theta_1^{N} \dteo(u,u')}{\nn H (u')} 
          + e^{T_0/(1-\theta_1)} \frac{T_0}{1-\theta_1}\, \dteo (u,u')\\
& \leq & \left(e^{T_0/(1-\theta_1)} \frac{\theta_1^{-\ell_0} }{\ep_4}\, |b|  + 2 e^{T_0/(1-\theta_1)} E |b|\, \theta_1^{N} 
        +  e^{T_0/(1-\theta_1)} \frac{T_0}{(1-\theta_1) }\,\right)\,  \dteo (u,u')
 \leq E \, |b|\, \dteo(u,u') , 
\end{eqnarray*}
since by (\ref{eq:E-cond}),  the choice of $N_0$ in (\ref{eq:No-cond}) and $N \geq N_0$, we have
$$e^{T_0/(1-\theta_1)}  \theta_1^{-\ell_0}/\ep_4  \leq \frac{E}{3}  \quad, \quad
2 e^{T_0/(1-\theta_1)} \, \theta_1^{N} \leq \frac{1}{3} \quad , \quad
e^{T_0/(1-\theta_1)} \frac{T_0}{(1-\theta_1) } \leq \frac{E}{3} .$$
Hence $\nn H \in \kk_{E|b|}$.
\endofproof

\bs

\noindent
{\it Proof of Lemma} 7.2.
(a)  Let $u, u'\in \hddmt_j$ for some $m \leq m_0$, $t = 1,2$ and $j = 1, \ldots, j_m$. For any $i = 1,2$
for  $v = \vm_{i,j}(u)\in  \Xmt_{i,j}$ and $v' = \vm_{i,j}(u')\in \Xmt_{i,j}$, 
we have $\ell (v,v') \geq N$  and $u= \sigma^{N}(v), u' = \sigma^N (v') \in \hddmt_j$.
This, $H \in \kk_{E|b|}$, $u, u' \in \cc'_m$ and (\ref{eq:ccmo-cond}) imply
$$|\ln H(v) - \ln H(v')| \leq \frac{|H(v) - H(v')|}{\min\{H(v), H(v')\}} \leq 
E \, |b|\, \dteo (v,v')  = E  \, |b| \, \theta_1^{N} \dteo(u,u') \leq E \, C_6\, \theta_1^N < \ln 2 ,$$
using (\ref{eq:No-cond}) as well.
Hence $|\ln H(v) - \ln H(v')|\leq \ln 2$, so $\frac{1}{2} \leq \frac{H(v)}{H(v')} \leq 2$.

\ms

(b) Consider again $\hddmt_j$ for some $m \leq m_0$, $t = 1,2$ and let $j = 1, \ldots, j_m$.
Assume e.g. $H(\vm_{1,j}(u')) \geq H(\vm_{2,j}(u'))$ for some $u' \in \hddmt_j$. Then, 
for every $u \in \hddmt_j$, using part (a) twice, we get 
$$H(\vm_{1,j}(u)) \geq H(\vm_{1,j}(u'))/2 \geq H(\vm_{2,j}(u'))/2 \geq H(\vm_{2,j}(u))/4 .$$
Similarly, if $H(\vm_{2,j}(u')) \geq H(\vm_{1,j}(u'))$ for some $u' \in \hddmt_j$, then $H(\vm_{2,j}(u)) \geq H(\vm_{1,j}(u))/4$
for all $u \in \hddmt_j$.

\ms

(c) Assume that $u,u' \in \hddmt_j$ for some $m \leq m_0$, some $t = 1,2$, and $j = 1, \ldots, j_m$.
Let $i  \in \{1,2\}$.
Consider the case when for some $v\in \Xmt_{i,j}$ we have $|h(v)|\geq \frac{3}{4}H(v)$.
Fix $v$ with this property and consider an arbitrary $v'\in \Xmt_{i,j}$. It follows from (7.1) and (7.2) that
$$ |h(v') -h(v)| \leq E |b|\,\theta^{N} H(v)\, \diamte (\tddmt_j)  \leq H(v)\, \frac{\ep_3}{32} .$$
Using $2H(v) \geq H(v')$ which follows from (a), one obtains
\begin{eqnarray*}
|h(v')| 
& \geq & |h(v)| - (\ep_3/32)  H(v) \geq (3/4  - \ep_3/32) H(v) \geq \frac{1}{4} \, H(v') .
\end{eqnarray*}
Thus, in this case the second alternative in (c) holds for all $v' \in \Xmt_{i,j}$.

In the same way one shows that if $|h(v)| \leq \frac{1}{4}H(v)$ for some
$v\in \Xmt_{i,j}$, then the first alternative in (c) holds. 
\endofproof

\bs

\section{Appendix II}
\setcounter{equation}{0}

\noindent
{\it Proof of Lemma} 6.2. Let $\ep \in (0,\theta/3)$. Then $\ep$ is an expansivity constants for $\sigma_A: \saa \longrightarrow \saa$.
As mentioned in Remark 1 in Sect. 6.2, given $\varphi \in \ff_\theta(\saa)$ there exists a constant $K_\varphi$ so that 
\be
|\varphi_n(x) - \varphi_n(y)| \leq K_\varphi \:\: \mbox{\rm for all } \:\: x,y\in \saa \:\: \mbox{\rm with }\:\: 
\dte(\sigma^i(x), \sigma^i(y)) < \ep \: \forall\:i = 0, 1, \ldots, n-1.
\ee
We can take e.g. $K_\varphi =  \frac{|\varphi|_\theta}{1-\theta}$.

Since the matrix $\aa$ satisfies $\aa^{\tp_0} > 0$ for some integer $\tp_0 > 0$ (see Sect. 2), the shift has the
{\it specification property} (see e.g. p. 60 in \cite{Ba}): for every $\delta > 0$ there exists an integer $p = p(\delta) > 0$ so that 
for all points $x_1,x_2, \ldots, x_\ell \in \saa$ and all integers $n_1, n_2 ,\ldots n_\ell > 0$
and $p_1, p_2, \ldots, p_\ell \geq p$ there exists $z \in \saa$ such that 
$\dte(\sigma^{m(j-1)+ i} (z), \sigma^i(x_j)) \leq \delta$, 
where $m_j = n_1+ \ldots + n_j + p_1 + \ldots + p_j$, for all $i = 0,1, \ldots, n_j-1$ and all $j = 1,2,\ldots, \ell$.

We will use this property with $\delta = \ep$. {\bf Fix $p = p(\ep)$} with the property described in the above definition.

Given an integer $n \geq 1$ and $\delta > 0$,  a subset $E$ of $\saa$ is called 
{\it $(n,\delta)$-separated} if for every $x, y \in E$, $x \neq y$, there exists $j = 0,1 \ldots, n-1$ so that
$\dte(\sigma^j(x), \sigma^j(y)) > \delta$. 
We say that $E$ is a maximal $(n,\delta)$-separated set
if whenever $E \subset E'$ and $E'$ is an $(n,\delta)$-separated set, we have $E = E'$. 

For any $\delta \in (0,\theta)$, $x \in \saa$ and $n \geq 1$,  the open ball
$$B_n(x,\delta) = \{ y\in \saa : \dte(\sigma^i(x), \sigma^i(y)) < \delta \;, \; i = 0,1,\ldots, n-1\}$$
is contained in the cylinder $\cc[x_0,x_1, \ldots, x_{n-1}]$ in $\saa$, so by (4.2), 
$\nu(B_n(x, \delta)) \leq c_2 e^{g_n(x)}$. Moreover, if $\theta^{r} \leq \delta$ for some integer $r \geq 0$,
then $\cc_{n+r} [x_0,x_1, \ldots, x_{n+r-1}] \subset B_n(x,\delta)$, so
$\nu(B_n(x,\delta)) \geq c_1 e^{g_{n+r}(x)} \geq c_1(\delta) e^{g_n(x)}$. So we have something similar to (4.2):
$$c_1(\delta) e^{g_n(x)} \leq \nu(B_n(x, \delta)) \leq c_2 e^{g_n(x)} ,$$
where $c_1(\delta) \geq c_1 e^{r g_0}$. Set $c_1' = c_1(\ep) \leq c_1$.

For any finite $E \subset \saa$, $n\geq 1$ and $\varphi \in \ff_\theta(\saa)$ set
$$Z_n (\varphi, E) = \sum_{x\in E} e^{\varphi_n(x)} .$$ 
It is known (see e.g. \cite{B3} or \cite{HR}) that 
$\di \Pr_\sigma(\varphi) = \lim_{n\to\infty} \log Z_n (\varphi, E_n)$
for any choice of the $(n,\ep)$-separated sets $E_n$ in $\saa$. Moreover, we know from Lemma 3 in \cite{B3}
how fast this convergence is occuring. We will describe this now in details\footnote{Although Lemmas 1 - 3 in \cite{B3}
are proved for homeomorphisms $f: X \longrightarrow X$ on compact metric spaces $X$, the arguments in the 
proofs work  for the shift $\sigma_A : \saa \longrightarrow \saa$. }.

We state a bit more precise version of the first lemma which is what we need here for the one-sided shift
$\sigma_A : \saa \longrightarrow \saa$. We sketch the proof for completeness, although it is almost the same as
the one in \cite{B3}.

\bs

\noindent
{\bf Lemma 10.1.} (Lemma 1 in \cite{B3}) {\it Let $\ep > 0$ be as above. There exists
an integer $N  = N(\ep) \geq 1$ with the following property:  for every $\delta \in [\ep/2, 2\ep]$ there exists a constant
$$\di C_{\delta, \ep}  = C_{\delta, \ep} (\varphi) \leq M e^{K_\varphi}$$
for some global constant $M = M(\ep)> 0$ depending on $\ep$ only such that
$$Z_{n}(\varphi, \delta) \leq C_{\delta, \ep} \, Z_n(\varphi, \ep)$$
for all integers $n > N$.}

\bs

\noindent
{\it Sketch of Proof.}
Let $\ep/2 \leq \delta \leq \ep$.
Then choose the minimal integer $N \geq 0$ so that
$\dte(\sigma^k(x), \sigma^k(y)) \leq 2\ep$ for all $k \leq N$ implies $\dte (x,y) \leq \delta$. 
It is enough to choose $N$ so that $\theta^N \leq \delta$, i.e. we need $\theta^N \leq \ep/2$. Then we choose
$\alpha > 0$ so that $\dte(x,y) \leq \alpha$ implies $\dte(\sigma^k(x), \sigma^k(y)) \leq \delta$ for all $k \leq N$.
We can simply take $\alpha = \theta^N \delta$. 
Let $M = M(\ep)$ be the maximum number
of points in the product metric space $\saa \times \saa$ (with the maximum metric) so that any two points are at least a distance
$\alpha$ apart. 

Let $n > N$. Consider now an arbitrary $(n , \delta)$-separated set $E$ and a
maximal $(n, \ep)$-separated set $F$ in $\saa$. Then for any $x \in E$ there exists $a(x) \in F$ with 
$\dte(\sigma^i(x), \sigma^i(a(x))) \leq \ep$ for all $i < n$. 

Given any $a\in F$, set $E_a = \{ x \in E: a(x) = a\}$. We will show that $|E_a| \leq M$.
It is enough to show that $X = \{ (x, \sigma^{n-N}(x)) : x\in E_a\}$ has at most $M$ distinct elements.
 If $|X| > M$, then by the choice of $M$, there exist $x, y\in X$, $x\neq y$, with
either $\dte(x,y) \leq \alpha$ (and then $\dte(\sigma^i(x), \sigma^i(y)) \leq \delta$ for all $0 \leq i \leq N$), 
or $\dte(\sigma^{n-N}(x),\sigma^{n-N}(y)) \leq \alpha$ (and then $\dte(\sigma^{n-N+ i}(x), \sigma^{n-N+i}(y)) \leq \delta$ 
for all $0 \leq  i \leq N$).
Now $x,y \in E_a$ gives
$$\dte(\sigma^i(x), \sigma^i(y)) \leq \dte(\sigma^i(x), \sigma^i(a)) + \dte(\sigma^i(a), \sigma^i(y)) \leq 2\ep $$
for all $i < n$. Then the choice of $N$ implies that $\dte( \sigma^{i+N} (x) , \sigma^{i+N}(y)) \leq \delta$.
This shows that $\dte(\sigma^k(x), \sigma^k(y)) \leq \delta$ for all $k < n$, which is a contradiction, since
$E$ is $(n,\delta)$-separated. Thus, $|X| \leq M$, so we have $|E_a| \leq M$ for all $a\in F$.

Now using $|\varphi_{n}(x) - \varphi_{n}(a)| \leq K_\varphi$ for $x \in E_a$, we derive
$$\sum_{x\in E} e^{\varphi_{n}(x)} \leq \sum_{a\in F} \sum_{x \in E_a} e^{\varphi_{n}(x)}
\leq \sum_{a\in F} \sum_{x \in E_a} e^{K_\varphi} e^{\varphi_{n}(a)} \leq M e^{K_\varphi} Z_n(\varphi, \ep) .$$
This proves the lemma.
\endofproof

\bs

{\bf Fix an integer $N = N(\ep)$} with the property in Lemma 10.1.

\bs

\noindent
{\bf Lemma 10.2.} (Lemma 2 in \cite{B3}) {\it There exists constants $D_\ep = D_\ep(\varphi) > 0$
and $E_\ep = E_\ep(\varphi) > 0$ such that 
$$E_\ep\, \prod_{j=1}^k Z_{n_j}(\varphi, \ep) \leq Z_{n_j+ n_{j+1} + \ldots + n_k}(\varphi, \ep)
\leq D_\ep\, \prod_{j=1}^k Z_{n_j}(\varphi, \ep)$$
for all integers $n_1, n_2, \ldots, n_k \geq N$ and all $k \geq 1$. Moreover we can take
\be
D_\ep = D_\ep(\varphi) = e^{K_\varphi} \, C_{\frac{\ep}{2}, \ep}(\varphi)
\ee
and}
\be
E_\ep = E_\ep(\varphi) = \frac{e^{- p\, \|\varphi\|_0 + K_\varphi}}{D_\ep (\varphi) \, C_{\ep, 3\ep} (\varphi) \, \max \{ 1 , Z_p(\varphi, \ep)\}}
= \frac{e^{- p\, \|\varphi\|_0}}{C_{\ep/2,\ep} \, C_{\ep, 3\ep} (\varphi) \, \max \{ 1 , Z_p(\varphi, \ep)\}} .
\ee

\ms

The proof of this lemma is given in all details in \cite{B3},  and it works without change in our case using the above Lemma 10.1. 
As a consequence of the above, in the same way as in \cite{B3}, one obtains

\bs

\noindent
{\bf Lemma 10.3.} (Lemma 3 in \cite{B3}) {\it With $\ep > 0$ as above, we have
\be
\frac{1} {D_\ep(\varphi)} \, e^{n \Pr_\sigma(\varphi)} \leq Z_n(\varphi, \ep) \leq \frac{1} {E_\ep(\varphi)} \, e^{n \Pr_\sigma(\varphi)} 
\ee
for all integers $n \geq N$.}

\bs

We now continue with the proof of Lemma 6.2. We follow the argument in the first half of the proof of Lemma 4.1 in \cite{MV}.

Fix an arbitrary $q_0 < 0$ and consider $q\in [q_0, 0]$. Let $|b| \geq b_0$ and let the constant $K > 0$ satisfy the assumptions in Lemma 6.2.
We will use Lemmas 10.1, 10.2 and 10.3  with $\varphi = \varphi(q) = g + q \Psi$. Notice that for a sufficiently large global constant $C_1 = C_1(q_0) > 0$,
the constant $K_\varphi = C_1 K$ satisfies (10.1) for $\varphi$ (using here $\Pr(g) = 0$ and (4.2); see also Lemma 3.12 in \cite{MV}).
Moreover, since $p$ is a constant, we can choose $C_1$ so that $\|p \varphi\|_0 \leq C_1$ and 
$ \max \{ 1 , Z_p(\varphi, \ep)\} \leq C_1$
for any choice of $q\in [q_0,0]$. Thus, $|E_\ep(\varphi)| \leq C_2(\ep)$ for some global constant $C_2(\ep)$ depending only on $\ep$.

Let $E_n = \{y_1, \ldots, y_k\}$ be a maximal $(n,\ep)$-separated set in $\saa$. Then $\cup_{i=j}^k B_n(y_j,\ep) = \saa$.
Thus, for every $x\in \saa$ there exists $j \leq k$ with
$\dte(\sigma^i(x), \sigma^i(y_j)) \leq \ep$ for all $i = 0,1, \ldots, n-1$. This implies $|\Psi_n(x) - \Psi_n(y_j)| \leq K$.
Hence for any $q\in \R$ we have
\begin{eqnarray*}
e^{n f_n(q) }
 =  \int_{\saa} e^{q \Psi_n(x)}\, d\nu(x) \leq \sum_{j=1}^k \int_{B_n(y_j,\ep)} e^{q \Psi_n(x)}\, d\nu(x)
 \leq  \sum_{j=1}^k  e^{|q| K + q \Psi_n(y_j)}\, \nu(B_n(y_j, \ep)) .
\end{eqnarray*}
Since $\nu(B_n(y_j, \ep)) \leq c_2\, e^{g_n(y_j)}$ by (4.2), it follows that
\begin{eqnarray*}
e^{n f_n(q) }
& \leq & c_2 e^{|q| K} \sum_{j=1}^k e^{ q \Psi_n(y_j) + g_n(y_j)} = c_2 e^{|q| K} Z_n(g + q \Psi, E_n) .
\end{eqnarray*}
This and (10.4) with $\varphi = g + q \Psi$ imply
\begin{eqnarray*}
f_n(q) 
 \leq  \frac{\log c_2 + |q| K}{n} + \frac{1}{n} \log Z_n(g+ q \psi, E_n)
 \leq  \Pr_\sigma(g + q \Psi) +  \frac{\log c_2 + |q| K}{n} - \frac{1}{n} \log E_\ep(\varphi) .
\end{eqnarray*}
Now $|E_\ep(\varphi )| \leq C_2(\ep)$ implies
$$\di f_n(q)  \leq  \Pr_\sigma(g + q \Psi) +  \frac{|q| K}{n} + \frac{C(\ep)}{n}$$
for some global constant $C(\ep) > 0$ depending only on $\ep$. This proves the lemma.
\endofproof

\bs

\end{document}